%% file: lmr-paper-1-num-tests-arxiv.tex
\documentclass[a4paper, 10pt]{article}
\pdfoutput=1

\usepackage{amssymb}
\usepackage{graphicx}
\usepackage{colortbl}
\usepackage[makeroom]{cancel}
\usepackage{mathtools}

\definecolor{gainsboro}{rgb}{0.86, 0.86, 0.86}
\definecolor{lightgray}{rgb}{0.83, 0.83, 0.83}
\definecolor{silver}{rgb}{0.75, 0.75, 0.75}
\definecolor{ashgrey}{rgb}{0.7, 0.75, 0.71}
\definecolor{battleshipgrey}{rgb}{0.52, 0.52, 0.51}

\input{customized_packages_and_macros}
\usepackage{url}
\urldef{\mailsa}\path|{smatculevich,ulanger}@ricam.oeaw.ac.at |
\urldef{\mailsb}\path|serepin@jyu.fi, | 
\urldef{\mailsc}\path|repin@pdmi.ras.ru |

\begin{document}


\title{Guaranteed error 
bounds for 
stabilised space-time IgA approximations to parabolic problems}

%
\author{
Ulrich Langer 
\thanks{RICAM, Austrian Academy of Sciences, Linz, Austria, \texttt{ulrich.langer@ricam.oeaw.ac.at}} , 
Svetlana Matculevich 
\thanks{RICAM, Austrian Academy of Sciences, Linz, Austria,       \texttt{svetlana.matculevich@ricam.oeaw.ac.at}} ,
and
Sergey Repin \thanks{University of Jyvaskyla, Finland;
St. Petersburg Department of V.A. Steklov Institute of Mathematics RAS, 
\texttt{serepin@jyu.fi, repin@pdmi.ras.ru}}}
%


\maketitle
	
\begin{abstract}
The paper is concerned with space-time IgA approximations of  parabolic initial-boundary value problems. 
We deduce guaranteed and fully computable error bounds adapted to special features of IgA 
approximations and investigate their applicability. The derivation method is based on 
the analysis of respective integral identities and purely functional arguments. Therefore, the estimates do 
not contain mesh-dependent constants and are valid for any approximation from the admissible (energy) 
class. In particular, they provide computable error bounds for norms associated with {stabilised space--time} 
IgA approximations as well as imply efficient error indicators enhancing the performance of fully adaptive 
solvers. {The last section of the paper contains a series of numerical examples where approximate 
solutions are recovered by IgA techniques. The mesh refinement algorithm 
is governed by a local error indicator generated by the error majorant. 
Numerical results discussed in the last section illustrate both reliability, as well as the quantitative 
efficiency of the error estimates presented.}


%
%
\end{abstract}

\section{Introduction}
\label{sec:introduction}
Time-dependent systems governed by parabolic partial differential equations (PDEs)
are typical models in scientific and engineering applications.
This triggers their active investigation in {modelling, mathematical analysis and numerical solution.}
By virtue of fast development of parallel computers, treating 
time as yet another dimension in space in evolutionary equations became quite natural. 
The \emph{space-time approach} is not affected by the disadvantages of time-marching schemes. 
Its various versions can be useful in combination with parallelisation methods, e.g., those 
discussed in \cite{LMR:Gander:2015,LMR:GanderNeumueller:2016a,LMR:LangerMooreNeumueller:2016a}.

Investigation of effective adaptive refinement methods is crucial for constructing fast and efficient 
solvers for PDEs. At the same time, scheme localisation is strongly linked with reliable and 
quantitatively efficient a posteriori error estimation tools. These tools are intended to identify the areas 
of a computational domain with relatively high discretisation errors and by that provide a fully 
automated refinement strategy in order to reach the desired accuracy level for the current approximation. 
Local refinement tools {in IgA such as T-splines, THB-splines, and LR-splines} were
combined with various a posteriori error estimation techniques, e.g., error estimates using hierarchical 
bases \cite{LMR:DorfelJuttlerSimeon2010,LMR:Vuongetall2011}, 
residual-based 
\cite{LMR:Johannessen2009,LMR:Wangetall2011,LMR:BuffaGiannelli2016,
LMR:KumarKvamsdalJohannessen2015,LMR:GantnerHaberlikPraetorius2017}, 
and goal-oriented error estimates \cite{LMR:ZeeVerhoosel2011,LMR:DedeSantos2012,
LMR:Kuru:2013,LMR:Kuruetall2014}.
Below we use a different (functional) method providing fully guaranteed error estimates
in various weighted norms equivalent to the global energy norm. These estimates include only global 
constants (independent of the mesh characteristic $h$) and are valid for any approximation from the 
admissible functional space.
%
Functional type error estimates (so-called majorants and minorants of deviation from the exact 
solution) were introduced in \cite{LMR:Repin:1997,LMR:Repin:1999} and later applied to different 
mathematical models \cite{LMR:RepinDeGruyterMonograph:2008,Malietall2014}. 
They provide guaranteed, sharp, and fully computable upper and lower bounds of errors. 
This approach, in combination with the IgA approximations generated by tensor-product splines,
was proposed and investigated in \cite{LMR:KleissTomar2015} for elliptic boundary value problems (BVP).
%

{ 
In this paper, we derive {new} functional type a posteriori error estimates for time-dependent problems 
(cf. also \cite{LMR:Repin:2002}) 
in the context of the stabilised space-time IgA scheme introduced 
in \cite{LMR:LangerMooreNeumueller:2016a}.} This scheme exploits the time-upwind test function 
based on the space-time streamline diffusion method (see, e.g., 
\cite{LMR:Hansbo:1994,LMR:Johnson:1987,LMR:JohnsonSaranen:1986}) and the approximations provided 
by the IgA framework. 
{ The obtained functional estimates, in turn, provide bounds for the error measured in new 
stronger norm induced by alternative stabilised space-time variational formulation of the parabolic problem.}
By exploiting the universality and efficiency of considered error estimates as well as the smoothness 
of IgA approximations, we aim to construct fully-adaptive, fast and efficient parallel 
space-time methods that could tackle complicated problems inspired by industrial applications. 
{We also 
study the numerical properties of newly derived error bounds and compare their performance to the 
behaviour of the bounds known from \cite{LMR:Repin:2002} on an extensive set of examples.}

The outline of the paper is as follows: Section~\ref{sec:model-problem} states the problem, discusses its 
solvability and provides an overview of the existing error control tools for 
initial BVPs (I--BVPs). 
In Section~\ref{sec:general-error-estimate}, we deduce new functional type a posteriori error estimates 
using a stabilised formulation of parabolic I--BVPs. Our analysis is based 
{  on a series of transformations performed 
on a stabilised variational setting; the result of these transformations defines respective 
generalised solutions}. Section~\ref{sec:discretization-non-moving-domain} presents 
a stabilised space-time IgA scheme with its main properties along with an overview of main ideas 
and definitions of the IgA framework. 
{Section~\ref{sec:algorithms} is devoted to the algorithmic realisation of 
an adaptive procedure based on the a posteriori error estimates discussed above.
Finally, in Section~\ref{sec:numerical-examples} we present and discuss obtained 
numerical results that demonstrate the efficiency of several majorants and the error identity for a 
comprehensive range of examples.}

\section{Parabolic model problem}
\label{sec:model-problem}

Let {$\overline{Q} := Q \cup \partial Q$, $Q := \Omega \times (0, T)$}, denote a space-time cylinder, 
where $\Omega \subset \Rd$, $d \in \{1, 2, 3\}$, is a bounded Lipschitz domain with a boundary 
$\partial \Omega$, and $(0, T)$ is a given time interval with the final time $T$, $0 < T < +\infty$. 
Here,  {the boundary $\partial Q$ of the space-time cylinder $Q$} is defined as 
$\partial Q := \Sigma \, \cup \, \overline{\Sigma}_{0} \, \cup \, \overline{\Sigma}_{T}$ with 
$\Sigma = \partial \Omega \times (0, T)$,  $\Sigma_{0} =  \Omega \times \{0\}$, and 
$\Sigma_{T} = \Omega \times \{T\}$. 
We discuss our approach to guaranteed error control of space-time approximations within the paradigm 
of a classical {\em linear parabolic I--BVP}: find $u: \overline{Q} \rightarrow \mathds{R}$ satisfying 
{  the parabolic PDE, the boundary condition, and the initial condition }
\begin{equation}
\partial_t u - \laplace_x u = f \;\;\; {\rm in} \;\;\; Q, \quad
u = 0 \;\;\; {\rm on} \;\;\; \Sigma, \quad
\mbox{and}\quad
u = u_0 \;\;\; {\rm on} \;\;\; { \overline{\Sigma}_0}, \label{eq:equation}
\end{equation}
respectively, where $\partial_t$ is a time derivative, $\laplace_x$ denotes the Laplace operator in space, 
{$f \in \L{2}(Q)$ is a given source function, and $u_0 \in \HD{1}{0}(\Sigma_0)$ is 
prescribed initial data.}
%
Here, $\L{2}(Q)$ is a space of square-integrable functions over $Q$ {equipped with the usual norm 
and the scalar product} denoted respectively by 
$$\| \, v \, \|_{Q} := \| \, v \, \|_{\L{2}(Q)} = (v, v)^{\rfrac{1}{2}}_{\L{2}(Q)} 
\quad \mbox{and} \quad
(v, w)_Q = (v, w)_{\L{2}(Q)} := \int_Q v(x, t) w(x, t) \dxt, \quad \forall v, w \in \L{2} (Q).$$ 
By $\H{k}(Q)$, $k \geq 1$, we denote standard Sobolev spaces of functions having derivatives of the 
order $k$ with respect to (w.r.t.) space and time. Next, we introduce the Sobolev 
spaces
\begin{equation}
\begin{alignedat}{3}
\HD{1, 0}{0}(Q) & := \big\{ u \in \L{2}(Q)\!\! && : 
\! \nabla_x u \in [\L{2}(\Omega)]^d, u \big|_{\Sigma} = 0 \big\}, \\
V^1_{0} := \HD{1}{0}(Q) & := \big\{ u \in \H{1}(Q) && : 
u \big|_{\Sigma} = 0 \big\}, \\ 
\HD{1}{0, \overline{0}}(Q) & := \big\{ u \in V^1_{0} && : 
u \big|_{\Sigma_T} = 0 \big\}, \\
V^1_{0,\underline{0}} := H^{1}_{0, \underline{0}} (Q) & := \big\{ u \in V^1_{0} && : 
u\big|_{\Sigma_0} = 0 \big\}, \\
V^{\Delta_x}_{0} := \HD{\Delta_x}{0}(Q) & := \big\{u \in V^1_{0} && :
\, \Delta_x u \in \L{2}(Q) \big\}, \\
V^{\Delta_x}_{0, \underline{0}} := \HD{\Delta_x}{0, \underline{0}}(Q) & := \big\{u \in V^{\Delta_x}_{0} && :
\, u\big|_{\Sigma_0} = 0 \big\}, \quad \mbox{with the norm} \quad 
\| w \|^2_{V^{\Delta_x}_{0, \underline{0}}} 
:= \| \Delta_x w \|^2_Q + \| \partial_t w \|^2_Q, \\
V^{\nabla_x \partial_t}_{0, \underline{0}} = H^{\nabla_x \partial_t}_{0, \underline{0}} (Q) & 
:= \{ w \in V^{\Delta_x}_{0, \underline{0}} && : \nabla_x \partial_t w \in \L{2}(Q) \}. 
\end{alignedat}
\label{eq:spaces}
\end{equation}
%
%
For the vector-valued functions, we additionally introduce the Hilbert spaces 
\begin{alignat*}{2}
	H^{\dvrg_x, 0}(Q) & := 
	\big \{  \flux \in [\L{2}(Q)]^d \;:\; 
	         \dvrg_x \flux \in \L{2} (Q) 
	\big \} \quad \mbox{and} \quad \\
	H^{\dvrg_x, 1}(Q) & := \big\{\flux \in H^{\dvrg_x, 0}(Q) \;:\; 
	         \partial_t \flux \in [\L{2}(Q)]^d \big\}.
	\label{eq:y-set-div-1}
\end{alignat*} 	
%
%
%

{ 
It follows from \cite{LMR:Ladyzhenskaya:1985} that,} 
if $f \in \L{2}(Q)$ and $u_0 \in \HD{1}{0}(\Sigma_0)$, then problem~\eqref{eq:equation} is uniquely 
solvable in $V^{\Delta_x}_{0}$, and the solution $u$ depends continuously on $t$ in the 
$\HD{1}{0}(\Omega)$-norm. Moreover, according to \cite[Remark 2.2]{LMR:Ladyzhenskaya:1985}, 
the norm $\| \, \nabla_{x} u (\cdot, t) \, \|^2_{\Omega}$ is an absolutely continuous function of 
$t \in [0, T]$ for any $u \in V^{\Delta_x}_0$. If $u_0 \in \L{2}(\Sigma_0)$, then the problem is 
uniquely  solvable in a wider class {$H^{1, 0}_0(Q)$}, and meets the modified variational formulation
\begin{equation}
(\nabla_x {u}, \nabla_x {w})_Q - (u, \partial_t w)_Q 
=: a(u, w) = l(w) 
:= (f, {w})_Q + (u_0, {w})_{\Sigma_0}
\label{eq:variational-formulation}
\end{equation}
for all $w \in V^1_{0, \overline{0}}(Q)$, where
$(u_0, w)_{\Sigma_0} := \int_{\Sigma_0} u_0(x) \, {w}(x,0) dx = \int_{\Omega} u_0(x) \,{w}(x,0) dx$.
{  According to well-established arguments (see 
\cite{LMR:Ladyzhenskaya:1954,LMR:Ladyzhenskaya:1985,LMR:Zeidler:1990a,DautrayLions2000Vol5,DautrayLions2000Vol6}), 
without loss of generality, we can `homogenise'}
the problem, i.e., consider \eqref{eq:variational-formulation} with $u_0 = 0$.

Our main goal is to establish fully computable estimates for space-time approximations of this class of 
problems. For this purpose, we use a functional approach to derive a posteriori error estimates. 
The first and the simplest forms of such estimates have been derived in \cite{LMR:Repin:2002} for 
\eqref{eq:equation}. The paper \cite{LMR:Repin:2002} provides the upper bound of the norm 
\begin{equation}
|\!|\!| e |\!|\!|^2_{(\nu_1, \nu_2)} 
:= \nu_1 \, \| \nabla_x e \|^2_Q + \nu_2\,\| e \|^2_{\Sigma_T}, \quad \nu_i \geq 0,
\label{eq:e-energy-norm}
\end{equation}
where $e = u - v$ is the difference between the exact solution $u$ and any approximation $v$ in the 
respective energy class $V^1_{0}$.
{ 
Assuming for simplicity that the initial condition is satisfied exactly, it is shown that for any 
$v \in V^1_{0}$ approximating  $u \in V^1_{0}$ and any
$\flux \in H^{\dvrg_x, 1}(Q)$, we have the following inequality:}
\begin{alignat}{2}
|\!|\!| e |\!|\!|^2_{(2 - \nu, 1)} & := 
(2 - \nu) \, \| e \|^2_Q + \,\| e \|^2_{\Sigma_T} \nonumber\\
& \leq \overline{\rm M}^{\rm I}(v, \flux; \beta) := 
 \tfrac{1}{\nu}\, \Big( (1 + \beta) \, \| \flux - \nabla_x v \|^2_Q
 + (1 + \tfrac{1}{\beta}) \, \CFriedrichs^2 \, \| \dvrg_x \flux + f  - \partial_t v \|^2_Q \Big), 
\label{eq:majorant-heat-eq} 
\end{alignat}
where $\nu \in (0, 2]$ is an auxiliary parameter, and 
$\CFriedrichs$ stands for the constant in the Friedrichs inequality \cite{Friedrichs1937}
\begin{equation}
\| w \|_Q \! \leq \! \CFriedrichs \| \nabla_x w \|_Q, \quad \forall w \in \HD{1, 0}{0}(Q). 
\label{eq:friedrichs}
\end{equation}
The numerical properties of 
\eqref{eq:majorant-heat-eq} w.r.t. the time-marching and space-time methods are discussed in 
\cite{LMR:MatculevichRepin:2014,LMR:MatculevichNeitaanmakiRepin:2015,LMR:HolmMatculevich:2017}
in the framework of finite-difference and finite-element schemes.
The advanced upper error-bound $\overline{\rm M}^{\rm I\!I}(v, \flux, \eta)$ (valid for the same error 
norm) introduced in \cite{LMR:Repin:2002} contains an additional auxiliary function $\eta \in V^1_0$. 
For the same $v$ and $\flux$, as well as any $\eta \in V^1_0$, any fixed parameters 
$\nu \in (0, 2]$ and $\gamma \in [1, + \infty)$, an alternative majorant has 
the form 
\begin{alignat}{2}
(2 - \nu) \, \| \nabla_x e \|^2_Q 
+ (1  & - \tfrac{1}{\gamma}) \, \| e \|_{\Sigma_T} 
\leq \overline{\rm M}^{\rm I\!I} (v, \flux, \eta)
:= \gamma \, \| \eta \|^2_{\Sigma_T} 
+ \| u_0 - v \|^2_{\Sigma_0} - 2 \, (\eta, u_0 - u_h)_{\Sigma_0}
+ 2 \mathcal{F}(v, \eta) \nonumber\\
& + \tfrac{1}{\nu} \, 
\Big \{ (1 +\beta) \, \| \flux - \nabla_x v + \nabla_x \eta \|^2_Q 
           + \CFriedrichs^2 \, (1 + \tfrac{1}{\beta}) \, 
              \| \dvrg_x \flux + f - \partial_t v - \partial_t w \|^2_Q \Big\},
\label{eq:majorant-II}
\end{alignat}
where
$$\mathcal{F}(v, \eta) := (\nabla_x v, \nabla_x \eta) + (\partial_t v - f, \eta).$$
The optimal parameters $\nu$, $\beta$, and $\gamma$ are easily defined in each particular case by 
minimisation of the respective majorant.
%

Finally, we note that for the case where $u, v \in V^{\Delta_x}_0$, the heat equation 
\eqref{eq:equation} imposes the error identity (see \cite{LMR:AnjamPauly:2016}):
\begin{alignat}{2}
\| \Delta_x e \|^2_Q +  \| \partial_t e \|^2_Q +  \| \nabla_x e \|^2_{\Sigma_T} 
=: |\!|\!|  e |\!|\!|^2_{\mathcal{L}} 
\equiv \EI^2 (v)
:= \| \nabla_x (u_0 - v) \|^2_{\Sigma_0} + \|  \Delta_x v + f - \partial_t v\|^2_Q.
\label{eq:strong-norm-error-identity}
\end{alignat}
The numerical performance of estimates $\overline{\rm M}^{\rm I}$ and
$\overline{\rm M}^{\rm I\!I}$, and of the identity 
\eqref{eq:strong-norm-error-identity} is studied in 
Section \ref{sec:numerical-examples}.

\section{Error majorants}
\label{sec:general-error-estimate}

In this section, we derive error majorants of the functional type for a stabilised weak formulation of 
parabolic I--BVPs. They provide guaranteed and fully computable upper 
bounds of the distance between 
{  the exact solution $u$ and some approximation $v$.} 
Functional nature of the majorants allows 
us to obtain
a posteriori error estimates
for any conforming 
approximation
$v \in V^{\Delta_x}_{0, \underline{0}}$.

We begin by testing \eqref{eq:equation} with the time-upwind test function
\begin{equation}
\lambda \, w + \mu \, \partial_t w, \quad 
w \in V^{\nabla_x \partial_t}_{0, \underline{0}}, \quad 
\lambda, \mu \geq 0,
\label{eq:upwind-test}
\end{equation}
and arrive at the stabilised weak formulation for $u \in V^1_{0, \underline{0}}$, i.e., 
\begin{equation}
\big(\partial_t u,  \lambda \, w + \mu \, \partial_t w \big)_Q 
+ \big(\nabla_x {u}, \nabla_x (\lambda \, w + \mu \, \partial_t w) \big)_Q =: a_s (u, w) = l_s (w)
:= (f, \lambda \, w + \mu \, \partial_t w)_Q, \quad \forall w \in V^{\nabla_x \partial_t}_{0, \underline{0}}.
\label{eq:stabilized-bilinear-form}
\end{equation}
%
{Then, the error $e = u - v$, $u, v \in V^{\nabla_x \partial_t}_{0, \underline{0}}$ 
(  this condition 
is required to ensure the existence of the term $\| \nabla_x e \|^2_{\Sigma_T}$}), 
is measured in terms of the norm generated by the bilinear form 
$a_s (u, w)$, i.e., 
\begin{equation}
| \! |\! | e | \! |\! |^2_{s, \nu_i}
:= \nu_{1} \, \| \nabla_x {e}\|^2_{Q} +  \nu_{2} \, \| \partial_t e \|^2_{Q} 
+ \nu_{3} \, \| \nabla_x e \|^2_{\Sigma_T} + \nu_{4} \, \| e \|^2_{\Sigma_T},
\label{eq:error-general}
\end{equation}
where $\{\nu_{i} \}_{i = 1, ..., 4}$ are the positive weights introduced in the derivation process. 

To obtain guaranteed error bounds of $| \! |\! | e | \! |\! |^2_{s, \nu_{i}}$, we apply the method similar 
to the one developed in \cite{LMR:Repin:2002,LMR:MatculevichRepin:2014} for parabolic I--BVPs. 
{  
As a starting point, }we consider the space of functions $V^{\nabla_x \partial_t}_{0, \underline{0}}$ 
(see \eqref{eq:spaces}) equipped with the norm 
$$\| w \|^2_{V^{\nabla_x \partial_t}_{0, \underline{0}}} 
:=  {\sup\limits_{t \in [0, T]} \| \nabla_x w(\cdot, t)\|^2_\Omega}
    + \| w \|^2_{V^{\Delta_x}_{0, \underline{0}}},
$$
where
$\| w \|^2_{V^{\Delta_x}_{0, \underline{0}}} 
:= \| \Delta_x w \|^2_Q + \| \partial_t w \|^2_Q$, 
which is dense in $V^{\Delta_x}_{0, \underline{0}}$. 
According to \cite[Remark 2.2]{LMR:Ladyzhenskaya:1985}, the norms 
$\| \cdot \|_{V^{\nabla_x \partial_t}_{0, \underline{0}}}$
and
$\| \cdot \|_{V^{\Delta_x}_{0, \underline{0}}}$
are equivalent.
%
Below, we exploit the density this property 
to derive the majorants \eqref{eq:error-general}
in Theorems \ref{th:theorem-majorant-general-1} and \ref{th:theorem-majorant-general-2}.
%
%
\begin{theorem}
\label{th:theorem-majorant-general-1}
%
For any 
$v \!\in\! V^{\Delta_x}_{0, \underline{0}}$ and ${\flux} \!\in\! H^{\dvrg_x, 0}(Q)$, 
the following estimate holds:
{
 
\begin{alignat}{2}
(2 - \tfrac{1}{\gamma}) \, 
(\lambda \,\| \nabla_x e \|^2_{Q}  & + \mu \, \| \partial_t e \|^2_{Q}) 
+ \lambda \, \| e \|^2_{\Sigma_T} + \mu \, \| \nabla_x{e} \|^2_{\Sigma_T} \nonumber\\
& \qquad\qquad 
=: | \! |\! |  e | \! |\! |^{{\rm I}, 2}_{s} 
\leq \overline{\rm M}^{\rm I}_{s, h} (v, \flux; \gamma, \beta, \alpha) 
:= \gamma\Big\{ \lambda \, \overline{\rm M}^{\rm I}(v, \flux; \beta)
		       + \mu \, \widetilde{\rm M}^{\rm I}(v, \flux; \alpha) \Big\},
\label{eq:estimate-1-non-moving}
\end{alignat}
where $\overline{\rm M}^{\rm I}(v, \flux; \beta)$ is the majorant defined in 
\eqref{eq:majorant-heat-eq} with $\nu= 1$, i.e.,  
$$
\overline{\rm M}^{\rm I}(v, \flux; \beta) := 
 (1 + \beta) \, \| \R^{\rm I}_{\rm d} \|^2_Q
 + (1 + \tfrac{1}{\beta}) \, \CFriedrichs^2 \, \| \R^{\rm I}_{\rm eq} \|^2_Q
$$
and
$$
\widetilde{\rm M}^{\rm I}(v, \flux; \alpha) :=
(1 + \alpha) \, \| \dvrg_x \R^{\rm I}_{\rm d}\|^2_{Q} 
+ (1 + \tfrac{1}{\alpha})\,\| \R^{\rm I}_{\rm eq} \|^2_{Q}.$$
Here, 
$\CFriedrichs$ is {the} Friedrichs constant \eqref{eq:friedrichs},
$\R^{\rm I}_{\rm eq}$ and $\R^{\rm I}_{\rm d}$ are the residuals defined by relations 
\begin{alignat}{2}
	\R^{\rm I}_{\rm eq}  (v, \flux) & := f - \partial_t v + \dvrg_x \: \flux \quad \mbox{and} \quad 
	\R^{\rm I}_{\rm d}  (v, \flux) := \flux - \nabla_x {v},
\label{eq:residuals}
\end{alignat}
$\lambda, \mu >0$ are the weights introduced in \eqref{eq:upwind-test},  
$\gamma \in \big[\tfrac{1}{2}, +\infty)$, and $\alpha, \beta > 0.$ 
}
\end{theorem}
%
%
\noindent{\bf Proof:}  
{ 
Let $\{ u_n \}_{n = 1}^{\infty}$ be a sequence in $V^{\nabla_x \partial_t}_{0, \underline{0}}$
such that 
$u_n  \rightarrow u$ in  $V^{\Delta_x}_{0, \underline{0}}$
for $n \rightarrow \infty$.
It satisfies the identity (cf. \eqref{eq:stabilized-bilinear-form})
}
\begin{equation}
a_s (u_n, w) = (f_n, \lambda \, w + \mu \, \partial_t w)_Q, \quad \mbox{where} \quad
f_n = {(u_n)}_t - \Delta_x u_n \in \L{2}(Q).
\label{eq:stabilized-bilinear-form-short}
\end{equation}
%
{ 
Next, we consider a sequence 
$\{ v_n \}_{n = 1}^{\infty} \in V^{\nabla_x \partial_t}_{0, \underline{0}}$
approximating 
$\{ u_n\}_{n = 1}^{\infty}$
in the sence that 
$v_n \rightarrow   u$ in $V^{\Delta_x}_{0, \underline{0}}$ 
for $n \rightarrow \infty$.
}
By subtracting $a_s(v_n, w)$ from the left- and right-hand side (LHS and 
RHS, respectively) of \eqref{eq:stabilized-bilinear-form-short} and by setting the test function $w$ to be the following 
difference $e_n = u_n - v_n \in V^{\nabla_x \partial_t}_{0, \underline{0}}$, we arrive at the error 
identity
\begin{multline}
	\lambda \, \| \nabla_x {e_n} \|^2_Q  + \mu \, \| \, \partial_t e_n \|^2_Q  
	+ \tfrac12 \, (\mu \, \| \nabla_x{e_n} \|^2_{\Sigma_T} + \lambda \| e_n \|^2_{\Sigma_T}) \\
	= \lambda \Big( (f_n - \partial_t v_n,  e_n)_Q - (\nabla_x {v_n}, \nabla_x {e_n})_{Q} \! \Big) 
	+ \mu \Big( (f_n - \partial_t  v_n, \partial_t  e_n)_Q 
	                  - (\nabla_x {v_n}, \nabla_x \, \partial_t e_n)_{Q} \!\Big).
	\label{eq:energy-balance-equation-general-estimate}
\end{multline}

{First, we modify the RHS of \eqref{eq:energy-balance-equation-general-estimate} by means of the relation} 
$$(\dvrg_x \flux, \lambda \, e_n  \!+\! \mu \, \partial_t e_n)_Q 
+ (\flux, \nabla_x (\lambda \, e_n  \!+\! \mu \, \partial_t e_n))_Q = 0.$$
The obtained result can be presented as follows:
\begin{alignat}{2}
\lambda \, \| \nabla_x {e_n} \|^2_Q + \mu \, \| \partial_t \, e_n \|^2_Q
& + \tfrac12 \, ( \mu \, \| \nabla_x{e_n} \|^2_{\Sigma_T} + \lambda \, \| e_n \|^2_{\Sigma_T}) \nonumber\\
	& = \lambda \, \big( (f_n - \partial_t v_n + \dvrg_x \: \flux,  e_n)_Q + ( \flux - \nabla_x {v_n}, \nabla_x {e_n})_{Q} \big) \nonumber\\
	& \qquad \qquad \qquad \qquad \qquad + \mu \, \big( (f_n - \partial_t v_n + \dvrg_x \: \flux, \partial_t {e_n})_Q + ( \flux - \nabla_x {v}, \nabla_x \partial_t {e_n})_{Q}\big) \nonumber\\
	& = \lambda \, \big( ( \R^{\rm I}_{\rm eq}  (v_n, \flux),  e_n)_Q + (\R^{\rm I}_{\rm d}  (v_n, \flux), \nabla_x {e_n})_{Q} \big) \nonumber\\
	 & \qquad \qquad \qquad \qquad \qquad+ \mu \, \big( (\R^{\rm I}_{\rm eq}  (v_n, \flux), \partial_t {e_n})_Q + (\R^{\rm I}_{\rm d}  (v_n, \flux), \nabla_x \partial_t {e_n})_{Q}\big).    
	\label{eq:energy-balance-equation-with-flux}
\end{alignat}
%
{ 
Integrating by parts {
in} the term $(\R^{\rm I}_{\rm d}(v_n, \flux), \nabla_x \partial_t {e_n})_{Q}$ leads to the 
identity}
\begin{equation*}
\mu\, \big(\R^{\rm I}_{\rm d}, \nabla_x (\partial_t e_n)\big)_Q 
= 
- \mu \, (\dvrg_x ( \flux - \nabla_x {v_n}), \partial_t e_n)_{Q}
= - \mu \, (\dvrg_x \flux - \Delta_x v_n, \partial_t e_n)_{Q}.
\end{equation*}
Using density arguments, i.e.,  
$u_n \rightarrow u$,  
$v_n \rightarrow v$ in $ V^{\Delta_x}_{0, \underline{0}}$, and 
$f_n \rightarrow$ in $\L{2}(Q)$, as ${n \rightarrow \infty}$,
we arrive at the identity formulated for $e = u - v$ with $u, v \in V^{\Delta_x}_{0, \underline{0}}$, i.e., 
\begin{alignat}{2}
	\lambda \, \| \nabla_x {e} \|^2_Q 
	 + \mu \, \| \partial_t \, e \|^2_Q
	& + \tfrac12 \, ( \mu \, \| \nabla_x{e} \|^2_{\Sigma_T} + \lambda \, \| e \|^2_{\Sigma_T}) \nonumber\\
	& = \lambda \, \big(\left(\R^{\rm I}_{\rm eq},  e \right)_Q + (\R^{\rm I}_{\rm d}, \nabla_x {e})_{Q} \big)
	      + \mu \, \big( \left(\R^{\rm I}_{\rm eq}, \partial_t {e} \right)_Q 
	      - \mu \, (\dvrg_x \R^{\rm I}_{\rm d}, \partial_t e)_{Q} \big).
	\label{eq:energy-balance-equation-with-flux-new}
\end{alignat}
By means of the H\"{o}lder, Friedrichs, and Young inequalities with positive scalar-valued parameters
$\gamma$, $\beta$, and $\alpha$, we deduce estimate \eqref{eq:estimate-1-non-moving}.
\hfill $\square$
\vskip 10pt

The next theorem assumes higher regularity on the approximation $v$ and the auxiliary function $\flux$. 

\begin{theorem}
\label{th:theorem-majorant-general-2} 
{ 
For any $v \!\in\! V^{\nabla_x \partial_t}_{0, \underline{0}}$ and ${\flux} \!\in\! H^{\dvrg_x, 1}(Q)$, 
we have the inequality
\begin{alignat*}{2}
(2 - \tfrac{1}{\zeta}) (\lambda \,\| \nabla_x {e}\|^2_{Q} & + \mu \, \| \partial_t e \|^2_{Q}) 
+ \mu \, (1 - \tfrac{1}{\epsilon}) \| \nabla_x e \|^2_{\Sigma_T} 
+ \lambda\, \| e \|^2_{\Sigma_T}  =: | \! |\! | e | \! |\! |^{{\rm I\!I}, 2}_{s} 
\leq \overline{\rm M}^{\rm I\!I}_{s, h}(v, \flux; 
\zeta, \alpha, \epsilon, \beta) \\
& \quad
:= \epsilon \, \mu \| \R^{\rm I}_{\rm d} \|^2_{\Sigma_T} 
+ \zeta \, \Big( \lambda \big( (1 + \alpha) \, 
   \overline{\rm M}^{\rm I}(v, \flux; \beta)
+ (1 + \tfrac{1}{ \alpha}) \, \tfrac{\mu^2}{\lambda^2} \, 
\| \partial_t \R^{\rm I}_{\rm d} \|^2_Q\big)
  + \mu \, \| \R^{\rm I}_{\rm eq} \|_{Q}^2 \Big),
\end{alignat*}
where $\overline{\rm M}^{\rm I}(v, \flux; \beta)$ is the majorant defined in 
\eqref{eq:majorant-heat-eq} with $\nu= 1$, 
$\CFriedrichs$ is the Friedrichs constant in \eqref{eq:friedrichs}, 
$\R^{\rm I}_{\rm eq} (v, y)$ and $\R^{\rm I}_{\rm d}  (v, y)$ are the residuals defined in 
\eqref{eq:residuals}, 
$\lambda, \mu > 0 $ are the parameters from \eqref{eq:upwind-test}, 
$\zeta \in \big[\tfrac{1}{2}, +\infty)$, $\epsilon \in [1, +\infty)$, 
and $\beta, \alpha >0$.
}
\end{theorem}
\section{Stabilised formulation  and its IgA discretisation}
\label{sec:discretization-non-moving-domain}

{   For the reader convenience, we recall the general concept of the IgA approach, the definitions
of B-splines (NURBS), and their use in geometrical representation of the space-time cylinder $Q$
as well as  
in the construction of the IgA trial spaces,
{ 
which are used to approximate the solution of 
the variational problem \eqref{eq:variational-formulation}.
}
}

Throughout the paper, $p \geq 2$ denotes the degree of polynomials used for the IgA approximations,
whereas $n$ denotes the number of basis functions used to construct a $B$-spline curve. A {\em knot-vector} 
is a non-decreasing set of coordinates in the parameter domain, written as 
$\Xi = \{ \xi_1, ..., \xi_{n+p+1}\}$, $\xi_i \in \mathds{R}$, where $\xi_1 = 0$ and $\xi_{n+p+1} = 1$. 
The knots can be repeated, and multiplicity of the $i$-th knot is indicated by $m_i$. In what follows, 
we consider only open knot vectors, i.e., $m_1 = m_{n+p+1} = p+1$. 
For the one-dimensional parametric domain $\Qhat := (0, 1)$, ${\mathcal{\Khat}}_h := \{ \Khat \}$ 
denotes a locally quasi-uniform mesh, where each element $\Khat \in \mathcal{\Khat}_h$ is constructed 
by distinct neighbouring knots. The global size of $ \mathcal{\Khat}_h$ is denoted by 
$\hhat := \max_{\Khat \in \mathcal{\Khat}_h} \{ \hhat_{\Khat}\}$, 
\mbox{where} 
$\hhat_{\Khat} := {\rm diam} (\Khat).$

The \emph{univariate B-spline basis functions $\Bhat_{i, p}: \Qhat \rightarrow \mathds{R}$} are defined 
by means of Cox-de Boor recursion formula
\begin{alignat}{2}
\Bhat_{i, p} (\xi) := \tfrac{\xi - \xi_i}{\xi_{i+p} - \xi_i} \, \Bhat_{i, p-1} (\xi)
                         + \tfrac{\xi_{i+p+1} - \xi}{\xi_{i+p+1} - \xi_{i+1}} \Bhat_{i+1, p-1} (\xi), 
                         \quad 
\Bhat_{i, 0} (\xi) \, := 
\begin{cases} 
1 & \mbox{if} \quad \xi_i \leq \xi \leq \xi_{i+1}  \\
0 & \mbox{otherwise}
\end{cases}
,
\end{alignat}
%
and are $(p-m_i)$-times continuously differentiable across the $i$-th knot with multiplicity $m_i$. 
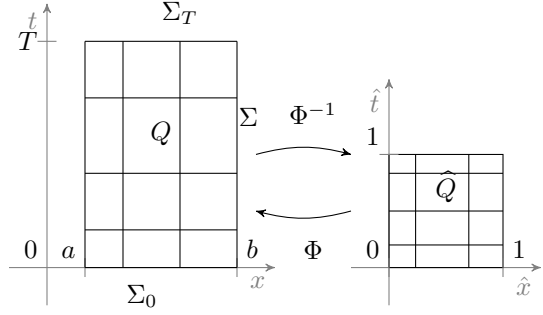
\begin{figure}
\centering
\begin{tikzpicture}[scale=0.5]
\def \ticksize {0.25};
\def \T {6.0};
\coordinate [label={above left:$a$}] (a) at (1.0, 0.0);
\coordinate [label={above right:$b$}] (b) at (5.0, 0.0);
\coordinate [label={above right:$$}] (aT) at (1.0, \T);
\coordinate [label={above right:$$}] (bT) at (5.0, \T);
\coordinate [label={below:$\Sigma_0$}] (Sigma_0) at (2.5, -0.2);
\coordinate [label={above:$\Sigma_T$}] (Sigma_T) at (3.5, \T + 0.2);
\coordinate [label={left:$T$}] (T) at (0.0, \T);
\coordinate [label={above left:$0$}] (O) at (0.0, 0.0);
\coordinate [label={above:$Q$}] (Q) at (3.0, \T/2);
\coordinate [label={above:$\Sigma$}] (sigma) at (5.3, 3.5);

\coordinate [label={above left:$$}] (xi10) at (2.0, 0.0);
\coordinate [label={above left:$$}] (xi20) at (3.5, 0.0);

\coordinate [label={above left:$$}] (xi1T) at (2.0, \T);
\coordinate [label={above left:$$}] (xi2T) at (3.5, \T);

\coordinate [label={above left:$$}] (yi10) at (1.0, 1.0);
\coordinate [label={above left:$$}] (yi20) at (5.0, 1.0);

\coordinate [label={above left:$$}] (yi12) at (1.0, 2.5);
\coordinate [label={above left:$$}] (yi22) at (5.0, 2.5);

\coordinate [label={above left:$$}] (yi13) at (1.0, 4.5);
\coordinate [label={above left:$$}] (yi23) at (5.0, 4.5);

\draw[->,thin,gray] (-1,0) --++(7,0)node[below left]{$x$};
\draw[->,thin,gray] (0,-1) --++(0,8)node[below left]{$t$};

\draw[thin, gray] (1, \ticksize) -- (1, -\ticksize);
\draw[thin, gray] (5.0, \ticksize) -- (5.0, -\ticksize);
\draw[thin, gray] (0.0-\ticksize, 6.0) -- (0.0 +\ticksize, 6.0);

\draw[black] (a) -- (b);
\draw[black] (aT) -- (bT);
\draw[black] (a) to [bend left=0] (aT);
\draw[black] (b) to [bend left=0] (bT);

\draw[black] (xi10) to [bend left=-0] (xi1T);
\draw[black] (xi20) to [bend left=0] (xi2T);

\draw[black] (yi10) -- (yi20);
\draw[black] (yi12) -- (yi22);
\draw[black] (yi13) -- (yi23);

\draw[->,thin,gray] (8,0) --(13,0)node[below left]{$\hat{x}$};
\draw[->,thin,gray] (9,-1) --(9,5)node[below left]{$\hat{t}$};

\draw[black] (9, 0) -- (12, 0) -- (12, 3) -- (9, 3) -- (9, 0);
\coordinate [label={above:$\Qhat$}] (Q) at (10.5, \T/4);

\coordinate [label={above left:$$}] (xi10hat) at (9.7, 0.0);
\coordinate [label={above left:$$}] (xi20hat) at (9.7, 3.0);

\coordinate [label={above left:$$}] (xi1That) at (11.1, 0.0);
\coordinate [label={above left:$$}] (xi2That) at (11.1, 3.0);

\coordinate [label={above left:$$}] (yi10hat) at (9.00, 0.6);
\coordinate [label={above left:$$}] (yi20hat) at (12.00, 0.6);

\coordinate [label={above left:$$}] (yi12hat) at (9.00, 1.5);
\coordinate [label={above left:$$}] (yi22hat) at (12.00, 1.5);

\coordinate [label={above left:$$}] (yi13hat) at (9.00, 2.5);
\coordinate [label={above left:$$}] (yi23hat) at (12.00, 2.5);

\draw[black] (xi10hat) -- (xi20hat);
\draw[black] (xi1That) -- (xi2That);

\draw[black] (yi10hat) -- (yi20hat);
\draw[black] (yi12hat) -- (yi22hat);
\draw[black] (yi13hat) -- (yi23hat);

\draw[black, ->] (5.5, \T/2) to [bend left=15] (8.0, \T/2);
\draw[black, <-] (5.5, \T/4) to [bend left=-15] (8.0, \T/4);

\draw (7, \T/2 + 0.3) node[label={above:{$\Phi^{-1}$}}] (Phiinverse) {};  
\draw (7, \T/4 - 0.3) node[label={below:{$\Phi$}}] (Phi) {};  

\coordinate [label={above left:$0$}] (Ohat) at (9.0, 0.0);
\coordinate [label={above right:$1$}] (1hat) at (12.0, 0.0);
\coordinate [label={above left:$1$}] (1hat) at (9.0, 3.0);

\draw[thin, gray] (9, \ticksize) -- (9, -\ticksize);
\draw[thin, gray] (12, \ticksize) -- (12, -\ticksize);
\draw[thin, gray] (-\ticksize + 9.0, 3.0) -- (9.0+\ticksize, 3.0);

\end{tikzpicture}
\caption{Mapping of the single-patch reference (parameter) domain $\Qhat$ 
to the physical single-patch space-time cylinder $Q$.}
\end{figure}
The {\emph{multivariate B-splines}} on the parameter domain $\Qhat := (0, 1)^{d+1}$, $d = \{1, 2, 3\}$, 
are defined as a tensor-product of the corresponding {univariate} ones. In the multidimensional case, 
we define the knot-vector dependent on the coordinate direction 
$\Xi^\alpha = \{ \xi^\alpha_1, ..., \xi^\alpha_{n^\alpha+p^\alpha+1}\}$, 
$\xi^\alpha_i \in \mathds{R}$, where $\alpha = 1, ..., d+1$ indicates the direction (in space or in time). 
Furthermore, we introduce a set of multi-indices
${\mathcal{I}} = \big\{\, i = (i_1,  ..., i_{d+1}): i_\alpha = 1, ..., n_\alpha$, 
$\,\alpha = 1, ..., d+1 \big\}$ and a multi-index $p := (p_1, ..., p_{d+1})$ indicating the order of 
polynomials. The tensor-product of univariate B-spline basis functions generates multivariate B-spline 
basis functions 
$$\Bhat_{i, p} ({\boldsymbol \xi}) := \prod_{\alpha = 1}^{d+1} \Bhat_{i_\alpha, p_\alpha} (\xi^\alpha), 
\quad 
{\boldsymbol \xi} = (\xi^1, ..., \xi^{d+1}) \in \Qhat.$$
The \emph{univariate and multivariate NURBS basis} functions are defined in the parametric domain by 
means of B-spine basis functions, i.e.,  for the given $p$ and any $i \in {\mathcal{I}}$, the NURBS basis 
functions $\Rhat_{i, p}: \Qhat \rightarrow \mathds{R}$ are defined as 
\begin{equation}
\Rhat_{i, p} ({\boldsymbol \xi}) 
:= \tfrac{w_i \, \Bhat_{i, p} ({\boldsymbol \xi})}{W({\boldsymbol \xi})}.
\end{equation}
Here, $W({\boldsymbol \xi})$ is the weighting function
%
$W({\boldsymbol \xi}) := \sum_{i \in {\mathcal{I}}} w_i \, \Bhat_{i, p} ({\boldsymbol \xi}),$
%
where $w_i \in \mathds{R}^+$. 

The physical space-time domain $Q \subset \mathds{R}^{d+1}$ is defined by the geometrical mapping  
$\Phi: \Qhat \rightarrow Q$ of the parametric domain $\Qhat := (0, 1)^{d+1}$:
\begin{equation}
Q := \Phi(\Qhat ) \subset \mathds{R}^{d+1}, \quad 
\Phi({\boldsymbol \xi}) := \sum_{i \in {\mathcal{I}}} \Rhat_{i, p}({\boldsymbol \xi}) \, {\bf P}_i,
\label{eq:geom-mapping}
\end{equation}
where $\{{\bf P}_i\}_{i \in \mathcal{I}} \in \mathds{R}^{d+1}$ are the control points. For simplicity, 
we assume the same polynomial degree for all coordinate directions, i.e., $p_{\alpha} = p$ for all 
$\alpha = 1, ... , d+1$. %
By means of geometrical mapping \eqref{eq:geom-mapping}, the mesh $\mathcal{K}_h$ 
discretising $Q$ is defined as 
$\mathcal{K}_h := \big\{K = \Phi(\Khat) : \Khat \in \mathcal{\Khat}_h \big\}$.
The global mesh size is denoted by 
\begin{equation}
h := \max\limits_{K \in \mathcal{K}_h} \{ \, h_{K}\,\}, \quad  
h_{K} := \| \nabla \Phi \|_{\L{\infty} (K)} \, \hhat_{\Khat}.
\label{eq:global-mesh-size}
\end{equation}
%
Moreover, we assume that $\mathcal{K}_h$ is a quasi-uniform mesh, i.e., there exists a positive constant 
$C_u$ independent of $h$, such that 
%
$h_{K} \leq h \leq C_u \, h_{K}.$

The finite dimensional spaces $V_h$ on $Q$ are constructed by a push-forward of the NURBS basis 
functions, i.e., 
\begin{equation}	
V_h := {\rm span} \,\big\{ \phi_{h,i} := \Rhat_{i, p} \circ \Phi^{-1} \big\}_{i \in \mathcal{I}},
\label{eq:vh-v0h}
\end{equation}
where  the geometrical mapping $\Phi$ is invertible in $Q$, with smooth inverse on each element 
$K \in \mathcal{K}_h$, see, e.g., see \cite{LMR:TagliabueDedeQuarteroni:2014,LMR:Bazilevsetal2006}. 
The subspace
$$V_{0h} := V_h \cap V^{\Delta_x, 1}_{0, \underline{0}}$$
is introduced for the functions satisfying homogeneous {  boundary and initial conditions.}

A stable space-time IgA scheme for \eqref{eq:equation} was presented and analysed in 
\cite{LMR:LangerMooreNeumueller:2016a}, where the authors proved its efficiency for fixed and 
moving spatial computational domains. In our analysis, we use spline bases of sufficiently high order, 
so that {$v_h \in V_{0h} \subset V^{\Delta_x, 1}_{0, \underline{0}}$.}
In order to provide an efficient discretisation method, we consider \eqref{eq:stabilized-bilinear-form}, 
where $\lambda = 1$ and $\mu = \delta_h = \theta h$ in \eqref{eq:upwind-test} with { 
some positive parameter $\theta$ and the global mesh-size $h$ 
defined in \eqref{eq:global-mesh-size}.
}
%
It implies the {  stabilised space-time IgA scheme: find $u_h \in V_{0h}$ satisfying 
\begin{equation}
(\partial_t u_h, v_h + \delta_h \partial_t v_h)_Q
+ \big (\nabla_x {u_h}, \nabla_x ({v_h} + \delta_h \partial_t v_h) \big)_Q
=: a_{s,h} (u_h, v_h) = l_{s,h}(v_h) 
:= (f, v_h + \delta_h \, \partial_t v_h)_Q.
\label{eq:discrete-scheme}
\end{equation}
for all $v_h \in V_{0h}$.
}
The $V_{0h}$-coercivity of $a_h(\cdot, \cdot): V_{0h} \times V_{0h} \rightarrow \mathds{R}$ 
w.r.t. the norm 
\begin{equation}
|\!|\!| v_h |\!|\!|^2_{s, h} 
:= \| \nabla_x {v_h }\|^2_{Q} 
+ \delta_h \, \| \partial_t v_h  \|^2_{Q} + 
\| v_h \|^2_{\Sigma_T} 
+ \delta_h \, \| \nabla_x v_h  \|^2_{\Sigma_T} 
\label{eq:error-non-moving}
\end{equation}
{follows from \cite[Lemma 1]{LMR:LangerMooreNeumueller:2016a}}.
%
{ 
As 
was noted in \cite{LMR:LangerMooreNeumueller:2016a}, coercivity 
implies
that the 
IgA
solution 
$u_h \in V_{0h}$ of \eqref{eq:discrete-scheme} is unique. Moreover, since 
the
IgA 
scheme 
\eqref{eq:discrete-scheme}
is posed in 
the finite dimensional space $V_{0h}$, uniqueness 
yields
existence of the solution in $V_{0h}$. 
}
{Moreover, following \cite{LMR:LangerMooreNeumueller:2016a,UL:LangerMatculevichRepin:2017a},
{  we can show boundedness of the bilinear form in \eqref{eq:discrete-scheme}}
{  with respect to}
appropriately chosen norms. Combining coercivity and boundedness properties of  
$a_{s, h}(\cdot, \cdot)$ with{ the }consistency of the scheme and 
approximation results for IgA spaces implies a corresponding a priori error estimate 
presented in Theorem \ref{th:theorem-8} below.
}

\begin{theorem}
\label{th:theorem-8}
Let $u \in H^{s}_{0}(Q) := H^{s}(Q) \cap H^{1, 0}_{0}(Q)$, $s \in \mathds{N}$, $s \geq 2$, 
be the exact solution to \eqref{eq:variational-formulation}, and let $u_h \in V_{0h}$ be the solution to 
\eqref{eq:discrete-scheme} with some fixed parameter $\theta$. Then, the following a priori error estimate 
%
\begin{equation}
\| u - u_h\|_{s, h} \leq C \, h^{r-1}\, \| u \|_{H^r(Q)} 
\label{eq:error-in-h}
\end{equation}
holds, 
{where $r = \min \{ s, p+1 \}$, and 
$C >0$ is a generic constant independent of $h$.
}
\end{theorem}

\ProofBegin
See, e.g., \cite[Theorem 8]{LMR:LangerMooreNeumueller:2016a}.
\ProofEnd

\vskip 5pt
\noindent
{ 
Corollary \ref{cor:majorant-1} presents a posteriori error majorants for $\lambda = 1$ and 
$\mu = \delta_h$, where $\delta_h = \theta \, h$, $\theta >0$.
\begin{corollary}
\label{cor:majorant-1}
(i) 
If $v \!\in\! V^{\Delta_x}_{0, \underline{0}}$ and ${\flux} \!\in\! H^{\dvrg_x, 0}(Q)$, 
Theorem \ref{th:theorem-majorant-general-1} yields the estimate
\begin{alignat}{2}
(2 - \tfrac{1}{\gamma}) \, 
(\| \nabla_x e \|^2_{Q}  + \delta_h  \, \| \partial_t e \|^2_{Q}) 
& + \| e \|^2_{\Sigma_T} + \delta_h \, \| \nabla_x{e} \|^2_{\Sigma_T}  
=: | \! |\! | e | \! |\! |^{{\rm I }, 2}_{s, h} \nonumber \\
& \leq \overline{\rm M}^{\rm I}_{s, h}(v, \flux; \gamma, \beta, \alpha) 
:= \gamma \, \Big( \overline{\rm M}^{\rm I}(v, \flux; \gamma, \beta) 
			+ \delta_h\, \widetilde{\rm M}^{\rm I}(v, \flux; \alpha) \Big),
\label{eq:estimate-deltah-1-non-moving}
\end{alignat}
where $\overline{\rm M}^{\rm I}$ and $\widetilde{\rm M}^{\rm I}$ are defined in 
Theorem \ref{th:theorem-majorant-general-1} and the best $\beta$ and $\alpha$ 
are given by relations 
$$\beta = \tfrac{\CFriedrichs \, \|\R^{\rm I}_{\rm eq}\|_{Q}}{\| \R^{\rm I}_{\rm d} \|_{Q}} 
\quad \mbox{ and} \quad 
\alpha = \tfrac{\|\R^{\rm I}_{\rm eq}\|_{Q}}{\| \dvrg_x \R^{\rm I}_{\rm d} \|_{Q}}.$$

\noindent
A particularly useful form of \eqref{eq:estimate-deltah-1-non-moving} follows $\gamma = 1$, i.e., 
\begin{alignat*}{2}
\,\| \nabla_x {e}\|^2_{Q} & + \delta_h \, \| \partial_t e \|^2_{Q} 
+ \| e \|^2_{\Sigma_T} + \delta_h \, \| \nabla_x e \|^2_{\Sigma_T} =: | \! |\! | e | \! |\! |^2_{s, h}
\leq \overline{\rm M}^{\rm I}_{s, h}(v, \flux; \alpha, \beta) 
:= \overline{\rm M}^{\rm I}(v, \flux; \beta) 
       + \delta_h\, \widetilde{\rm M}^{\rm I}(v, \flux; \alpha).
\end{alignat*}
\noindent
(ii) 
If $v \in V^{\nabla_x \partial_t}_{0, \underline{0}}$ and ${\flux} \in H^{\dvrg_x, 1}(Q)$, then
Theorem \ref{th:theorem-majorant-general-2}
yields 
\begin{alignat}{2}
(2 - \tfrac{1}{\zeta}) (\| \nabla_x {e}\|^2_{Q} & + \delta_h\, \| \partial_t e \|^2_{Q}) 
+  (1 - \tfrac{1}{\epsilon}) \| \nabla_x e \|^2_{\Sigma_T} 
+ \delta_h \, \| e \|^2_{\Sigma_T}  
=: | \! |\! | e | \! |\! |^{{\rm I\!I}, 2}_{s, h} 
\leq \overline{\rm M}^{\rm I\!I}_{s, h}
(v, \flux; \zeta, \beta, \alpha, \epsilon^{\rm I}) \nonumber \\
& := \epsilon \, \delta_h \, \| \R^{\rm I}_{\rm d} \|^2_{\Sigma_T} 
+ \zeta \, \Big( 
(1 + \alpha) \, \overline{\rm M}^{\rm I}(v, \flux; \gamma, \beta)  
+ (1 + \tfrac{1}{\alpha}) \delta_h^2 \, \| \partial_t \R^{\rm I}_{\rm d} \|^2_Q
+ \delta_h \, \| \R^{\rm I}_{\rm eq} \|_{Q}^2 \Big),
\label{eq:estimate-2-non-moving-domain-2} 
\end{alignat}
where the optimal parameters are given by relations 
$$\beta = \tfrac{\CFriedrichs \| \R^{\rm I}_{\rm eq}\|_{Q}}{\| \R^{\rm I}_{\rm d} \|_{Q}} 
\quad {and} \quad 
\alpha = 
\tfrac{\delta_h \, \| \partial_t \R^{\rm I}_{\rm d} \|_{Q}}{\sqrt{(1 + \beta) \, \| \R^{\rm I}_{\rm d} \|^2_{Q} 
				+ \big(1 + \tfrac{1}{\beta}\big) \CFriedrichs^2 \, \|\R^{\rm I}_{\rm eq}\|^2_{Q}}}.$$
\noindent
For $\zeta = 1$ and $\epsilon = 2$, we obtain
\begin{alignat}{2}
\,\| \nabla_x {e} & \|^2_{Q}  + \delta_h\, \| \partial_t e \|^2_{Q} +
\| e \|^2_{\Sigma_T} + \tfrac{\delta_h}{2} \, \| \nabla_x e \|^2_{\Sigma_T} 
\leq \overline{\rm M}^{\rm I\!I}_{s, h}
(v, \flux; \beta, \alpha) \nonumber \\[5pt]
& := \epsilon \, \delta_h \, \| \R^{\rm I}_{\rm d} \|^2_{\Sigma_T} 
+ \zeta \, \Big( (1 + \alpha) \, 
				\overline{\rm M}^{\rm I}(v, \flux; \beta)  
			          + (1 + \tfrac{1}{\alpha})  
			          + \delta_h^2 \, \| \partial_t \R^{\rm I}_{\rm d} \|^2_Q
+ \delta_h \, \| \R^{\rm I}_{\rm eq} \|_{Q}^2 \Big).
\label{eq:estimate-2}
\end{alignat}
In both (i) and (ii), $\R^{\rm I}_{\rm d}$ and $\R^{\rm I}_{\rm eq}$ are defined in \eqref{eq:residuals}, 
$\CFriedrichs$ is the Friedrichs constant in \eqref{eq:friedrichs}, $\delta_h$ is the discretisation parameter, 
$\gamma, \zeta \in \big[\tfrac{1}{2}, +\infty)$, 
$\epsilon \in [1, +\infty)$, and $\beta, \alpha >0$. 
\end{corollary}
}

\section{Numerical realisation}
\label{sec:algorithms}
In this section, we discuss the IgA {  discretisation of the variational formulation presented 
above as well as the estimates that control the reconstructed approximations quality. We also suggest 
efficient algorithms for the reconstruction of a posteriori error bounds. The numerical examples presented in 
Section \ref{sec:numerical-examples} demonstrate computational properties of the majorants that follow 
from \cite{LMR:Repin:2002}, of the error identity $\EI$, and of the error bounds exposed in Section 
\ref{sec:general-error-estimate}.}

{
\subsection{Computation of the majorants in the IgA framework}
\label{subsec:discretisation}
}
We consider the approximations
$$u_h \in V_{0h} := V_h \cap V^{\Delta_x, 1}_{0, \underline{0}},$$ 
{  where $V_{h}$ is defined in \eqref{eq:vh-v0h}} and we consider 
NURBS of degree $p = 2$. Due to the restriction on knots-multiplicity of 
$\hat{\mathcal{S}}^{p}_h$ in the framework of one-patch domains, $u_h \in C^{p-1}$ is 
automatically provided. 
It is important to note that the scope of this paper is limited to a single-patch domain since it is important 
to first fully analyse {the} behaviour of the error-estimation tool in a simplified setting. The extension of this 
simpler setting to a widely used in practice multi-patch case, in which the physical domain is decomposed 
into several simple patches, will be a focus of the subsequent paper.

{Then approximation has the form}
$$u_h(x, t) = u_h(x_1, . . . , x_{d+1}) 
:= \sum_{i \in \mathcal{I}} \underline{ \rm u}_{h,i} \, \phi_{h, i} (x_1, . . . , x_{d+1}) ,
$$
where $\underline{\rm u}_h 
:= \big[ \underline{ \rm u}_{h,i}\big]_{i \in \mathcal{I}} \in {\mathds{R}}^{|\mathcal{I}|}$
{  is the vector of degrees of freedom} 
{ 
(also called  control points in the IgA community) 
defined }
by the linear system 
\begin{equation}
{\rm K}_h \, \underline{ \rm  u}_h = {\rm f}_h, \quad 
{\rm K}_h := \big[a_{s, h}(\phi_{h,i},\phi_{h,j}) \big]_{i, j \in \mathcal{I}},
\quad
{\rm f}_h :=\big[l_{s, h}(\phi_{h,i}) \big]_{i \in \mathcal{I}}.
\label{eq:system-uh}
\end{equation}
{ 
In the numerical tests presented in Section \ref{sec:numerical-examples}, we analyse the approximation 
properties of $u_h$ by looking at the convergence of the error $e = u - u_h$ measured in terms of the 
following three norms earlier defined in \eqref{eq:e-energy-norm} (with $\nu= 1$), 
\eqref{eq:strong-norm-error-identity}, and \eqref{eq:error-non-moving}, i.e., 
}
\begin{equation}
\begin{alignedat}{2}
{|\!|\!|  e |\!|\!|^2_{(1, 1)}} =  {|\!|\!|  e |\!|\!|^2} 
& := \| \nabla_x e \|^2_Q + \| e \|^2_{\Sigma_T}, \\ 
|\!|\!|  e |\!|\!|^2_{\mathcal{L}} 
& := \| \Delta_x e \|^2_Q + \| \partial_t e \|^2_Q +  \| \nabla_x e \|^2_{\Sigma_T}, \quad \mbox{and} \\
|\!|\!|  e |\!|\!|^2_{s, h}             
& := \| \nabla_x e \|^2_Q + \delta_h \, \| \partial_t  e \|^2_Q + \| e \|^2_{\Sigma_T} +  \delta_h \, \| \nabla_x e \|^2_{\Sigma_T}. 
\end{alignedat}
\label{eq:errors}
\end{equation}
The majorant for $|\!|\!|  e |\!|\!|^2$ (defined in \eqref{eq:majorant-heat-eq} with $\nu= 1$) has 
the form  
\begin{equation}
\overline{\rm M}^{\rm I}(u_h, \flux_h) 
:= (1 + \beta) \, \| \flux_h - \nabla_x u_h \|^2_Q
    + (1 + \tfrac{1}{\beta}) \, \CFriedrichs^2 \, \| \dvrg_x \flux_h + f  - \partial_t u_h\|^2_Q
= (1 + \beta) \,  \overline{\mathrm m}^{{\rm I}}_{\mathrm{d}} 
   + (1 + \tfrac{1}{\beta}) \, \CFriedrichs^2 \,  \overline{\mathrm m}^{{\rm I}}_{\mathrm{eq}}, 
\label{eq:majorant-elliptic} 
\end{equation}
where $\beta > 0$ and ${\flux}_h \in Y_h \subset H^{\dvrg_x, 0}(Q)$. The space 
$Y_h \equiv \mathcal{S}^{q}_h 
:= \big\{ {{\boldsymbol \psi}}_{h,i} := \oplus^{d+1} \hat{\mathcal{S}}^{q}_h \circ \Phi^{-1} \big\}$ 
is generated by the push-forward of $\oplus^{d+1} \hat{\mathcal{S}}^{q}_h$, where 
$\hat{\mathcal{S}}^{q}_h$ is the space of NURBS of degree $q$ {  approximating }each of $d+1$ 
components of $\flux_h = \big(y_h^{(1)}, \ldots ,y_h^{(d+1)}\big)^{\rm T}$. 
The best estimate follows from the minimisation of $\overline{\rm M}^{\rm I}(u_h, \flux_h)$ w.r.t.
$$\flux_h(x, t) = \flux_h(x_1, . . . , x_{d+1}) 
= \sum_{i \in \mathcal{I} \times (d+1)} 
\underline{ \rm \bf y}_{h,i} \, {{\boldsymbol \psi}}_{h,i}(x_1, . . . , x_{d+1}).$$ 
Here, ${{\boldsymbol \psi}}_{h,i}$ are the basis functions generating the space $Y_h$, and 
$\underline{ \rm \bf y}_{h} 
:= \big[ \underline{ \rm \bf y}_{h,i}\big]_{i \in \mathcal{I}} \in {\mathds{R}}^{(d+1)|\mathcal{I}|}$ 
is defined by the linear system
\begin{equation}
\left( {\CFriedrichs^2} \, {\rm Div}_h + {\beta} \, {\rm M}_h \right)\, \underline{ \rm \bf y}_{h} 
= - {\CFriedrichs^2} \, {\rm z}_h + {\beta} \, {\rm g}_h, 
\label{eq:system-fluxh}
\end{equation}
%
where 
\begin{equation}
\begin{array}{r@{$\;$}l l} 
	{ {\rm Div}_h} & 
	:= 
	\big[ (\dvrg_x {{\boldsymbol \psi}}_{h,i}, \dvrg_x {{\boldsymbol \psi}}_{h,j})_Q
	\big]_{i, j=1}^{(d+1)|\mathcal{I}|} 
, \qquad \,
	{\rm z}_{h} := 
	\big[\big(f - v_t,  \dvrg_x {{\boldsymbol \psi}}_{h,j} \big)_Q  
	\big]_{j=1}^{(d+1)|\mathcal{I}|}  
	,  
	\\[5pt]
	{{\rm M}_h} & := 
	\big[ ({{\boldsymbol \psi}}_{h,i}, {{\boldsymbol \psi}}_{h,j})_Q
	\big]_{i, j=1}^{(d+1)|\mathcal{I}|}
	, 
	\qquad \qquad \qquad
	{\rm g}_h := 
	\big[ \big(\nabla_x v, {{\boldsymbol \psi}}_{h,j}\big)_Q 
	\big]_{j=1}^{(d+1)|\mathcal{I}|}. 
\end{array}
\label{eq:system-matrix-components-definitions}
\end{equation}
%

Next, {  we consider {  a }
discretisation of} the second form of the majorant 
$\overline{\rm M}^{\rm I\!I}(u_h, \flux_h, \eta_h)$. For simplicity of exposition, we assume that 
the initial condition on $\Sigma_0$ is satisfied exactly, and
{  parameters $\delta^{\rm I\!I}$ and $\gamma$ are set to $1$.}
In order to make the reconstruction of $\eta_h$ transparent and overcome minimisation of 
$\overline{\rm M}^{\rm I\!I}(u_h, \flux_h, \eta_h)$ w.r.t. $\eta_h$, we represent $\eta_h$ as  
$\eta_h = w_h - u_h$. Here, $u_h$ is the approximation at hand obtained by solving 
\eqref{eq:system-uh} and $w_h$ is the solution to the same variational problem 
\eqref{eq:stabilized-bilinear-form-short} using wider approximation space 
$${  
{W}_{0h} := W_h \cap H^1_{0}(Q),  \quad 
\mbox{with} \quad 
W_h \equiv {\mathcal{S}}^{r}_h 
:=\big\{ \chi_{h,i} := \hat{\mathcal{S}}^{r}_h \circ {\upchi}^{-1} \big\},}$$ 
where ${\mathcal{S}}^{r}_h$ is the space of NURBS of degree $r$. 
{  As a result, the function $w_h$ can be represented 
in the form}
$$w_h(x, t) = w_h(x_1, . . . , x_{d+1}) 
:= \sum_{i \in \mathcal{I}} \underline{ \rm w}_{h,i} \, \chi_{h,i}.
$$
Here, 
$\underline{\rm w}_h 
:= \big[ \underline{ \rm w}_{h,i}\big]_{i \in \mathcal{I}} \in {\mathds{R}}^{|\mathcal{I}|}$
is the vector of control points of $w_h$ defined by the linear system
\begin{equation}
{\rm K}^{(r)}_h \, \underline{ \rm  w}_h = {\rm f}^{(r)}_h, 
\label{eq:system-wh}
\end{equation}
where 
${\rm K}^{(r)}_h := \big[a_{s, h}(\chi_{h, i}, \chi_{h, j}) \big]_{i, j \in \mathcal{I}}$,
${\rm f}^{(r)}_h :=\big[l_{s, h}(\chi_{h, i}) \big]_{i \in \mathcal{I}}$,
and $r$ indicates the degree of splines used to construct the basis $\chi_{h,i}$.
Taking the new representation of $\eta_h$ into account, \eqref{eq:majorant-II} can be reformulated as 
follows:
\begin{alignat}{2}
\| \nabla_x e \|^2_Q 
\leq \overline{\rm M}^{\rm I\!I} (u_h, w_h)
& = \|  w_h - u_h \|^2_{\Sigma_T} + 2 \, \mathcal{F}(u_h, w_h - u_h) 
+ (1 + {\beta} ) \, \big\| {\R}^{\rm I\!I} _{\rm d}\big\|^2_Q
+ \CFriedrichs^2 \, (1 + \tfrac{1}{{\beta}} ) \, \big\| {\R}^{\rm I\!I}_{\rm eq}\big\|^2_Q, 
\label{eq:majorant-II-w}
\end{alignat}
where 
$${\R}^{\rm I\!I} _{\rm d}(u_h, \flux_h, w_h) 
:= \flux_h + \nabla_x w_h - 2 \, \nabla_x u_h \quad \mbox{and} \quad
{\R}^{\rm I\!I}_{\rm eq}(\flux_h, w_h) 
:= \dvrg_x \flux_h + f - \partial_t w_h.$$
%
Since $\partial_t w_h$ is approximated by a richer space, the term
$\big\| {\R}^{\rm I\!I}_{\rm eq}(\flux_h, w_h) \big\|^2_Q$ is expected to be smaller than 
$\| {\R}_{\rm eq}(\flux_h, u_h) \|^2_Q$. Therefore, the value of the error bound 
$\overline{\rm M}^{\rm I\!I}$ must be improved. The optimal parameter ${\beta}$ is 
calculated by ${\beta} 
:= \CFriedrichs \| {\R}^{\rm I\!I}_{\rm eq}\|_Q / \| {\R}^{\rm I\!I}_{\rm d}\|_Q$.

In \cite{LMR:Repin:2002}, it was shown that if $w_h = u$ and $\flux_h = \nabla_x u$,  inequality
\eqref{eq:majorant-II-w} can be reformulated as follows:
\begin{alignat*}{2}
\| \nabla_x e \|^2_Q \leq \overline{\rm M}^{\rm I\!I} (u_h, u)
& := \| u - u_h \|^2_{\Sigma_T} + 2 \, \mathcal{F}(u_h, u - u_h) 
+ 4 \, (1 + {\beta} ) \, \| \nabla_x (u - u_h) \|^2_Q. 
\end{alignat*}
Moreover, after rearranging the terms of 
\begin{alignat*}{2}
\mathcal{F}(u_h, u - u_h) 
=  \Big( (\nabla_x u, \nabla_x (u - u_h )) + (\partial_t u - f, u - u_h) \Big) 
+ (\nabla_x (u_h - u), \nabla_x (u - u_h )) + (\partial_t (u_h - u), u - u_h),
\end{alignat*}
it is easy to see that the first scalar product on the RHS of $\mathcal{F}(u_h, u - u_h)$ vanishes.
As a result, we obtain
\begin{alignat}{2}
\| \nabla_x e \|^2_Q 
\leq \overline{\rm M}^{\rm I\!I} (u_h, u) & := \| u - u_h \|^2_{\Sigma_T} 
+ (4 \, (1 + {\beta}) - 2) \, \| \nabla_x (u - u_h) \|^2_Q
- 2 (\partial_t (u - u_h), u - u_h) \nonumber\\
& \; = (4 \, (1 + {\beta}) - 2) \, \| \nabla_x (u - u_h) \|^2_Q. 
\label{eq:equivalence-2}
\end{alignat}
Thus, we have the following double inequality
$$\| \nabla_x e \|^2_Q 
\leq \overline{\rm M}^{\rm I\!I} (u_h, u) 
\leq C_{\overline{\rm M}^{\rm I\!I}\rm gap} \| \nabla_x e \|^2_Q, \quad 
C_{\overline{\rm M}^{\rm I\!I}\rm gap} := (4 \, (1 + {\beta}) - 2),$$
and therefore ${C^{-1}_{\overline{\rm M}^{\rm I\!I}\rm gap}} \, \overline{\rm M}^{\rm I\!I} (u_h, u)$ 
can be used for more efficient error indication.

\subsection{Algorithms}
\label{subsec:algorithms}

Next, we concentrate on the algorithms providing an adaptive procedure based on the a 
posteriori error estimates presented above. A reliable $u_h$-approximation procedure is summarised 
in Algorithm \ref{alg:reliable-uh-reconstruction}. We assume that $f$, $u_0$, and $Q$ in 
\eqref{eq:equation} are given. As an input to Algorithm \ref{alg:reliable-uh-reconstruction}, the initial 
(or obtained on a previous refinement step) mesh $\mathcal{K}_h$ discretising the space-time cylinder $Q$ is provided. 
As an output, Algorithm \ref{alg:reliable-uh-reconstruction} returns a refined version of the mesh denoted by $\mathcal{K}_{h_{\rm ref}}$. 
Overall, the algorithm is structured according to the classic block-chain 
$${\rm APPROXIMATE} 
\rightarrow {\rm ESTIMATE} 
\rightarrow {\rm MARK} 
\rightarrow {\rm REFINE}.$$

The {\rm APPROXIMATE} step involves assembling of the system of the IgA solution $u_h$,
i.e., the matrix ${{\rm {K}}_h}$ and RHS ${{\rm {f}}_h}$ in \eqref{eq:system-uh}, and 
solving it with sparse direct ${\rm LU}$ factorisations (like Eigen SparseLU \cite{eigenweb} 
that is used in our numerical example). 
{  Such a choice of a solver 
is made in order to provide a fair comparison of time spent on solving \eqref{eq:system-uh}, 
\eqref{eq:system-fluxh}, and \eqref{eq:system-wh}.}
{ 
On coarser grids, sparse direct solvers are quite efficient. 
However, on finer grids, iterative solvers like multigrid become more and more efficient 
in a nested iteration setting, where one can use the interpolated coarse grid 
solution as an initial guess on the next, adaptively refined grid,
see, e.g., \cite{LMR:Bakhvalov:1966,LMR:Brandt:1973,LMR:Hackbusch:1985}.
}
The time spent on assembling and solving sub-procedures 
is tracked and saved in the vectors ${t_{\rm as}(u_h)}$ and ${t_{\rm sol}(u_h)}$, respectively. This notation 
is used in the upcoming examples to analyse the efficiency of Algorithm \ref{alg:reliable-uh-reconstruction} 
and to compare the computational costs for its subroutines. After the {\rm APPROXIMATE} step, the error 
contained in $u_h$ is evaluated in terms of several norms defined in \eqref{eq:errors}, i.e., 
$|\!|\!|  e |\!|\!|$, $|\!|\!|  e |\!|\!|_{s, h}$, and $|\!|\!|  e |\!|\!|_{\mathcal{L}}$. 
To measure the time for element-wise (e/w) assembling of the latter quantities, 
we use $t_{\rm e/w}(|\!|\!|  e |\!|\!|)$, 
$t_{\rm e/w}(|\!|\!|  e |\!|\!|_{s, h})$, and $t_{\rm e/w}(|\!|\!|  e |\!|\!|_{\mathcal{L}})$, respectively.

The next {\rm ESTIMATE} step focuses on the reconstruction of the global estimates 
$\overline{\rm M}^{\rm I} (u_h, \flux_h)$, $\overline{\rm M}^{\rm I\!I} (u_h, w_h)$, and 
$\overline{\rm M}_{s, h} (u_h, \flux_h)$, as well as the error identity $\EI$. The time spent on each of 
the error estimators is measured in the same way, for instance, ${t_{\rm as}(\flux_h)}$, 
${t_{\rm sol}(\flux_h)}$, and $t_{\rm e/w}(\overline{\rm M}^{\rm I})$ correspond to the times required 
to assemble system \eqref{eq:system-fluxh}, solve it, and evaluate e/w contributions 
of $\overline{\rm M}^{\rm I} (u_h, \flux_h)$.
Analogously, since $\overline{\rm M}^{\rm I\!I} (u_h, \flux_h)$ depends on $w_h$, we store 
in ${t_{\rm as}(w_h)}$ the time corresponding to the assembling of system \eqref{eq:system-wh} 
and in ${t_{\rm sol}(\flux_h)}$ the time spent on \eqref{eq:system-wh}. 
Element-wise evaluation costs are tracked in $t_{\rm e/w}(\overline{\rm M}^{\rm I\!I} (u_h, \flux_h))$. 
The reconstruction of $\overline{\rm M}_{s, h} (u_h, \flux_h)$ as well as $\EI$ narrows down to their 
e/w assembly since they do not have to be optimised and can be directly computed. 
Therefore, the time-expenses are saved in 
$t_{\rm e/w}(\overline{\rm M}_{s, h} (u_h, \flux_h))$ as well as $t_{\rm e/w}(\EI)$.
%
A detailed description of the majorant $\overline{\rm M}^{\rm I} (u_h, \flux_h)$ calculation procedure 
is presented in Algorithm \ref{alg:estimate-step-maj}, whereas the steps of
$\overline{\rm M}^{\rm I\!I} (w_h)$-reconstruction are described in Algorithm 
\ref{alg:estimate-step-maj-II}.

\begin{algorithm}[!t]
\begin{algorithmic} 
\small
\STATE {{\bf Input:}}
$\mathcal{K}_h$ \COMMENT{discretisation of $Q$} \\[4pt]
\qquad \quad \, ${\rm span} \,\big\{ \phi_{h,i}(x_1, . . . , x_{d+1}) \big\}$, $i = 1, ..., |\mathcal{I}|$ 
\COMMENT{$V_{h}$-basis} \\[4pt]
\STATE {{\bf APPROXIMATE}}: \\[4pt]
\begin{itemize}
\item ASSEMBLE the matrix ${{\rm {K}}_h}$ and RHS ${{\rm {f}}_h}$ \hfill :${t_{\rm as}(u_h)}$ \\[4pt]
\item SOLVE ${\rm K}_h \, \underline{ \rm  u}_h = {\rm f}_h$ \hfill :${t_{\rm sol}(u_h)}$ \\[4pt]
\item Reconstruct 
$u_h(x, t) = u_h(x_1, . . . , x_{d+1}) := \sum_{i \in \mathcal{I}} \underline{ \rm u}_{h,i} \, \phi_{h,i}$ \\[4pt]
%
\end{itemize}
%
{  Compute} the error $e = u - u_h$ measured in terms of $|\!|\!|  e |\!|\!|$, $|\!|\!|  e |\!|\!|_{s, h}$, and $|\!|\!|  e |\!|\!|_{\mathcal{L}}$
\hfill :${t_{\rm e/w}(|\!|\!|  e |\!|\!|) + t_{\rm e/w}(|\!|\!|  e |\!|\!|_{s, h}) + t_{\rm e/w}(|\!|\!|  e |\!|\!|_{\mathcal{L}})}$ \\[4pt]
\STATE {{\bf {ESTIMATE}}}: 
\begin{itemize}
\item {  compute} $\overline{\rm M}^{\rm I} (u_h, \boldsymbol{{y}}_h)$ 
\hfill :${t_{\rm as}(\boldsymbol{{y}}_h) + t_{\rm sol}(\boldsymbol{{y}}_h) + t_{\rm e/w}(\overline{\rm M}^{\rm I} )}$\\[4pt]
\item {  compute} $\overline{\rm M}^{\rm I\!I} (u_h, w_h)$ 
\hfill :${t_{\rm as}(w_h) + t_{\rm sol}(w_h) + t_{\rm e/w}(\overline{\rm M}^{\rm I\!I}) }$ \\[4pt]
\item {  compute} $\overline{\rm M}^{\rm I}_{s, h} (u_h, \boldsymbol{{y}}_h)$ 
\hfill :$t_{\rm e/w}(\overline{\rm M}^{\rm I}_{s, h})$ \\[4pt]
\item {  compute} $\EI (u_h)$ 
\hfill :${t_{\rm e/w}(\EI)}$ \\[4pt]
\end{itemize}
\STATE{{\bf MARK}}: Using the marking criteria $\mathds{M}_*(\sigma)$, select elements $K$ of the mesh 
$\mathcal{K}_h$ that must be refined \\[4pt]
%
\STATE {{\bf REFINE}}: Execute the refinement strategy:  
$\mathcal{K}_{h_{\rm ref}} = \mathcal{R}(\mathcal{K}_h)$ \\[4pt]
\STATE {{\bf Output:}} $\mathcal{K}_{h_{\rm ref}}$ \COMMENT{refined discretisation of $\Omega$}
\end{algorithmic}
\caption{\small Reliable reconstruction of $u_h$ (a single refinement step)}
\label{alg:reliable-uh-reconstruction}
\end{algorithm}

In the third chain-block {\rm MARK}, we use a marking criterion denoted by ${\mathds{M}}_{*}({\sigma})$. 
It provides an algorithm for defining the threshold 
$\mathfrak{S}_*$ for selecting those  $K \in \mathcal{K}_h$ 
for further refinement that satisfies the criterion
$${\mdIK \geq {\mathfrak{S}}_* ({\mathds{M}}_{*}({\sigma})), \quad K \in \mathcal{K}_h}.$$
Having reconstructed $\EI (u_h)$ in addition to $\mdInosq{} (u_h, \boldsymbol{{y}}_h)$, which is defined 
by one term of $\overline{\rm M}^{\rm I} (u_h, \boldsymbol{{y}}_h)$, we have a variety of different 
error indicators to base the mesh refinement strategy on.
In the open source C++ library G+Smo \cite{gismoweb} used for carrying out the numerical 
examples presented further, 
several marking strategies are considered. In particular, the marking based on `absolute 
threshold' is denoted as {\rm GARU} (an abbreviation for `greatest appearing residual utilisation'), the 
{`relative threshold'} is denoted as {\rm PUCA} (which stands for `percent-utilising cutoff ascertainment'), 
and the most widely used bulk marking (also known as the D\"orfler marking \cite{Dorfler1996}) is denoted 
by {\rm BULK}. In further examples, we mainly use the latter marking criterion.
In the case of uniform refinement, all elements of $\mathcal{K}_h$ are marked for 
refinement (i.e, $\sigma =0$). If the numerical IgA scheme is implemented correctly, the error is supposed 
to decrease at least as $O(h^p)$, which is verified throughout the numerical tests in Section 
\ref{sec:numerical-examples}.

Finally, on the last {\rm REFINE} step, we apply the refinement algorithm $\mathcal{R}$ to those 
elements that have been selected on the {\rm MARK} level. Since the THB-splines are based on 
subdomains of different hierarchical levels, the procedure $\mathcal{R}$ increases  the level of 
subdomains by applying the dyadic cell refinement. 

In the following, we concentrate on the structure of Algorithm \ref{alg:estimate-step-maj}, which clarifies the 
{\rm ESTIMATE} step of Algorithm \ref{alg:reliable-uh-reconstruction} in the context of functional type 
error estimates. On the Input step, the algorithm receives the approximate solution $u_h$ 
reconstructed by the IgA scheme. Moreover, since the majorant is minimised w.r.t. the vector-valued 
variable $\flux_h \in Y_{h}$, the collection of basis functions generating the space 
$Y_{h} := {\rm span} \,\big\{ {\boldsymbol \psi}_{h, i} \big\}$, $i = 1, ..., (d+1 ) |\mathcal{I}|$ 
is provided. The last input parameter $N^{\rm it}_{\rm maj}$ defines the number of the optimisation 
loops executed to obtain a good enough minimiser of $\overline{\rm M}^{\rm I}$. According to the tests 
performed in \cite{RepinSauterSmolianski2003,LMR:MatculevichRepin:2014,LMR:HolmMatculevich:2017}, 
one or two iterations are usually sufficient to achieve the reasonable accuracy of error majorant. 
Another criterion to exit the cycle earlier and, 
therefore, minimise the computational costs of the error control, can be the condition that the ratio 
$(1 + \tfrac{1}{\beta}) \, \CFriedrichs^2 \, \mfI / (1 + \beta) \, \mdI$ is small enough.
In this case, the efficiency index is automatically close to one. When the calculation of 
$\overline{\rm M}^{\rm I}$ is followed up by the reconstruction of $\overline{\rm M}^{\rm I\!I}$, 
we consider only $N^{\rm it}_{\rm maj} = 1$ iteration. In addition to $\overline{\rm M}^{\rm I}$ and 
$\overline{\rm M}^{\rm I\!I}$, 
we evaluate the majorant $\overline{\rm M}^{\rm I}_{s, h}$ specifically derived in Theorem 
\ref{th:theorem-majorant-general-1} for {the} 
stabilised scheme \eqref{eq:stabilized-bilinear-form-short} 
and the control of the error $|\!|\!|  e |\!|\!|_{s, h}$. 

We emphasise that both matrices ${{\rm {Div}}_h}$ and ${{\rm {M}}_h}$, as well as vectors 
${\rm z}_h$ and ${\rm g}_h$, are assembled only once. The loop iterates $N^{\rm it}_{\rm maj}$ times 
such that each time the optimal $\underline{ \rm \bf y}^{(n)}_{h}$ and $\beta^{{\rm I}, (n)}$ are 
reconstructed. In our implementation, the optimality system for the flux (see \eqref{eq:system-fluxh}) is 
solved by the sparse direct ${\rm LDL^{\rm T}}$ Cholesky factorisations.
The time spent on {\rm ASSEMBLE} and {\rm SOLVE} steps w.r.t. the system 
\eqref{eq:system-fluxh} is measured by ${t_{\rm as}(\flux_h)}$ and ${t_{\rm sol}(\flux_h)}$ 
respectively and compared to ${t_{\rm as}}(u_h)$ and ${t_{\rm sol}}(u_h)$ in forthcoming 
numerical examples. It is crucial to note that the matrices ${{\rm Div}_h}$ and ${{\rm M}_h}$ have block 
structure (of $(d+1) \times (d+1)$ blocks) due to the properties of the approximation spaces $V_h$ and 
$Y_h$.
Moreover, since ${{\rm Div}_h}$, ${{\rm M}_h}$, ${\rm r}_h$ and ${\rm g}_h$ are generated 
by the scalar product of the derivatives or divergence w.r.t. spatial coordinates only, $(d+1, d+1)$-th block  
of ${{\rm Div}_h}$ is zero as well as the $(d+1)$-th block of the RHS of \eqref{eq:system-fluxh}, i.e., 
\begin{alignat*}{2}
\left( {\CFriedrichs^2} \, 
\begin{bmatrix} 
{\rm Div}_h^{(d)} & 0 \\ 
0 & 0 
\end{bmatrix}
+ 
{\beta} \, 
\begin{bmatrix} 
{\rm M}_h^{(d)} & 0 \\ 
0 & {M}_h^{(1)} 
\end{bmatrix}
\right)
\,
\cdot
\,
\begin{bmatrix} 
{\flux}_h^{(d)}\\ 
{\flux}_h^{(1)} 
\end{bmatrix}
& =
- {\CFriedrichs^2} \, 
\begin{bmatrix} 
{\rm z}_h^{(d)}\\ 
0 
\end{bmatrix}
+
{\beta} \, 
\begin{bmatrix} 
{\rm g}_h^{(d)}\\ 
0 
\end{bmatrix},
\end{alignat*}
{ 
where $(1)$-block corresponds to the time variable. This resolves into the vector $Y_h$ with zero $(d+1)$-th 
block, which in turn allows us to solve the system composed only of spatial blocks}. Besides the computational 
costs related to the assembling and solving of \eqref{eq:system-uh} and \eqref{eq:system-fluxh}, we measure 
the time spent on the e/w evaluation of all the majorants.

Analogously to the selection of $q$ for the space $Y_h$, we let $r = p + l$, $l \in \mathds{N}^+$. 
At the same time, we use a coarser mesh $\mathcal{K}_{Lh}$,  $L \in \mathds{N}^+$ {  
in order to recover $w_h$}. For the reader convenience, we collect the notation related to the spaces 
parametrisation in Table \ref{tab:table-of-notation}. 
The sequence of steps of the $w_h$-approximation, as well as $\overline{\rm M}^{\rm I\!I}$-reconstruction corresponding to it, are presented in Algorithm \ref{alg:estimate-step-maj-II}. 
Its structure is similar to the structure of Algorithm 
\ref{alg:estimate-step-maj} with the exception that the free variable of 
$\overline{\rm M}^{\rm I\!I}(v_h, \flux_h, w_h)$ is a scalar
function and we solve system \eqref{eq:system-wh} to reconstruct the degrees of freedom (d.o.f.) of 
$w_h$ only once.

{  Evaluation of the error identity does not require any optimisation techniques. 
Therefore, it can be computed straightforwardly by using}
$$\EI^2 (u_h) := \| \nabla_x (u_0 - u_h) \|^2_{\Sigma_0} + \|  \Delta_x u_h + f - \partial_t u_h\|^2_Q$$
without any overhead in time performance. Time spent on the element-wise assembly of $\EI$ is tracked 
in $t_{\rm e/w}(\EI)$.

\begin{algorithm}[!t]
\small
\caption{\small \quad \small {\bf ESTIMATE} step (majorant $\overline{\rm M}^{\rm I}$ minimisation)}
\begin{algorithmic} 
\STATE  {{\bf Input:}} 
 $u_h$ \COMMENT{approximation}\\[4pt] 
\qquad \quad \, $\mathcal{K}_h$ \COMMENT{discretisation of $\Omega$} \\[4pt]
\qquad \quad \, ${\rm span} \,\big\{ {\boldsymbol \psi}_{h, i} \big\}$, $i = 1, ..., (d+1) |\mathcal{I}|$ \COMMENT{$Y_{h}$-basis} \\[4pt]
\qquad \quad \, $N^{\rm it}_{\rm maj}$ \COMMENT{number of optimisation iterations} \\[4pt]
\STATE  {ASSEMBLE}  
${{\rm {Div}}_h}, {{\rm {M}}_h} \in \mathds{R}^{(d+1) |\mathcal{I}| \times (d+1) |\mathcal{I}|}$ 
and ${\rm z}_h$, ${\rm g}_h \in {\mathds R}^{(d+1) |\mathcal{I}|}$
\hfill :${t_{\rm as}(\boldsymbol{{y}}_h)}$\\[4pt]
\STATE  {Set} $\beta^{(0)} = 1$ \\[4pt]
\FOR {$n = 1$ {\bf to} $N^{\rm it}_{\rm maj}$} 
\STATE {SOLVE} \quad 
$\big( \rfrac{\CFriedrichs^2}{\beta^{(n-1)} \, } {{\rm {Div}}_h} + {{\rm {M}}_h} \big)\, \underline{ \rm \bf y}^{(n)}_{h} = 
	-\rfrac{\CFriedrichs^2}{\beta^{(n-1)} \,} {\rm z}_h + {\rm g}_h$ 
\hfill :${t_{\rm sol}(\boldsymbol{{y}}_h)}$ \\[4pt]
\STATE  {Reconstruct}
$\flux^{(n)}_h 
:= \sum_{i \in \mathcal{I} \times (d+1)} \underline{ \rm \bf y}^{(n)}_{h,i} \, {{\boldsymbol \psi}}_{h,i}$ \\[4pt]
\STATE  {{  Compute} } 
$\overline{\mathrm m}^{{\rm I}, {(n)}}_{\mathrm{eq}} := \| \, \dvrg_x \flux^{(n)}_h + f  - \partial_t u_h \, \|_{\Omega}$ and 
$\overline{\mathrm m}^{{\rm I}, {(n)}}_{\mathrm{d}} := \| \flux^{(n)}_h - \nabla_x u_h \,\|_{\Omega}$ 
\hfill :${\bf t_{\rm e/w}(\overline{\rm M}^{{\rm I}})}$\\[4pt]
\STATE  {{  Compute} }
$
\beta^{(n)} = 
\tfrac{\CFriedrichs \,  \overline{\mathrm m}^{{\rm I}, {(n)}}_{\mathrm{eq}} }{  \overline{\mathrm m}^{{\rm I}, {(n)}}_{\mathrm{d}} }$ \\[4pt]
\ENDFOR 
%
\STATE Assign 
$\flux_h = \flux^{(n)}_h$,
${\mdInosq} = \overline{\mathrm m}^{{\rm I}, {(n)}}_{\mathrm{d}}$, 
${\mfInosq} = \overline{\mathrm m}^{{\rm I}, {(n)}}_{\mathrm{eq}}$ \\[4pt]
\STATE {{  Compute}}
$\overline{\rm M}^{{\rm I}, 2} (u_h, \flux^{(n)}_h; \beta) 
	:= (1 + \beta) \, {\mfI} + 
   	    (1 + \tfrac{1}{\beta}) \, {\CFriedrichs^2} \, {\mdI}$ \\[4pt]
\STATE {{  Compute}}
$\alpha = \tfrac{ \mfInosq{} }{\| \dvrg_x (\flux_h - \nabla_x u_h) \|_{Q}}$ \\[4pt]
\STATE {{  Compute}}
$\overline{\rm M}^{{\rm I}, 2}_{s, h} (u_h, \flux^{(n)}_h; \beta) 
:= \overline{\rm M}^{{\rm I}, 2} 
+ \delta_h\, \big( (1 + \alpha)\, \| \dvrg_x (\flux_h - \nabla_x u_h) \|^2_{Q} + (1 + \tfrac{1}{\alpha}) \, {\mfI} \big)$ \\[2pt]
%
\STATE {{\bf Output:}} $\overline{\rm M}^{\rm I}$, $\overline{\rm M}^{\rm I}_{s, h}$ 
\COMMENT{total error majorants on $\Omega$} \\[4pt]
\qquad \quad \;\;\;\, ${\mdInosq}$ 
\COMMENT{indicator of the error distribution over $\mathcal{K}_h$}
\end{algorithmic}
\label{alg:estimate-step-maj}
\end{algorithm}

\begin{table}[htbp]
\small
\begin{center}
\begin{tabular}{r c p{12cm} }
\toprule
$p$ & $$ & the degree of the splines used for the approximation of $u_h$ \\[2pt]
$q$ & $$ & the degree of the splines used for the approximation of $\flux_h$ \\[2pt]
$r$ & $$ & the degree of the splines used for the approximation of $w_h$ \\[2pt]
$m$ & $ $ & $q - p$ \\[2pt]
$l$ & $ $ & $r - p$ \\[2pt]
$S_h^{p}$ & $$ & the approximation space for the scalar-functions generated by splines \\[2pt]
$\oplus^d S_h^{q}$ & $$ & the approximation space for the $d$-dimensional vector-functions generated by splines \\[2pt]
$S_h^{q} \oplus S_h^{q}$ & $$ & the approximation space for the two-dimensional vector-functions generated by splines \\[2pt]
$M$  & $ $ & the coarsening ratio of the global mesh size for reconstruction of $\flux_h$ to the global mesh-size 
for {the} approximation of $u_h$ \\[2pt]
$L$  & $ $ & the coarsening ratio of the global mesh size for reconstruction of $w _h$ to the global mesh-size 
for {the} approximation of $u_h$ \\[2pt]
$\mathcal{K}_{h}$ ($\mathcal{K}^{u_h}_h$) & $ $ & the mesh used for the approximation of $u_h$ \\[2pt]
$\mathcal{K}_{Mh}$  ($\mathcal{K}^{\flux_h}_h$, $M = 1$)& $ $ & the mesh used for the approximation of $\flux_h$ \\[2pt] 
$\mathcal{K}_{Lh}$  ($\mathcal{K}^{w_h}_h$, $L = 1$)& $ $ & the mesh used for the approximation of $w_h$ \\[2pt] 
$N_{\rm ref}$ & $ $ & the number of uniform or adaptive refinement steps \\[2pt]
$N_{\rm ref, 0}$ & $ $ & the number of initial refinement steps performed before testing \\[2pt]
${\mathds{M}}_{*}({\sigma})$ & $ $ & the marking criterion $*$ with the parameter $\sigma$\\[2pt]
\bottomrule
\end{tabular}
\end{center}
\caption{Summary of some notations introduced in the text.}
\label{tab:table-of-notation}
\end{table}


\begin{algorithm}[!t]
\begin{algorithmic} 
\small
\caption{\small \quad \small {\bf ESTIMATE} step (advanced majorant $\overline{\rm M}^{\rm I\!I}$ minimisation)}
\label{alg:estimate-step-maj-II}
\STATE  {{\bf Input:}} 
$u_h$ \COMMENT{approximation}\\[2pt] 
\qquad \quad \, $\flux_h$ \COMMENT{auxiliary vector-function reconstructed by Algorithm \ref{alg:estimate-step-maj}}\\[4pt] 
\qquad \quad \, $\mathcal{K}_h$ \COMMENT{disctretisation of $\Omega$}, \\[4pt]
\qquad \quad \, ${\rm span} \,\big\{ {\boldsymbol \chi}_{h, i} \big\}$, $i = 1, ..., |\mathcal{I}|$ \COMMENT{$W_{h}$-basis}, \\[4pt]
\STATE  {\bf ASSEMBLE}  
${{\rm {K}}^{(r)}_h} \in \mathds{R}^{|\mathcal{I}| \times |\mathcal{I}|}$ 
and ${\rm f}^{(r)}_h \in {\mathds R}^{|\mathcal{I}|}$
\hfill :${t_{\rm as}({{w}}_h)}$ \\[4pt]
\item {\bf SOLVE} ${\rm K}^{(r)}_h \, \underline{ \rm  w}_h = {\rm f}^{(r)}_h$
\hfill :${t_{\rm sol}({{w}}_h)}$ \\[4pt]
\STATE {Reconstruct}
$w_h := \sum_{i \in \mathcal{I}} \underline{ \rm w}_{h,i} \, {{\boldsymbol \chi}}_{h,i}$ \\[4pt]
\STATE {{  Compute}} \hfill :${t_{\rm e/w}(\overline{\rm M}^{{\rm I\!I}})}$
\\[-28pt]
\begin{alignat*}{2}
\overline{\rm m}^{\rm I\!I}_{\rm eq}(\flux_h, w_h) & := \| \, \dvrg_x \flux_h + f - \partial_t w_h \, \|^2_{\Omega},  \\
\overline{\rm m}^{\rm I\!I} _{\rm d}(u_h, \flux_h, w_h)  & := \| \, \flux_h + \nabla_x w_h - 2 \, \nabla_x u_h \,\|^2_{\Omega}, \quad \mbox{and} \\
\mathcal{F}(u_h, w_h - u_h) & := \big(\nabla_x u_h, \nabla_x (w_h - u_h)\big) + (\partial_t u_h - f, w_h - u_h)
\end{alignat*}
\STATE {{  Compute}}
$\beta = \tfrac{ \CFriedrichs \| {\R}^{\rm I\!I}_{\rm eq}\|_Q}{\| {\R}^{\rm I\!I}_{\rm d}\|_Q}$ \\[4pt]
\STATE {{  Compute}}
%
$\overline{\rm M}^{\rm I\!I} (u_h, \flux_h, w_h)
	:= \|  w_h - u_h \|^2_{\Sigma_T} + 2 \, \mathcal{F}(u_h, w_h - u_h) 
+ (1 + {\beta} ) \, \big\| {\R}^{\rm I\!I} _{\rm d}\big\|^2_Q
+ \CFriedrichs^2 \, (1 + \tfrac{1}{{\beta}} ) \,  \big\| {\R}^{\rm I\!I}_{\rm eq}\big\|^2_Q$ \\[4pt]
\STATE {{\bf Output:}} $\overline{\rm M}^{\rm I\!I}$ \COMMENT{total error majorant on $\Omega$} \\
\qquad \quad \;\;\;\, $\overline{\rm m}^{\rm I\!I} _{\rm d}$ \COMMENT{indicator of the error distribution over $\mathcal{K}_h$}
\end{algorithmic}
\end{algorithm}


\section{Numerical examples}
\label{sec:numerical-examples}

{In the last section, we study the numerical behaviour of the error control tools discussed above on 
a series of benchmark examples. 
We start with a simple example to make the implementation of the majorants
clear to the reader, and to provide some important properties of these a posteriori error estimators.
}
The complexity of numerical tests will increase by the end of the section, where we add local drastic 
changes to the exact solutions, and consider domains with {a} more complicated shape.


\subsection{Example 1}
\label{ex:unit-domain-example-2}
\rm
As a starting point, we consider a simple example where we take the solution
%
$$u(x,t) = (1 - x) \, x^2 \,  (1 - t) \, t, \quad (x, t) \in \overline{Q} := [0, 1]^2, $$ 
and compute the RHS 
$$f(x,t) = - (1 - x)\, x^2 \,(1 - 2\,t) - (2 - 6\,x) \, (1 - t)\,t, \quad (x, t) \in Q := (0, 1)^2.$$ 
The solution $u(x, t)$ obviously satisfies homogeneous Dirichlet boundary conditions on 
$\Sigma = \partial \Omega \times (0, 1)$ 
{ 
and homogeneous initial conditions on ${\overline{\Sigma}}_0$.
}

First of all, we test {  the behaviour of a posteriori}
error estimates by executing the uniform refinement strategy. We start with the
initial mesh obtained by one global refinement ($N_{\rm ref, 0} = 1$), and we proceed further with 
{  further}
eight 
{  uniform}
refinement steps ($N_{\rm ref} = 8$). 
{  The approximation spaces considered are the  following:} 
$u_h \in S_{h}^{2}$, $\flux_h \in S_{7h}^{3} \oplus S_{7h}^{3}$, and $w_h \in S_{7h}^{3}$, 
where the coarsening parameter is 
{  given}
by $M = L = 7$. Table 
\ref{tab:unit-domain-example-2-error-majorant-v-2-y-3-uniform-ref} describes the performance of 
each error estimate (with optimal functions reconstructed according to Algorithms  
\ref{alg:estimate-step-maj} and \ref{alg:estimate-step-maj-II}). Here, the values of the error-norms
$|\!|\!| e |\!|\!|_Q$, $\|\!|\!|  e |\!|\!|_{s, h}$, and $|\!|\!|  e |\!|\!|_{\mathcal{L}}$ 
are followed by the efficiency indices of $\overline{\rm M}^{\rm I}$ ($\overline{\rm M}^{\rm I\!I}$), 
$\overline{\rm M}^{\rm I}_{s, h}$, and the identity ${\EI}$, respectively, i.e., %
\begin{alignat}{2}
\Ieff (\overline{\rm M}^{\rm I}) 
:= \tfrac{\overline{\rm M}^{\rm I}}{|\!|\!| e |\!|\!|_Q} = 
 \tfrac{\overline{\rm M}^{\rm I}}{\| \nabla_x e \|_Q}, \quad
\Ieff (\overline{\rm M}^{\rm I\!I}) 
:= \tfrac{\overline{\rm M}^{\rm I\!I}}{  C_{\overline{\rm M}^{\rm I\!I}\rm gap} \, \| \nabla_x e \|_Q}, 
\quad 
\Ieff (\overline{\rm M}^{\rm I}_{s, h}) 
:= \tfrac{\overline{\rm M}^{\rm I}_{s, h}}{|\!|\!|  e |\!|\!|_{s, h}}, \quad 
\Ieff ({\EI}) := \tfrac{{\EI}}{|\!|\!|  e |\!|\!|_{\mathcal{L}}} = 1.
\end{alignat}
Even though the definition of the last efficiency index seems trivial, we expose it in order to control 
the accuracy of the numerical integration procedures. 
From Table \ref{tab:unit-domain-example-2-error-majorant-v-2-y-3-uniform-ref}, it is obvious that for 
this rather smooth example a posteriori error estimates maintain very high efficiency since we can 
reconstruct optimal $\flux_h$ and $w_h$ with very low costs. By analysing Table 
\ref{tab:unit-domain-example-2-error-majorant-v-2-y-3-uniform-ref}, it is easy to see that 
$\overline{\rm M}^{\rm I\!I}$ improves the performance of $\overline{\rm M}^{\rm I}$ for about 
$9$--$10\%$, whereas the time for assembling and solving \eqref{eq:system-wh} is a thousand times 
smaller than the time spent on \eqref{eq:system-uh}, see the last row of Table 
\ref{tab:unit-domain-example-2-times-v-2-y-3-uniform-ref} with corresponding ratios. 
$\overline{\rm M}^{\rm I}_{s, h}$ performs similarly to $\overline{\rm M}^{\rm I}$. However, if the 
parameter $\theta$ in the space-time IgA scheme \eqref{eq:discrete-scheme} is independent of $h$, 
$\overline{\rm M}^{\rm I}$ does converge slower than $|\!|\!| e |\!|\!|_{s, h}$ for the uniform 
refinement case. As expected, the sharpest error indication is 
provided by the error identity ${\EI}$, its efficiency index stays equal to $1$ on all refinement levels. 
When it comes to the time performance of ${\EI}$, it does not require any computational overhead 
w.r.t. the element-wise evaluation of the error 
$|\!|\!|  e |\!|\!|_{\mathcal{L}}$ since it depends solely on the approximation $u_h$ at hand. However, 
we should emphasise that in order to use ${\EI}$, the solution and its approximation must satisfy higher 
regularity assumptions, i.e., $u, v \in V^{\Delta_x}_0$. Such regularity is easy to provide in problems 
similar to this example but it has to be weakened in more complicated cases. 

The time spent on assembling and solving the systems
for defining the functions minimising the error functionals
is illustrated in Table \ref{tab:unit-domain-example-2-times-v-2-y-3-uniform-ref}. The last row 
demonstrates dimensionless ratios of such time spent on the variables $u_h$, $\flux_h$, $w_h$. 
We see that the minimum time is required on the reconstruction of $w_h$. 
{ 
The time effort spent on 
$\flux_h$ also stays low due to the relatively small number of d.o.f. we keep for the flux variable. 
The last column of Table \ref{tab:unit-domain-example-2-times-v-2-y-3-uniform-ref} provides the ratio of 
the total time $t_{\rm appr.}$ spent on reconstruction of the approximation, 
which includes time for assembling and solving of system \eqref{eq:system-uh}, 
i.e., $t_{\rm as}(\flux_h) + t_{\rm sol}(\flux_h)$, 
to the time $t_{\rm er.est.}$ spent on the error estimates. 
The latter is summarised from $t_{\rm as}(\flux_h)$, $t_{\rm as}(w_h)$, 
$t_{\rm sol}(\flux_h)$, and $t_{\rm sol}(w_h)$. We can see that 
this ratio grows with the increase of iterations
as well, and reaches a quite substantial value at the last step, i.e.,
} 
$$ \tfrac{t_{\rm appr.}}{t_{\rm er.est.}} := 
\tfrac{t_{\rm sol}(u_h) + t_{\rm as}(u_h)}{
t_{\rm sol}(\flux_h) + t_{\rm as}(\flux_h) + t_{\rm sol}(w_h) + t_{\rm as}(w_h)} 
= 5616.1.$$
\begin{table}[!t]
\scriptsize
\centering
\newcolumntype{g}{>{\columncolor{gainsboro}}c} 	
\newcolumntype{k}{>{\columncolor{lightgray}}c} 	
\newcolumntype{s}{>{\columncolor{silver}}c} 
\newcolumntype{a}{>{\columncolor{ashgrey}}c}
\newcolumntype{b}{>{\columncolor{battleshipgrey}}c}
\begin{tabular}{c|cga|ck|cb|cc}
\parbox[c]{0.8cm}{\centering \# ref. } & 
\parbox[c]{1.4cm}{\centering  $\| \nabla_x e \|_Q$}   & 	  
\parbox[c]{1.0cm}{\centering $\Ieff (\overline{\rm M}^{\rm I})$ } & 
\parbox[c]{1.4cm}{\centering $\Ieff (\overline{\rm M}^{\rm I\!I})$ } & 
\parbox[c]{1.0cm}{\centering  $|\!|\!|  e |\!|\!|_{s, h}$ }   & 	  
\parbox[c]{1.4cm}{\centering $\Ieff (\overline{\rm M}^{\rm I}_{s, h})$ } & 
\parbox[c]{1.0cm}{\centering  $|\!|\!|  e |\!|\!|_{\mathcal{L}}$ }   & 	  
\parbox[c]{1.4cm}{\centering$\Ieff ({\EI})$ } & 
\parbox[c]{1.2cm}{\centering e.o.c. ($|\!|\!|  e |\!|\!|_{s, h}$)} & 
\parbox[c]{1.2cm}{\centering e.o.c. ($|\!|\!|  e |\!|\!|_{\mathcal{L}}$)} \\
\midrule
     2 &     2.5516e-03 &         1.15 &         1.06 &     2.5516e-03 &         1.25 &     7.9057e-02 &         1.00 &     3.42 &     1.71 \\
   4 &     1.5947e-04 &         1.39 &         1.20 &     1.5947e-04 &         1.70 &     1.9764e-02 &         1.00 &     2.36 &     1.18 \\
   6 &     9.9670e-06 &         1.31 &         1.15 &     9.9670e-06 &         2.36 &     4.9411e-03 &         1.00 &     2.09 &     1.05 \\
   8 &     6.2294e-07 &         1.06 &         1.03 &     6.2294e-07 &         4.06 &     1.2353e-03 &         1.00 &     2.02 &     1.01 \\
\end{tabular}
\caption{{\em Example 1}. 
Efficiency of $\overline{\rm M}^{\rm I}$, $\overline{\rm M}^{\rm I\!I}$, 
$\overline{\rm M}^{\rm I}_{s, h}$, and ${\EI}$ for 
$u_h \in S^{2}_{h}$, 
$\flux_h \in S^{3}_{7h} \oplus S^{3}_{7h}$, and 
$w_h \in S^{3}_{7h}$, w.r.t. uniform refinements.}
\label{tab:unit-domain-example-2-error-majorant-v-2-y-3-uniform-ref}
\end{table}
\begin{table}[!t]
\scriptsize
\centering
\newcolumntype{g}{>{\columncolor{gainsboro}}c} 	
\begin{tabular}{c|ccc|cgg|cgg|c}
& \multicolumn{3}{c|}{ d.o.f. } 
& \multicolumn{3}{c|}{ $t_{\rm as}$ }
& \multicolumn{3}{c|}{ $t_{\rm sol}$ } 
& $\tfrac{t_{\rm appr.}}{t_{\rm er.est.}}$ \\
\midrule
\parbox[c]{0.8cm}{\centering \# ref. } & 
\parbox[c]{0.8cm}{\centering $u_h$ } &  
\parbox[c]{0.6cm}{\centering $\flux_h$ } &  
\parbox[c]{0.6cm}{\centering $w_h$ } & 
\parbox[c]{1.4cm}{\centering $u_h$ } & 
\parbox[c]{1.4cm}{\centering $\flux_h$ } & 
\parbox[c]{1.4cm}{\centering $w_h$ } & 
\parbox[c]{1.4cm}{\centering $u_h$ } & 
\parbox[c]{1.4cm}{\centering $\flux_h$ } & 
\parbox[c]{1.4cm}{\centering $w_h$ } &
\\
\midrule
   2 &         36 &         50 &         25 &   1.74e-03 &   2.28e-03 &   1.12e-03 &         2.61e-04 &         1.59e-04 &         1.11e-04 & 0.54 \\
   4 &        324 &         50 &         25 &   2.09e-02 &   1.36e-03 &   9.95e-04 &         8.05e-03 &         8.00e-05 &         6.20e-05 & 11.59\\
   6 &       4356 &         50 &         25 &   3.30e-01 &   1.07e-03 &   7.85e-04 &         6.89e-01 &         6.30e-05 &         2.67e-04 & 466.36\\
   8 &      66564 &         98 &         49 &   3.05e+00 &   4.94e-03 &   1.59e-03 &         3.61e+01 &         1.50e-04 &         2.91e-04 & 5616.1\\
    \midrule
    &       &         &    &
    \multicolumn{3}{c|}{ $t_{\rm as} (u_h)$ \quad : \quad $t_{\rm as} (\flux_h)$ \quad : \qquad $t_{\rm as} (w_h)$ } &      
    \multicolumn{3}{c|}{\; $t_{\rm sol} (u_h)$ \, : \quad $t_{\rm sol} (\flux_h)$  \quad:  \qquad  $t_{\rm sol} (w_h)$\;} & \\
 \midrule
         &       &       &       &    1923.28 &       3.11 &       1.00 &        124016.90 &             0.52 &             1.00 & \\
\end{tabular}
\caption{{\em Example 1}. 
Assembling and solving time (in seconds) spent for the systems defining d.o.f. of 
$u_h \in S^{2}_{h}$, $\flux_h \in S^{3}_{7h} \oplus S^{3}_{7h}$, and 
$w_h \in S^{3}_{7h}$ w.r.t. uniform refinements.}
\label{tab:unit-domain-example-2-times-v-2-y-3-uniform-ref}
\end{table}
%
\begin{table}[!t]
\scriptsize
\centering
\newcolumntype{g}{>{\columncolor{gainsboro}}c} 	
\newcolumntype{k}{>{\columncolor{lightgray}}c} 	
\newcolumntype{s}{>{\columncolor{silver}}c} 
\newcolumntype{a}{>{\columncolor{ashgrey}}c}
\newcolumntype{b}{>{\columncolor{battleshipgrey}}c}
\begin{tabular}{c|cga|ck|cb|cc}
\parbox[c]{0.8cm}{\centering \# ref. } & 
\parbox[c]{1.4cm}{\centering  $\| \nabla_x e \|_Q$}   & 	  
\parbox[c]{1.0cm}{\centering $\Ieff (\overline{\rm M}^{\rm I})$ } & 
\parbox[c]{1.4cm}{\centering $\Ieff (\overline{\rm M}^{\rm I\!I})$ } & 
\parbox[c]{1.0cm}{\centering  $|\!|\!|  e |\!|\!|_{s, h}$ }   & 	  
\parbox[c]{1.4cm}{\centering $\Ieff (\overline{\rm M}^{\rm I}_{s, h})$ } & 
\parbox[c]{1.0cm}{\centering  $|\!|\!|  e |\!|\!|_{\mathcal{L}}$ }   & 	  
\parbox[c]{1.4cm}{\centering$\Ieff ({\EI})$ } & 
\parbox[c]{1.2cm}{\centering e.o.c. ($|\!|\!|  e |\!|\!|_{s, h}$)} & 
\parbox[c]{1.2cm}{\centering e.o.c. ($|\!|\!|  e |\!|\!|_{\mathcal{L}}$)} \\
\midrule
   3 &     6.3789e-04 &         1.02 &         1.00 &     6.3789e-04 &         1.02 &     3.9528e-02 &         1.00 &     2.71 &     1.36 \\
   5 &     9.2549e-05 &         1.09 &         1.00 &     9.2549e-05 &         1.09 &     1.3474e-02 &         1.00 &     1.47 &     0.78 \\
   7 &     9.3372e-06 &         1.05 &         1.00 &     9.3372e-06 &         1.05 &     4.6582e-03 &         1.00 &     3.00 &     1.32 \\
   8 &     4.3954e-06 &         1.05 &         1.00 &     4.3954e-06 &         1.05 &     2.9007e-03 &         1.00 &     1.88 &     1.18 \\
\end{tabular}
\caption{{\em Example 1}. 
Efficiency of $\overline{\rm M}^{\rm I}$, $\overline{\rm M}^{\rm I\!I}$, $\overline{\rm M}^{\rm I}_{s, h}$, and ${\EI}$ for 
$u_h \in S^{2}_{h}$, $\flux_h \in S^{3}_{7h} \oplus S^{3}_{7h}$, and $w_h \in S^{3}_{7h}$, w.r.t. adaptive refinements 
(with the marking criterion ${\mathds{M}}_{\rm BULK}(0.4)$).}
\label{tab:unit-domain-example-2-error-majorant-v-2-y-3-adaptive-ref}
\end{table}
\begin{table}[!t]
\scriptsize
\centering
\newcolumntype{g}{>{\columncolor{gainsboro}}c} 	
\begin{tabular}{c|ccc|cgg|cgg|c}
& \multicolumn{3}{c|}{ d.o.f. } 
& \multicolumn{3}{c|}{ $t_{\rm as}$ }
& \multicolumn{3}{c|}{ $t_{\rm sol}$ } 
& $\tfrac{t_{\rm appr.}}{t_{\rm er.est.}}$ \\
\midrule
\parbox[c]{0.8cm}{\centering \# ref. } & 
\parbox[c]{0.8cm}{\centering $u_h$ } &  
\parbox[c]{0.6cm}{\centering $\flux_h$ } &  
\parbox[c]{0.6cm}{\centering $w_h$ } & 
\parbox[c]{1.4cm}{\centering $u_h$ } & 
\parbox[c]{1.4cm}{\centering $\flux_h$ } & 
\parbox[c]{1.4cm}{\centering $w_h$ } & 
\parbox[c]{1.4cm}{\centering $u_h$ } & 
\parbox[c]{1.4cm}{\centering $\flux_h$ } & 
\parbox[c]{1.4cm}{\centering $w_h$ } \\
\midrule
   2 &         36 &         50 &         25 &   1.66e-02 &   1.82e-02 &   1.53e-02 &         1.92e-04 &         5.80e-05 &         1.52e-04 & 0.49 \\
   4 &        240 &         50 &         25 &   2.45e-01 &   1.71e-02 &   1.39e-02 &         6.60e-03 &         6.30e-05 &         5.60e-05 & 8.08 \\
   6 &       2451 &         50 &         25 &   2.57e+00 &   1.85e-02 &   1.10e-02 &         2.86e-01 &         6.10e-05 &         1.21e-04 & 96.22 \\
   8 &      11422 &         50 &         25 &   1.19e+01 &   1.77e-02 &   8.98e-03 &         4.41e+00 &         6.70e-05 &         1.24e-04 & 606.97\\
    \midrule
    &       &         &    &
    \multicolumn{3}{c|}{ $t_{\rm as} (u_h)$ \quad : \quad $t_{\rm as} (\flux_h)$ \quad : \qquad $t_{\rm as} (w_h)$ } &      
    \multicolumn{3}{c|}{\; $t_{\rm sol} (u_h)$ \, : \quad $t_{\rm sol} (\flux_h)$  \quad:  \qquad  $t_{\rm sol} (w_h)$\;} & 
    \\
 \midrule

         &       &       &       &    1323.22 &       1.97 &       1.00 &         35548.15 &             0.54 &             1.00 &  \\
\end{tabular}
\caption{{\em Example 1}. Assembling and solving time (in seconds) spent for the systems generating
d.o.f. of $u_h \in S^{2}_{h}$, $\flux_h \in S^{3}_{7h} \oplus S^{3}_{7h}$, and 
$w_h \in S^{3}_{7h}$ w.r.t. adaptive refinements (with the marking criterion ${\mathds{M}}_{\rm BULK}(0.4)$).}
\label{tab:unit-domain-example-2-times-v-2-y-3-adaptive-ref}
\end{table}
%

\noindent
In case of the adaptive refinement strategy, which we execute{ in} nine refinement steps, i.e., 
$N_{\rm ref} = 9$, Table \ref{tab:unit-domain-example-2-error-majorant-v-2-y-3-adaptive-ref} 
demonstrates that the advanced form of the majorant $\overline{\rm M}^{\rm I\!I}$ (column four) is 
again $19$--$20\%$ sharper than the values of $\overline{\rm M}^{\rm I}$ (column three). 
The computing time for the reconstruction of error estimates is illustrated in 
Table \ref{tab:unit-domain-example-2-times-v-2-y-3-adaptive-ref}.

\begin{figure}[!t]
	\centering
	\captionsetup[subfigure]{oneside, labelformat=empty}
	%
	\subfloat[]{
	\includegraphics[width=4.1cm, trim={0cm 0.0cm 0cm 0cm}, clip]{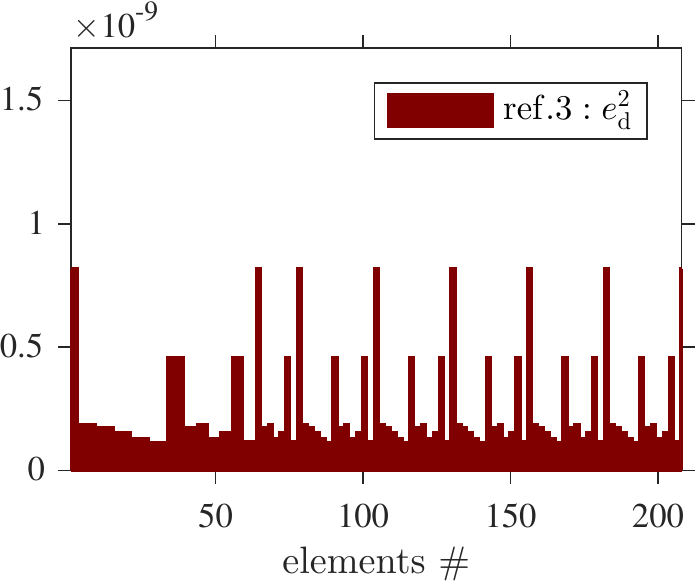}
	}
	\subfloat[]{
	\includegraphics[width=4.1cm, trim={0cm 0.0cm 0cm 0cm}, clip]{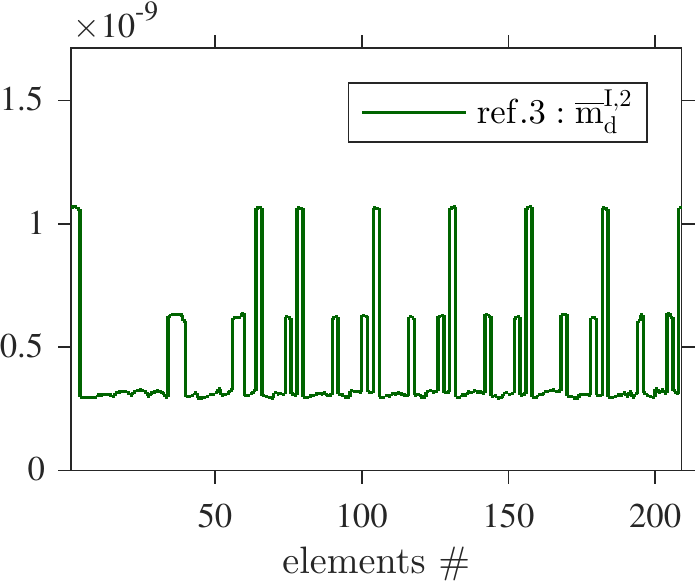}
	}
	\subfloat[]{
	\includegraphics[width=4.1cm, trim={0cm 0.0cm 0cm 0cm}, clip]{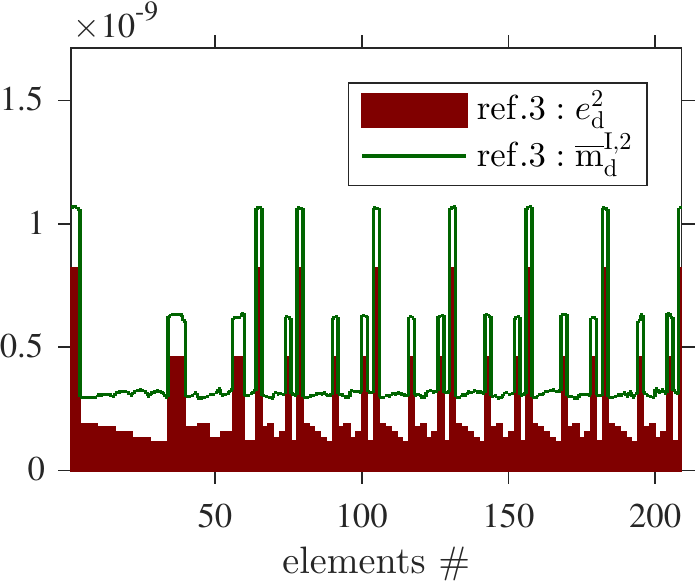}
	}
	\\[-20pt]
	\hskip 10pt
	\subfloat[]{
	\includegraphics[width=4.2cm, trim={0cm 0.0cm 0cm 0cm}, clip]{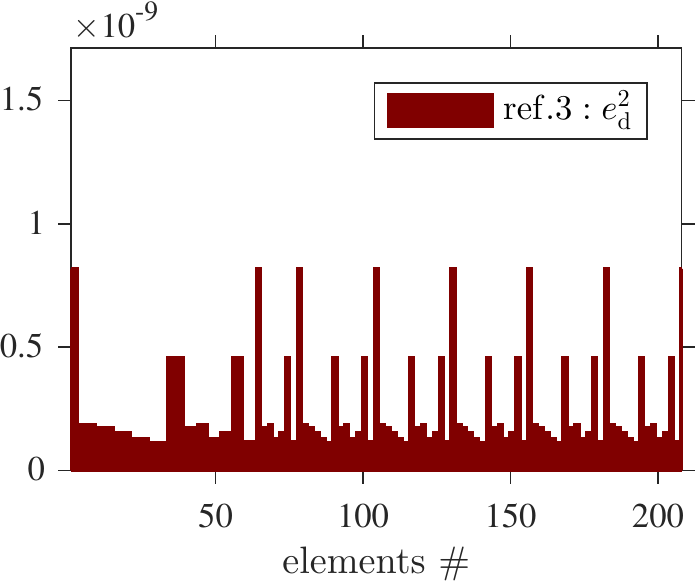}
	}
	\hskip -5pt
	\subfloat[]{
	\includegraphics[width=4.6cm, trim={0cm 0.0cm 0cm 0cm}, clip]{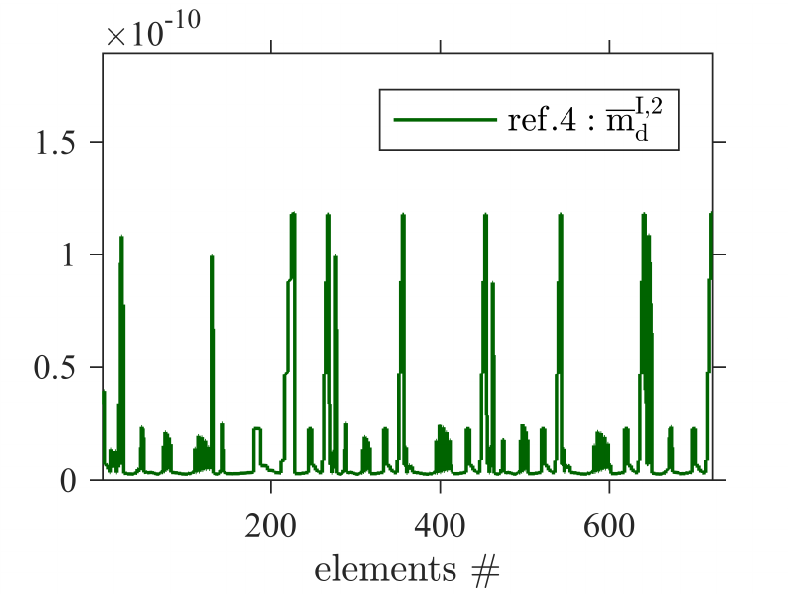}
	}
	\hskip -15pt
	\subfloat[]{
	\includegraphics[width=4.6cm, trim={0cm 0.0cm 0cm 0cm}, clip]{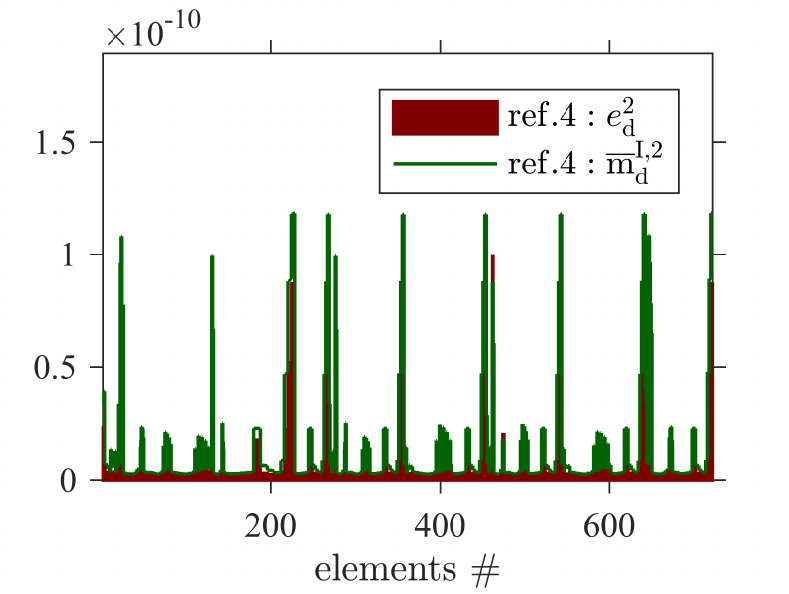}
	}
	\vskip -15pt
	\caption{{\em Example 1}. 
	Distribution of $e^2_{{\rm d}, K} := |\!| \nabla_x e |\!|^2_{K}$ and 
	$\overline{\rm m}^{\rm I, 2}_{{\rm d}, K} := \| \flux_h - \nabla_x u_h\|^2_K$ w.r.t. refinements 2 and 3.}
\label{fig:unit-domain-example-2-ed-md-distribution}
\end{figure}

\begin{figure}[!t]
	\centering
	\captionsetup[subfigure]{oneside, labelformat=empty}
	%
	\subfloat[]{
	\includegraphics[width=4.1cm, trim={0cm 0.0cm 0cm 0cm}, clip]{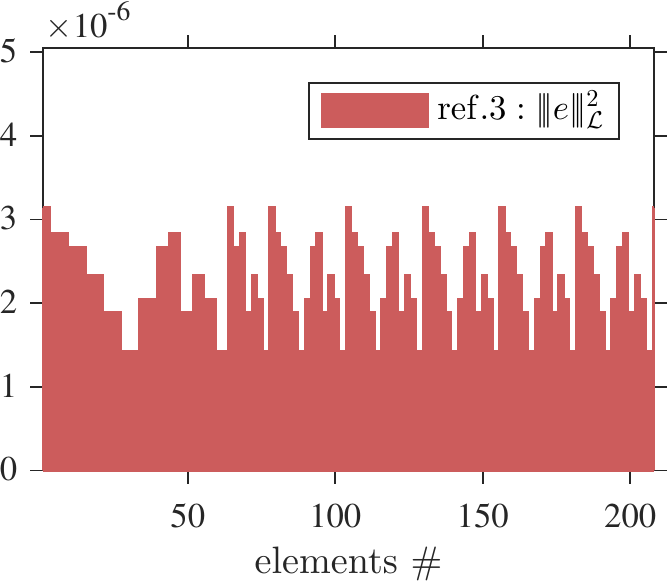}
	}
	\subfloat[]{
	\includegraphics[width=4.1cm, trim={0cm 0.0cm 0cm 0cm}, clip]{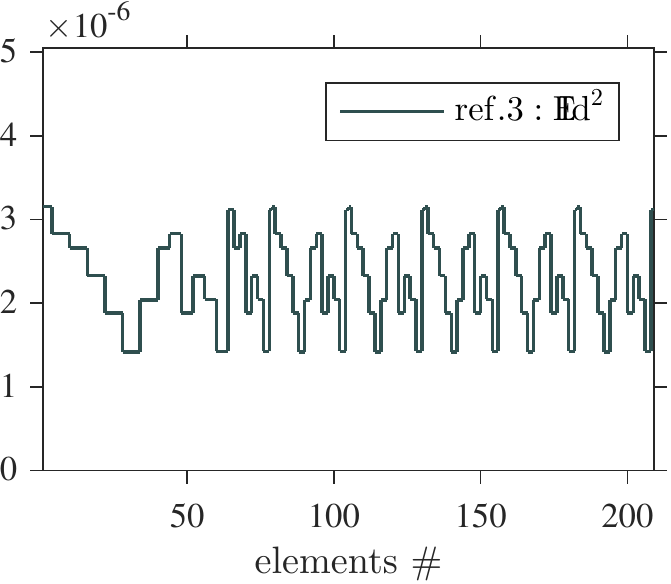}
	}
	\subfloat[]{
	\includegraphics[width=4.1cm, trim={0cm 0.0cm 0cm 0cm}, clip]{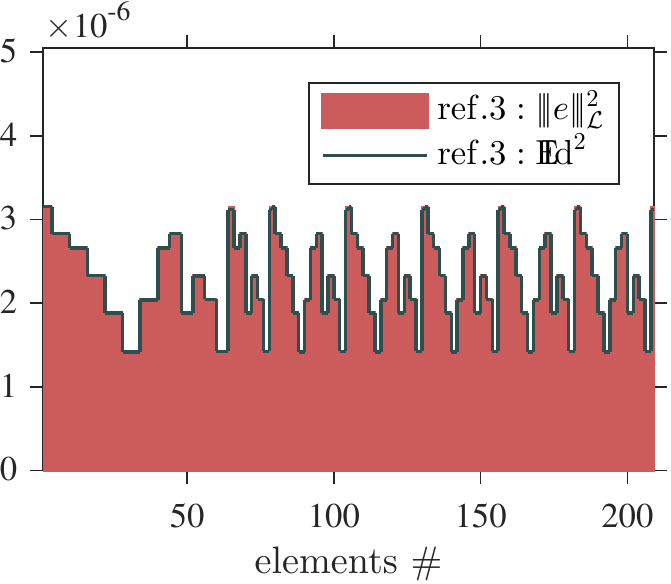}
	}
	\\[-20pt]
	\hskip 5pt
	\subfloat[]{
	\includegraphics[width=4.2cm, trim={0cm 0.0cm 0cm 0cm}, clip]{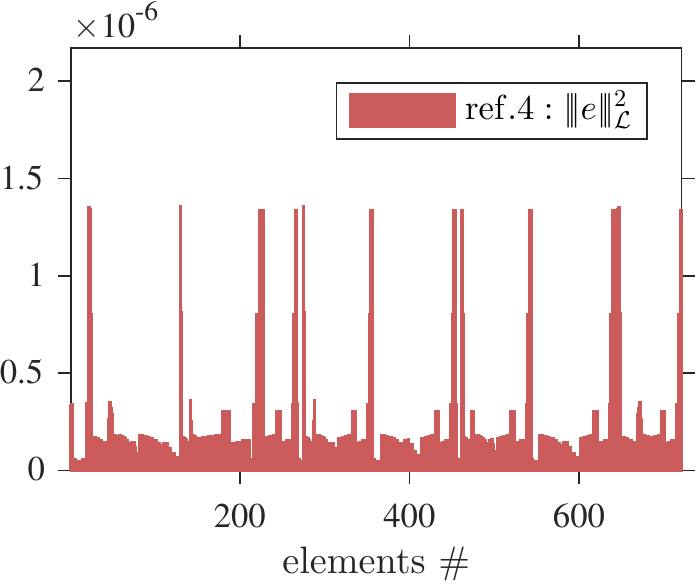}
	}
	\hskip -5pt
	\subfloat[]{
	\includegraphics[width=4.7cm, trim={0cm 0.0cm 0cm 0cm}, clip]{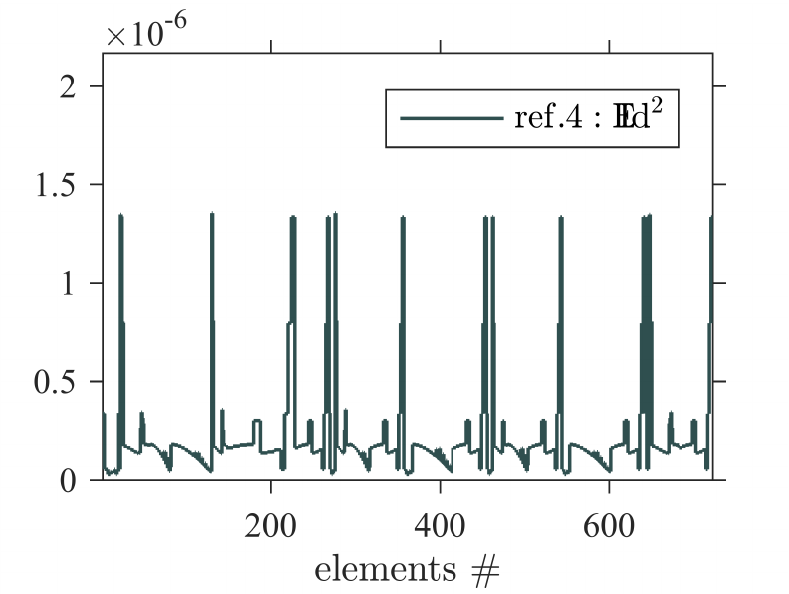}
	}
	\hskip -15pt
	\subfloat[]{
	\includegraphics[width=4.7cm, trim={0cm 0.0cm 0cm 0cm}, clip]{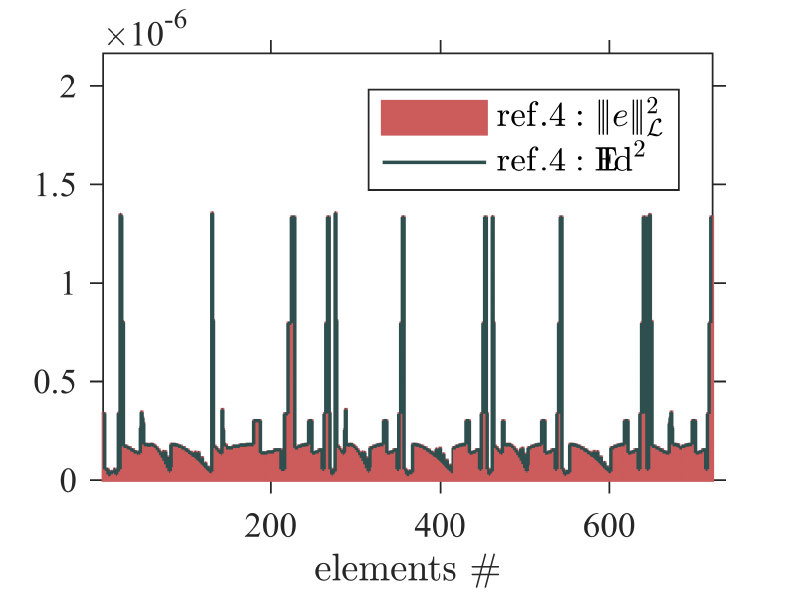}
	}
	\vskip -15pt
	\caption{{\em Example 1}. 
	Distribution of $|\!|\!| e|\!|\!|^2_{\mathcal{L}, K}$ and $\EI_K$ w.r.t. refinements 2 and 3.}
	\label{fig:unit-domain-example-2-e-sol-e-id-distribution}
\end{figure}

We also analyse the quantitative sharpness of the error indication provided 
by $\overline{\rm M}^{\rm I}$ and the local indicator generated by ${\EI}$. 
This can be done by comparing local error contributions and local indicators w.r.t. 
numbered elements $K \in \mathcal{K}_{h}$ (see Figure \ref{fig:unit-domain-example-2-ed-md-distribution} 
and \ref{fig:unit-domain-example-2-e-sol-e-id-distribution}). We can see that { the}
local error contributions 
$e^{2}_{{\rm d}, K}= \| \nabla_x e\|^2_K$ (red bars in the first column of Figure 
\ref{fig:unit-domain-example-2-ed-md-distribution}) are efficiently mimicked by the error indicators
$\overline{\rm m}^{{\rm I}, 2}_{{\rm d}, K}$ (green contour line in the second column). 
The resemblance of these distributions is even stronger emphasised in the third column of 
Figure \ref{fig:unit-domain-example-2-ed-md-distribution}, where plots from the first and the second 
columns overlap. The quantitative sharpness of the ${\EI}$-distribution is analogously confirmed by 
Figure \ref{fig:unit-domain-example-2-e-sol-e-id-distribution}.

\subsection{Example 2}
\label{ex:unit-domain-example-3}
\begin{figure}[!t]
	\centering
	\subfloat[$k_1 = k_2 = 1$]{
	\includegraphics[scale=0.6]{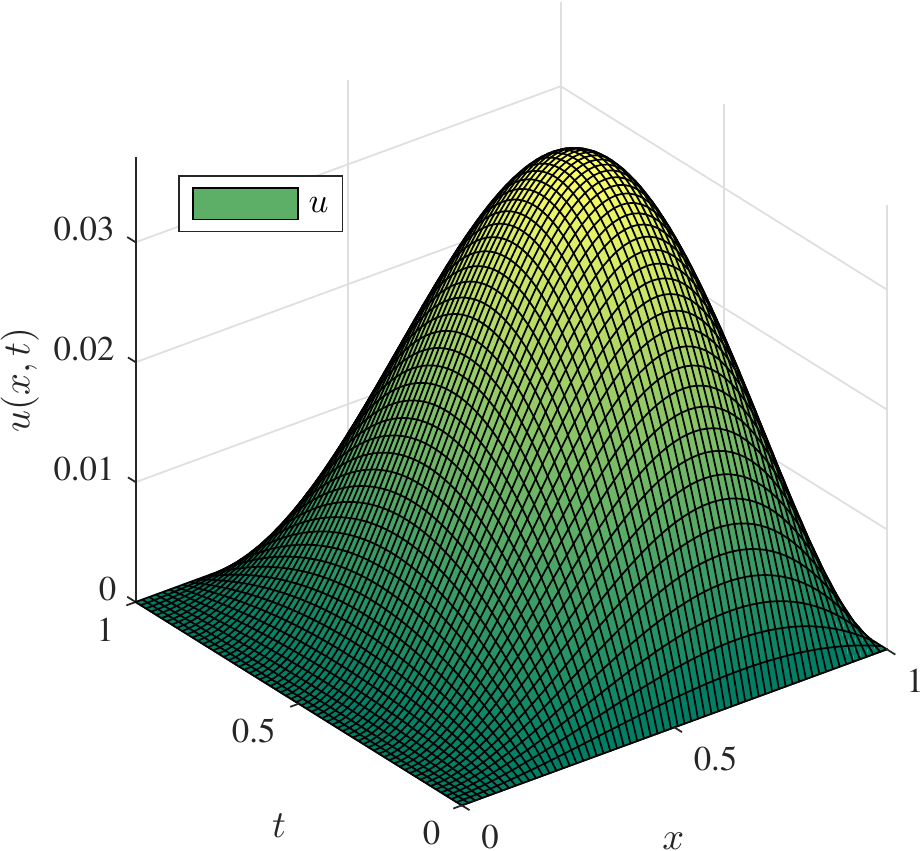}
	\label{fig:example-3-1-exact-solution-a}}
	\quad
	\subfloat[$k_1 = 6, k_2 = 3$]{
	\includegraphics[scale=0.6]{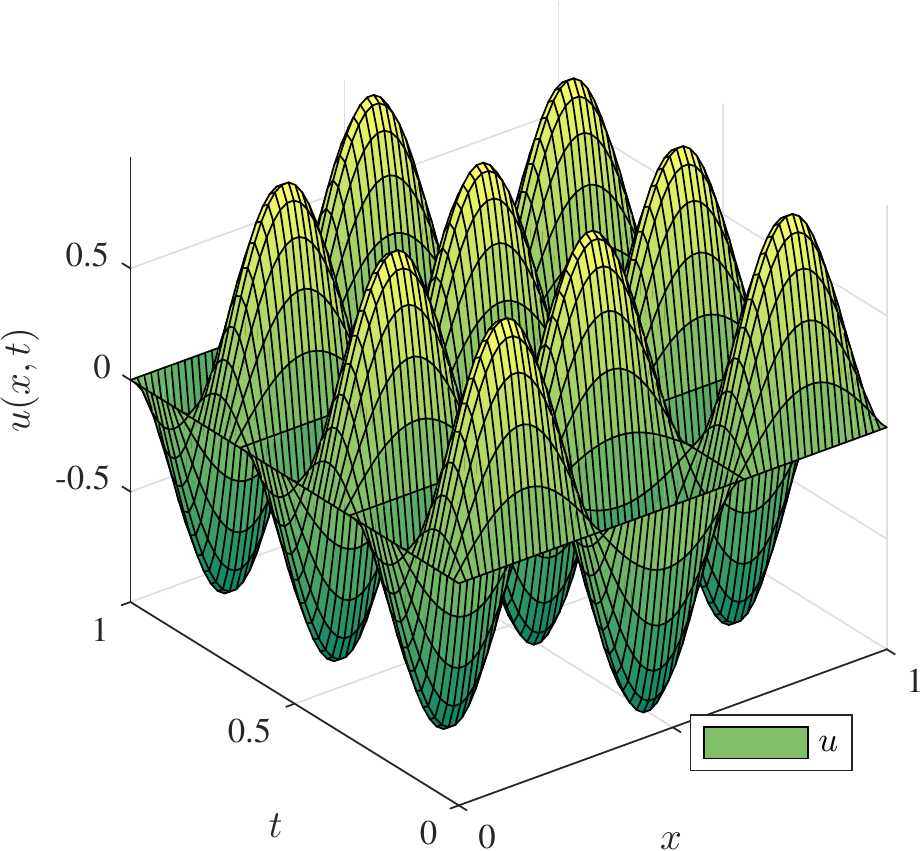}
	\label{fig:example-3-2-exact-solution-a}}
	\caption{\small {\em Example 2}. 
	Exact solution $u = \sin (k_1\,\pi\,x)\, \sin(k_2\,\pi\,t)$.}
	\label{fig:example-3-1-exact-solution}
\end{figure}

\rm
Next, we consider the example with the exact solution such that the change of the gradient depends on 
user-defined parameters. 
Let
$Q = (0, 1)^2$ be the unit square, and let the exact solution, 
the RHS, and the Dirichlet boundary conditions be chosen as follows:
\begin{alignat*}{4}
u(x, t) 	& = \sin (k_1\,\pi\,x)\, \sin(k_2\,\pi\,t)  			
\quad && (x, t)  && \in \overline{Q} = [0, 1]^2, \\
f(x, t) 		& = \sin (k_1\,\pi\,x)\,( k_2\,\pi\,\cos(k_2 \, \pi \, t) + k_1^2 \, \pi^2 \,  \sin(k_2 \, \pi \, t)) 		
\quad && (x, t)  && \in {Q} = (0, 1)^2, \\
u_0(x, t) 	& = 0,	     								
\quad && (x, t)  && \in 
{\overline{\Sigma}}_0, \\
u_D(x, t) 	& = 0,							        		
\quad && (x, t)  && \in \Sigma := \partial \Omega \times (0, 1).
\end{alignat*}
%
%
In the first part of the example (referred to {\em Example 2-1}), we chose the parameters as 
$k_1 = k_2 = 1$. For such $k_1$ and $k_2$, the exact solution is illustrated in Figure 
\ref{fig:example-3-1-exact-solution-a}. The function $u_h$ is approximated by $S_h^{2}$, whereas 
$\flux_h \in S_{7h}^{4} \oplus S_{7h}^{4}$ and $w_h \in S^{4}_{7h}$. We consider eight adaptive 
refinement steps ($N_{\rm ref} = 8$) preceded by three global refinements ($N_{\rm ref, 0} = 3$) 
to generate the initial mesh. For the marking criterion, we use bulk marking with the parameter 
$\sigma = 0.6$. 

\begin{table}[!t]
\scriptsize
\centering
\newcolumntype{g}{>{\columncolor{gainsboro}}c} 	
\newcolumntype{k}{>{\columncolor{lightgray}}c} 	
\newcolumntype{s}{>{\columncolor{silver}}c} 
\newcolumntype{a}{>{\columncolor{ashgrey}}c}
\newcolumntype{b}{>{\columncolor{battleshipgrey}}c}
\begin{tabular}{c|cga|ck|cb|cc}
\parbox[c]{0.8cm}{\centering \# ref. } & 
\parbox[c]{1.4cm}{\centering  $\| \nabla_x e \|_Q$}   & 	  
\parbox[c]{1.0cm}{\centering $\Ieff (\overline{\rm M}^{\rm I})$ } & 
\parbox[c]{1.4cm}{\centering $\Ieff (\overline{\rm M}^{\rm I\!I})$ } & 
\parbox[c]{1.0cm}{\centering  $|\!|\!|  e |\!|\!|_{s, h}$ }   & 	  
\parbox[c]{1.4cm}{\centering $\Ieff (\overline{\rm M}^{\rm I}_{s, h})$ } & 
\parbox[c]{1.0cm}{\centering  $|\!|\!|  e |\!|\!|_{\mathcal{L}}$ }   & 	  
\parbox[c]{1.4cm}{\centering$\Ieff ({\EI})$ } & 
\parbox[c]{1.2cm}{\centering e.o.c. ($|\!|\!|  e |\!|\!|_{s, h}$)} & 
\parbox[c]{1.2cm}{\centering e.o.c. ($|\!|\!|  e |\!|\!|_{\mathcal{L}}$)} \\
\midrule
%
  2 &     3.8987e-03 &         1.80 &         1.00 &     3.8988e-03 &         1.80 &     3.1934e-01 &         1.00 &     2.68 &     1.76 \\
   4 &     9.5149e-04 &         1.77 &         1.02 &     9.5149e-04 &         1.77 &     1.6041e-01 &         1.00 &     4.18 &     2.59 \\
   6 &     2.6155e-04 &         2.82 &         1.07 &     2.6170e-04 &         2.82 &     8.2550e-02 &         1.00 &     2.74 &     1.59 \\
   8 &     8.2844e-05 &         2.79 &         1.21 &     8.2854e-05 &         2.79 &     4.7088e-02 &         1.00 &     2.49 &     1.79 \\
\end{tabular}
\caption{{\em Example 2-1}. 
Efficiency of $\overline{\rm M}^{\rm I}$, $\overline{\rm M}^{\rm I\!I}$, 
$\overline{\rm M}^{\rm I}_{s, h}$, and ${\EI}$ for $u_h \in S^{2}_{h}$,
$\flux_h \in S_{7h}^{4} \oplus S_{7h}^{4}$, and $w_h \in S^{4}_{7h}$,
w.r.t. adaptive refinements (with the marking criterion ${\mathds{M}}_{\rm BULK}(0.6)$).}
\label{tab:unit-domain-example-3-1-error-majorant-v-2-y-4-adaptive-ref}
\end{table}

\begin{table}[!t]
\scriptsize
\centering
\newcolumntype{g}{>{\columncolor{gainsboro}}c} 	
\begin{tabular}{c|ccc|cgg|cgg|c}
& \multicolumn{3}{c|}{ d.o.f. } 
& \multicolumn{3}{c|}{ $t_{\rm as}$ }
& \multicolumn{3}{c|}{ $t_{\rm sol}$ } 
& $\tfrac{t_{\rm appr.}}{t_{\rm er.est.}}$ \\
\midrule
\parbox[c]{0.8cm}{\centering \# ref. } & 
\parbox[c]{0.8cm}{\centering $u_h$ } &  
\parbox[c]{0.6cm}{\centering $\flux_h$ } &  
\parbox[c]{0.6cm}{\centering $w_h$ } & 
\parbox[c]{1.4cm}{\centering $u_h$ } & 
\parbox[c]{1.4cm}{\centering $\flux_h$ } & 
\parbox[c]{1.4cm}{\centering $w_h$ } & 
\parbox[c]{1.4cm}{\centering $u_h$ } & 
\parbox[c]{1.4cm}{\centering $\flux_h$ } & 
\parbox[c]{1.4cm}{\centering $w_h$ } & 
\\
\midrule
   2 &        190 &        288 &        144 &   1.58e-01 &   6.92e-01 &   3.32e-01 &         1.51e-03 &         1.76e-03 &         3.94e-03 & 0.15\\
   4 &        716 &        288 &        144 &   7.13e-01 &   7.20e-01 &   3.78e-01 &         2.24e-02 &         1.40e-03 &         2.34e-03 & 0. 65\\
   6 &       2588 &        288 &        144 &   2.09e+00 &   7.06e-01 &   2.97e-01 &         1.73e-01 &         1.01e-03 &         2.20e-03  & 2.24\\
   8 &       8303 &        288 &        144 &   8.27e+00 &   4.08e-01 &   3.08e-01 &         8.17e-01 &         8.65e-04 &         1.74e-03 & 5.92\\
     \midrule
    &       &         &    &
    \multicolumn{3}{c|}{ $t_{\rm as} (u_h)$ \quad : \quad $t_{\rm as} (\flux_h)$ \quad : \qquad $t_{\rm as} (w_h)$ } &      
    \multicolumn{3}{c|}{\; $t_{\rm sol} (u_h)$ \, : \quad $t_{\rm sol} (\flux_h)$  \quad:  \qquad  $t_{\rm sol} (w_h)$\;} & 
    \\
 \midrule
         &       &       &       &      26.89 &       1.33 &       1.00 &           470.49 &             0.50 &             1.00 & \\
\end{tabular}
\caption{{\em Example 2-1}. 
Assembling and solving time (in seconds) spent for the systems generating
d.o.f. of $u_h \in S^{2}_{h}$,
$\flux_h \in S^{4}_{7h} \oplus S^{4}_{7h}$, and $w_h \in S^{4}_{7h}$
w.r.t. adaptive refinements (with the marking criterion ${\mathds{M}}_{\rm BULK}(0.6)$).}
\label{tab:unit-domain-example-3-times-v-2-y-4-adaptive-ref}
\end{table}

The resulting performance of the majorants and the error identity is presented in 
Table \ref{tab:unit-domain-example-3-1-error-majorant-v-2-y-4-adaptive-ref}. 
It is again clear that $\overline{\rm M}^{\rm I\!I}$ is $1.8$--$2.8$ times sharper than 
$\overline{\rm M}^{\rm I}$. The performance of ${\EI}$ remains sharp even 
though its values as well as the values of $|\!|\!| e |\!|\!|_{\mathcal{L}}$ decrease one order 
slower than $\overline{\rm M}^{\rm I}$ and $|\!|\!| e |\!|\!|$.
At the same time, if we compare the effort spent on the reconstruction of $\flux_h$ and 
$w_h$, we see from Table \ref{tab:unit-domain-example-3-times-v-2-y-4-adaptive-ref} that the 
approximation of $u_h$ takes longer. Total time expenses invested in $u_h$ are
again compared to the costs of error-control in the last column of Table 
\ref{tab:unit-domain-example-3-times-v-2-y-4-adaptive-ref}, where it is shown that ratios of such 
expenses reach $5.92$ on the last refinement step.

\begin{figure}[!t]
	\centering
	\captionsetup[subfigure]{oneside, margin={0.7cm,0cm}}
	\subfloat[{REF 5: \newline
	ref. based on $|\!|\!|  e |\!|\!|_{h, s}$}]{
	\spacetimeaxis{\includegraphics[width=4.0cm, trim={8.1cm 2cm 6cm 2cm}, clip]{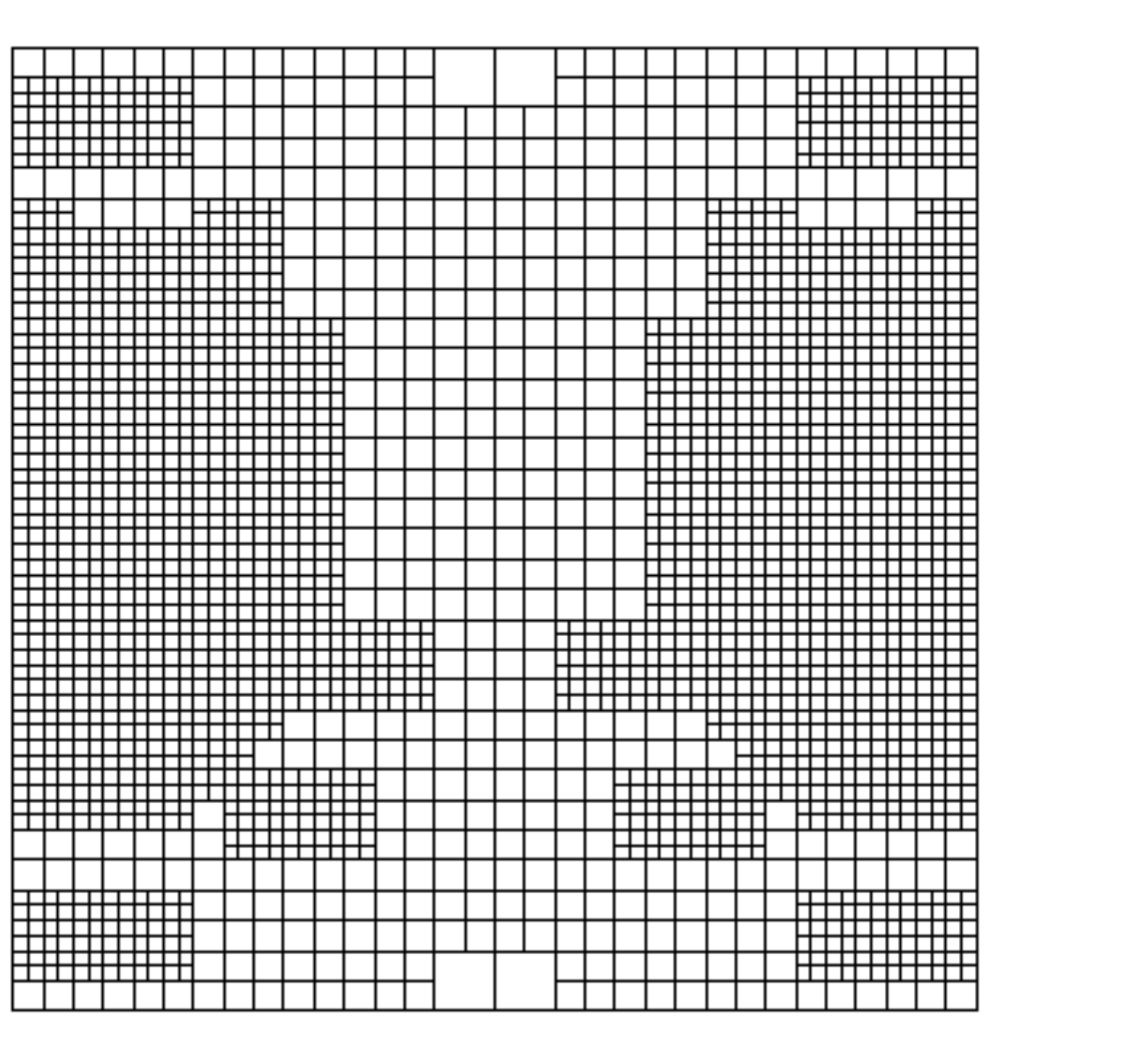}
	}
	}
	\hskip -28pt
	\subfloat[{REF 5: \newline
	ref. based on $\overline{\rm m}^{\rm I}_{\rm d}$ 
	}]{
	\spacetimeaxis{\includegraphics[width=4.0cm, trim={8.1cm 2cm 6cm 2cm}, clip]{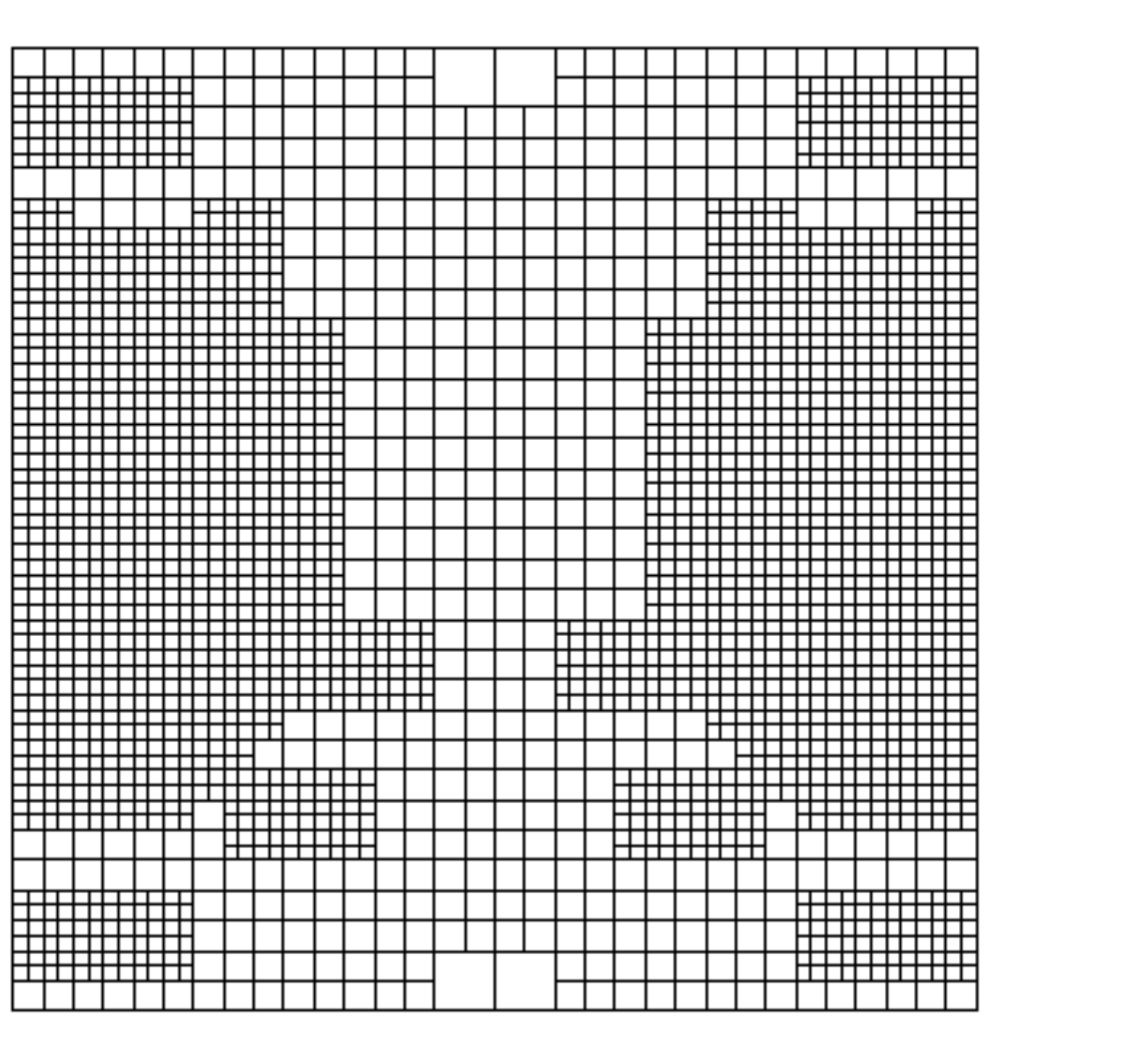}
	}
	}
	\hskip -28pt
	\subfloat[{REF 5: \newline
	ref. based on $|\!|\!|  e |\!|\!|_{\mathcal{L}}$}]{
	\spacetimeaxis{\includegraphics[width=4.0cm, trim={8.1cm 2cm 6cm 2cm}, clip]{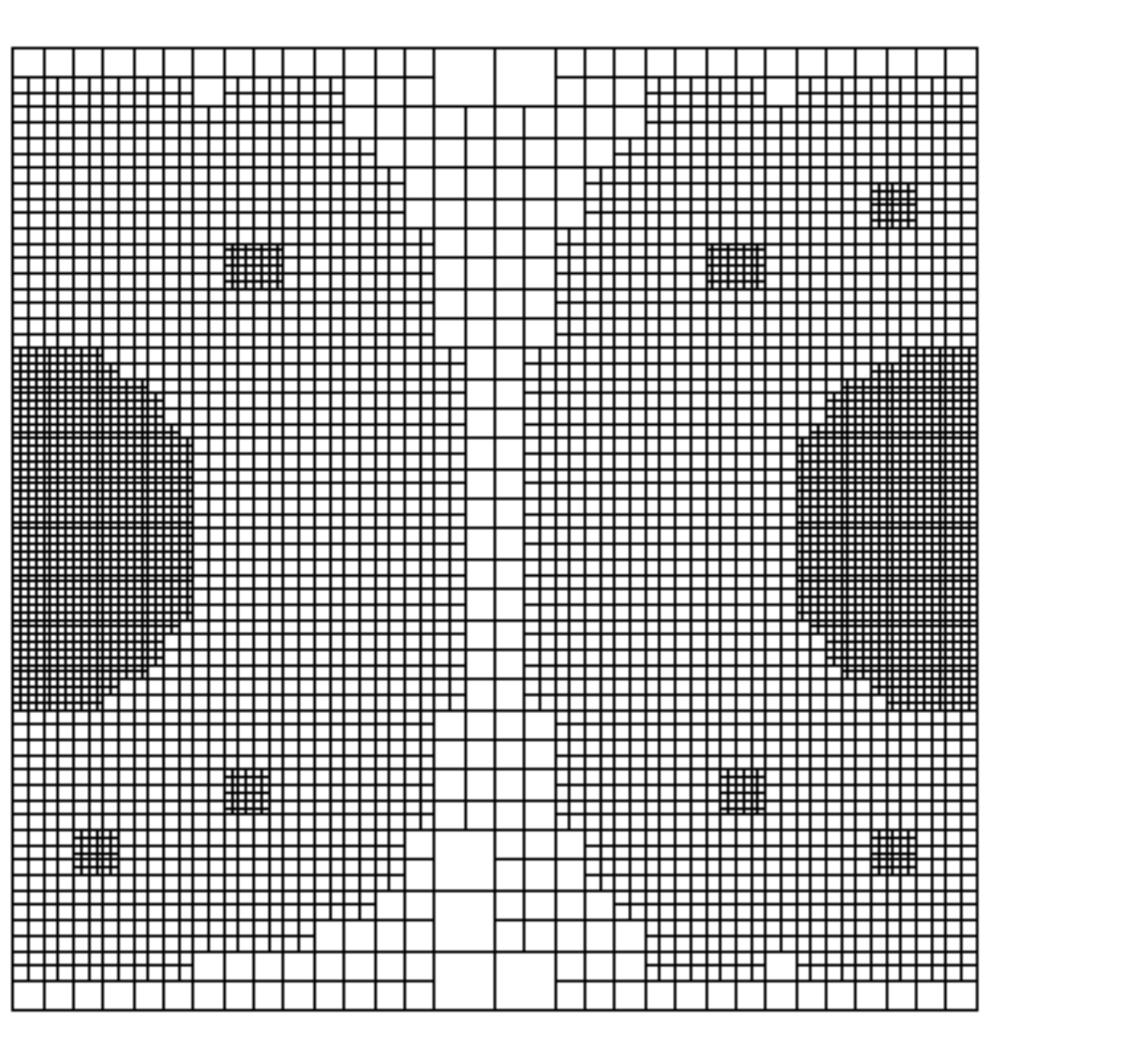}
	}
	}
	\hskip -28pt
	\subfloat[{REF 5: \newline
	ref. based on ${\EI}$}]{
	\spacetimeaxis{\includegraphics[width=4.0cm, trim={8.1cm 2cm 6cm 2cm}, clip]{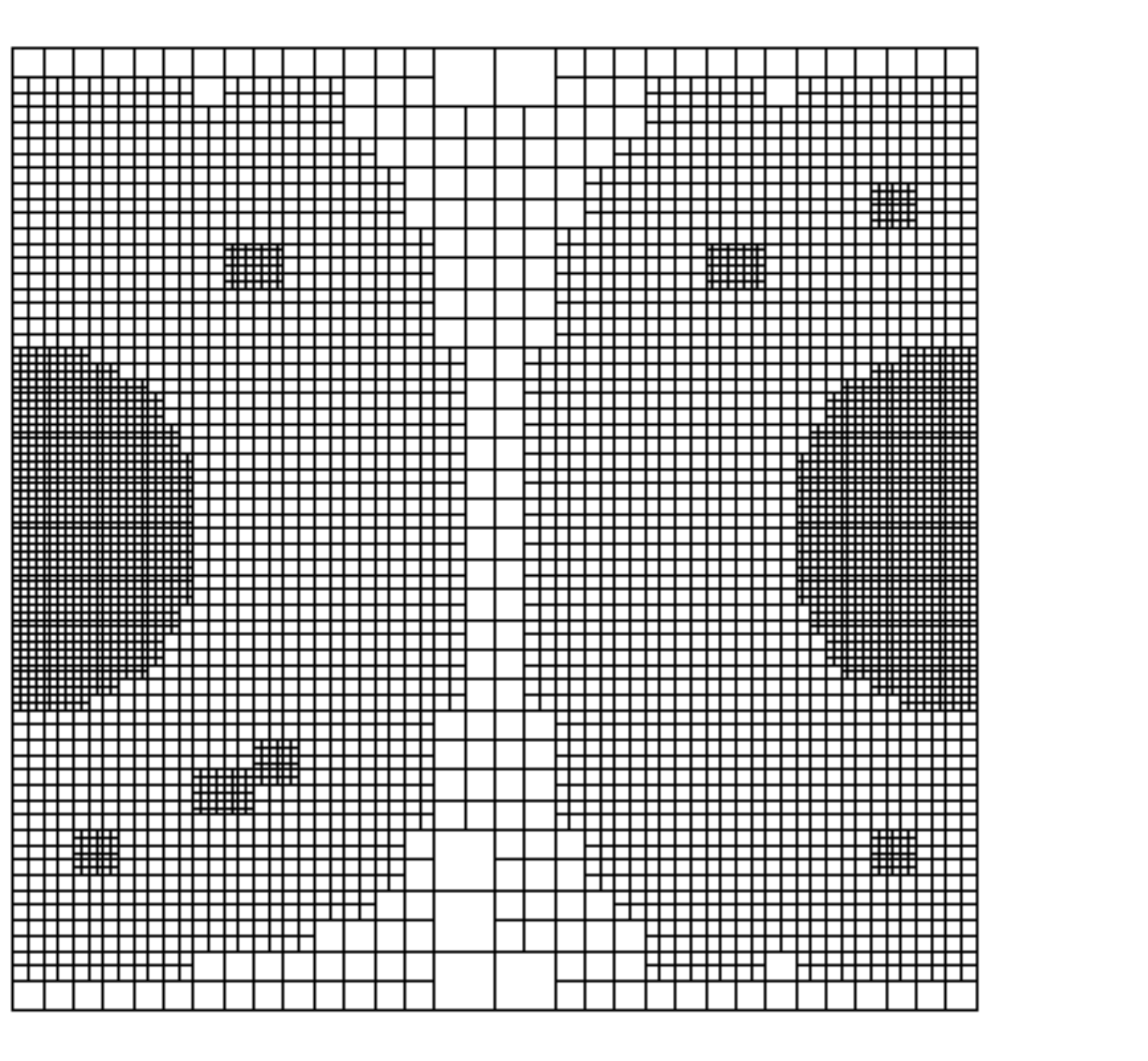}
	}
	}
	\\[-0.5pt]
	\subfloat[{REF 7: \newline
	ref. based on $|\!|\!|  e |\!|\!|_{h, s}$}]{
	\spacetimeaxis{\includegraphics[width=4.0cm, trim={8.1cm 2cm 6cm 2cm}, clip]{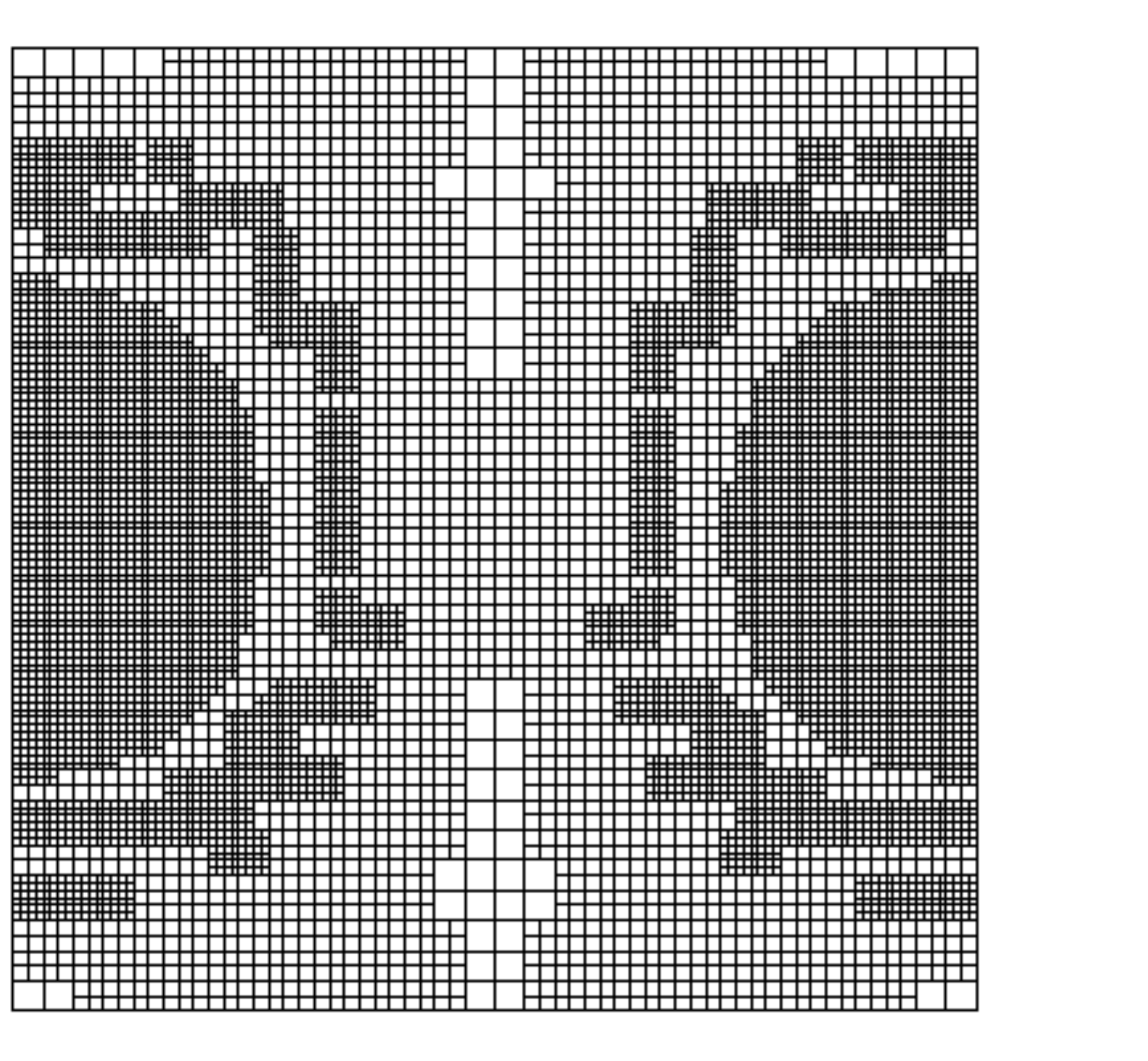}
	}
	}
	\hskip -28pt
	\subfloat[{REF 7: \newline
	ref. based on $\overline{\rm m}^{\rm I}_{\rm d}$ 
	}]{
	\spacetimeaxis{\includegraphics[width=4.0cm, trim={8.1cm 2cm 6cm 2cm}, clip]{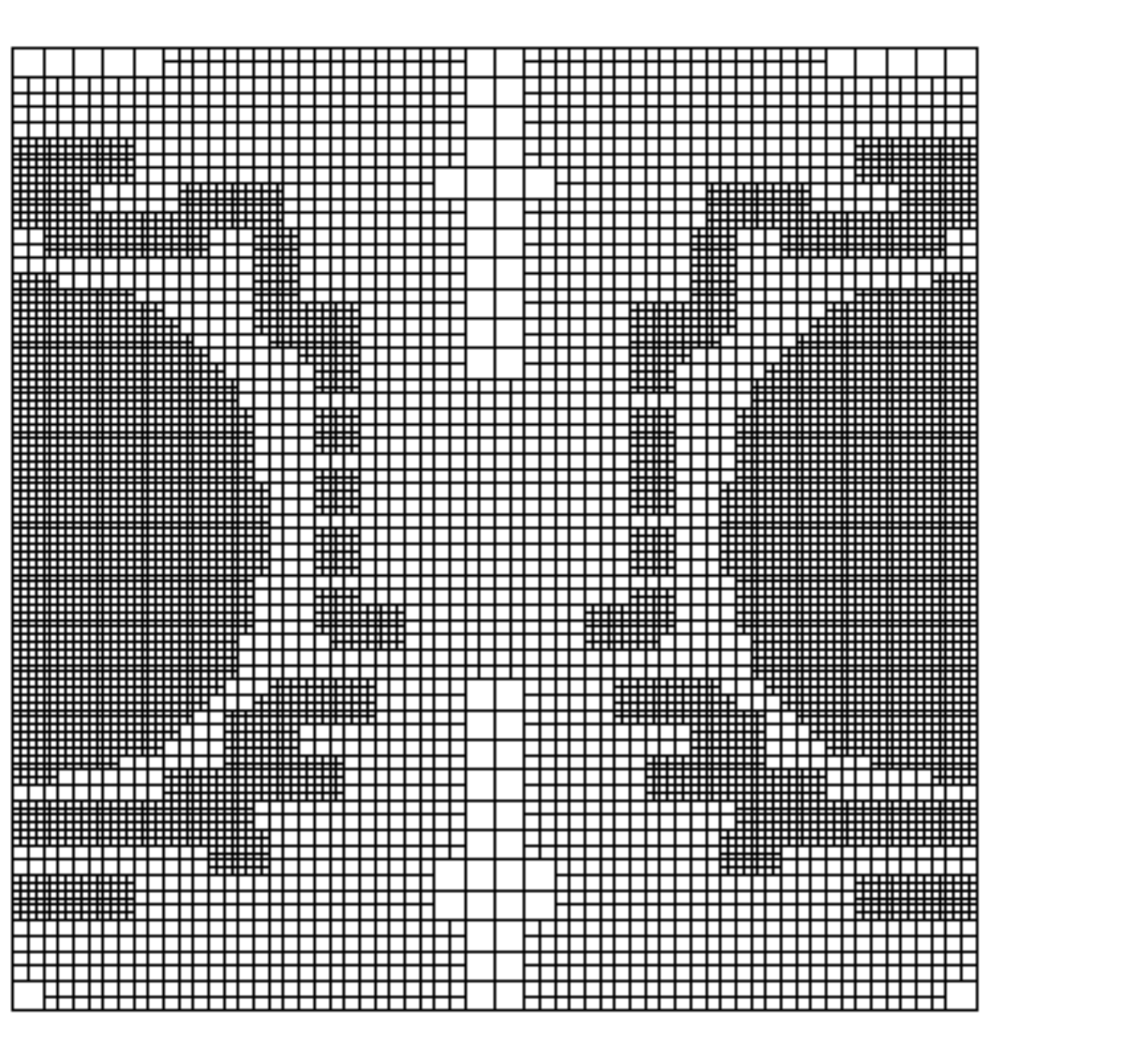}
	}
	}
	\hskip -28pt
	\subfloat[{REF 7: \newline
	ref. based on $|\!|\!|  e |\!|\!|_{\mathcal{L}}$}]{
	\spacetimeaxis{\includegraphics[width=4.0cm, trim={8.1cm 2cm 6cm 2cm}, clip]{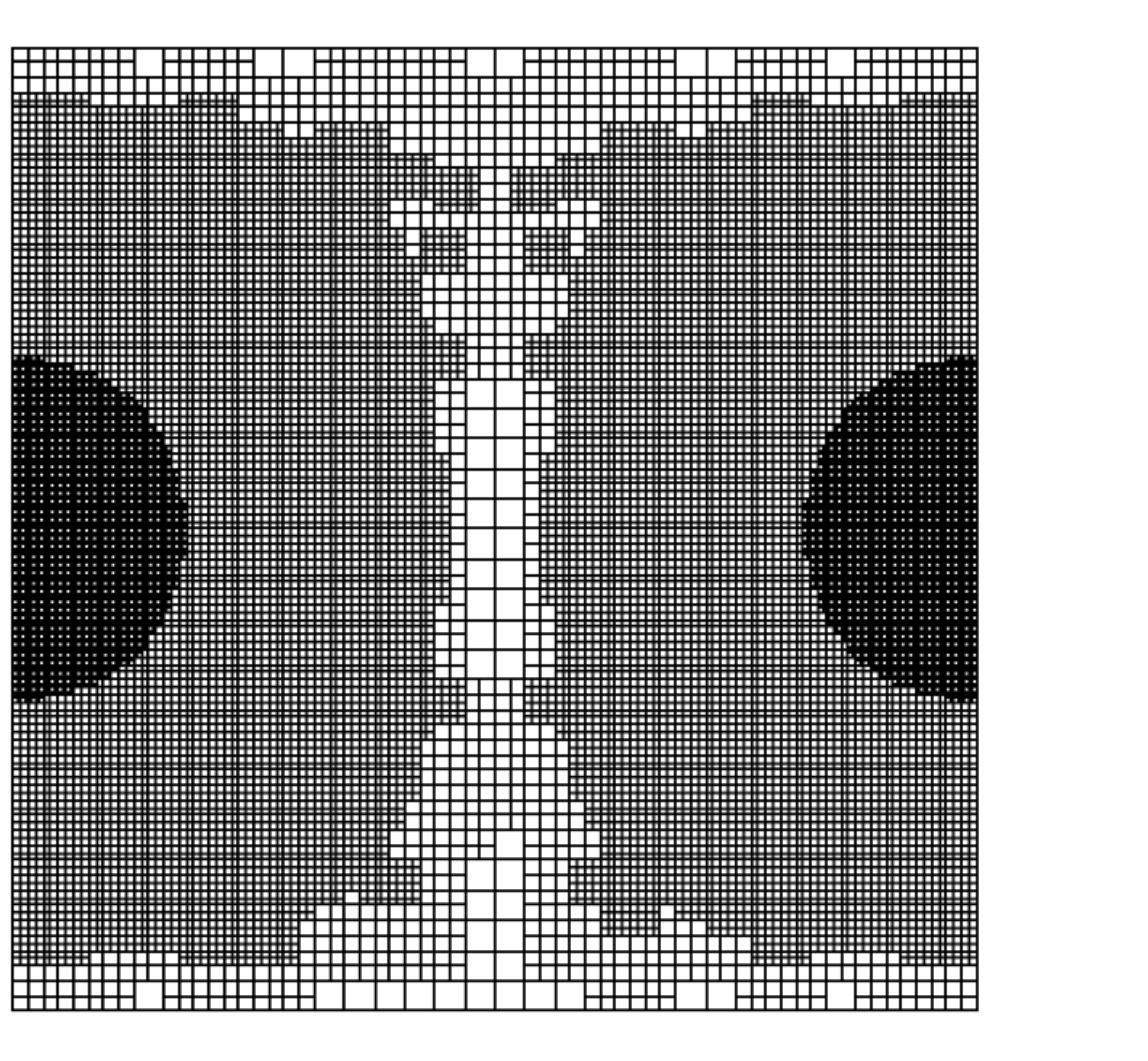}
	}
	}
	\hskip -28pt
	\subfloat[{REF 7: \newline
	ref. based on ${\EI}$}]{
	\spacetimeaxis{\includegraphics[width=4.0cm, trim={8.1cm 2cm 6cm 2cm}, clip]{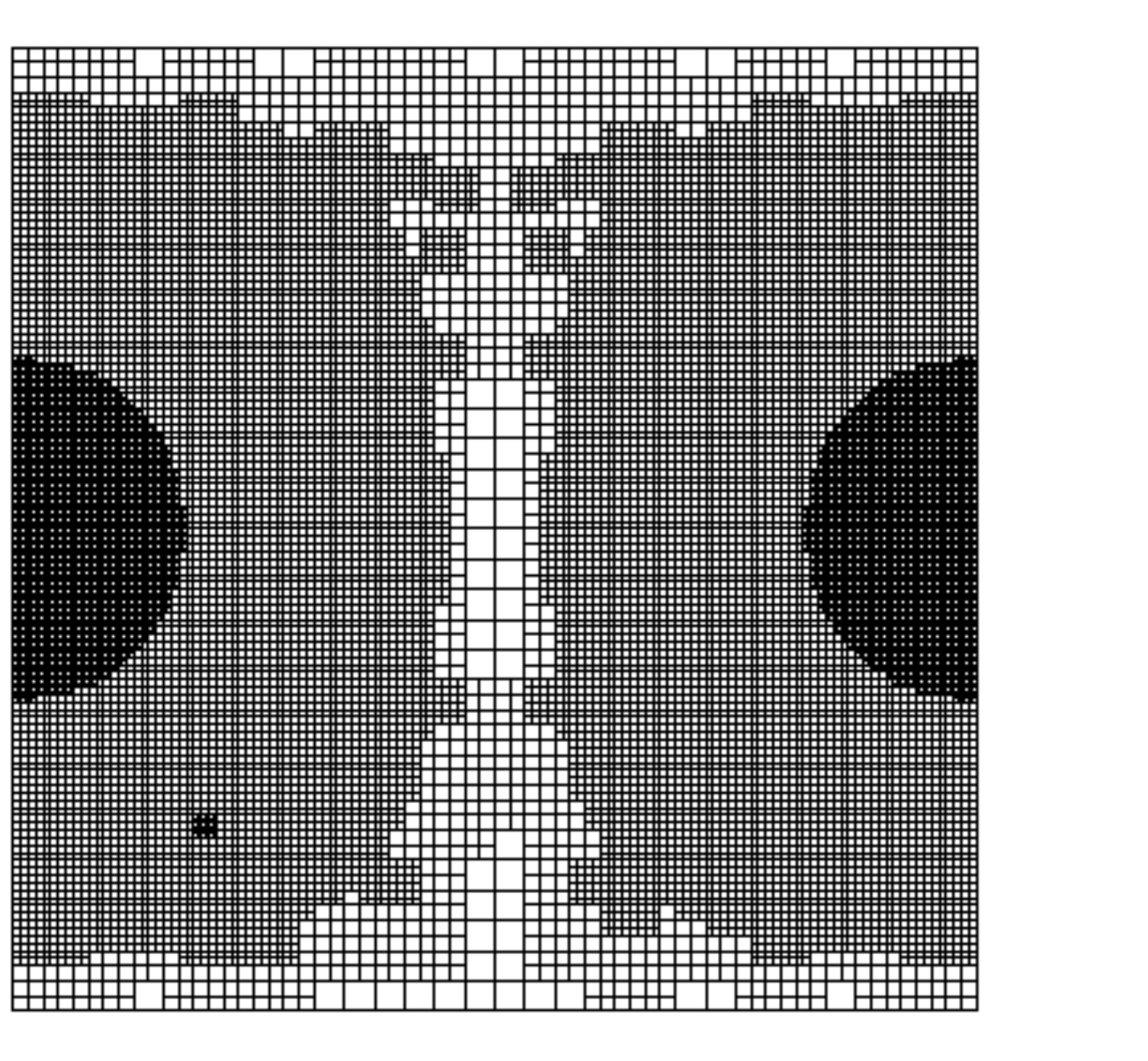}
	}
	}
	\caption{{\em Example 2-1}. 
	Comparison of meshes obtained by the refinement based on $|\!|\!| e |\!|\!|_{h, s}$, $\overline{\rm m}^{\rm I}_{\rm d}$, 
	$|\!|\!|  e |\!|\!|_{\mathcal{L}}$ and $|\!|\!|  e |\!|\!|_{\mathcal{L}}$, where
	$\flux_h \in S^{4}_{7h} \oplus S^{4}_{7h}$, and $w_h \in S^{4}_{7h}$,
	(with the marking criterion ${\mathds{M}}_{\rm BULK}(0.6)$).}
	\label{fig:unit-domain-example-3-meshes-v-2-y-4-adaptive-ref}
\end{figure}

\begin{figure}[!t]
	\centering
	\captionsetup[subfigure]{oneside, margin={0.7cm,0cm}, labelformat=empty}
	\subfloat[]{
	\includegraphics[width=4.2cm, trim={0cm 0cm 0cm 0cm}, clip]{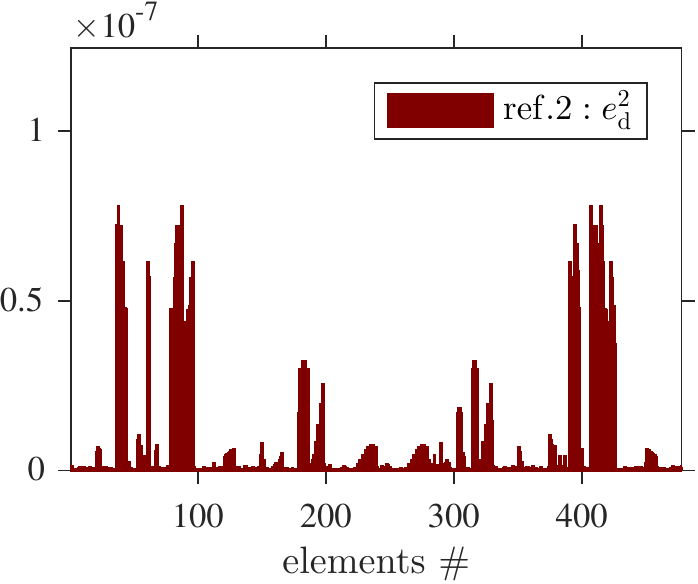}
	}
	\subfloat[]{
	\includegraphics[width=4.2cm, trim={0cm 0cm 0cm 0cm}, clip]{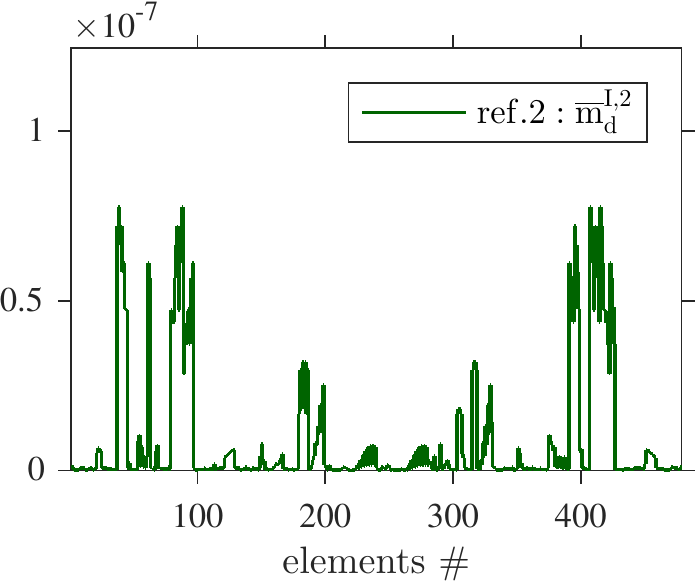}
	}
	\subfloat[]{
	\includegraphics[width=4.2cm, trim={0cm 0cm 0cm 0cm}, clip]{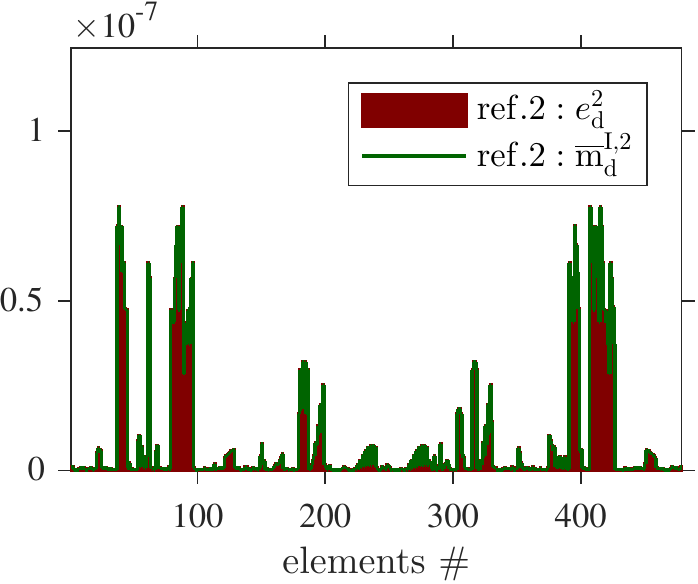}
	}
	\vskip -15pt
	\caption{{\em Example 2-1}. 
	Local distribution of $e_{{\rm d}, K}$ and $\overline{\rm m}_{{\rm d}, K}$ 
	on the refinement step $2$.}
	\label{fig:unit-domain-example-3-1-edmd-distribution}
\end{figure}

\begin{figure}[!t]
	\centering
	\captionsetup[subfigure]{oneside, margin={0.7cm,0cm}, labelformat=empty}
	\subfloat[]{
	\includegraphics[width=4.4cm, trim={0cm 0cm 0cm 0cm}, clip]{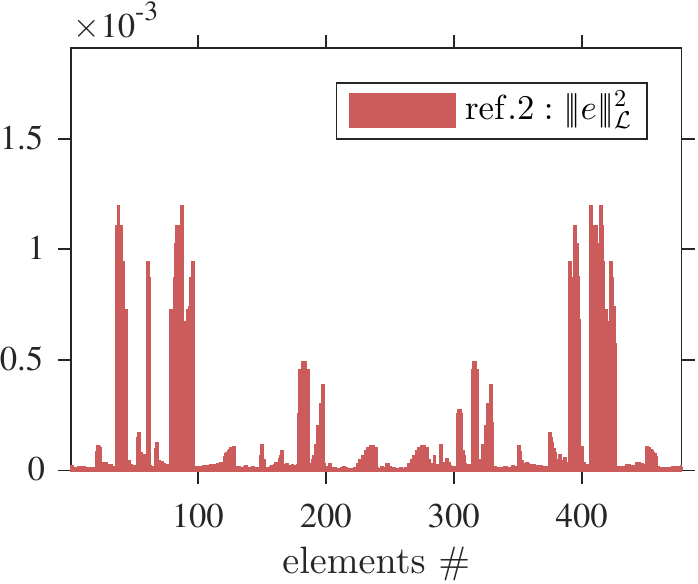}
	}
	\subfloat[]{
	\includegraphics[width=4.4cm, trim={0cm 0cm 0cm 0cm}, clip]{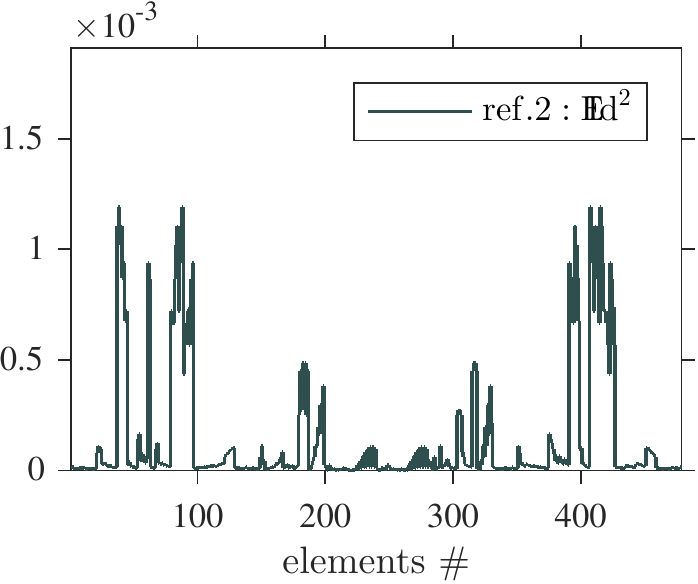}
	}
	\subfloat[]{
	\includegraphics[width=4.4 cm, trim={0cm 0cm 0cm 0cm}, clip]{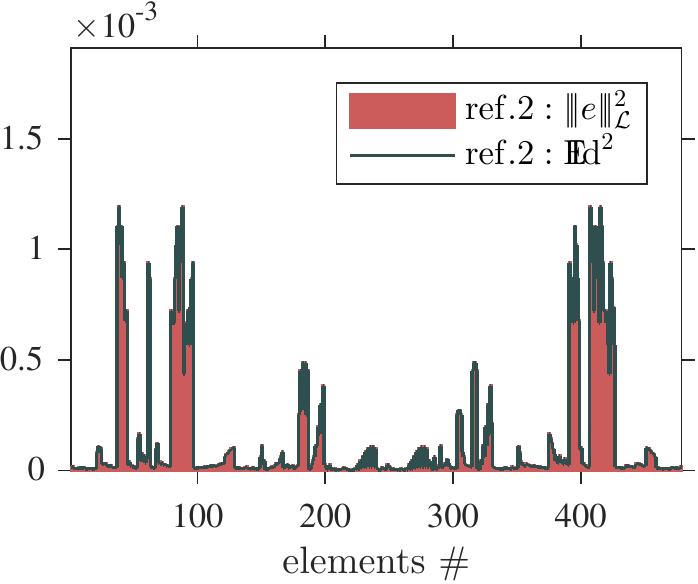}
	}
	\vskip -15pt
	\caption{{\em Example 2-1}. 
	Local distribution of $|\!|\!|  e |\!|\!|_{\mathcal{L}, K}$ and $\EI_K$ 
	on the refinement step $2$.}
	\label{fig:unit-domain-example-3-1-e-sol-e-id-distribution}
\end{figure}

The comparison of the meshes in Figure \ref{fig:unit-domain-example-3-meshes-v-2-y-4-adaptive-ref}
illustrates that the refinement based on $\| \nabla_x e \|_Q$ and the indicator $\overline{\rm m}^{\rm I}_{\rm d}$ 
(first and second columns) provide similar results. The same observation holds when we compare the 
meshes produced by refinement based on the distributions of $|\!|\!|  e |\!|\!|_{\mathcal{L}, K}$ 
and $\EI_K$. The meshes, which we obtain using the majorant, mimic the topology of the 
meshes in the second column. The similarity of the meshes in the third and fourth columns provides clear 
evidence on the sharpness and the efficiency of $\EI_K$, when the error indication is concerned. 
Moreover, the local error contribution $\| \nabla_x e \|_K$ and indicator 
$\overline{\rm m}^{\rm I}_{\rm d, K}$, as well as 
$|\!|\!|  e |\!|\!|_{\mathcal{L}, K}$ an $\EI_K$, are compared in Figures 
\ref{fig:unit-domain-example-3-1-edmd-distribution} and 
\ref{fig:unit-domain-example-3-1-e-sol-e-id-distribution}, respectively. We illustrate the values 
of $\| \nabla_x e \|_{K}$ and $\overline{\rm m}^{\rm I}_{{\rm d}, K}$ 
($|\!|\!|  e |\!|\!|_{\mathcal{L}, K}$ and $\EI_{K}$) on refinement step $2$.


{ 
To conclude the evaluation of the error indication properties for this test-case, we compare the performance 
of the local error indicator generated by the majorant and local residual 
$\| \Delta_x u_h + f - \partial_t u_h\|^2_K$, $K \in \mathcal{K}_h$. Figure 
\ref{fig:unit-domain-example-3-comparison-meshes-v-2-y-4-adaptive-ref} demonstrates the similarity of 
the meshes obtained when the refinement is performed by local true error distribution (first column) and 
by the  majorant error indicator (second column). At the same time, it emphasises that the refinement 
strategy provided by the local residual $\| \Delta_x u_h + f - \partial_t u_h\|^2_K$ is not as quantitatively 
exact.  The corresponding meshes are illustrated in the third column of Figure 
\ref{fig:unit-domain-example-3-comparison-meshes-v-2-y-4-adaptive-ref} and differ from the meshes in the 
first column. Moreover, such an indicator has only heuristic nature and does not provide reliable error 
estimation.

\begin{figure}[!t]
	\centering
	\captionsetup[subfigure]{oneside, margin={0.7cm,0cm}}
	\subfloat[{REF 5: \newline
	ref. based on $\| \nabla_x e\|^2_K$}]{
	\spacetimeaxis{\includegraphics[width=4.5cm, trim={2cm 3cm 2cm 5cm}, clip]{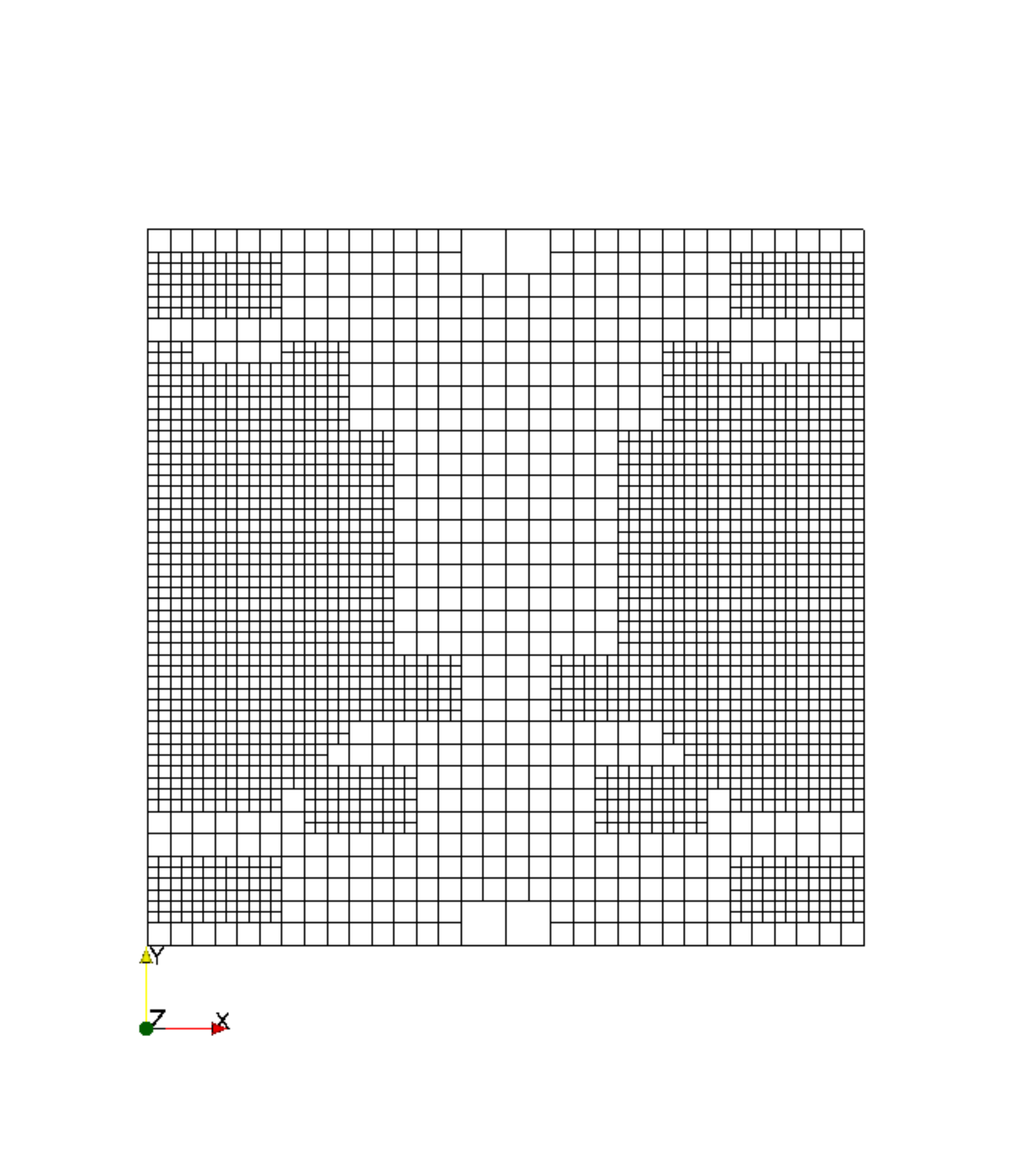}
	}
	}
	\subfloat[{REF 5: \newline
	ref. based on $\overline{\rm m}^{\rm I}_{{\rm d}, K}$ 
	}]{
	\spacetimeaxis{\includegraphics[width=4.5cm, trim={2cm 3cm 2cm 5cm}, clip]{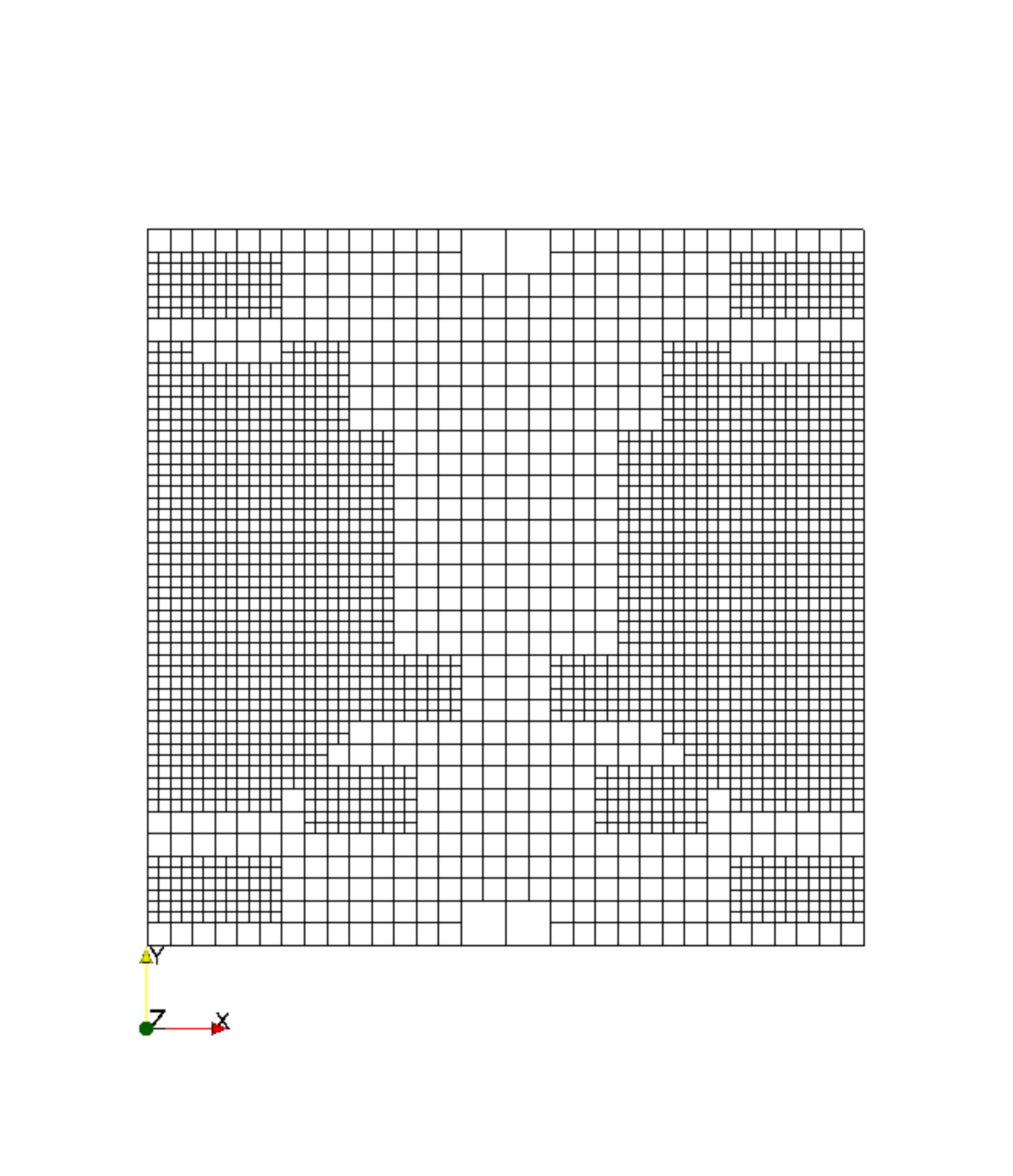}
	}
	}
	\subfloat[{REF 5: \newline
	ref. based on $\| \Delta_x u_h + f - \partial_t u_h\|^2_K$}]{
	\spacetimeaxis{\includegraphics[width=4.5cm, trim={2cm 3cm 2cm 5cm}, clip]{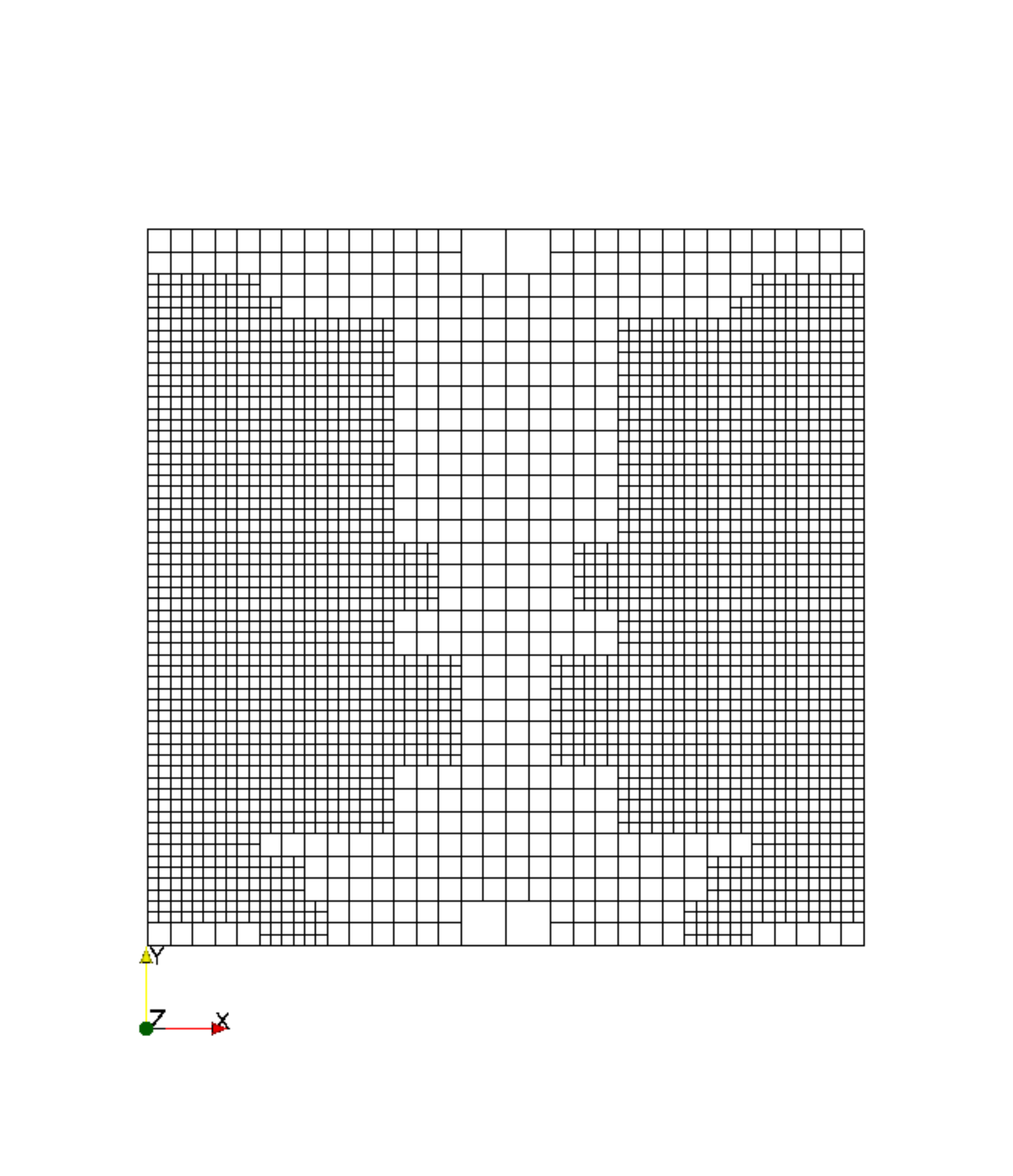}
	}
	}
	\\[-0.5pt]
	\captionsetup[subfigure]{oneside, margin={0.7cm,0cm}}
	\subfloat[{REF 6: \newline
	ref. based on $\| \nabla_x e\|^2_K$}]{
	\spacetimeaxis{\includegraphics[width=4.5cm, trim={2cm 3cm 2cm 5cm}, clip]{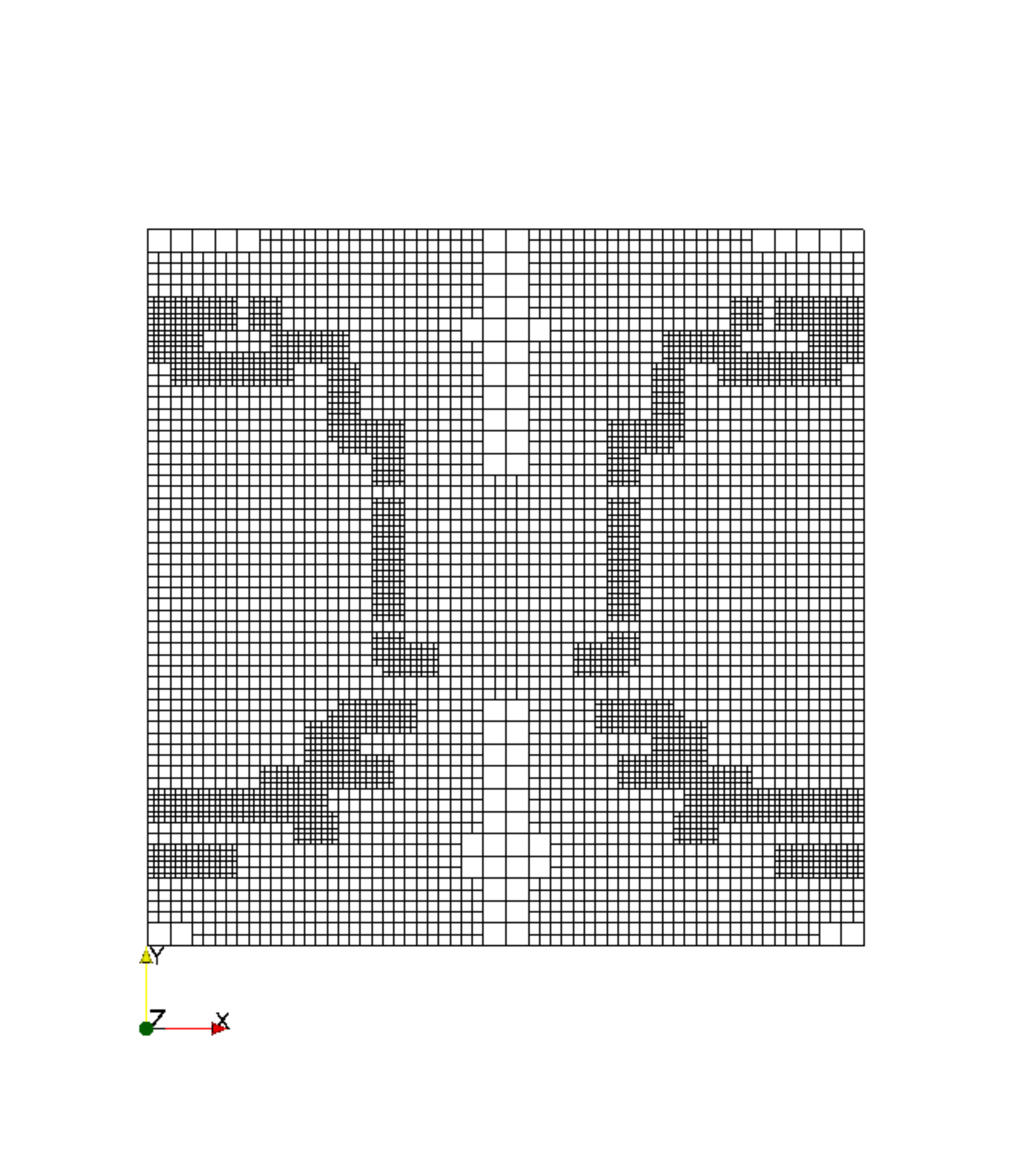}
	}
	}
	\subfloat[{REF 6: \newline
	ref. based on $\overline{\rm m}^{\rm I}_{{\rm d}, K}$ 
	}]{
	\spacetimeaxis{\includegraphics[width=4.5cm, trim={2cm 3cm 2cm 5cm}, clip]{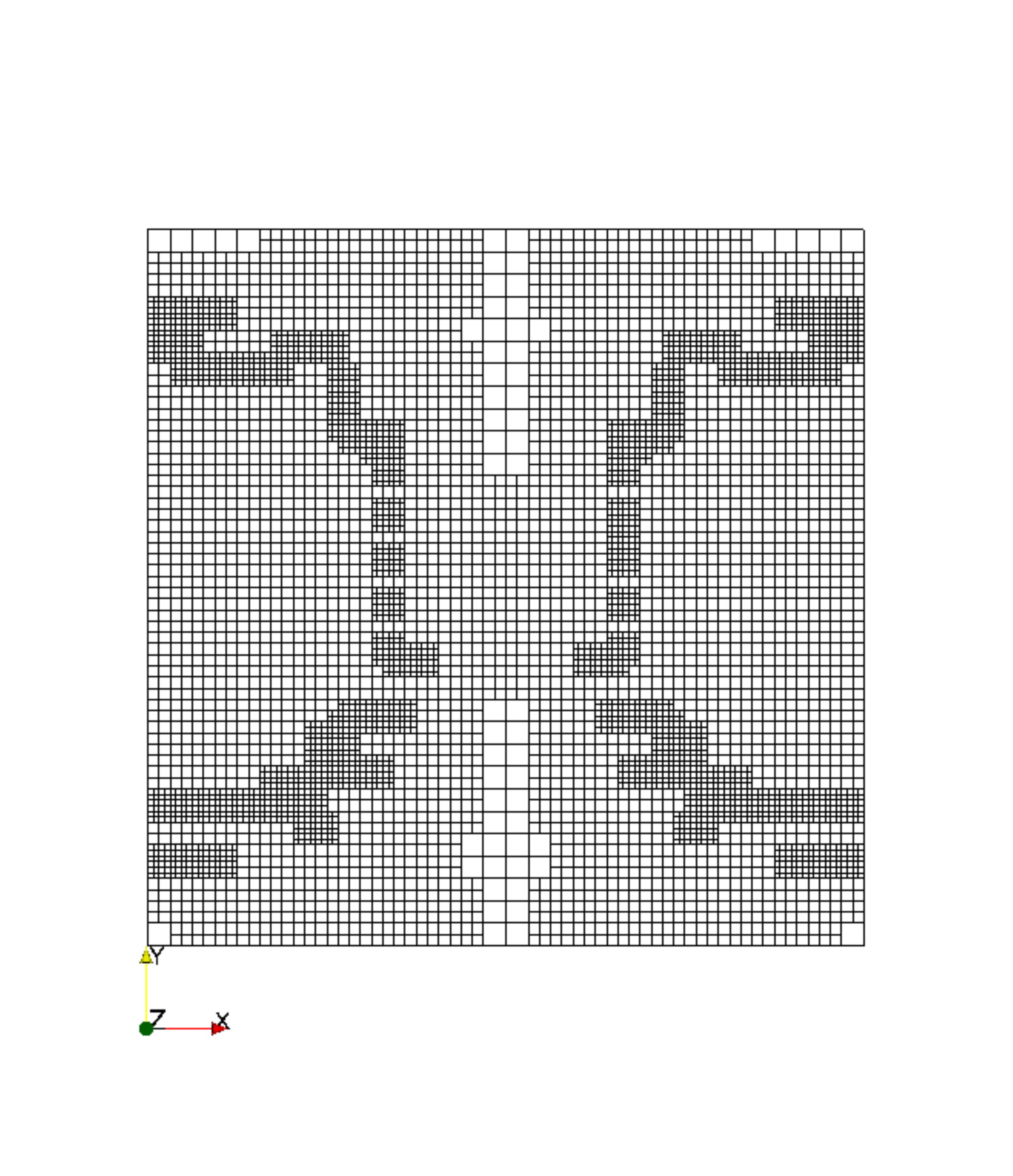}
	}
	}
	\subfloat[{REF 6: \newline
	ref. based on $\| \Delta_x u_h + f - \partial_t u_h\|^2_K$}]{
	\spacetimeaxis{\includegraphics[width=4.5cm, trim={2cm 3cm 2cm 5cm}, clip]{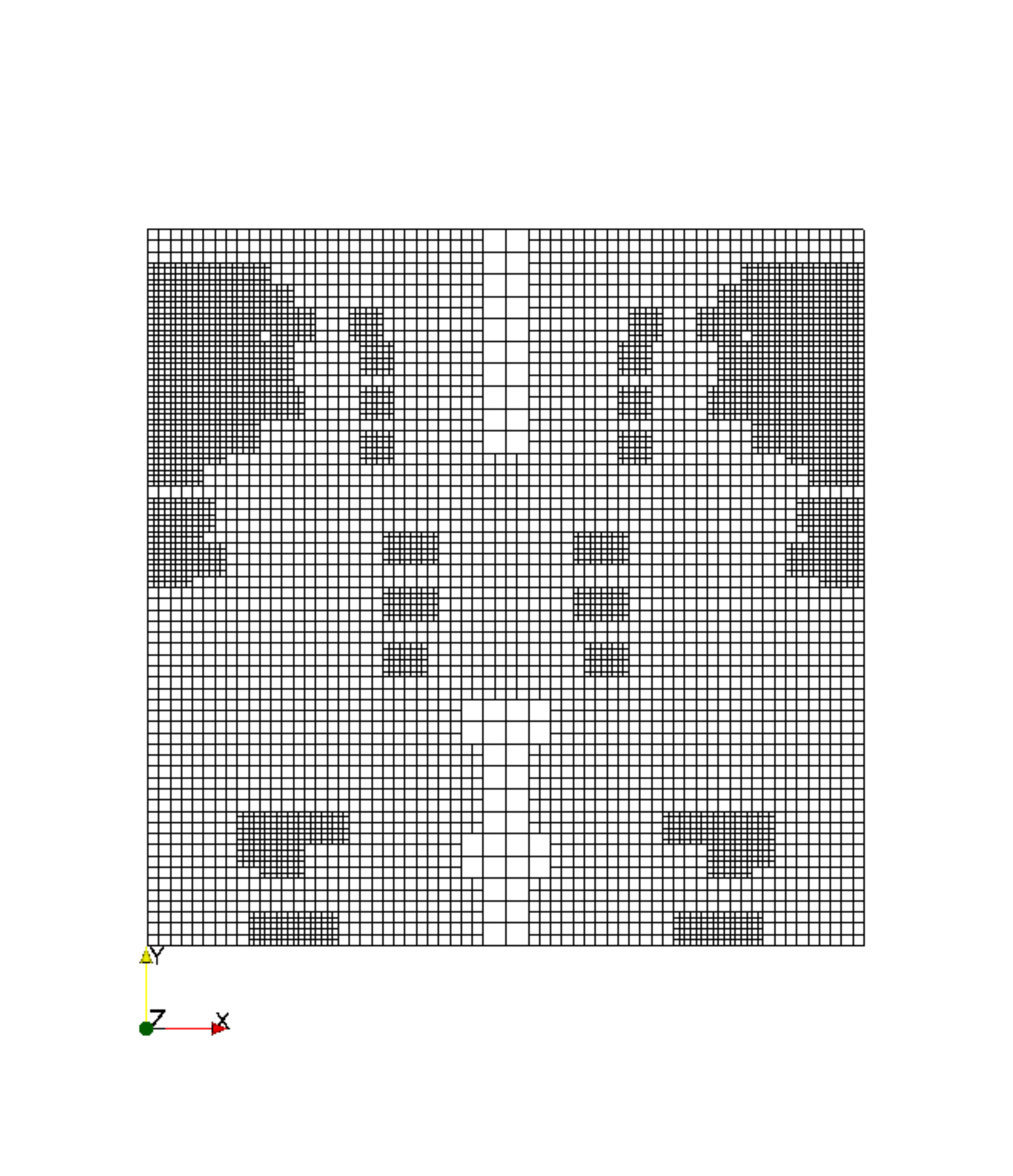}
	}
	}
	\\[-0.5pt]
	\captionsetup[subfigure]{oneside, margin={0.7cm,0cm}}
	\subfloat[{REF 7: \newline
	ref. based on $\| \nabla_x e\|^2_K$}]{
	\spacetimeaxis{\includegraphics[width=4.5cm, trim={2cm 3cm 2cm 5cm}, clip]{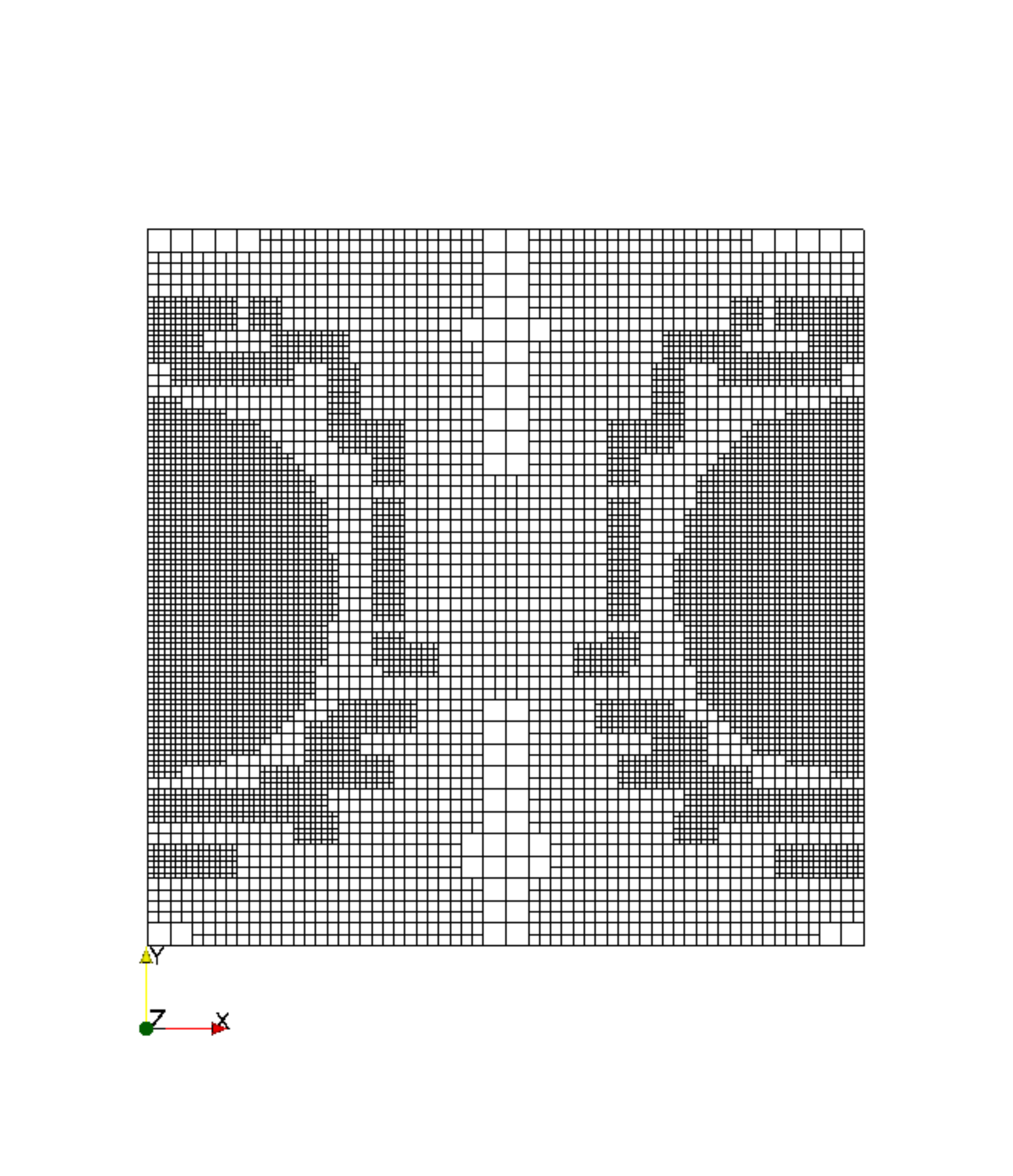}
	}
	}
	\subfloat[{REF 7: \newline
	ref. based on $\overline{\rm m}^{\rm I}_{{\rm d}, K}$ 
	}]{
	\spacetimeaxis{\includegraphics[width=4.5cm, trim={2cm 3cm 2cm 5cm}, clip]{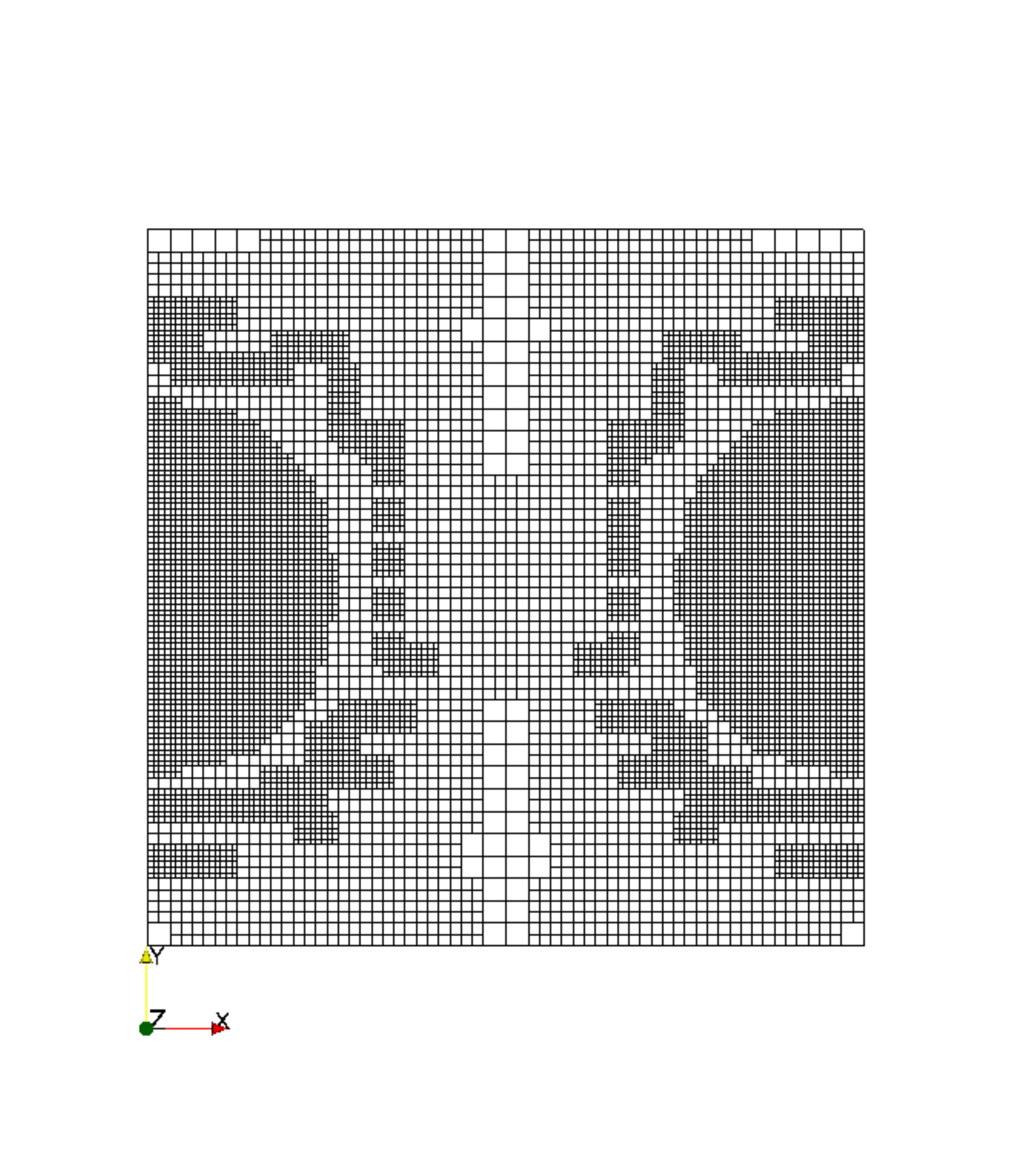}
	}
	}
	\subfloat[{REF 7: \newline
	ref. based on $\| \Delta_x u_h + f - \partial_t u_h\|^2_K$}]{
	\spacetimeaxis{\includegraphics[width=4.5cm, trim={2cm 3cm 2cm 5cm}, clip]{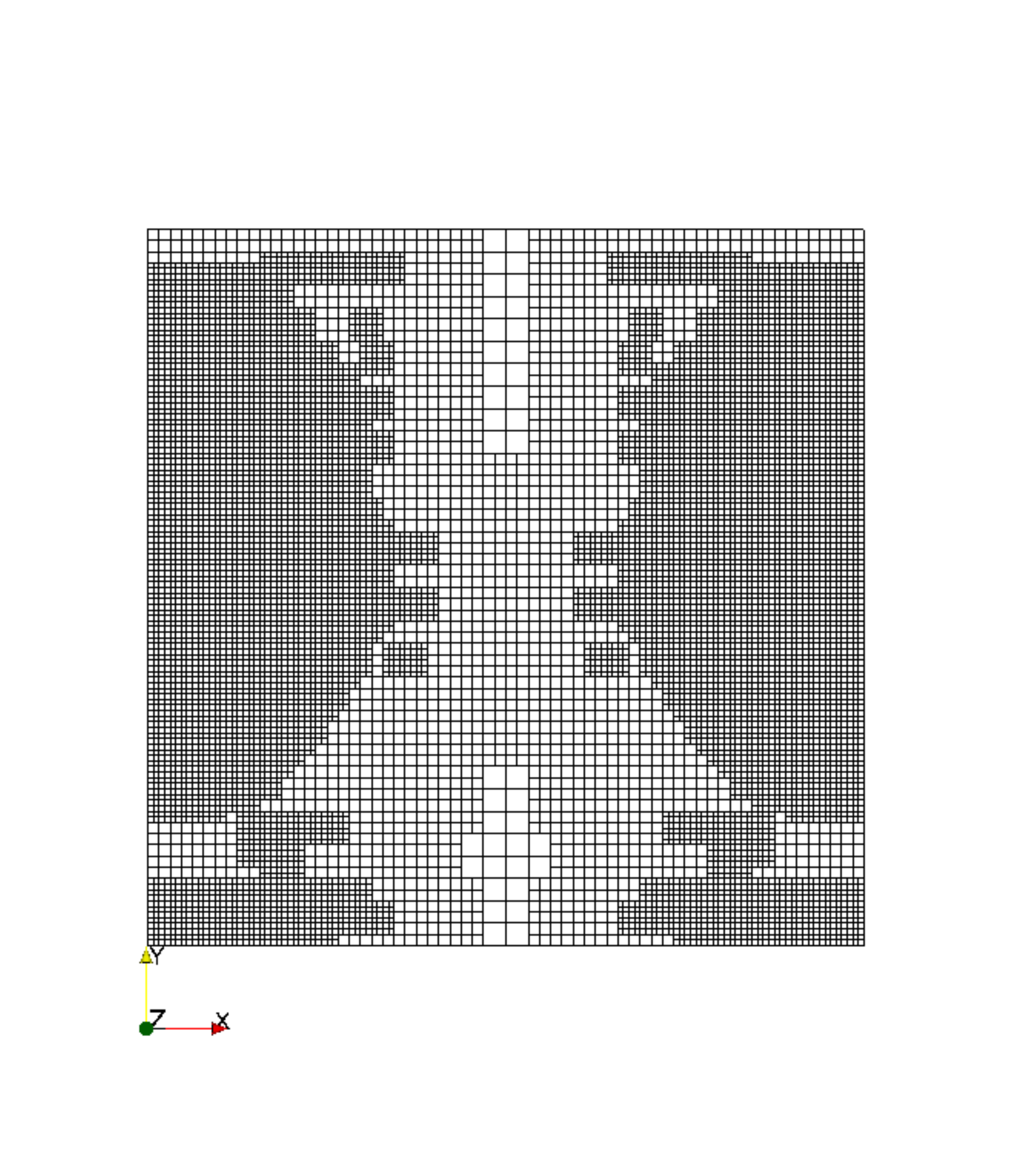}
	}
	}
	\caption{{\em Example 2-1}. 
	Comparison of meshes obtained by the refinement based on local contributions 
	$\| \nabla_x e\|^2_K$, $\overline{\rm m}^{\rm I}_{\rm d, K}$, and 
	residual $\| \Delta u_h + f - \partial_t u_h\|^2_K$, $K \in \mathcal{K}_h$, where
	$\flux_h \in S^{4}_{7h} \oplus S^{4}_{7h}$, and $w_h \in S^{4}_{7h}$,
	(with the marking criterion ${\mathds{M}}_{\rm BULK}(0.6)$).}
	\label{fig:unit-domain-example-3-comparison-meshes-v-2-y-4-adaptive-ref}
\end{figure}
}

Next, let us consider a more complicated case with $k_1 = 6$ and $k_2 = 3$ (see Figure 
\ref{fig:example-3-2-exact-solution-a}). We start with a configuration, where the initial mesh is obtained 
by four global refinements ($N^0_{\rm ref} = 4$), and we proceed with six adaptive steps 
($N_{\rm ref} = 6$) 
using the $\mathds{M}_{\rm BULK}(\sigma)$ marking criterion with $\sigma =0.6$. The obtained efficiency 
indices characterising all majorants and error identity are presented in Table 
\ref{tab:unit-domain-example-3-2-estimates-v-2-y-7-w-7-adaptive-ref}. 
{Here, the auxiliary functions $\flux_h$ and $w_h$ are taken from the approximation spaces 
$S^{7}_{7h} \oplus S^{7}_{7h}$ and $ S^{7}_{7h}$, respectively.} 
We see that the performance of $\overline{\rm M}^{\rm I}_{s, h}$ is identical to {that of}
the majorant $\overline{\rm M}^{\rm I}$, since $\theta$ from the space-time IgA scheme is set to zero 
in this example. 
The numerical performance of the majorant corresponding to the advanced discretisation scheme with parameter 
$\delta_h$ scaled proportionally to the local size of the element $h_K$ will be discussed in the follow-up 
report.

\begin{table}[!t]
\scriptsize
\centering
\newcolumntype{g}{>{\columncolor{gainsboro}}c} 	
\newcolumntype{k}{>{\columncolor{lightgray}}c} 	
\newcolumntype{s}{>{\columncolor{silver}}c} 
\newcolumntype{a}{>{\columncolor{ashgrey}}c}
\newcolumntype{b}{>{\columncolor{battleshipgrey}}c}
\begin{tabular}{c|cga|ck|cb|cc}
\parbox[c]{0.8cm}{\centering \# ref. } & 
\parbox[c]{1.4cm}{\centering  $\| \nabla_x e \|_Q$}   & 	  
\parbox[c]{1.0cm}{\centering $\Ieff (\overline{\rm M}^{\rm I})$ } & 
\parbox[c]{1.4cm}{\centering $\Ieff (\overline{\rm M}^{\rm I\!I})$ } & 
\parbox[c]{1.0cm}{\centering  $|\!|\!|  e |\!|\!|_{s, h}$ }   & 	  
\parbox[c]{1.4cm}{\centering $\Ieff (\overline{\rm M}^{\rm I}_{s, h})$ } & 
\parbox[c]{1.0cm}{\centering  $|\!|\!|  e |\!|\!|_{\mathcal{L}}$ }   & 	  
\parbox[c]{1.4cm}{\centering$\Ieff ({\EI})$ } & 
\parbox[c]{1.2cm}{\centering e.o.c. ($|\!|\!|  e |\!|\!|_{s, h}$)} & 
\parbox[c]{1.2cm}{\centering e.o.c. ($|\!|\!|  e |\!|\!|_{\mathcal{L}}$)} \\
\midrule
   2 &     1.7936e-01 &         1.15 &         1.01 &     1.7936e-01 &         1.15 &     3.2970e+01 &         1.00 &     2.18 &     1.20 \\
   4 &     2.8466e-02 &         1.15 &         1.03 &     2.8466e-02 &         1.15 &     1.3901e+01 &         1.00 &     2.41 &     0.83 \\
   6 &     7.3156e-03 &         1.29 &         1.11 &     7.3156e-03 &         1.29 &     7.2587e+00 &         1.00 &     1.61 &     0.52 \\
   8 &     1.9064e-03 &         2.02 &         1.32 &     1.9064e-03 &         2.02 &     3.6573e+00 &         1.00 &     2.21 &     0.82 \\
\end{tabular}
\caption{{\em Example 2-2}. 
Efficiency of $\overline{\rm M}^{\rm I}$, $\overline{\rm M}^{\rm I\!I}$, $\overline{\rm M}^{\rm I}_{s, h}$, and ${\EI}$ for 
$u_h \in S^{2}_{h}$, $\flux_h \in S^{7}_{7h} \oplus S^{7}_{7h}$ and $w_h \in S^{7}_{7h}$,
w.r.t. adaptive refinement steps (with the marking criterion ${\mathds{M}}_{\rm BULK}(0.6)$).}
\label{tab:unit-domain-example-3-2-estimates-v-2-y-7-w-7-adaptive-ref}
\end{table}

\begin{table}[!t]
\scriptsize
\centering
\newcolumntype{g}{>{\columncolor{gainsboro}}c} 	
\begin{tabular}{c|ccc|cgg|cgg|c}
& \multicolumn{3}{c|}{ d.o.f. } 
& \multicolumn{3}{c|}{ $t_{\rm as}$ }
& \multicolumn{3}{c|}{ $t_{\rm sol}$ } 
& $\tfrac{t_{\rm appr.}}{t_{\rm er.est.}}$ \\
\midrule
\parbox[c]{0.8cm}{\centering \# ref. } & 
\parbox[c]{0.8cm}{\centering $u_h$ } &  
\parbox[c]{0.6cm}{\centering $\flux_h$ } &  
\parbox[c]{0.6cm}{\centering $w_h$ } & 
\parbox[c]{1.4cm}{\centering $u_h$ } & 
\parbox[c]{1.4cm}{\centering $\flux_h$ } & 
\parbox[c]{1.4cm}{\centering $w_h$ } & 
\parbox[c]{1.4cm}{\centering $u_h$ } & 
\parbox[c]{1.4cm}{\centering $\flux_h$ } & 
\parbox[c]{1.4cm}{\centering $w_h$ } & \\
%
%
\midrule
   2 &        940 &       1058 &        529 &   7.96e-01 &   2.03e+01 &   1.02e+01 &         2.41e-02 &         5.16e-02 &         5.92e-02 & 0.03\\
   4 &       5590 &       1058 &        529 &   5.64e+00 &   1.80e+01 &   1.01e+01 &         4.56e-01 &         4.63e-02 &         5.65e-02 & 0.22\\
   6 &      21089 &       1058 &        529 &   3.14e+01 &   1.81e+01 &   1.01e+01 &         3.70e+00 &         3.12e-02 &         5.51e-02 & 1.24 \\
   8 &      67043 &       1058 &        529 &   2.06e+02 &   1.82e+01 &   1.23e+01 &         1.64e+01 &         4.59e-02 &         7.40e-02 & 7.26\\
     \midrule
    &       &         &    &
    \multicolumn{3}{c|}{ $t_{\rm as} (u_h)$ \quad : \quad $t_{\rm as} (\flux_h)$ \quad : \qquad $t_{\rm as} (w_h)$ } &      
    \multicolumn{3}{c|}{\; $t_{\rm sol} (u_h)$ \, : \quad $t_{\rm sol} (\flux_h)$  \quad:  \qquad  $t_{\rm sol} (w_h)$\;}& 
     \\
 \midrule
%
         &       &       &       &      16.76 &       1.48 &       1.00 &           221.53 &             0.62 &             1.00 & \\
\end{tabular}
\caption{{\em Example 2-2}. 
Assembling and solving time (in seconds) spent for the systems generating
d.o.f. of $u_h \in S^{2}_{h}$, $\flux_h \in S^{7}_{7h} \oplus S^{7}_{7h}$, and $w_h \in S^{7}_{7h}$
w.r.t. adaptive refinement steps (with the marking criterion ${\mathds{M}}_{\rm BULK}(0.6)$). }
\label{tab:unit-domain-example-3-2-times-v-2-y-9-adaptive-ref}
\end{table}
 
\begin{figure}[!t]
	\centering
	\captionsetup[subfigure]{oneside, margin={0.3cm,0cm}}

	\subfloat[{REF 3: \newline
	ref. based on ${\rm e}_{\rm d}$}]{
	\spacetimeaxis{\includegraphics[width=4.2cm, trim={8.1cm 2cm 5cm 2cm}, clip]{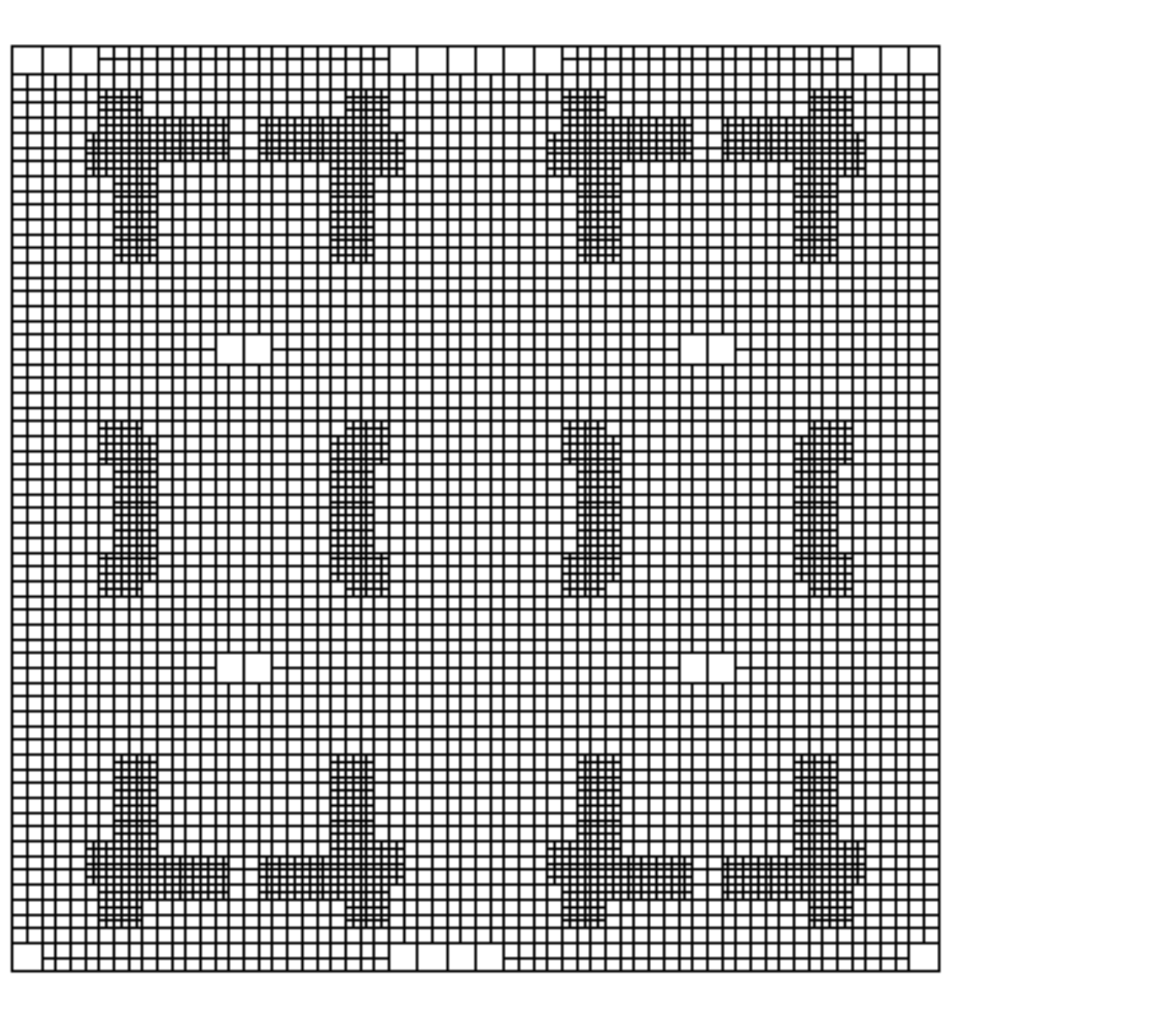}}
	}
	\hskip -30pt
	\subfloat[{REF 3: \newline
	ref. based on $\overline{\rm m}^{\rm I}_{\rm d}$}]{
	\spacetimeaxis{\includegraphics[width=4.2cm, trim={8.1cm 2cm 5cm 2cm}, clip]{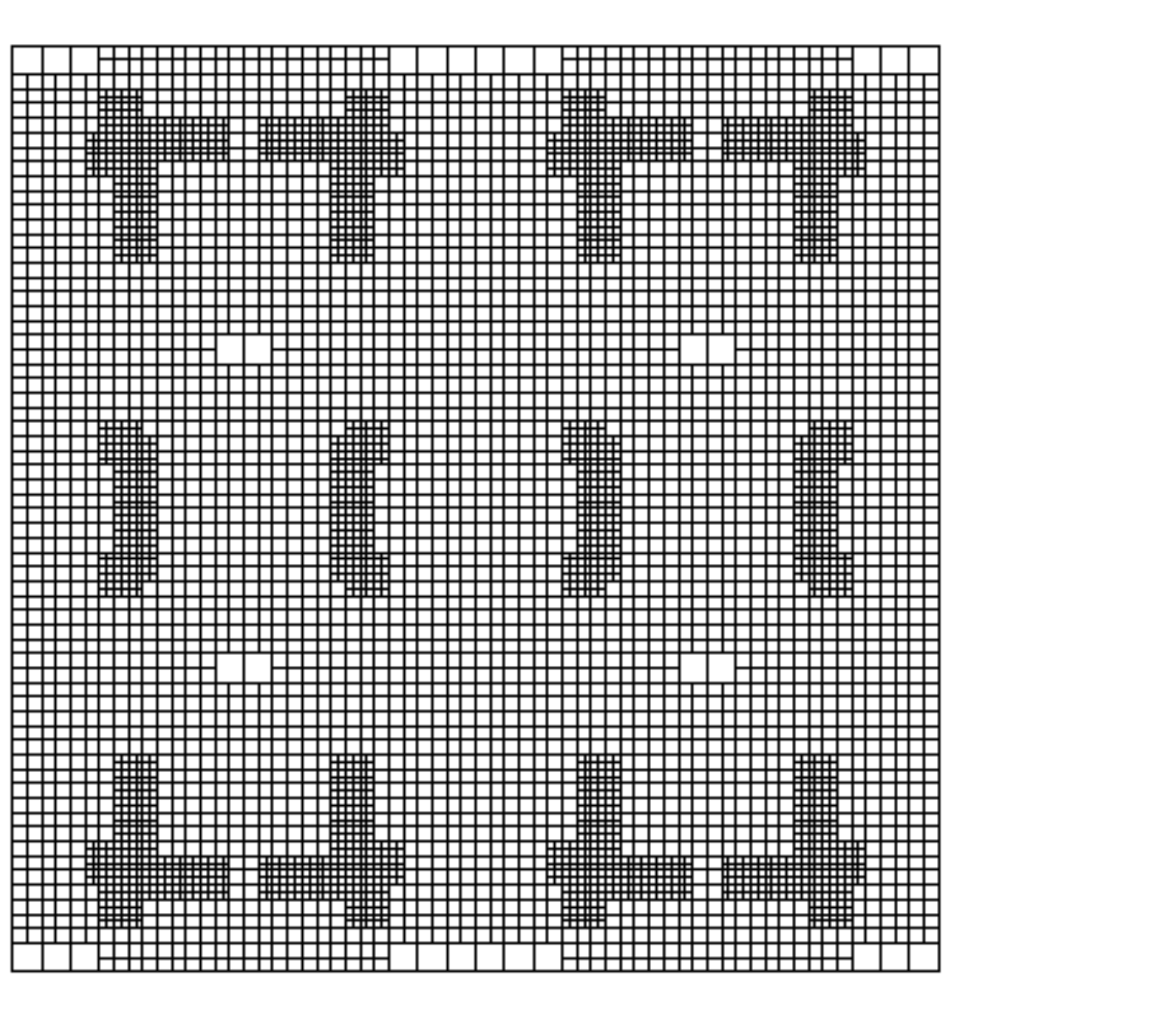}}
	}
	\hskip -30pt
	\subfloat[{REF 3: \newline
	ref. based on $|\!|\!|  e |\!|\!|_{\mathcal{L}}$}]{
	\spacetimeaxis{\includegraphics[width=4.2cm, trim={8.1cm 2cm 5cm 2cm}, clip]{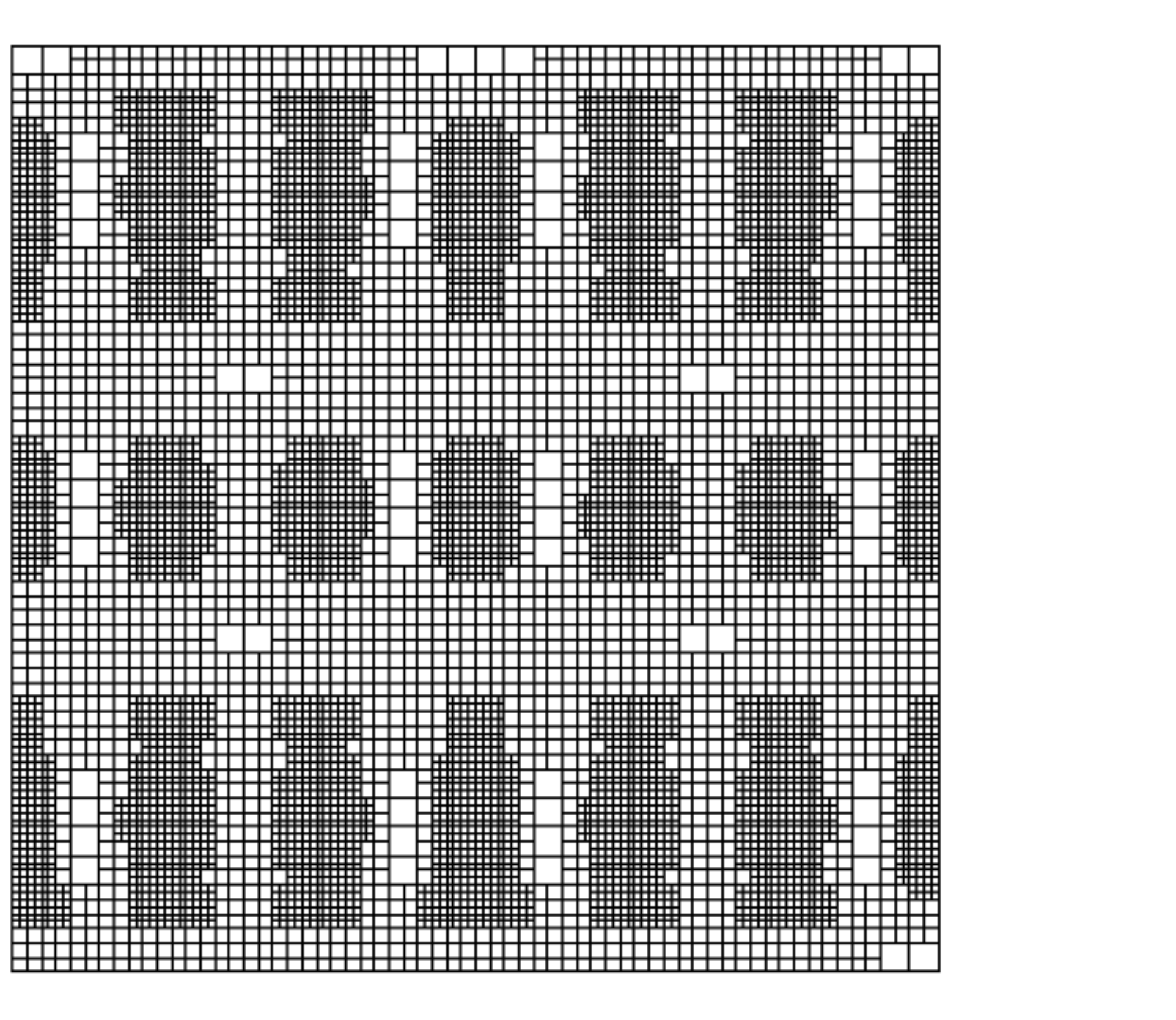}}
	}
	\hskip -30pt
	\subfloat[{REF 3: \newline
	ref. based on ${\EI}$}]{
	\spacetimeaxis{\includegraphics[width=4.2cm, trim={8.1cm 2cm 5cm 2cm}, clip]{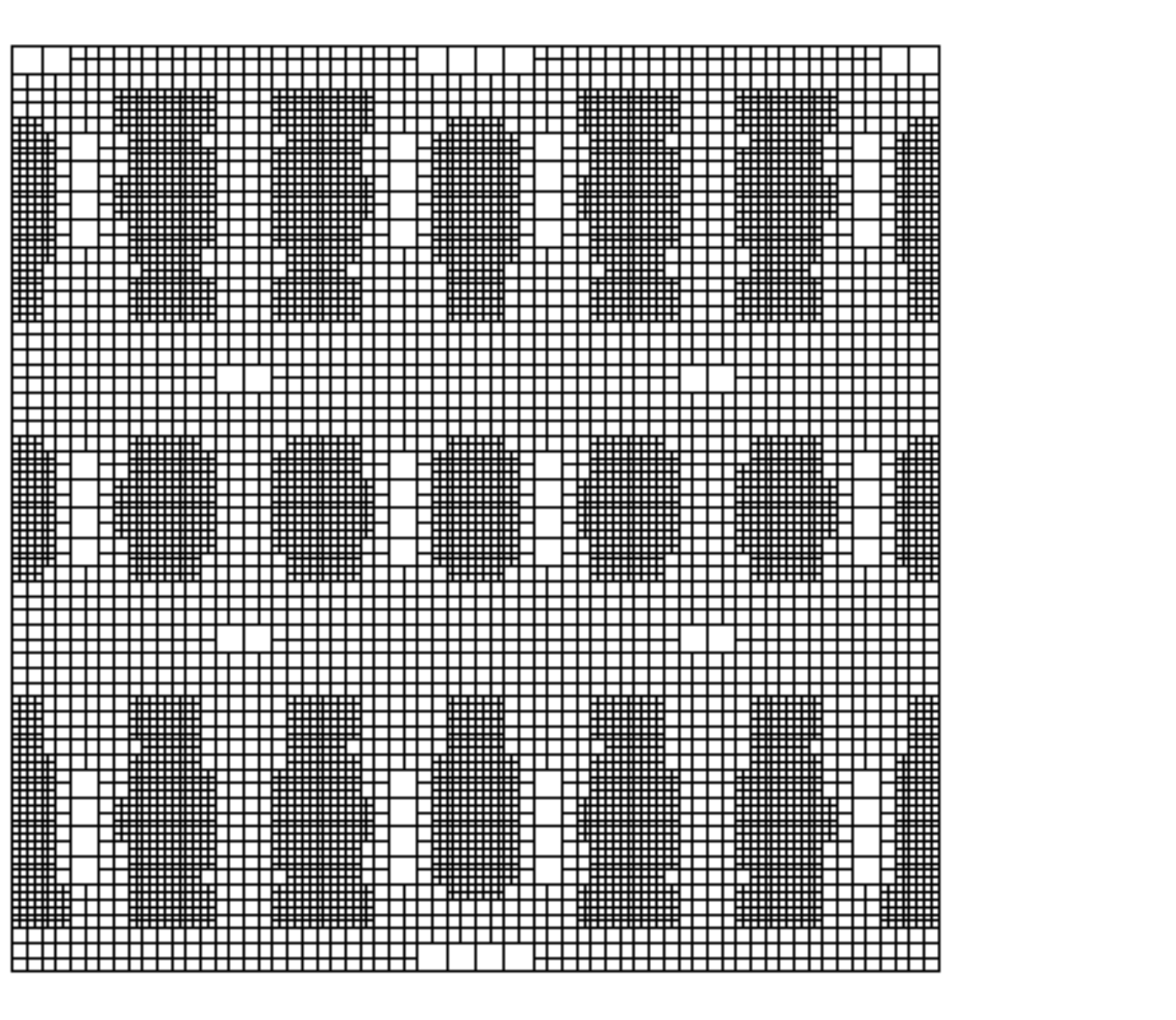}}
	}
	\\[-0.5pt]
	\subfloat[{REF 4: \newline
	ref. based on ${\rm e}_{\rm d}$}]{
	\spacetimeaxis{\includegraphics[width=4.2cm, trim={8.1cm 2cm 5cm 2cm}, clip]{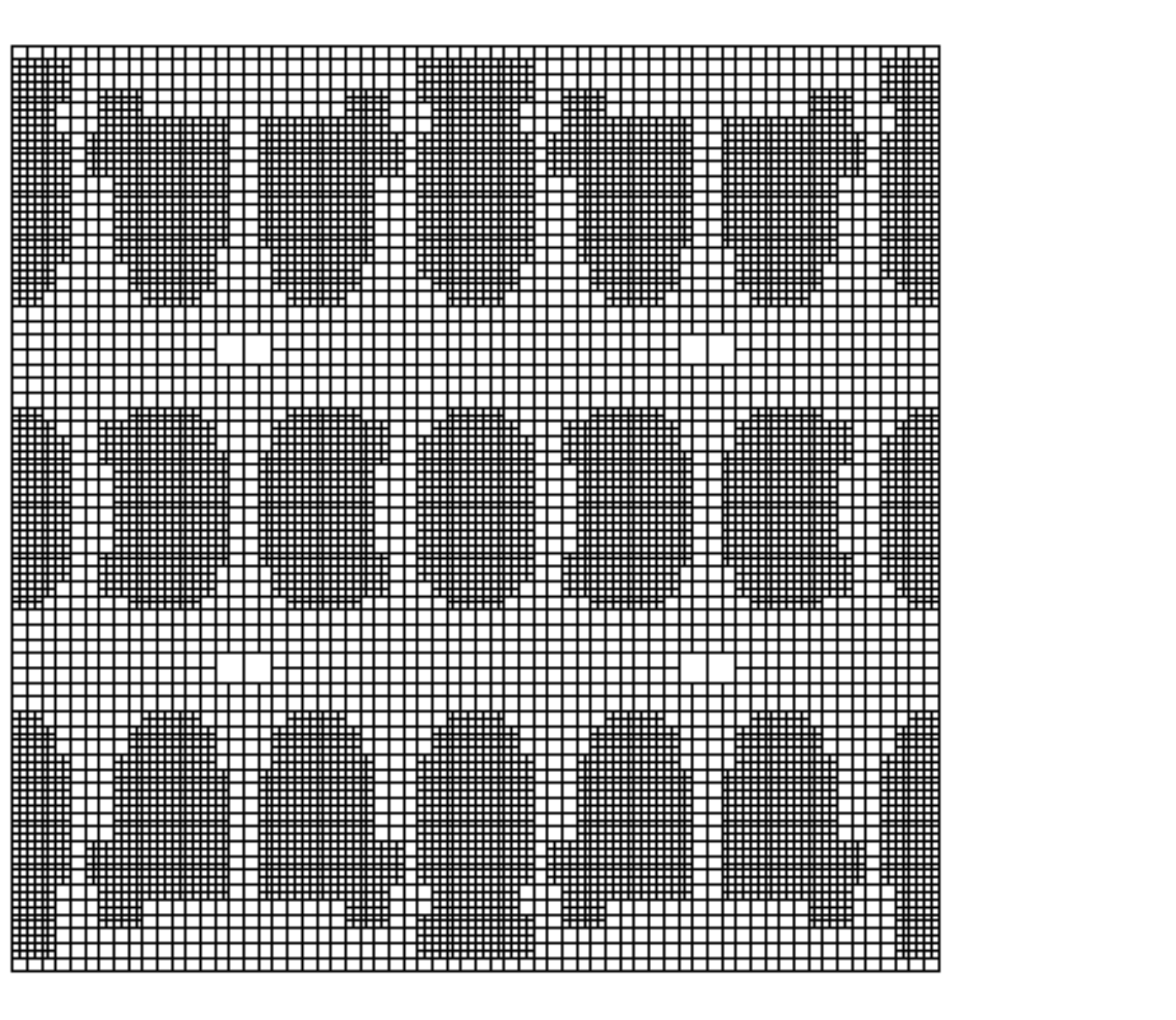}}
	}
	\hskip -30pt
	\subfloat[{REF 4: \newline
	ref. based on $\overline{\rm m}^{\rm I}_{\rm d}$}]{
	\spacetimeaxis{\includegraphics[width=4.2cm, trim={8.1cm 2cm 5cm 2cm}, clip]{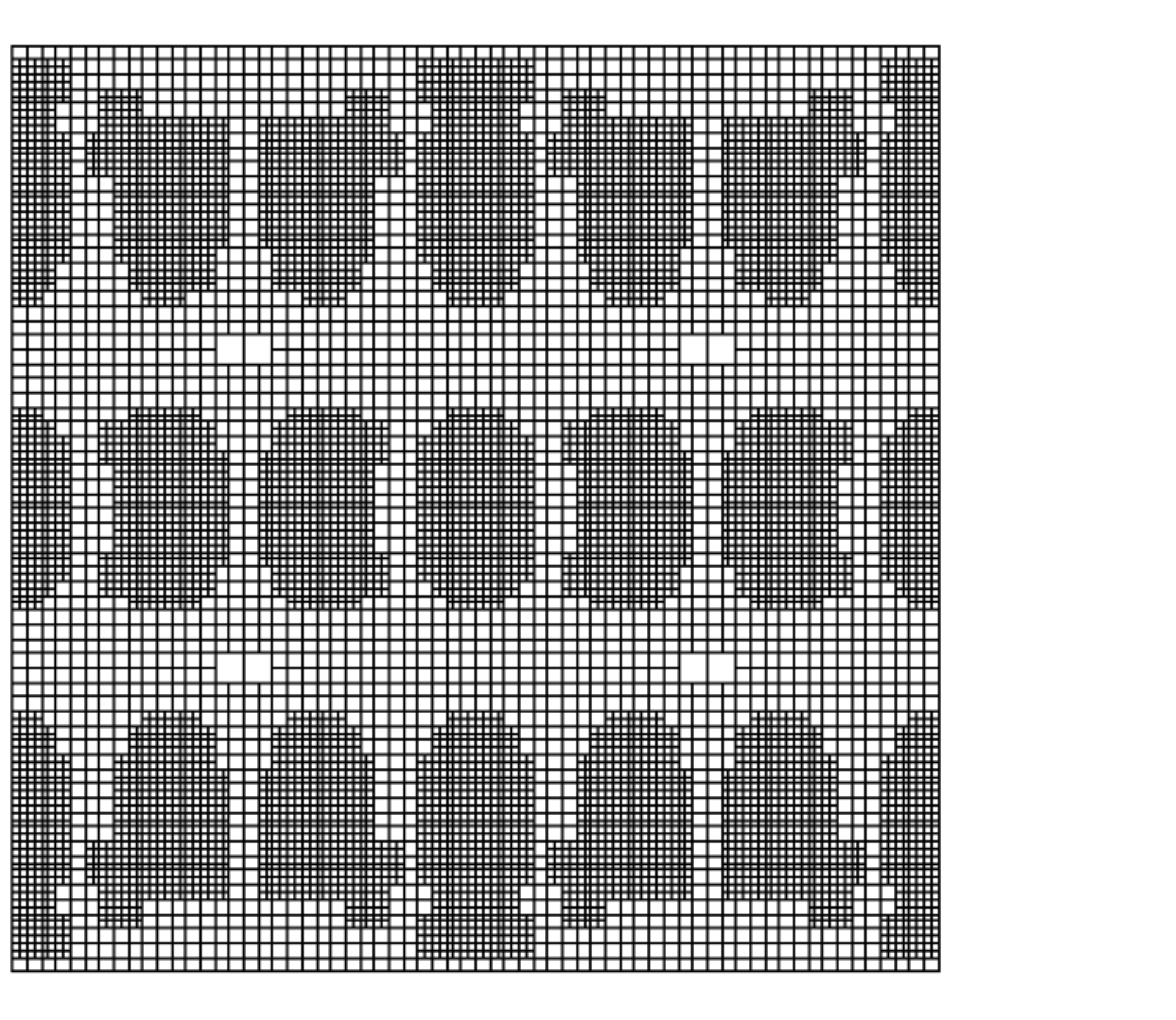}}
	}
	\hskip -30pt
	\subfloat[{REF 4: \newline ref. based on $|\!|\!|  e |\!|\!|_{\mathcal{L}}$}]{
	\spacetimeaxis{\includegraphics[width=4.2cm, trim={8.1cm 2cm 5cm 2cm}, clip]{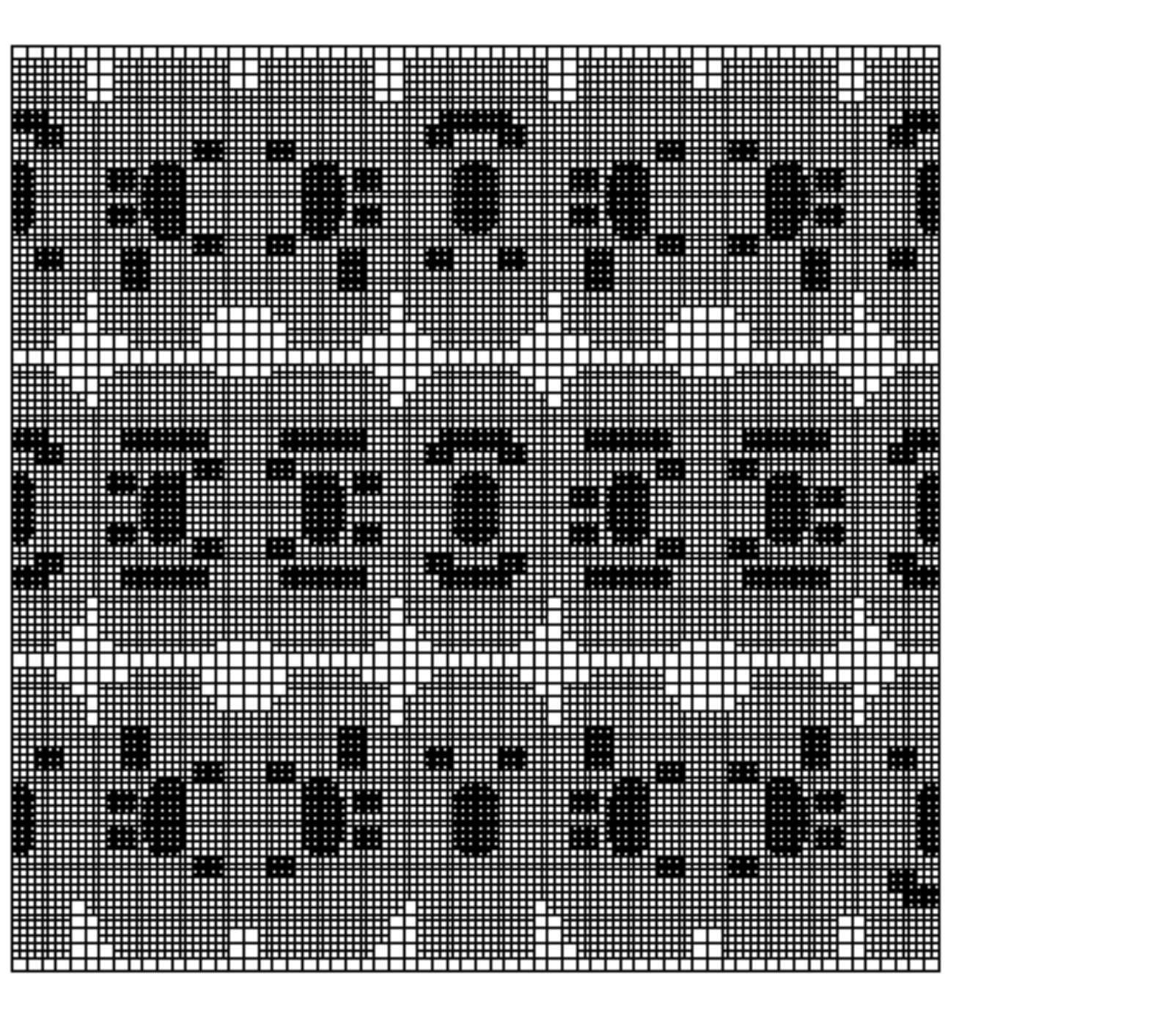}}
	}
	\hskip -30pt
	\subfloat[{REF 4: \newline
	ref. based on ${\EI}$}]{
	\spacetimeaxis{\includegraphics[width=4.2cm, trim={8.1cm 2cm 5cm 2cm}, clip]{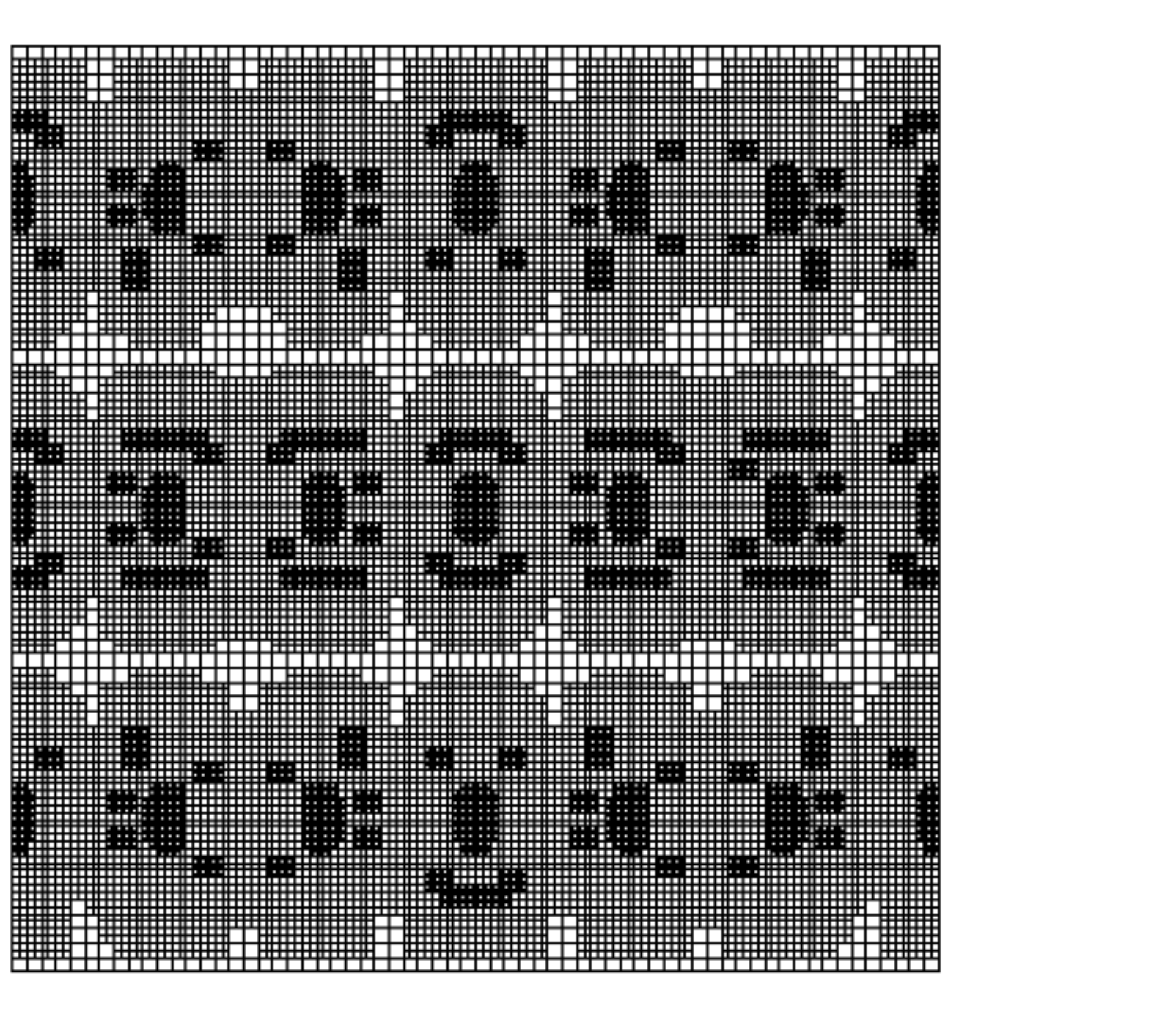}}
	}
	\caption{{\em Example 2-2}. 
	Comparison of the meshes obtained by refinement based on ${\rm e}_{\rm d}$, $\overline{\rm m}^{\rm I}_{\rm d}$, 
	$|\!|\!|  e |\!|\!|_{\mathcal{L}}$, and ${\EI}$ for 
	$u_h \in S^{2}_{h}$, $\flux_h \in S^{7}_{7h} \oplus S^{7}_{7h}$, and $w_h \in S^{7}_{7h}$,
	(with the marking criterion ${\mathds{M}}_{\rm BULK}(0.6)$).}
	\label{fig:unit-domain-example-3-2-meshes-v-2-y-9-adaptive-ref}
\end{figure}

\begin{figure}[!t]
	\centering
	\captionsetup[subfigure]{oneside, labelformat=empty}
	\subfloat[]{
	\includegraphics[width=3.8cm, trim={0cm 0cm 0cm 0cm}, clip]{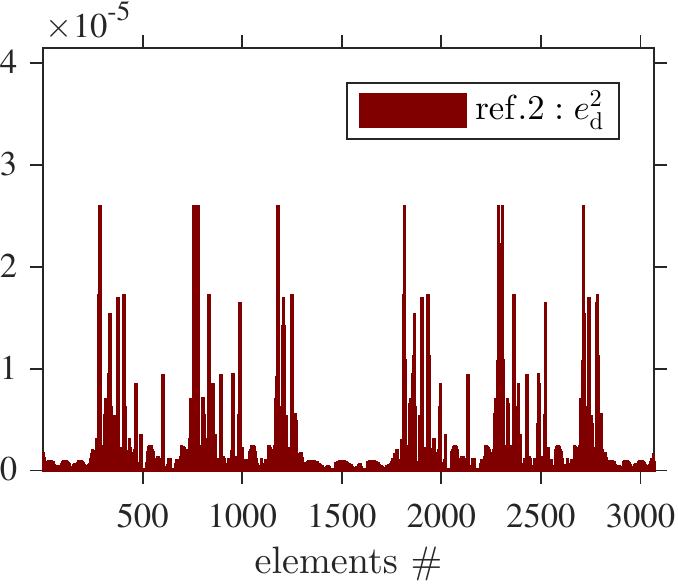}
	}
	\hskip -5pt
	\subfloat[]{
	\includegraphics[width=4.3cm, trim={0cm 0cm 0cm 0cm}, clip]{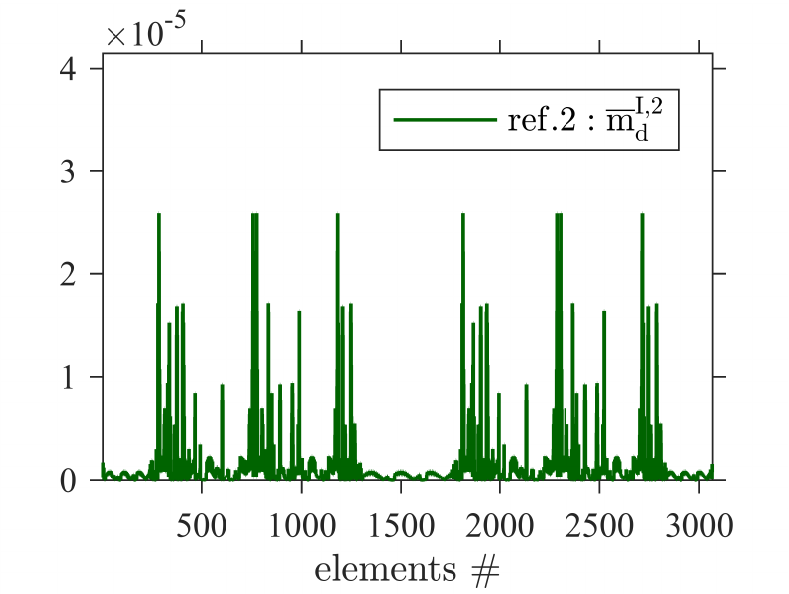}
	}
	\hskip -5pt
	\subfloat[]{
	\includegraphics[width=3.8cm, trim={0cm 0cm 0cm 0cm}, clip]{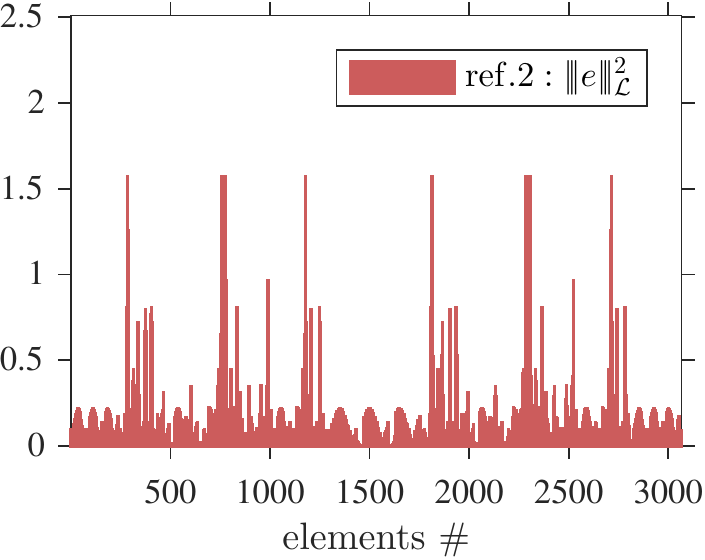}
	}
	\hskip -5pt
	\subfloat[]{
	\includegraphics[width=4.3cm, trim={0cm 0cm 0cm 0cm}, clip]{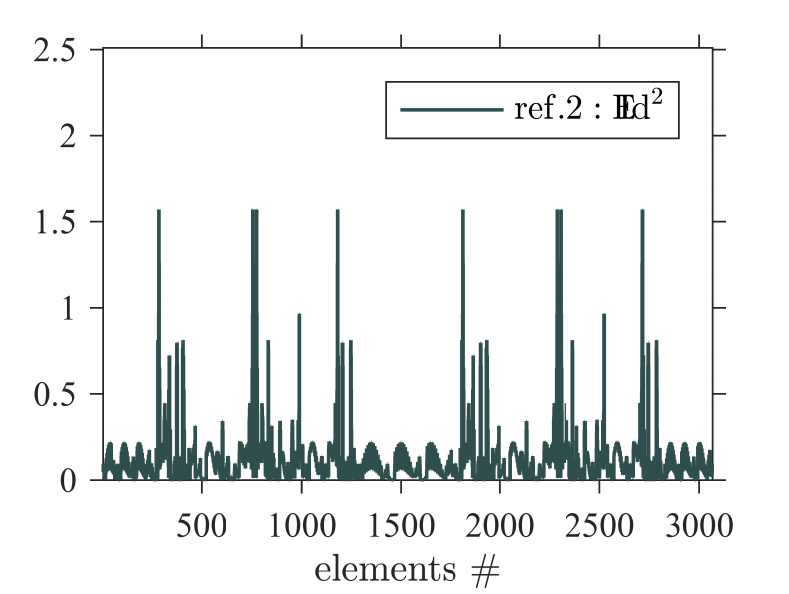}
	}
	\vskip -15pt
	\caption{{\em Example 2-2}. 
	Distribution of $e_{{\rm d}, K}$ and $\overline{\rm m}_{{\rm d}, K}$ as well as 
	$|\!|\!|  e |\!|\!|_{\mathcal{L}, K}$ and $\EI_K$ on the refinement step $2$.}
	\label{fig:unit-domain-example-3-2-edmd-e-sol-e-id-distribution}
\end{figure}

A comparison of the meshes corresponding to different refinement criteria is presented in
Figure \ref{fig:unit-domain-example-3-2-meshes-v-2-y-9-adaptive-ref}. 
The first two columns contain the meshes produced by the refinement based on ${\rm e}_{\rm d}$ and 
$\overline{\rm m}^{\rm I}_{\rm d}$, whereas the third and fourth columns correspond to the adaptive meshes { 
obtained on the steps 3 and 4 using local distributions of $|\!|\!|  e |\!|\!|_{\mathcal{L}, K}$ and ${\EI}_K$ 
for the refinement criterion.}
It is clear from the plots that the meshes related to ${\EI}$-based refinement are denser than the 
meshes in the  second column. 
{ 
Nevertheless, the error identity suggests similar areas of the mesh refinement. 
Therefore, it can be used as effectively as $\overline{\rm m}^{\rm I}_{\rm d}$ for mesh-adaptation. 
}
Moreover, Figure \ref{fig:unit-domain-example-3-2-edmd-e-sol-e-id-distribution} 
confirms that $\overline{\rm m}^{\rm I}_{{\rm d}, K}$ and ${\EI}_K$ are quantitatively sharp when it comes 
to estimating ${\rm e}_{{\rm d}, K}$ and  $|\!|\!|  e |\!|\!|_{\mathcal{L}, K}$, respectively.

\subsection{Example 3}
\label{ex:unit-domain-example-6}
\begin{figure}[!t]
	\centering
	\subfloat[]{
	\includegraphics[scale=0.6]{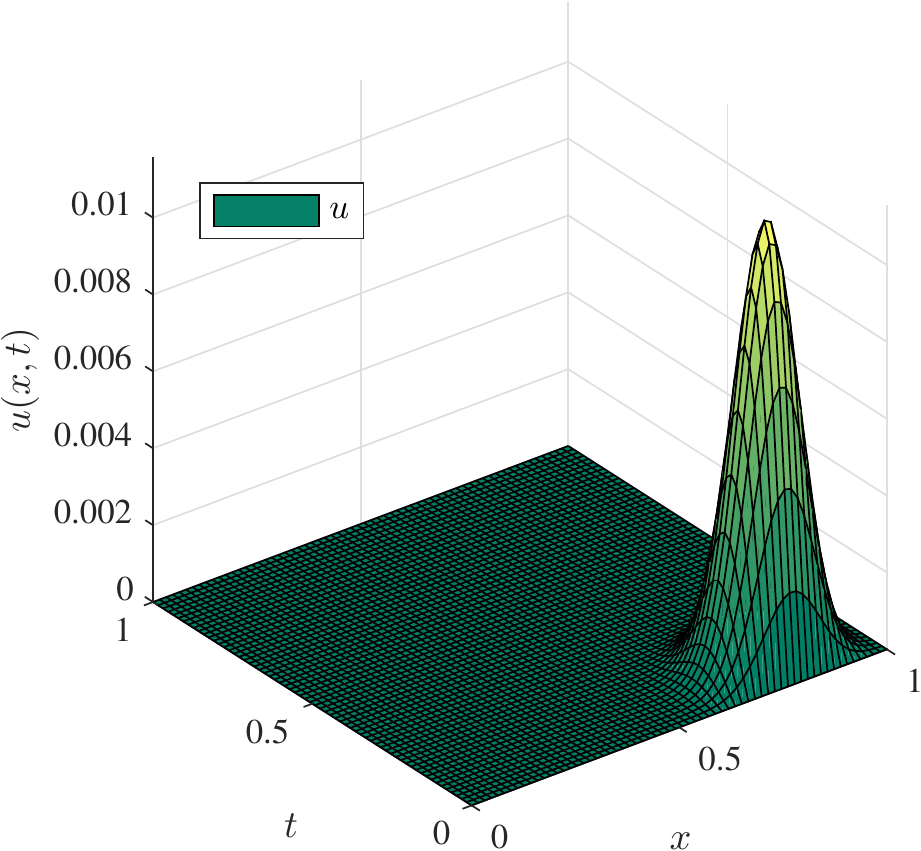}
	\label{fig:example-6-exact-solution-1}
	}
	\qquad
	\subfloat[]{
	\includegraphics[scale=0.6]{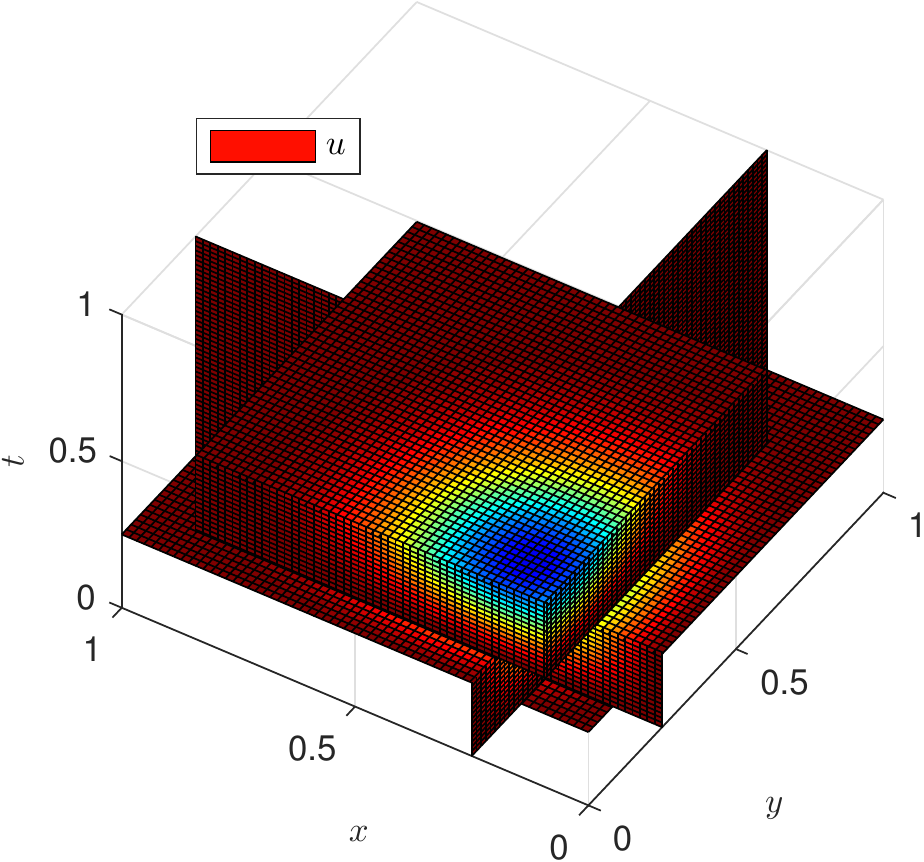}
	\label{fig:example-6-exact-solution-2}	}
	\caption{{\em Example 3}. 
	(a) Exact solution $u = x \, (x -  1) \, t \, (t - 1) \, e^{-100 \,|(x, t) - (0.8, 0.05)|}$. 
	{ 
	(b) Exact solution \\
	$u = x \, (x -  1) \, y \, (y -  1) \,  t \, (t - 1) \, e^{-100 \,|(x, y, t) - (0.25, 0.25, 0.25)|}$}.}
	\label{fig:example-6-exact-solution}
\end{figure}

\rm 
As another standard test case, we consider an example with a sharp local 
Gaussian {peak}
in the exact solution. Let $Q := (0, 1)^2$, and the solution to be defined by
\begin{alignat*}{3}
u(x, t) 	& = (x^2 - x) \, (t^2 - t) \, e^{-100 \,|(x, t) - (0.8, 0.05)|},  
& \quad  (x, t) \in \overline{Q},
\label{eq:example-8-exact-solution}
\end{alignat*}
where the  {  peak}
is located in the point $(x, t) = (0.8, 0.05)$, 
see Figure \ref{fig:example-6-exact-solution-1}. 
Then $f$ is computed by substituting $u$ into \eqref{eq:equation}. 
The Dirichlet boundary conditions are obviously homogeneous.
For this example, we consider only an adaptive refinement {procedure.}
For the discretisation spaces, we 
use our standard setting, i.e., $u_h \in S^{2}_h$ for the primal variable, as well as 
$\flux_h \in S^{3}_h \oplus S^{3}_h$ and $w_h \in S^{3}_h$ for the auxiliary functions. 
We start with four initial global refinements ($N^0_{\rm ref} = 4$), and continue with 
seven adaptive steps ($N_{\rm ref} = 7$). As marking criteria, we choose 
$\mathds{M}_{\rm BULK}(\sigma)$ with bulk parameter $\sigma =0.6$.
%

The analysis of the quantitative efficiency of the majorants and the error identity in terms of error 
estimation 
is provided in Table \ref{tab:unit-domain-example-6-estimates-v-2-y-3-adaptive-ref}, which confirms that 
$\overline{\rm M}^{\rm I\!I}$ is twice as sharp as that of $\overline{\rm M}^{\rm I}$, 
and the error identity sharply predicts 
the exact error $|\!|\!|  e |\!|\!|_{\mathcal{L}}$. When it comes to the time efficiency summarised 
in Table \ref{tab:unit-domain-example-6-times-v-2-y-3-adaptive-ref}, we see that the assembling of 
matrices for the $\flux_h$ and $w_h$ is 3 - 4 times more time-consuming. 
{ 
However, the solution of the corresponding systems is 1.4 - 3 times faster.
}
%
\begin{table}[!t]
\scriptsize
\centering
\newcolumntype{g}{>{\columncolor{gainsboro}}c} 	
\newcolumntype{k}{>{\columncolor{lightgray}}c} 	
\newcolumntype{s}{>{\columncolor{silver}}c} 
\newcolumntype{a}{>{\columncolor{ashgrey}}c}
\newcolumntype{b}{>{\columncolor{battleshipgrey}}c}
\begin{tabular}{c|cga|ck|cb|cc}
\parbox[c]{0.8cm}{\centering \# ref. } & 
\parbox[c]{1.4cm}{\centering  $\| \nabla_x e \|_Q$}   & 	  
\parbox[c]{1.0cm}{\centering $\Ieff (\overline{\rm M}^{\rm I})$ } & 
\parbox[c]{1.4cm}{\centering $\Ieff (\overline{\rm M}^{\rm I\!I})$ } & 
\parbox[c]{1.0cm}{\centering  $|\!|\!|  e |\!|\!|_{s, h}$ }   & 	  
\parbox[c]{1.4cm}{\centering $\Ieff (\overline{\rm M}^{\rm I}_{s, h})$ } & 
\parbox[c]{1.0cm}{\centering  $|\!|\!|  e |\!|\!|_{\mathcal{L}}$ }   & 	  
\parbox[c]{1.4cm}{\centering$\Ieff ({\EI})$ } & 
\parbox[c]{1.2cm}{\centering e.o.c. ($|\!|\!|  e |\!|\!|_{s, h}$)} & 
\parbox[c]{1.2cm}{\centering e.o.c. ($|\!|\!|  e |\!|\!|_{\mathcal{L}}$)} \\
\midrule
%
%
   3 &     6.1181e-02 &         1.58 &         1.24 &     6.1181e-02 &         1.58 &     1.8482e+01 &         1.00 &     2.00 &     1.07 \\
   5 &     1.4192e-02 &         1.32 &         1.12 &     1.4192e-02 &         1.32 &     9.0808e+00 &         1.00 &     2.77 &     1.75 \\
   7 &     4.1402e-03 &         1.27 &         1.07 &     4.1402e-03 &         1.27 &     4.8937e+00 &         1.00 &     2.46 &     1.71 \\
\end{tabular}
\caption{{\em Example 3}. 
Efficiency of $\overline{\rm M}^{\rm I}$, $\overline{\rm M}^{\rm I\!I}$, 
$\overline{\rm M}^{\rm I}_{s, h}$, and ${\EI}$ for $u_h \in S^{2}_{h}$,
$\flux_h \in S_{h}^{3} \oplus S_{h}^{3}$, and $w_h \in S^{3}_{h}$, 
w.r.t. adaptive refinements (with the marking criterion ${\mathds{M}}_{\rm BULK}(0.6)$).}
\label{tab:unit-domain-example-6-estimates-v-2-y-3-adaptive-ref}
\end{table}

\begin{table}[!t]
\scriptsize
\centering
\newcolumntype{g}{>{\columncolor{gainsboro}}c} 	
\begin{tabular}{c|ccc|cgg|cgg|c}
& \multicolumn{3}{c|}{ d.o.f. } 
& \multicolumn{3}{c|}{ $t_{\rm as}$ }
& \multicolumn{3}{c|}{ $t_{\rm sol}$ } 
& $\tfrac{t_{\rm appr.}}{t_{\rm er.est.}}$ \\
\midrule
\parbox[c]{0.8cm}{\centering \# ref. } & 
\parbox[c]{0.8cm}{\centering $u_h$ } &  
\parbox[c]{0.6cm}{\centering $\flux_h$ } &  
\parbox[c]{0.6cm}{\centering $w_h$ } & 
\parbox[c]{1.4cm}{\centering $u_h$ } & 
\parbox[c]{1.4cm}{\centering $\flux_h$ } & 
\parbox[c]{1.4cm}{\centering $w_h$ } & 
\parbox[c]{1.4cm}{\centering $u_h$ } & 
\parbox[c]{1.4cm}{\centering $\flux_h$ } & 
\parbox[c]{1.4cm}{\centering $w_h$ } \\
\midrule
  1 &        324 &        722 &        361 &   1.90e-01 &   1.01e+00 &   6.23e-01 &         3.43e-03 &         2.54e-03 &         1.03e-02 &             0.12 \\
   3 &        543 &       1142 &        571 &   4.31e-01 &   2.95e+00 &   1.61e+00 &         1.33e-02 &         1.65e-02 &         4.06e-02 &             0.10 \\
   5 &       1167 &       2342 &       1171 &   9.20e-01 &   7.79e+00 &   3.90e+00 &         5.62e-02 &         8.07e-02 &         1.71e-01 &             0.08 \\
   7 &       4006 &       7742 &       3871 &   5.43e+00 &   3.01e+01 &   1.83e+01 &         7.43e-01 &         5.12e-01 &         1.42e+00 &             0.12 \\
      \midrule
    &       &         &    &
    \multicolumn{3}{c|}{ $t_{\rm as} (u_h)$ \quad : \quad $t_{\rm as} (\flux_h)$ \quad : \qquad $t_{\rm as} (w_h)$ } &      
    \multicolumn{3}{c|}{\; $t_{\rm sol} (u_h)$ \, : \quad $t_{\rm sol} (\flux_h)$  \quad:  \qquad  $t_{\rm sol} (w_h)$\;} & 
    \\
 \midrule
%
         &       &       &       &       0.30 &       1.65 &       1.00 &             0.52 &             0.36 &             1.00 &             0.12 \\
\end{tabular}
\caption{{\em Example 3}. 
Assembling and solving time (in seconds) spent for the systems generating
d.o.f. of $u_h \in S^{2}_{h}$, $\flux_h \in S^{3}_{h} \oplus S^{3}_{h}$, and $w_h \in S^{3}_{h}$ 
w.r.t. adaptive refinements (with the marking criterion ${\mathds{M}}_{\rm BULK}(0.6)$).}
\label{tab:unit-domain-example-6-times-v-2-y-3-adaptive-ref}
\end{table}

\begin{figure}[!t]
	\centering
	\captionsetup[subfigure]{oneside, margin={0.3cm,0cm}}
	\subfloat[{REF 4: \newline
	ref. based on ${\rm e}_{\rm d}$}]{
	\spacetimeaxis{\includegraphics[width=4.2cm, trim={8.1cm 2cm 5cm 2cm}, clip]{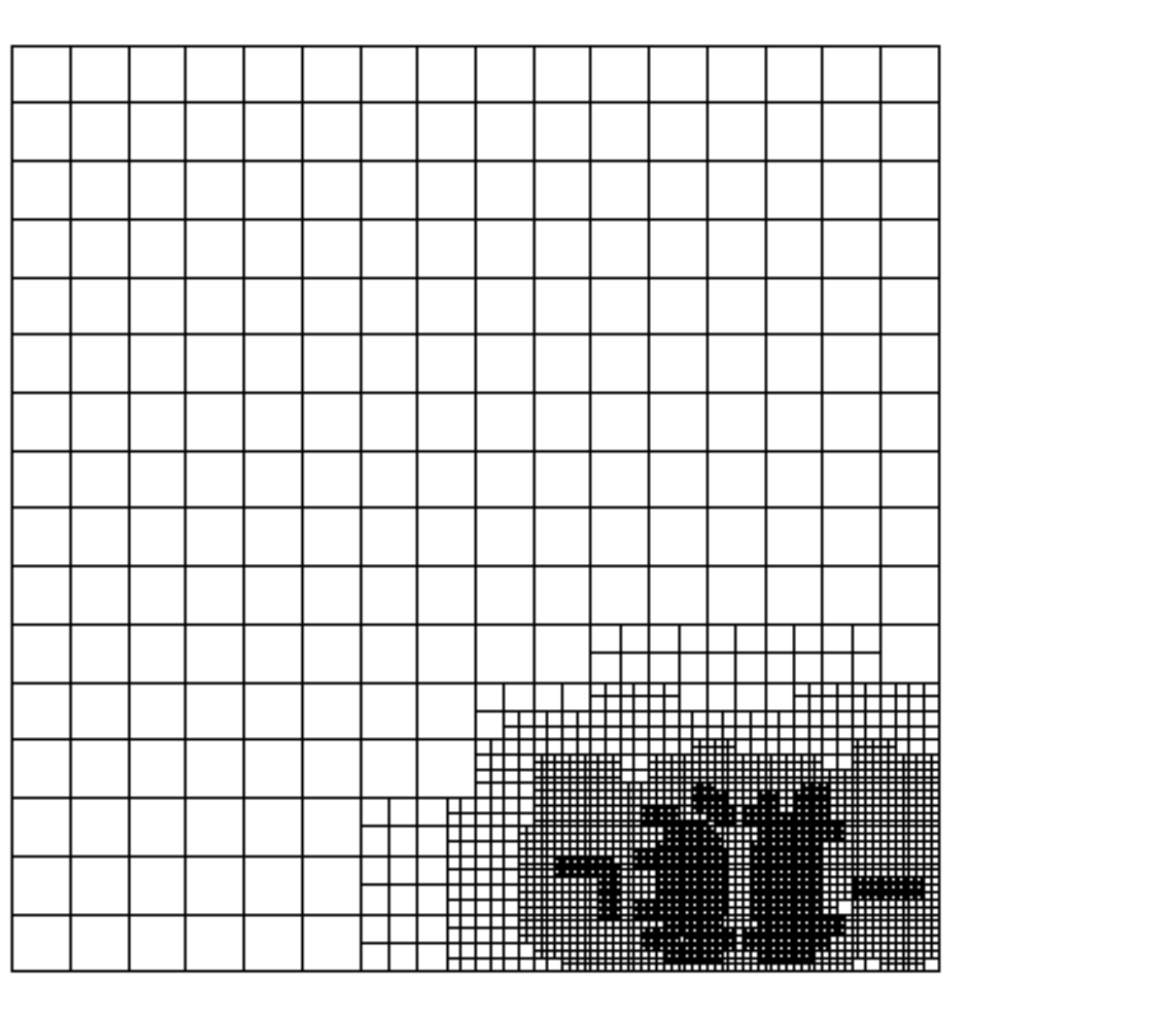}}
	}
	\hskip -30pt
	\subfloat[{REF 4: \newline
	ref. based on $\overline{\rm m}^{\rm I}_{\rm d}$ 
	}]{
	\spacetimeaxis{\includegraphics[width=4.2cm, trim={8.1cm 2cm 5cm 2cm}, clip]{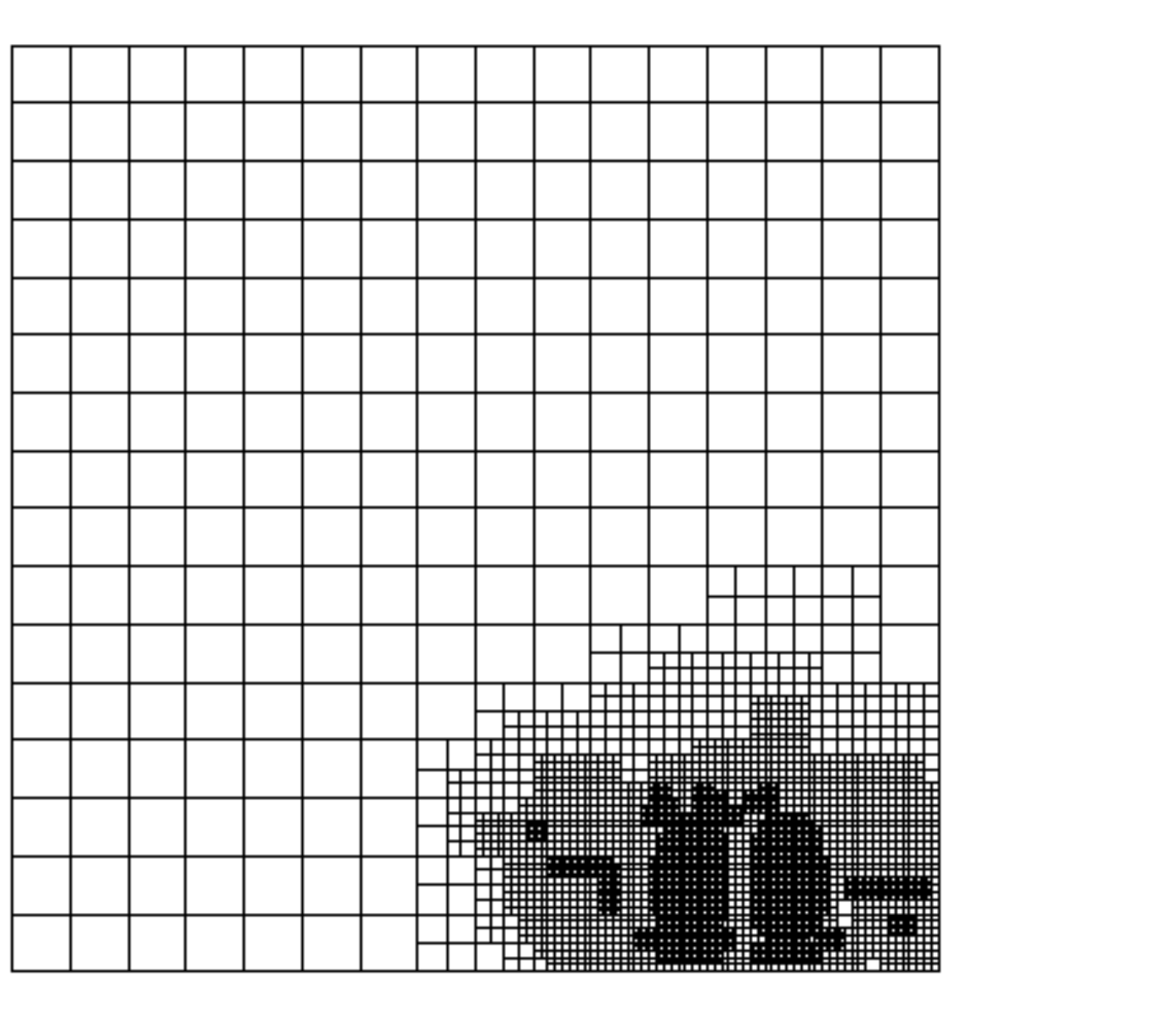}}
	}
	\hskip -30pt
	\subfloat[{REF 4: \newline
	ref. based on $|\!|\!|  e |\!|\!|_{\mathcal{L}}$}]{
	\spacetimeaxis{\includegraphics[width=4.2cm, trim={8.1cm 2cm 5cm 2cm}, clip]{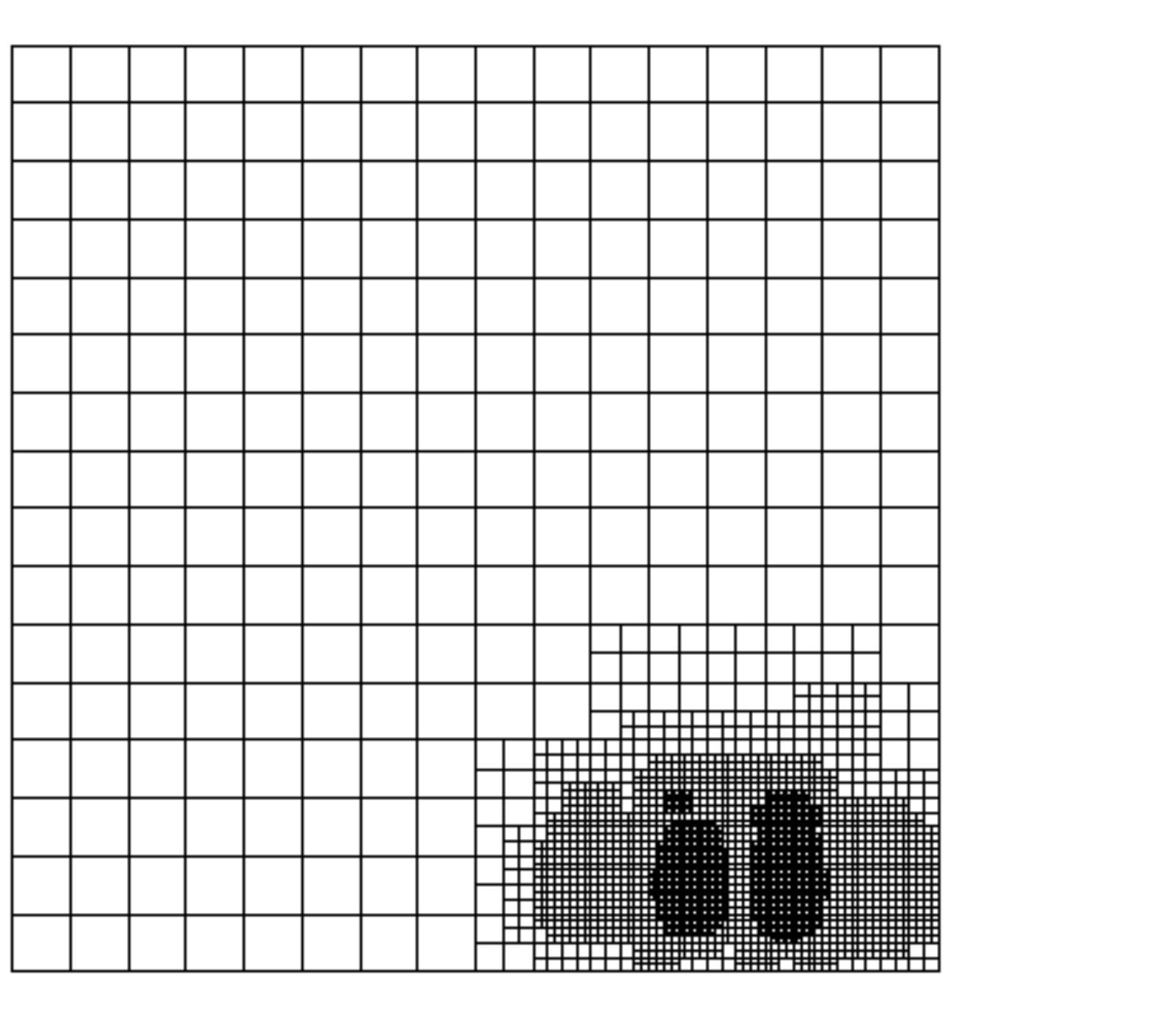}}
	}
	\hskip -30pt
	\subfloat[{REF 4: \newline
	ref. based on ${\EI}$}]{
	\spacetimeaxis{\includegraphics[width=4.2cm, trim={8.1cm 2cm 5cm 2cm}, clip]{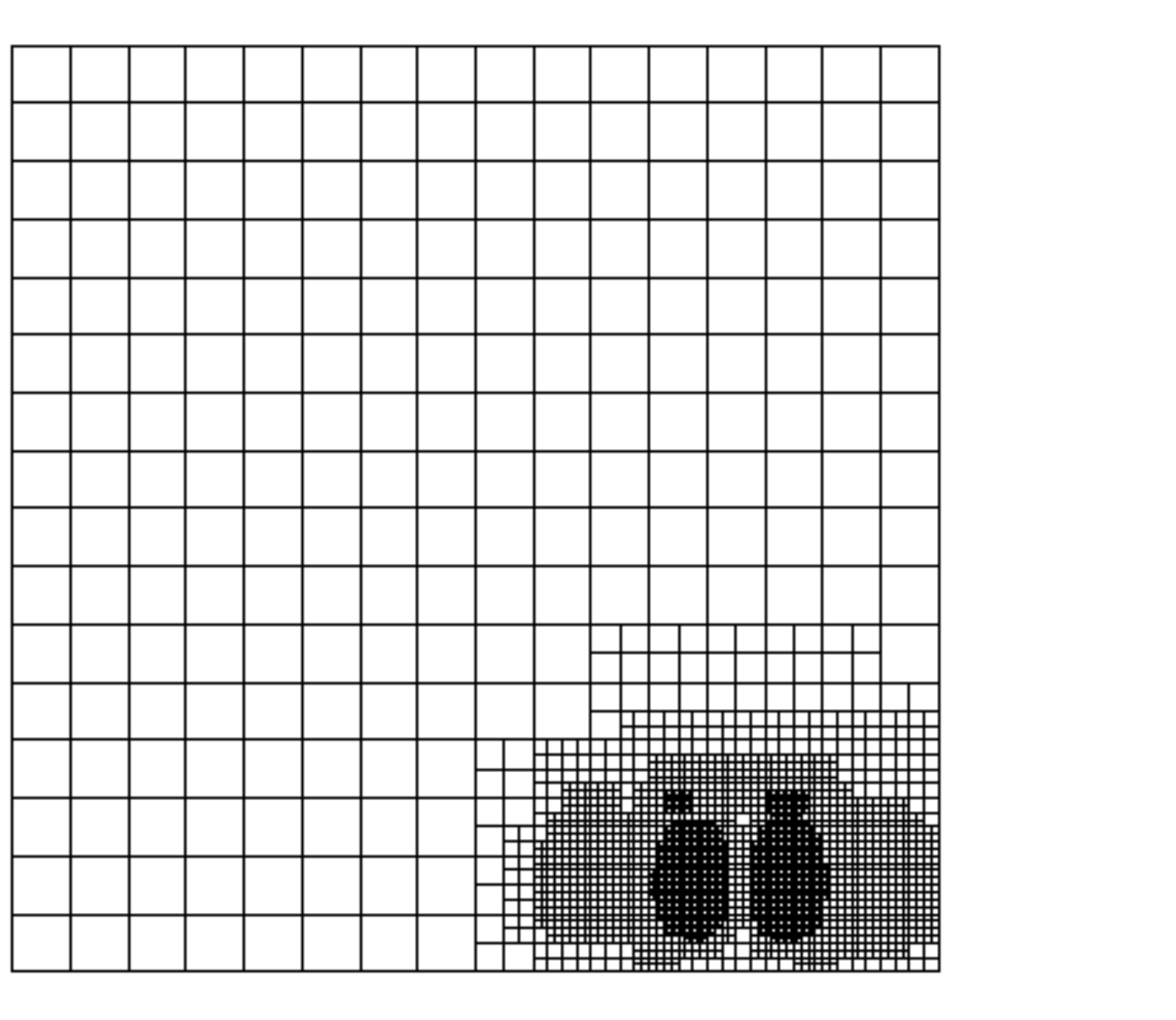}}
	}
	\\[-0.5pt]
	\subfloat[{REF 5: \newline
	ref. based on ${\rm e}_{\rm d}$}]{
	\spacetimeaxis{\includegraphics[width=4.2cm, trim={8.1cm 2cm 5cm 2cm}, clip]{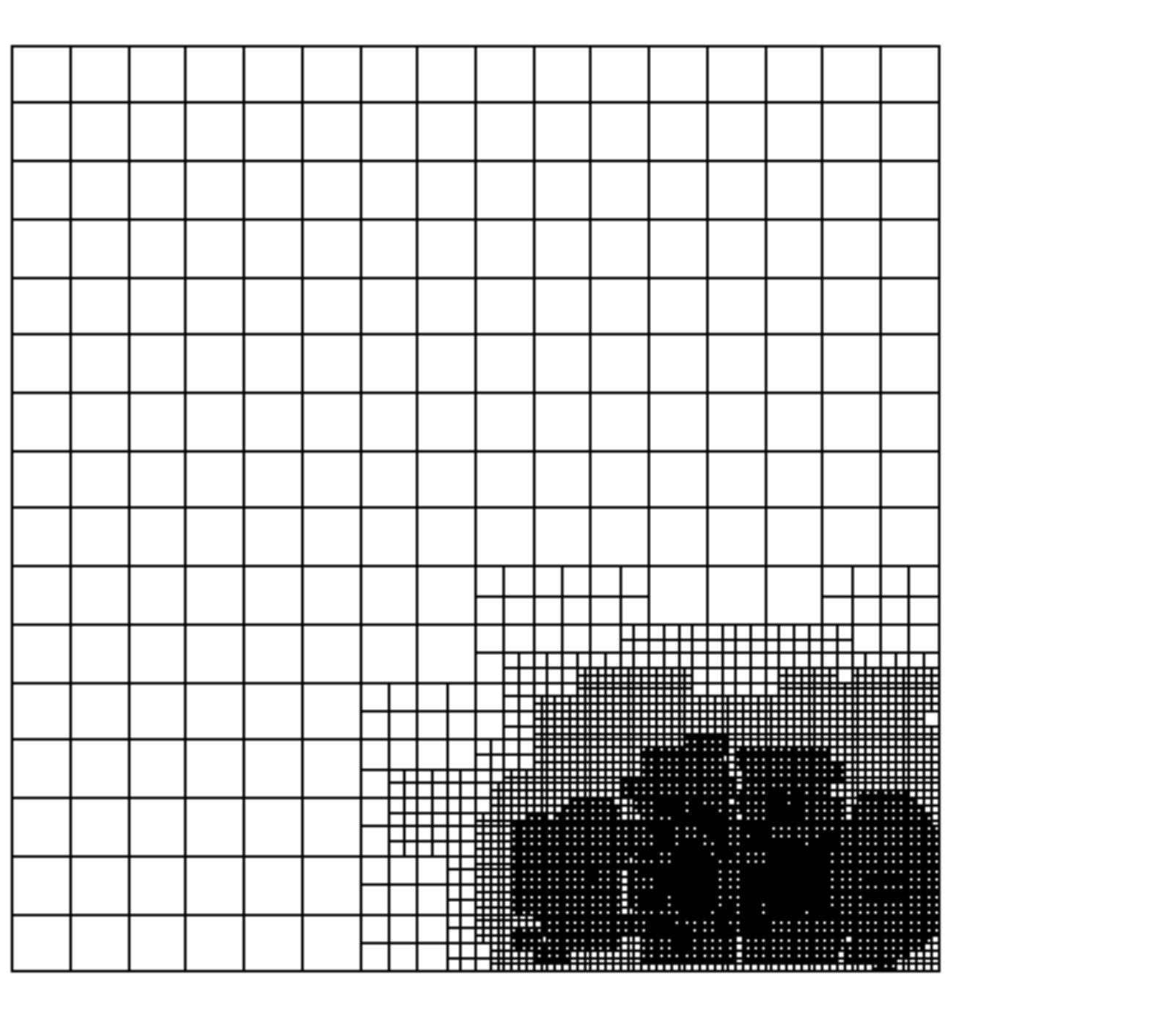}}
	}
	\hskip -30pt
	\subfloat[{REF 5: \newline
	ref. based on $\overline{\rm m}^{\rm I}_{\rm d}$
	}]{
	\spacetimeaxis{\includegraphics[width=4.2cm, trim={8.1cm 2cm 5cm 2cm}, clip]{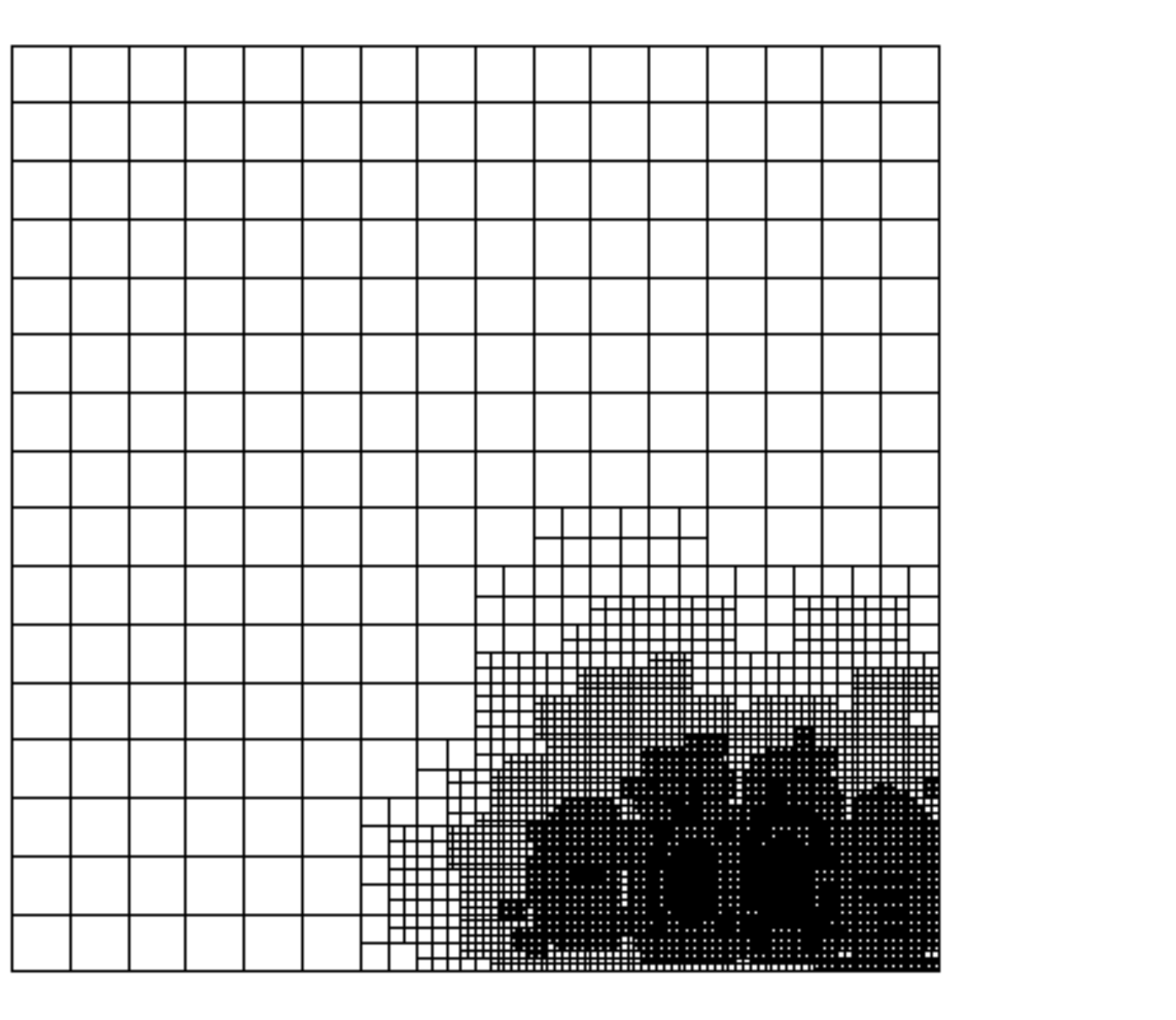}}
	}
	\hskip -30pt
	\subfloat[{REF 5: \newline
	ref. based on $|\!|\!|  e |\!|\!|_{\mathcal{L}}$
	}]{
	\spacetimeaxis{\includegraphics[width=4.2cm, trim={8.1cm 2cm 5cm 2cm}, clip]{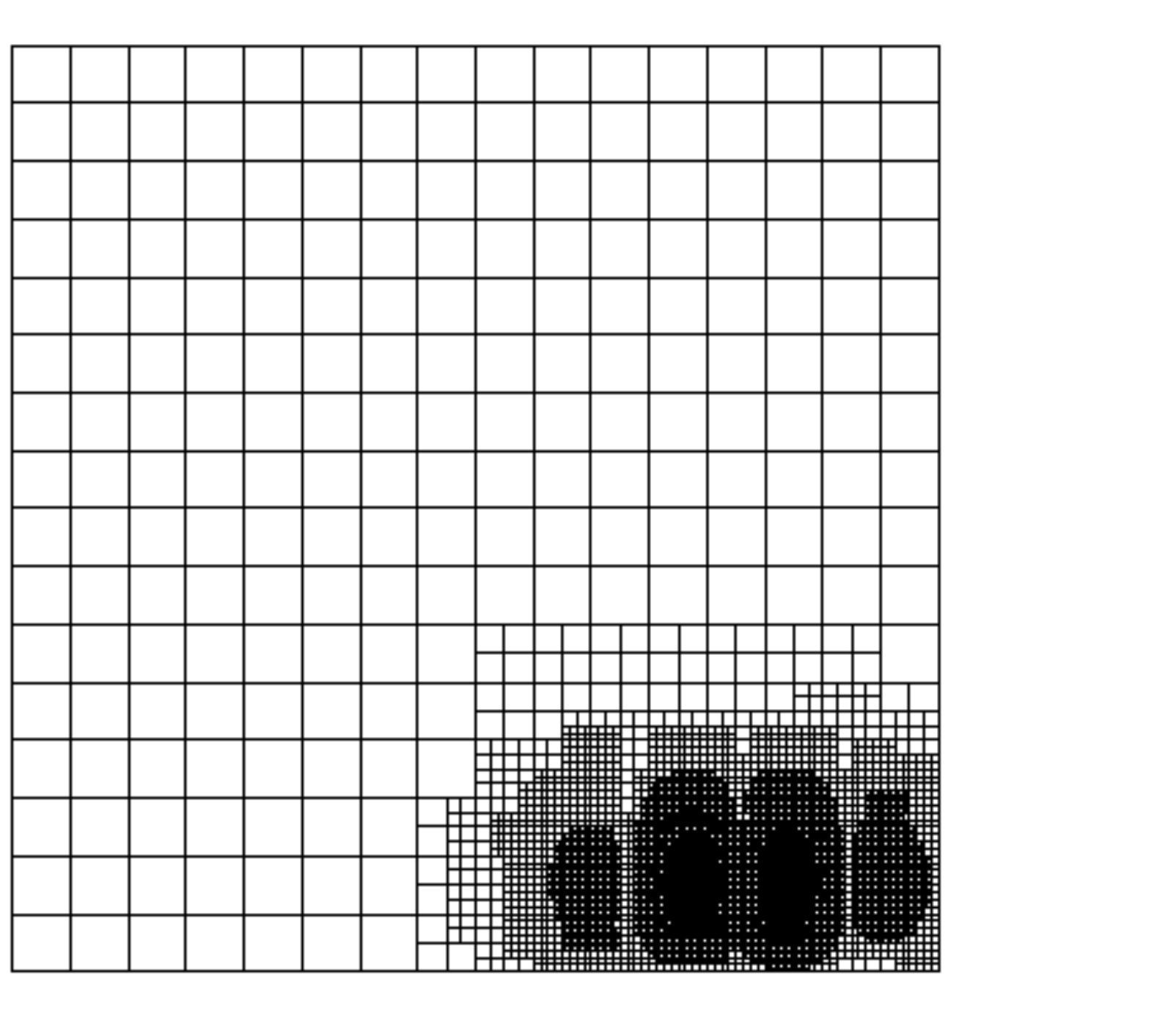}}
	}
	\hskip -30pt
	\subfloat[{REF 5: \newline
	ref. based on ${\EI}$}]{
	\spacetimeaxis{\includegraphics[width=4.2cm, trim={8.1cm 2cm 5cm 2cm}, clip]{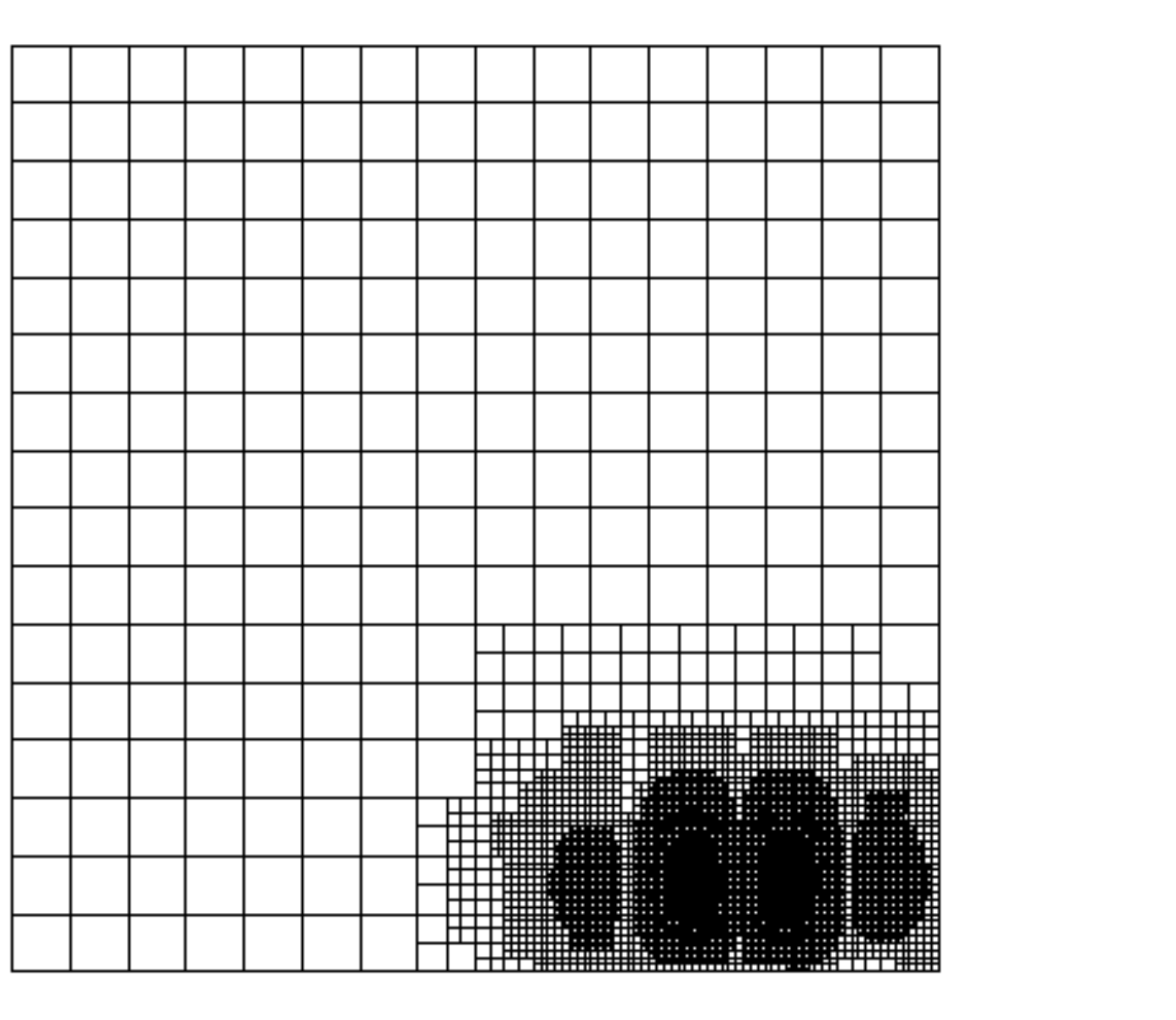}}
	}
	\caption{{\em Example 3}. 
	Comparison of the meshes obtained by refinement based on ${\rm e}_{\rm d}$, $\overline{\rm m}^{\rm I}_{\rm d}$, $|\!|\!|  e |\!|\!|_{\mathcal{L}}$, and ${\EI}$ for 
	$u_h \in S^{2}_{h}$, $\flux_h \in S_{h}^{3} \oplus S_{h}^{3}$, and $w_h \in S_{h}^{3}$, w.r.t. refinement steps
	4 and 5.}
	\label{fig:unit-domain-example-6-meshes-v-2-y-3-w-3-adaptive-ref}
\end{figure}

Coming back to the error indication properties of the majorant and the error identity, 
we analyse the meshes
presented in Figure \ref{fig:unit-domain-example-6-meshes-v-2-y-3-w-3-adaptive-ref}. 
The first two columns compare the meshes produced by the refinement based on ${\rm e}_{\rm d}$ and 
$\overline{\rm m}^{\rm I}_{\rm d}$. 
The third and the fourth columns illustrate practically matching meshes 
produced by values of $|\!|\!|  e |\!|\!|_{\mathcal{L}}$ and ${\EI}$. This test demonstrates that both 
$\overline{\rm m}^{\rm I}_{{\rm d}}$ and ${\EI}$ can be used for error indication and efficient 
mesh refinement algorithm. In particular, the error identity is suited for the cases when we can not afford
any time overhead for error analysis.

\begin{figure}[!t]
	\centering
	\captionsetup[subfigure]{oneside, labelformat=empty}
	\subfloat[]{
	\includegraphics[width=4.4cm, trim={0cm 0cm 0cm 0cm}, clip]{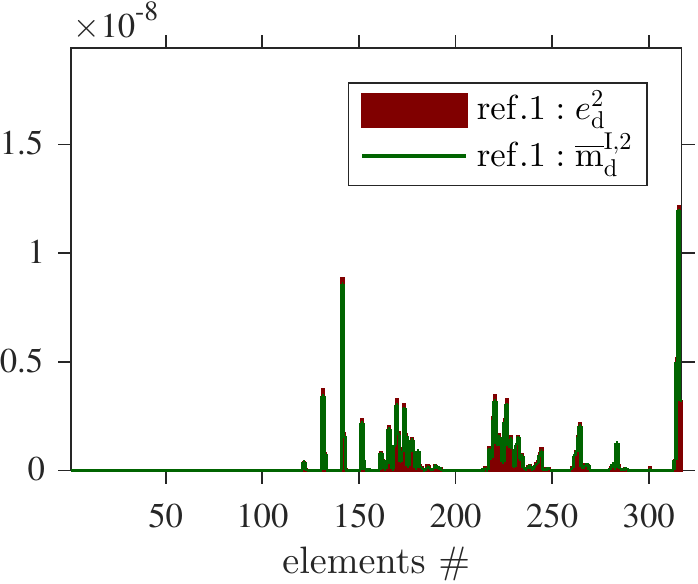}} 
	\hskip 10pt
	\subfloat[]{
	\includegraphics[width=4.2cm, trim={0cm 0cm 0cm 0cm}, clip]{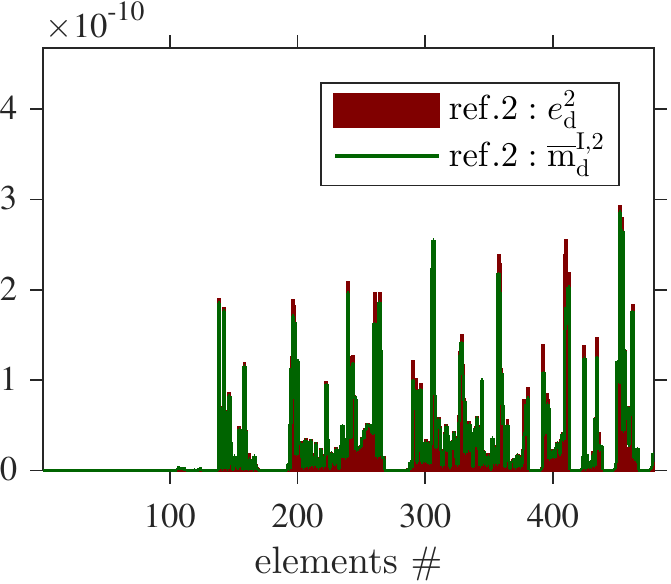}}
	\\[-20pt]
	\subfloat[]{
	\includegraphics[width=4.3cm, trim={0cm 0cm 0cm 0cm}, clip]{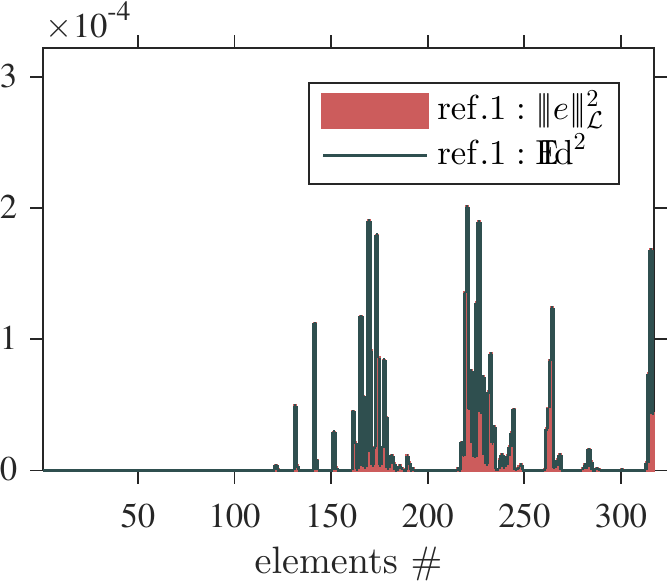}}
	\;
	\subfloat[]{
	\includegraphics[width=4.4cm, trim={0cm 0cm 0cm 0cm}, clip]{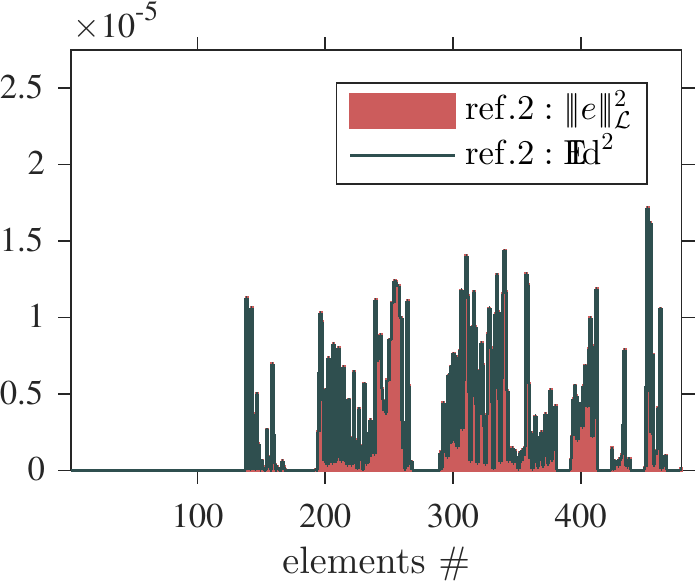}}
	\vskip -15pt
	\caption{{\em Example 3}. 
	Distribution of ${\rm e}_{{\rm d}, K}$, $\overline{\rm m}^{\rm I}_{{\rm d}, K}$, 
	as well as $|\!|\!|  e |\!|\!|_{\mathcal{L}, K}$, and $\EI_K$ on the refinement step 1 and 2.}
\end{figure}

{ 
To emphasise on the advantages of the space-time approach, 
we consider an analogous example for $d = 2$, i.e.,  $Q := (0, 1)^3$, 
such that the solution is defined by
\begin{alignat*}{3}
u(x, y, t) 	& = (x^2 - x) \,  (y^2 - y) \, (t^2 - t) \, e^{-100 \,|(x, y, t) - (0.25, 0.25, 0.25)|},  
& \quad  (x, y, t) \in \overline{Q}.
\end{alignat*}
Here, the Gaussian peak is located in the point $(x, y, t) = (0.25 0.25, 0.25)$. Analogously, 
$f$ is computed by substituting $u$ into \eqref{eq:equation}, 
as well as the Dirichlet boundary conditions remains  homogeneous. Figure 
\ref{fig:unit-domain-example-31-meshes-v-2-y-3-w-3-adaptive-ref} illustrates $2$-dimensional meshes
obtained by slicing $3$-dimensional meshes (in space and time) along the time variable. 
Figures demonstrate the advantage of the adaptive space-time approach over the time-stepping methods, 
which always have to consider refining or coarsening on subsequent time step. It is also worth 
emphasising on the universality of the applied functional error estimates w.r.t. to any discretisation. 
To our knowledge, it is the only approach, which allows fully unstructured discretisation in space and time
for evolutionary equations. Application of functional error estimates to FEM and their corresponding numerical 
properties are studied in \cite{LMR:HolmMatculevich:2017}.
\begin{figure}[!t]
	\centering
	\captionsetup[subfigure]{oneside, margin={0.3cm,0cm}}
	\subfloat[$t = 0$      ]{\includegraphics[width=4.5cm, trim={5cm 2cm 5cm 2cm}, clip]{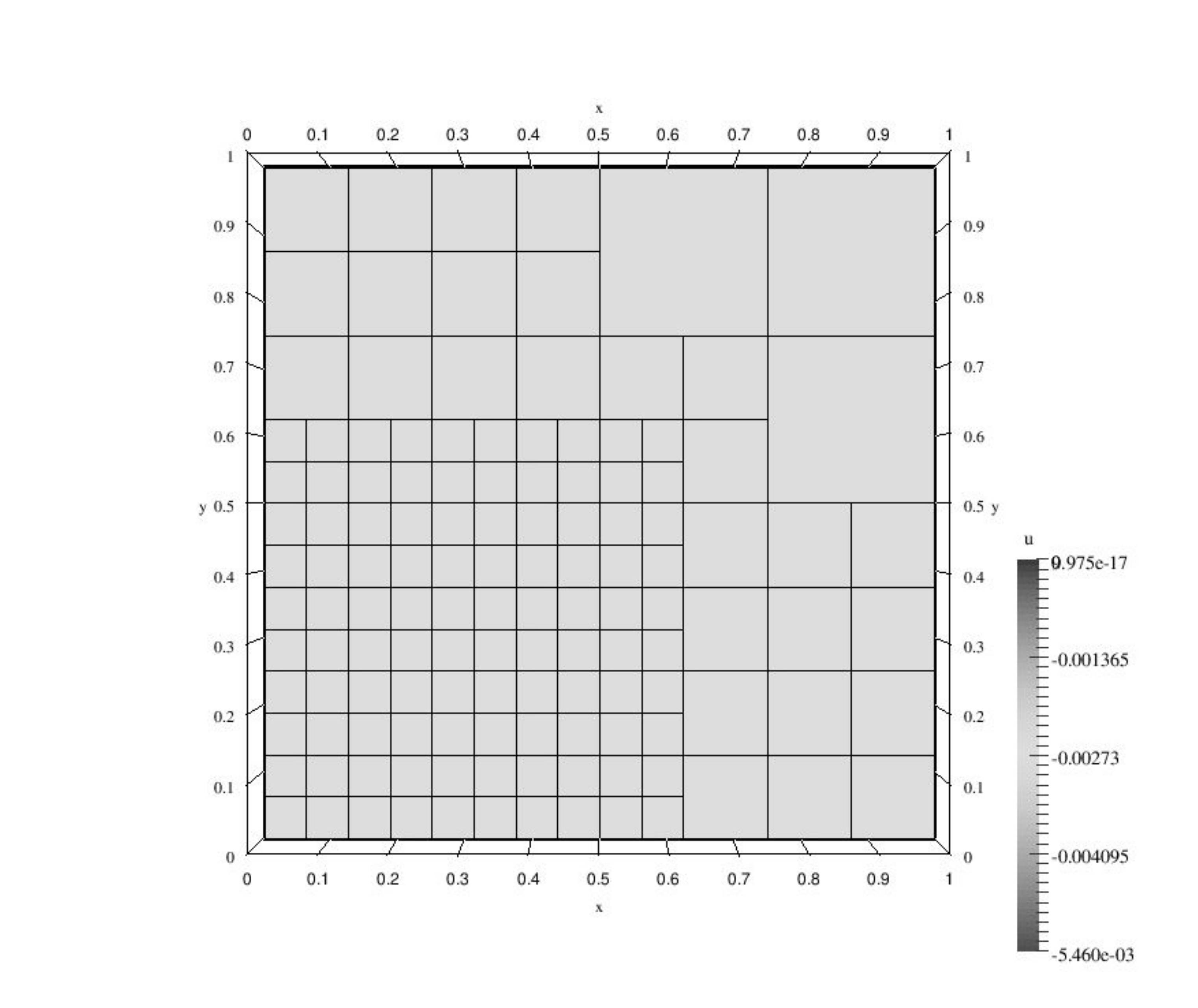}}\,
	\subfloat[$t = 0.125$]{\includegraphics[width=4.5cm, trim={5cm 2cm 5cm 2cm}, clip]{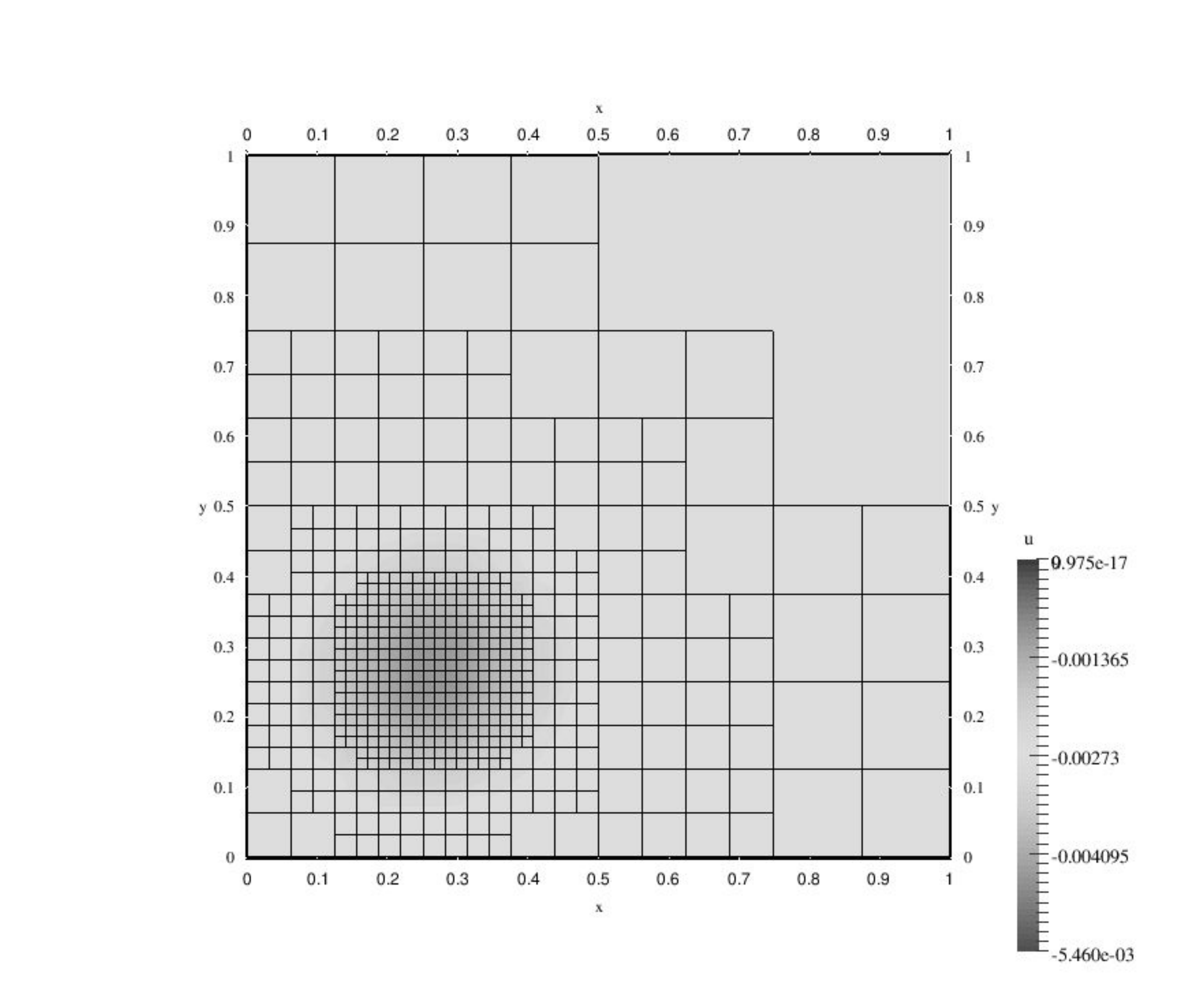}}\,
	\subfloat[$t = 0.25$  ]{\includegraphics[width=4.5cm, trim={5cm 2cm 5cm 2cm}, clip]{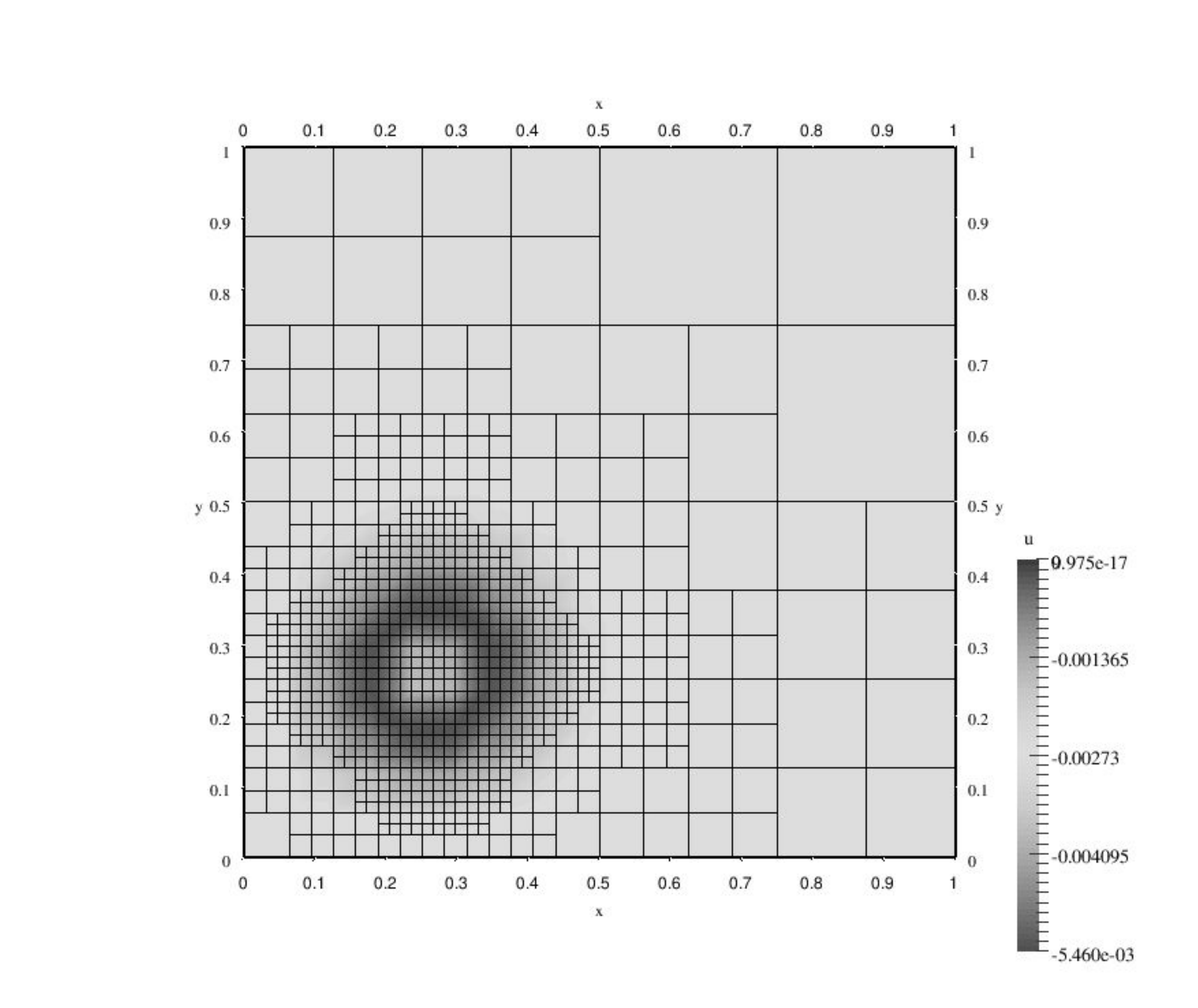}}
	\\[-10pt]
	\subfloat[$t = 0.5$    ]{\includegraphics[width=4.5cm, trim={5cm 2cm 5cm 2cm}, clip]{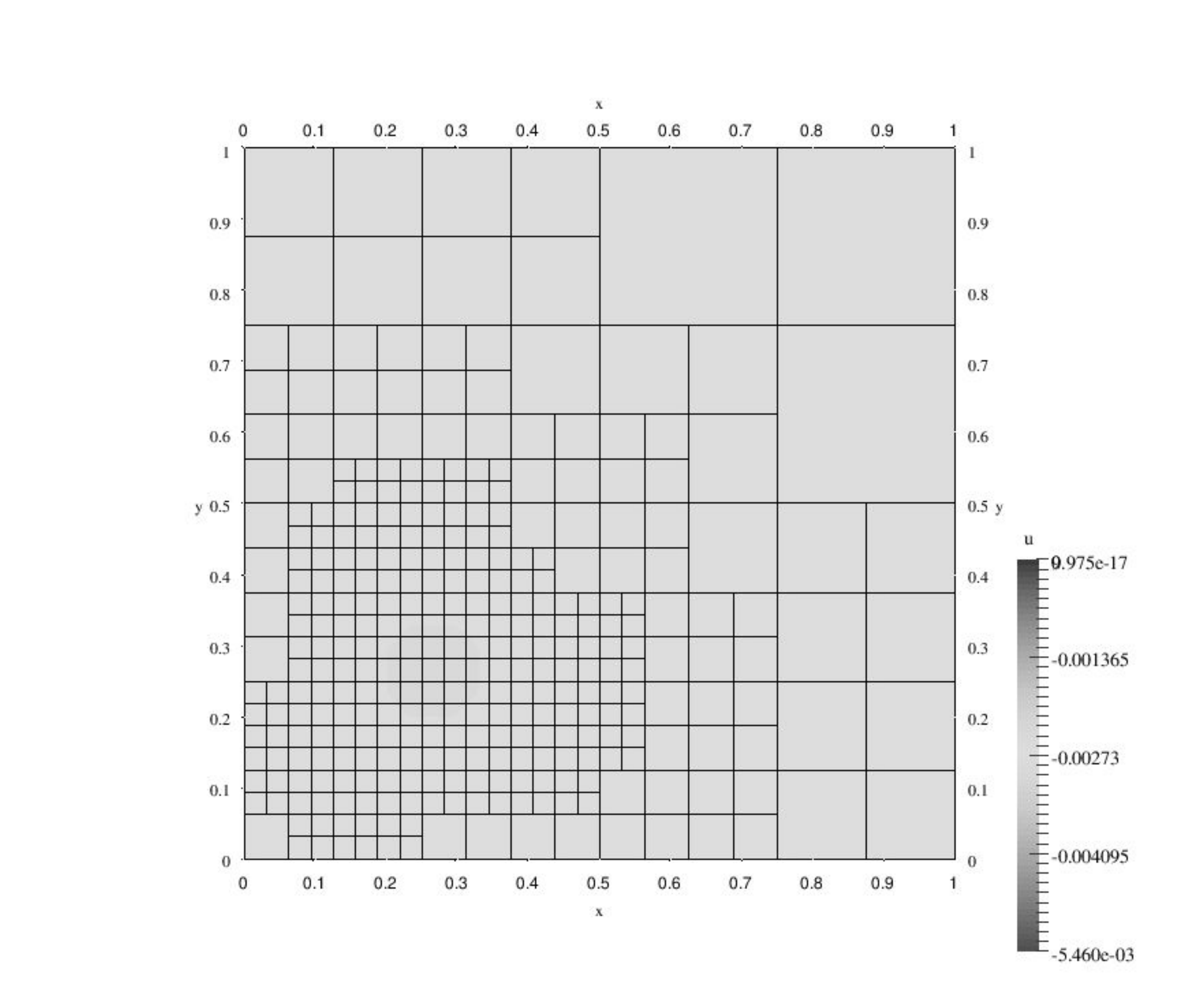}}
	\,
	\subfloat[$t = 0.75$   ]{\includegraphics[width=4.5cm, trim={5cm 2cm 5cm 2cm}, clip]{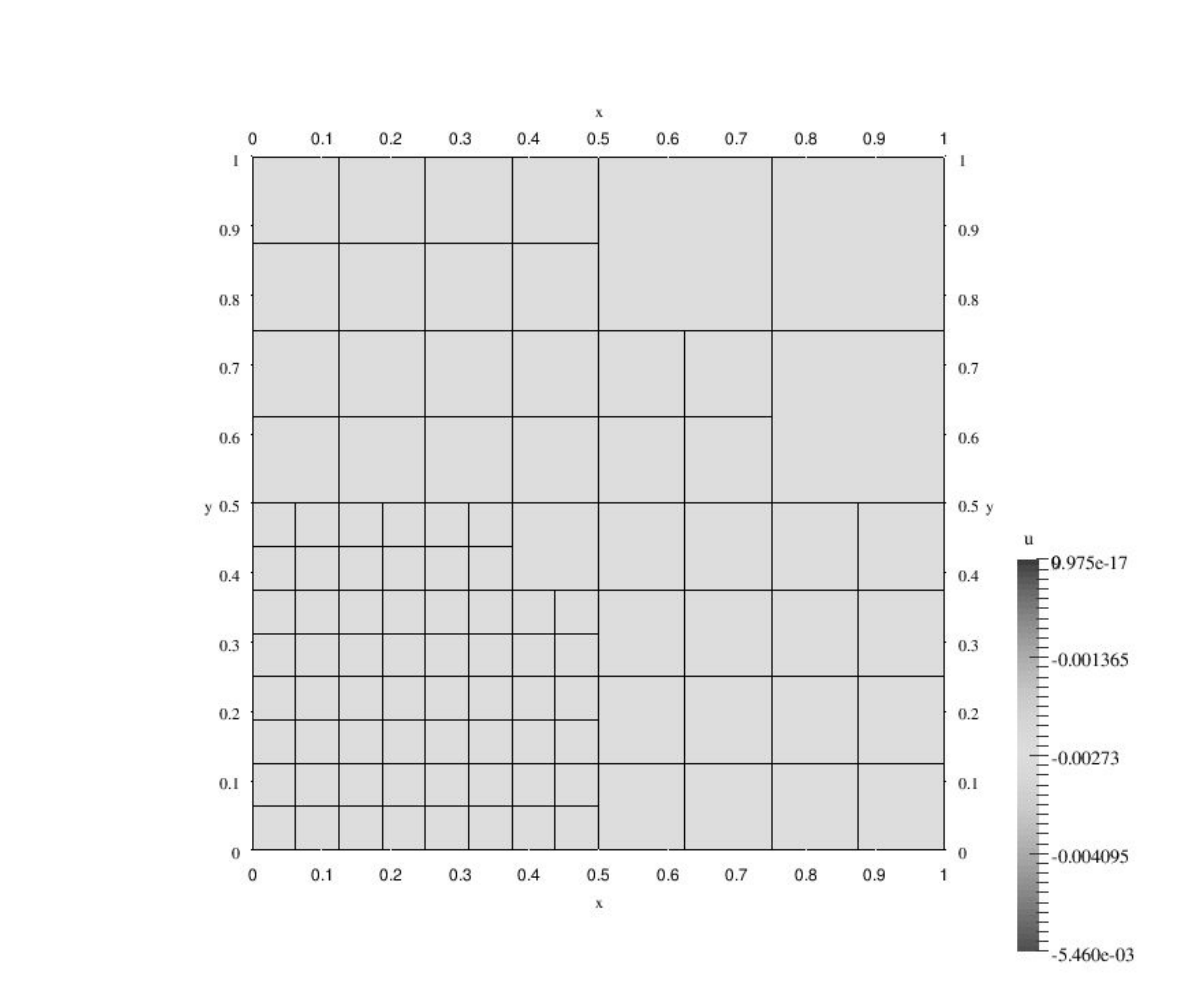}}\,
	\subfloat[$t = 1.0$     ]{\includegraphics[width=4.5cm, trim={5cm 2cm 5cm 2cm}, clip]{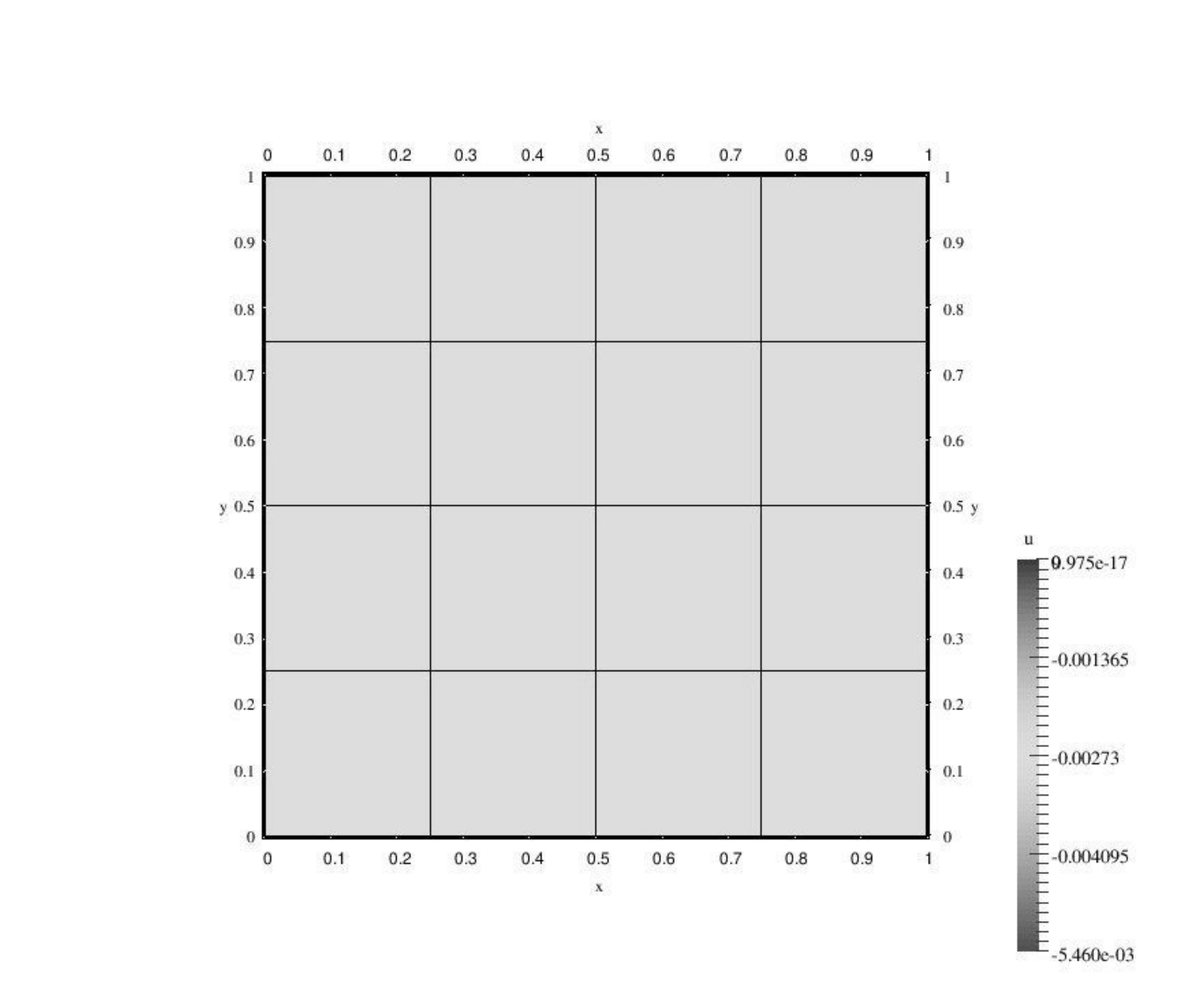}}
	\caption{{\em Example 3}. 
	Meshes obtained by the space-time refinement based on $\overline{\rm m}^{\rm I}_{\rm d}$ for 
	$u_h \in S^{2}_{h}$, $\flux_h \in S_{h}^{3} \oplus S_{h}^{3}$, and $w_h \in S_{h}^{3}$, w.r.t. 
	$t = \big\{0, 0.125, 0.25, 0.5, 0.75, 1.0 \big\}$.}
	\label{fig:unit-domain-example-31-meshes-v-2-y-3-w-3-adaptive-ref}
\end{figure}
}

\subsection{Example 4}
\label{ex:quarter-annulus-example-4}
\rm

\begin{figure}[!t]
	\centering
	\subfloat[]{
	\includegraphics[width=5cm, trim={0cm 0cm 0cm 2cm}, clip]{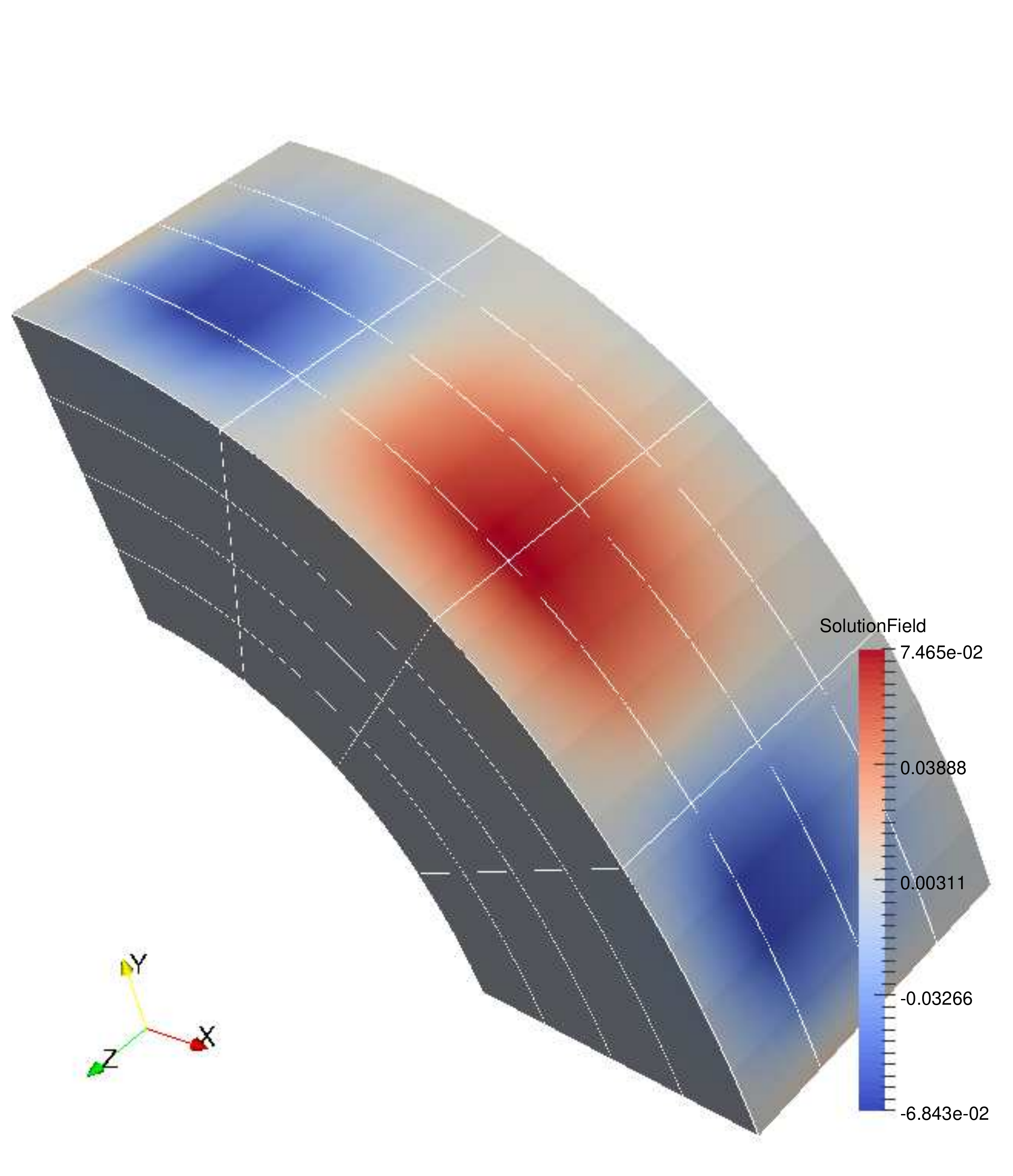}
	\label{fig:quarter-annulus-example-4-a}}
	\subfloat[]{
	\includegraphics[scale=0.6, trim={0cm 0cm 0cm 0cm}, clip]{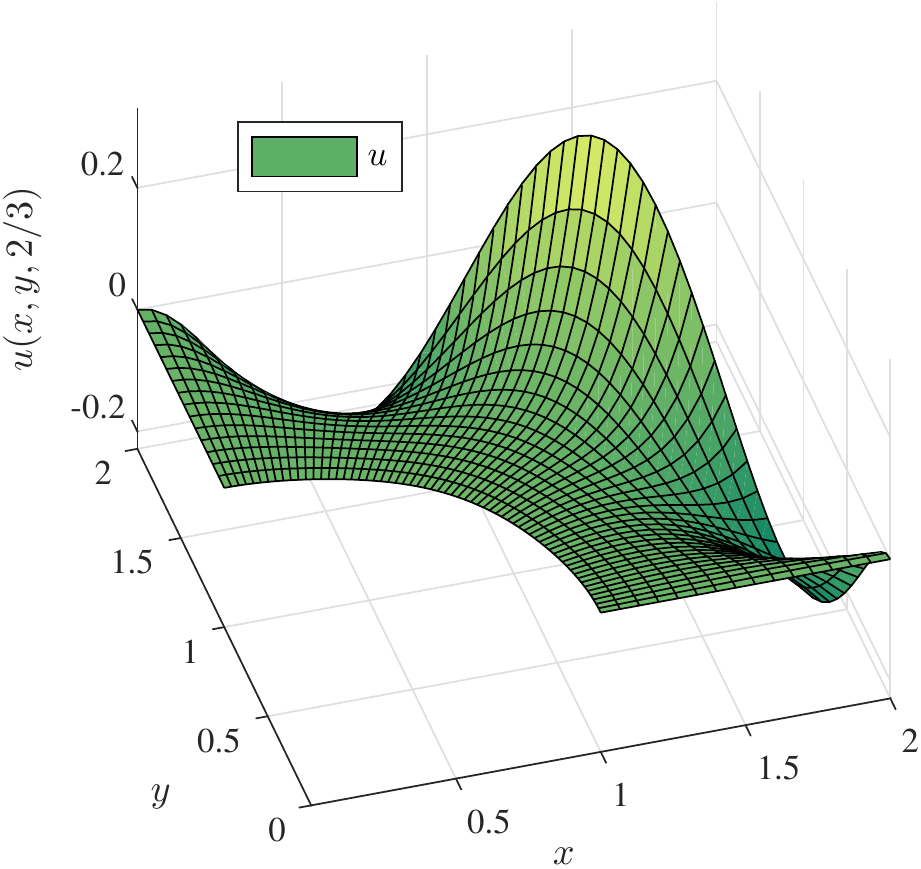}
	\label{fig:quarter-annulus-example-4-b}}
	\caption{\small {\em Example 4}. 
	(a) Exact solution $u = (1 - x) \, x^2 \, (1 - y) \, y^2 \, (1 - t) \, t^2$. 
	(b) $u$ at the time moment $t = \tfrac{2}{3}$.}
	\label{fig:quarter-annulus-example-4}
\end{figure}

Finally, in the last example, we test functional error estimates  on the three-dimensional 
{ 
space-time cylinder
}
$Q = \Omega_{ \,
\begin{tikzpicture}[scale=0.1]
\draw (0:1cm) -- (0:2cm)
arc (0:90:2cm) -- (90:1cm)
arc (90:0:1cm) -- cycle;
\end{tikzpicture}} \times (0, T)$, where  
$\Omega_{ 
\begin{tikzpicture}[scale=0.1]
\draw (0:1cm) -- (0:2cm)
arc (0:90:2cm) -- (90:1cm)
arc (90:0:1cm) -- cycle;
\end{tikzpicture}}
$ 
is of a quarter-annulus shape, and the final time of the time interval is $1$. 
The exact solution is defined by
\begin{alignat*}{4}
u(x, y, t) 	& = (1 - x) \, x^2 \, (1 - y) \, y^2 \, (1 - t) \, t^2, 
		& \quad  (x, y, t) & \in 
\overline{Q} := \overline{\Omega}_{ 
\begin{tikzpicture}[scale=0.1]
\draw (0:1cm) -- (0:2cm)
arc (0:90:2cm) -- (90:1cm)
arc (90:0:1cm) -- cycle;
\end{tikzpicture}}
\times [0, 1],
\end{alignat*}
{ 
see Figure \ref{fig:quarter-annulus-example-4}.
}
The RHS $f(x, y, t)$, $(x, y, t) \in 
Q := {\Omega}_{ 
\begin{tikzpicture}[scale=0.1]
\draw (0:1cm) -- (0:2cm)
arc (0:90:2cm) -- (90:1cm)
arc (90:0:1cm) -- cycle;
\end{tikzpicture}} \times (0, 1)$, 
is computed based on the substitution of $u$ into the equation 
\eqref{eq:equation} and the Dirichlet boundary conditions are defined as
{  $u_D = u$ on $\Sigma$.}
%

\noindent
The initial mesh for the test is generated by one uniform refinement $N^0_{\rm ref} = 1$.
We start the analysis from Table \ref{tab:quarter-annulus-example-4-estimates-v-2-y-3-uniform-ref}, 
where the performance of the studied error estimates is illustrated for both uniform and adaptive refinement 
strategies. It is easy to see that all majorants have adequate performance, taking into account 
that the auxiliary functions $\flux_h \in \oplus^3 S_{3h}^{3}$ and $w_h \in S^{3}_{3h}$ 
{ 
are from
spline spaces of just one order higher than 
the spline space for $u_h$. 
}
Such choice of the spaces 
is beneficial when the time expenditure on error estimation is concerned. Table 
\ref{tab:quarter-annulus-example-4-times-v-2-y-3-adaptive-ref} confirms that 
assembling and solving of the systems reconstructing d.o.f. of $u_h$ requires more time than assembling 
and solving routines for the systems generating $\flux_h$ and $w_h$.

\begin{table}[!t]
\scriptsize
\centering
\newcolumntype{g}{>{\columncolor{gainsboro}}c} 	
\newcolumntype{k}{>{\columncolor{lightgray}}c} 	
\newcolumntype{s}{>{\columncolor{silver}}c} 
\newcolumntype{a}{>{\columncolor{ashgrey}}c}
\newcolumntype{b}{>{\columncolor{battleshipgrey}}c}
\begin{tabular}{c|cga|ck|cb|cc}
\parbox[c]{0.8cm}{\centering \# ref. } & 
\parbox[c]{1.4cm}{\centering  $\| \nabla_x e \|_Q$}   & 	  
\parbox[c]{1.0cm}{\centering $\Ieff (\overline{\rm M}^{\rm I})$ } & 
\parbox[c]{1.4cm}{\centering $\Ieff (\overline{\rm M}^{\rm I\!I})$ } & 
\parbox[c]{1.0cm}{\centering  $|\!|\!|  e |\!|\!|_{s, h}$ }   & 	  
\parbox[c]{1.4cm}{\centering $\Ieff (\overline{\rm M}^{\rm I}_{s, h})$ } & 
\parbox[c]{1.0cm}{\centering  $|\!|\!|  e |\!|\!|_{\mathcal{L}}$ }   & 	  
\parbox[c]{1.4cm}{\centering$\Ieff ({\EI})$ } & 
\parbox[c]{1.2cm}{\centering e.o.c. ($|\!|\!|  e |\!|\!|_{s, h}$)} & 
\parbox[c]{1.2cm}{\centering e.o.c. ($|\!|\!|  e |\!|\!|_{\mathcal{L}}$)} \\
\bottomrule
\multicolumn{10}{l}{ \rule{0pt}{3ex} uniform refinement }\\[2pt]
\toprule
   3 &     9.6674e-03 &         6.17 &         5.37 &     9.6675e-03 &         6.18 &     1.9361e-01 &         1.00 &     3.30 &     1.46 \\
   4 &     2.1555e-03 &         6.89 &         5.04 &     2.1555e-03 &         6.90 &     9.0086e-02 &         1.00 &     2.55 &     1.30 \\
   5 &     5.2154e-04 &         3.35 &         2.67 &     5.2154e-04 &         3.39 &     4.4046e-02 &         1.00 &     2.23 &     1.13 \\
   6 &     1.2926e-04 &         2.04 &         1.67 &     1.2926e-04 &         2.15 &     2.1892e-02 &         1.00 &     2.10 &     1.05 \\
\bottomrule
\multicolumn{10}{l}{ \rule{0pt}{3ex}   
adaptive refinement using bulk marking criterion $\sigma =0.4$ }\\[2pt]
\toprule
%
   3 &     9.6919e-03 &         6.16 &         5.36 &     9.6919e-03 &         6.16 &     1.9385e-01 &         1.00 &     3.56 &     1.57 \\
   4 &     2.4997e-03 &         5.98 &         4.38 &     2.4997e-03 &         5.98 &     9.4309e-02 &         1.00 &     3.25 &     1.73 \\
   5 &     9.4272e-04 &         2.28 &         1.78 &     9.4272e-04 &         2.28 &     5.4120e-02 &         1.00 &     1.79 &     1.02 \\
   6 &     2.5862e-04 &         1.52 &         1.29 &     2.5862e-04 &         1.52 &     3.0451e-02 &         1.00 &     2.88 &     1.28 \\
\end{tabular}
\caption{{\em Example 4}. 
Efficiency of $\overline{\rm M}^{\rm I}$, $\overline{\rm M}^{\rm I\!I}$, 
$\overline{\rm M}^{\rm I}_{s, h}$, and ${\EI}$ for $u_h \in S^{2}_{h}$,
$\flux_h \in \oplus^3 S_{2h}^{3}$, and $w_h \in S^{3}_{2h}$, 
w.r.t. uniform refinement and adaptive refinement steps.}
\label{tab:quarter-annulus-example-4-estimates-v-2-y-3-uniform-ref}
\end{table}

\begin{table}[!t]
\scriptsize
\centering
\newcolumntype{g}{>{\columncolor{gainsboro}}c} 	
\begin{tabular}{c|ccc|cgg|cgg|c}
& \multicolumn{3}{c|}{ d.o.f. } 
& \multicolumn{3}{c|}{ $t_{\rm as}$ }
& \multicolumn{3}{c|}{ $t_{\rm sol}$ } 
& $\tfrac{t_{\rm appr.}}{t_{\rm er.est.}}$ \\
\midrule
\parbox[c]{0.8cm}{\centering \# ref. } & 
\parbox[c]{0.8cm}{\centering $u_h$ } &  
\parbox[c]{0.6cm}{\centering $\flux_h$ } &  
\parbox[c]{0.6cm}{\centering $w_h$ } & 
\parbox[c]{1.4cm}{\centering $u_h$ } & 
\parbox[c]{1.4cm}{\centering $\flux_h$ } & 
\parbox[c]{1.4cm}{\centering $w_h$ } & 
\parbox[c]{1.4cm}{\centering $u_h$ } & 
\parbox[c]{1.4cm}{\centering $\flux_h$ } & 
\parbox[c]{1.4cm}{\centering $w_h$ } \\
\bottomrule 
\multicolumn{10}{l}{  \rule{0pt}{3ex} uniform refinement} \\[2pt]
\toprule
   2 &        216 &        375 &        125 &   5.68e-02 &   9.14e-02 &   4.70e-02 &         1.12e-03 &         9.10e-03 &         2.15e-04 &             0.39 \\
   3 &       1000 &        375 &        125 &   3.87e-01 &   9.34e-02 &   4.51e-02 &         4.14e-02 &         1.40e-02 &         2.39e-04 &             2.81 \\
   4 &       5832 &       1029 &        343 &   2.97e+00 &   6.94e-01 &   2.62e-01 &         3.49e+00 &         1.79e-01 &         3.62e-03 &             5.67 \\
   5 &      39304 &       3993 &       1331 &   2.36e+01 &   5.74e+00 &   1.96e+00 &         1.09e+02 &         3.77e+00 &         1.67e-01 &            11.40 \\
   6 &     287496 &      20577 &       6859 &   1.80e+02 &   2.91e+01 &   9.53e+00 &         1.12e+04 &         5.17e+01 &         6.76e+00 &           117.73 \\
 \midrule
    &       &         &    &
    \multicolumn{3}{c|}{ $t_{\rm as} (u_h)$ \quad : \quad $t_{\rm as} (\flux_h)$ \quad : \qquad $t_{\rm as} (w_h)$ } &      
    \multicolumn{3}{c|}{\; $t_{\rm sol} (u_h)$ \, : \quad $t_{\rm sol} (\flux_h)$  \quad:  \qquad  $t_{\rm sol} (w_h)$\;} & \\
\midrule
 	 & 	 & 	 & 	 &      18.93 &       3.06 &       1.00 &          1664.71 &             7.64 &             1.00 &           \\
\bottomrule
\multicolumn{10}{l}{ \rule{0pt}{3ex}   
adaptive refinement using bulk marking criterion $\sigma =0.4$ ($N_{\rm ref, 0}$ = 1)} \\[2pt]
\toprule
   2 &        216 &        375 &        125 &   5.55e-01 &   5.01e-01 &   4.15e-01 &         8.41e-04 &         1.08e-02 &         3.65e-04 &             0.60 \\
   3 &        894 &        375 &        125 &   5.48e+00 &   5.79e-01 &   3.99e-01 &         2.28e-02 &         7.36e-03 &         2.38e-04 &             5.58 \\
   4 &       3127 &       1029 &        343 &   4.60e+01 &   5.15e+00 &   3.19e+00 &         5.16e-01 &         1.75e-01 &         2.68e-03 &             5.46 \\
   5 &      15990 &       3993 &       1331 &   2.25e+02 &   3.92e+01 &   2.29e+01 &         1.96e+01 &         3.41e+00 &         1.02e-01 &             3.72 \\
   6 &      61390 &      20577 &       6859 &   1.22e+03 &   2.33e+02 &   1.83e+02 &         3.85e+02 &         5.86e+01 &         6.63e+00 &             3.32 \\
 
 \midrule
    &       &         &    &
    \multicolumn{3}{c|}{ $t_{\rm as} (u_h)$ \quad : \quad $t_{\rm as} (\flux_h)$ \quad : \qquad $t_{\rm as} (w_h)$ } &      
    \multicolumn{3}{c|}{\; $t_{\rm sol} (u_h)$ \, : \quad $t_{\rm sol} (\flux_h)$  \quad:  \qquad  $t_{\rm sol} (w_h)$\;} & \\
\midrule
         &       &       &       &       6.63 &       1.27 &       1.00 &            58.07 &             8.84 &             1.00 &             \\
\end{tabular}
\caption{{\em Example 4}. 
Assembling and solving time (in seconds) spent for the systems generating
d.o.f. of $u_h \in S^{3}_{3h}$
$\flux_h \in \oplus^3 S^{3}_{3h}$, and 
$w_h \in S^{3}_{3h}$ w.r.t (a) uniform refinements and (b) 
adaptive refinements (using bulk marking criterion with $\sigma =0.4$).}
\label{tab:quarter-annulus-example-4-times-v-2-y-3-adaptive-ref}
\end{table}

\begin{figure}[!t]
	\centering
	\captionsetup[subfigure]{oneside, labelformat=empty}
	\subfloat[ref. \# 1: \quad $\hat{\Omega}$ and $\hat{\mathcal{K}}_h$]{
	{\includegraphics[width=4cm, trim={2cm 2cm 3cm 4cm}, clip]{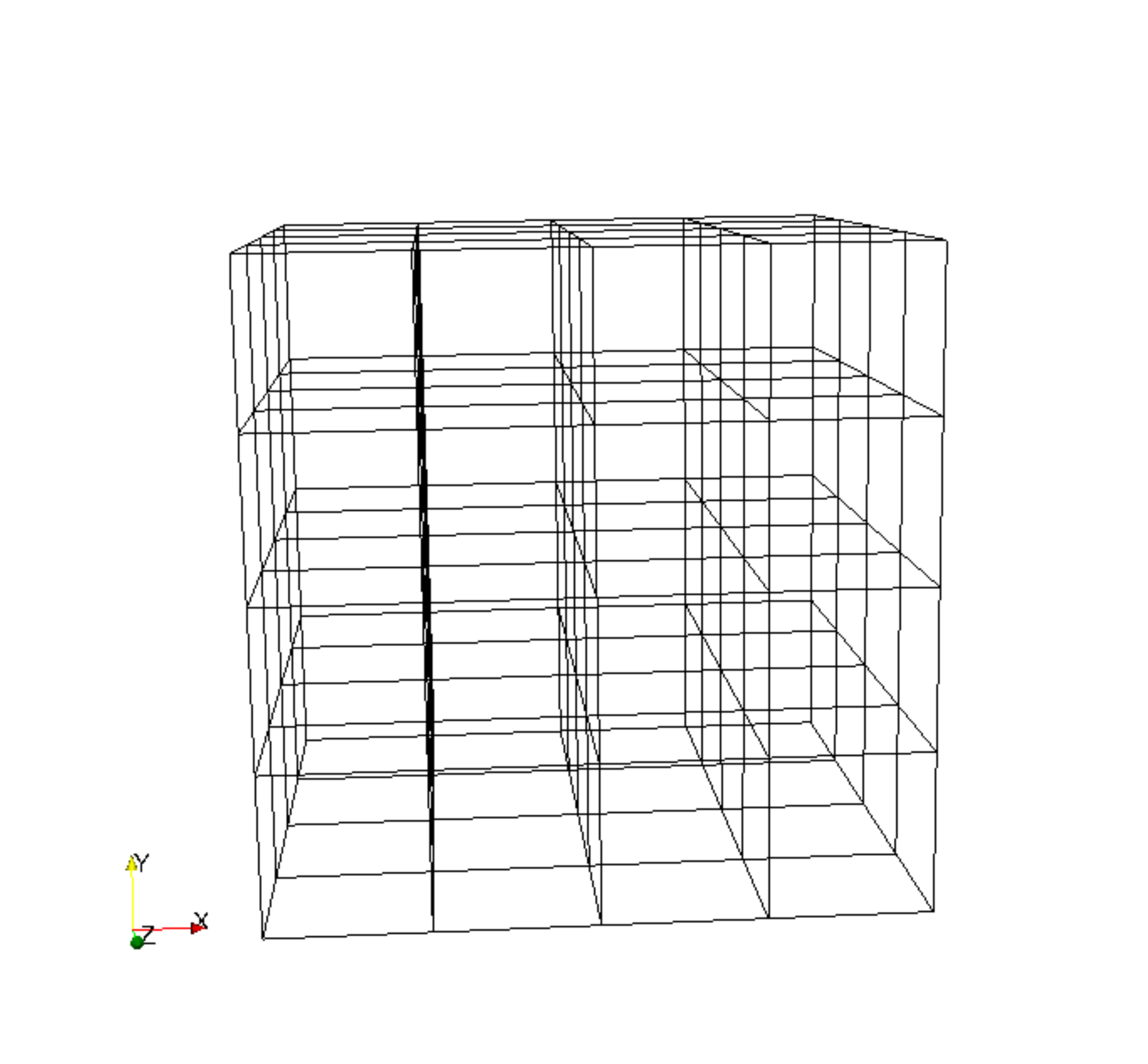}}}
	~
	\subfloat[ref. \# 1: \quad ${\Omega}$ and ${\mathcal{K}}_h$]{
	\includegraphics[width=4cm, trim={2cm 2cm 3cm 4cm}, clip]{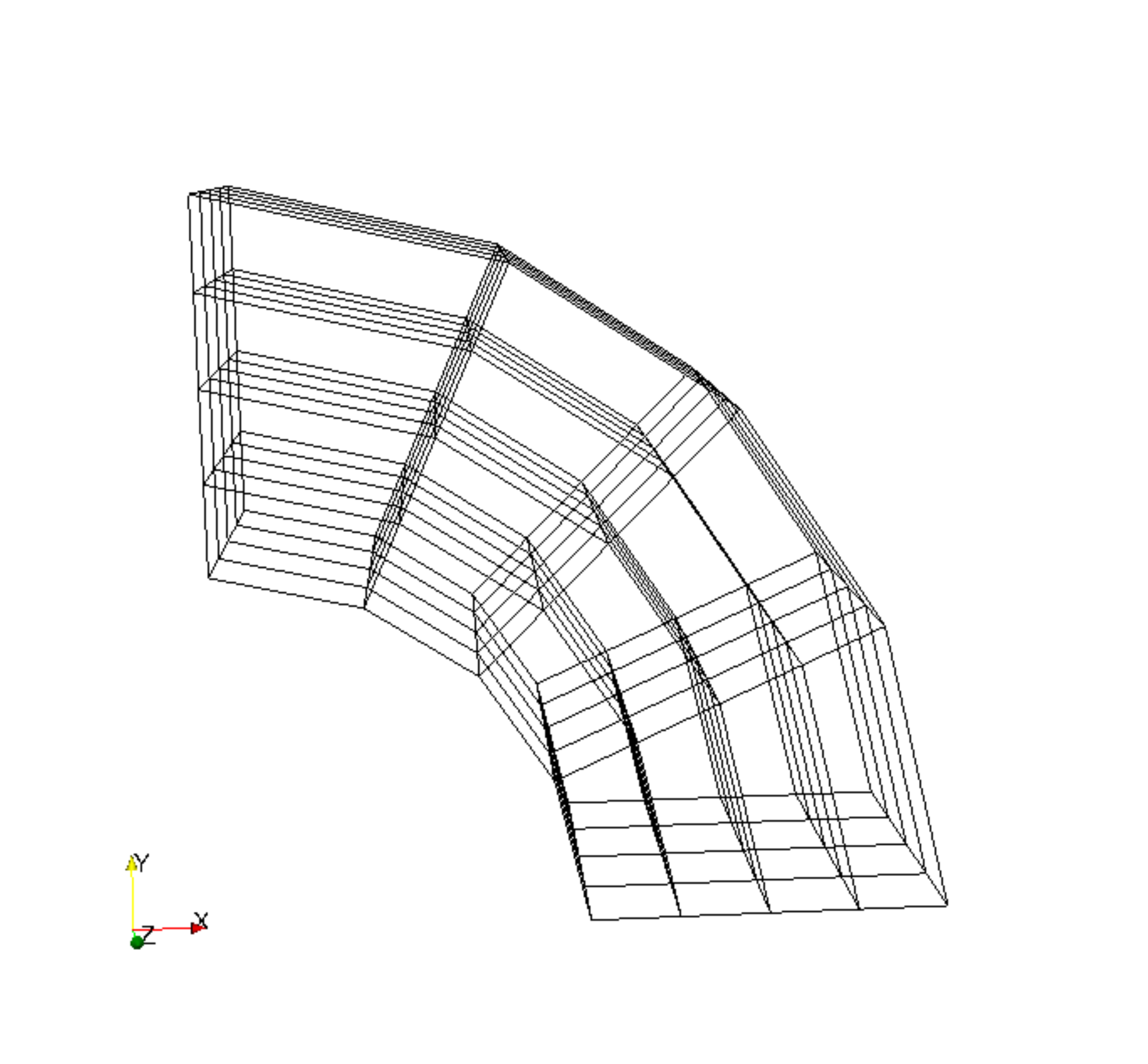}} 
	~
	%
	%
	\subfloat[ref. \# 3: \quad $\hat{\Omega}$ and $\hat{\mathcal{K}}_h$]{
	\includegraphics[width=4cm, trim={2cm 2cm 3cm 4cm}, clip]{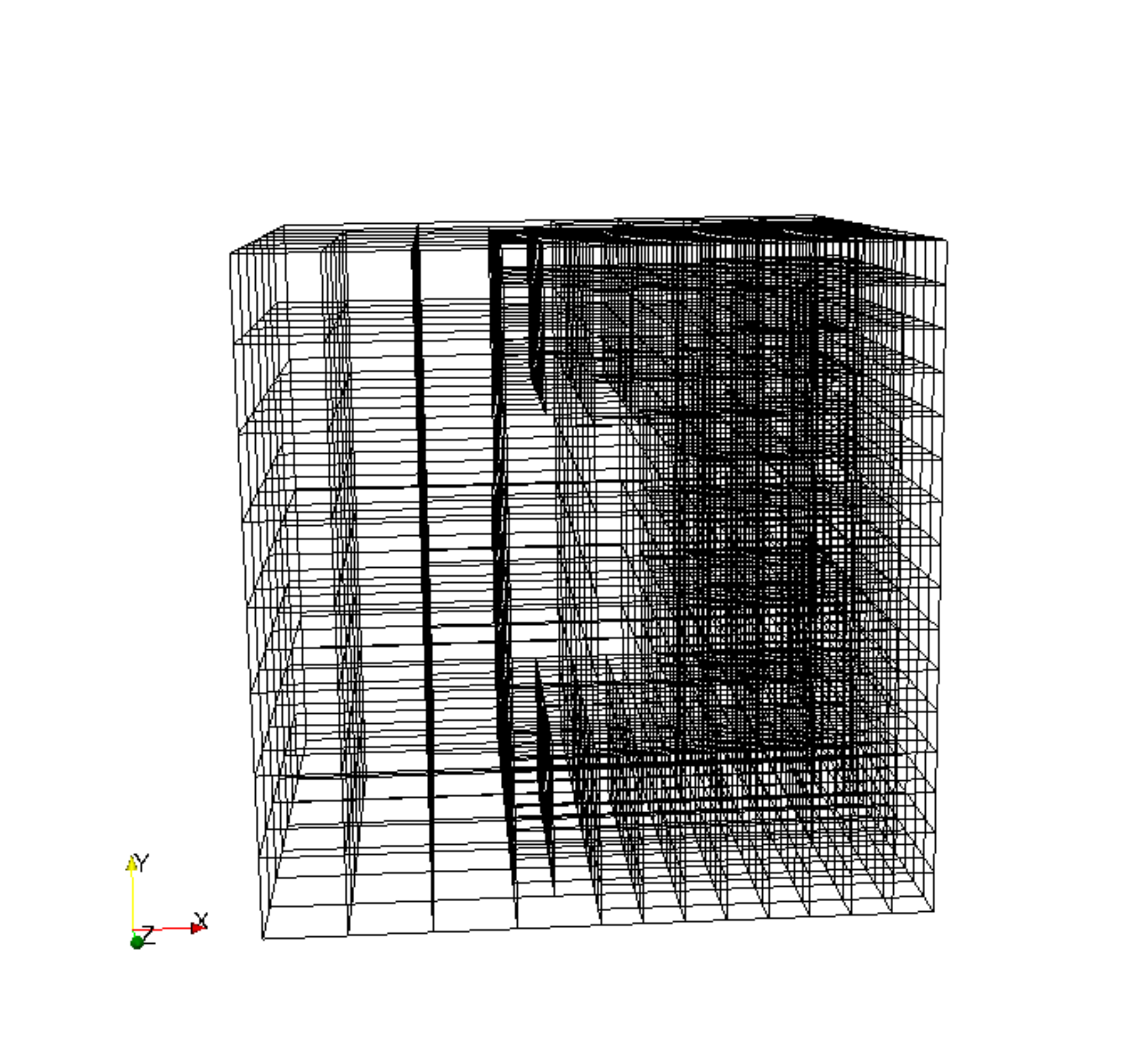}}
	~
	\subfloat[ref. \# 3: \quad ${\Omega}$ and ${\mathcal{K}}_h$]{
	\includegraphics[width=4cm, trim={2cm 2cm 3cm 4cm}, clip]{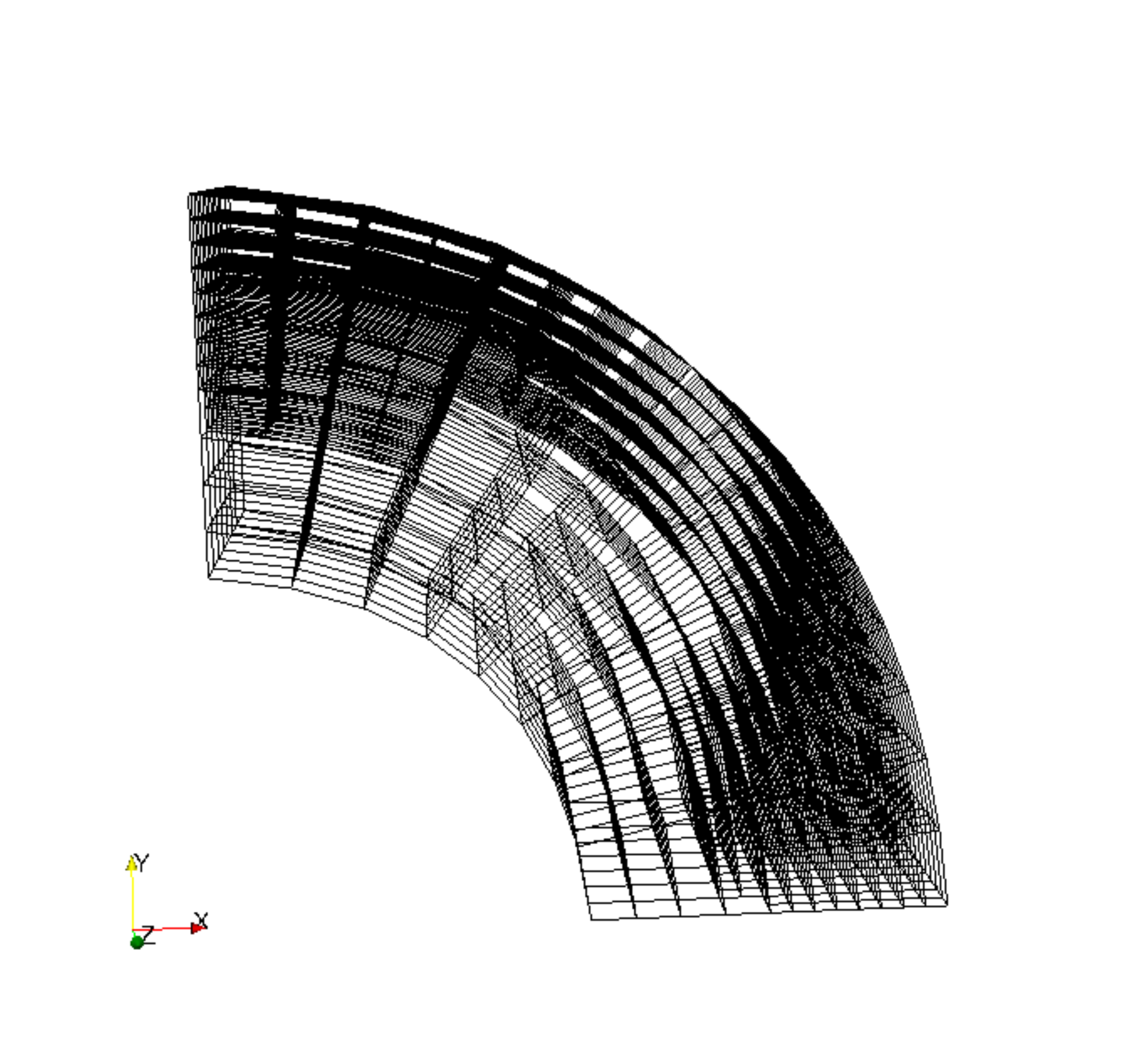}} 
%
	\caption{{\em Example 4}. Comparison of meshes on the physical and parametric domains w.r.t. adaptive refinement steps,
	 criterion ${\mathds{M}}_{{\rm \bf BULK}}(0.6)$.}
	\label{tab:quarter-annulus-example-4-meshes-on-param-and-phys-domain}
\end{figure}

{ 
All the numerical results presented below are obtained for the bulk marking criterion with the 
parameter $\sigma =0.4$.} 
Figure \ref{tab:quarter-annulus-example-4-meshes-on-param-and-phys-domain} presents an evolution 
of the adaptive meshes discretising the parametric space-time cylinder $\hat{Q}$ (left column) and 
the corresponding meshes discretising $Q$ (right column). 
{From the graphics presented,} we can see 
that the refinement is localised in the area close to the lateral surface of the quarter-annulus with the 
radius two. This can be explained by fast changes in the solution appearing close to this
`outer' surface, see $u$ at the time $t = \tfrac{2}{3}$ in Figure \ref{fig:quarter-annulus-example-4-b}. 

\begin{figure}[!t]
	\centering
	\captionsetup[subfigure]{oneside, labelformat=empty}
	\hskip15pt
	\subfloat[]{
	\includegraphics[width=4.6cm, trim={0cm 0cm 0cm 0cm}, clip]{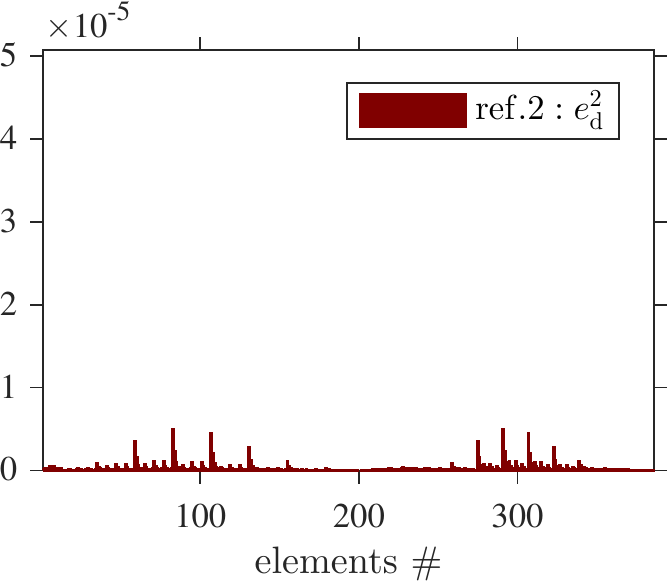}
	} 
	\hskip 5pt
	\subfloat[]{
	\includegraphics[width=4.6cm, trim={0cm 0cm 0cm 0cm}, clip]{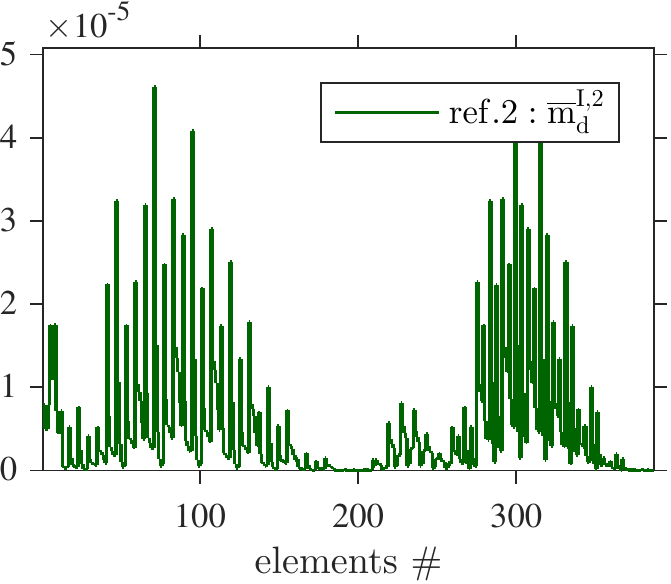}
	}
	\subfloat[]{
	\includegraphics[width=4.4cm, trim={0cm 0cm 0cm 0cm}, clip]{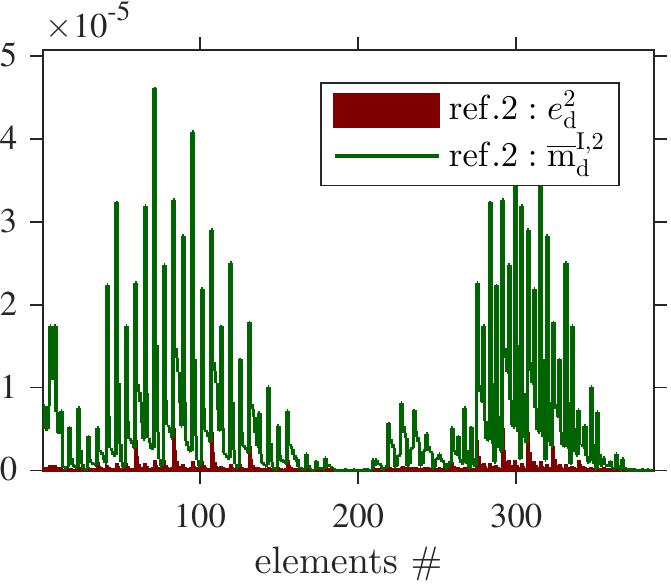}
	}
	\\[-20pt]
%
	%
	\subfloat[]{
	\includegraphics[width=4.8cm, trim={0cm 0cm 0cm 0cm}, clip]{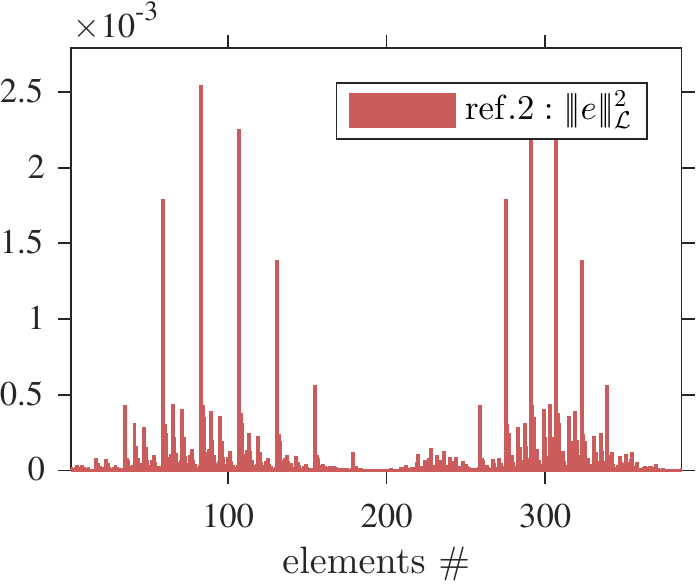}} 
	\hskip5pt
	\subfloat[]{
	\includegraphics[width=4.7cm, trim={0cm 0cm 0cm 0cm}, clip]{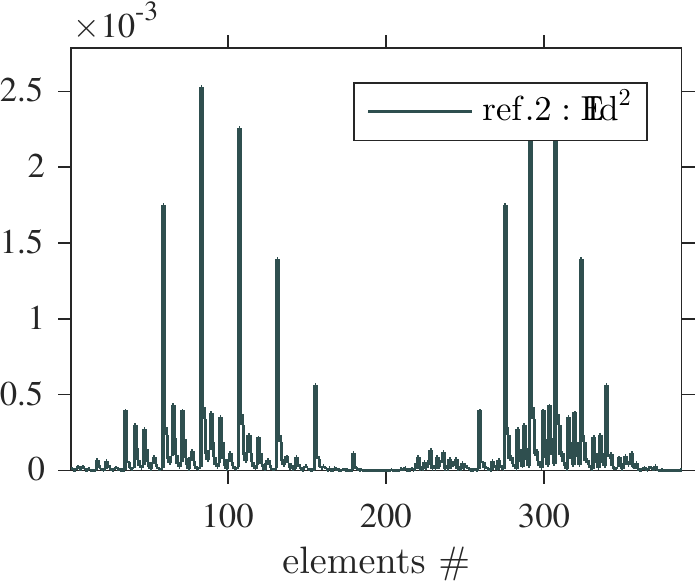}}
	\subfloat[]{
	\includegraphics[width=4.6cm, trim={0cm 0cm 0cm 0cm}, clip]{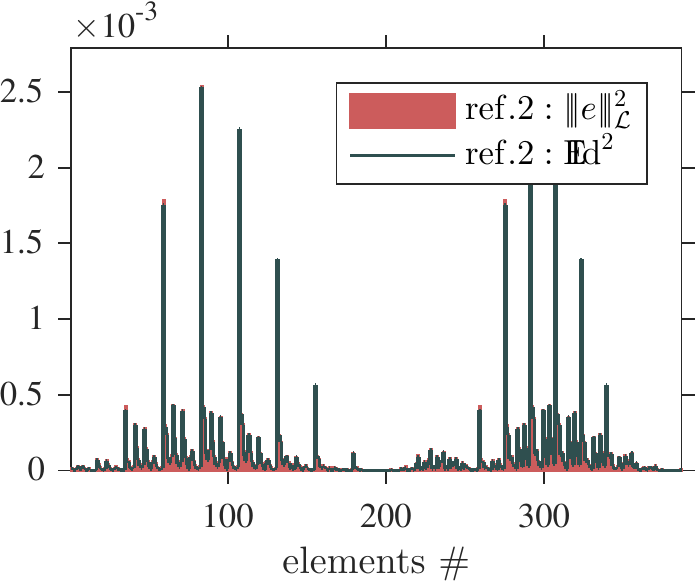}}
	\vskip -15pt
	\caption{{\em Example 4}. 
	Distribution of ${\rm e}_{{\rm d}, K}$, $\overline{\rm m}^{\rm I}_{{\rm d}, K}$ (first row)
	as well as $|\!|\!|  e |\!|\!|_{\mathcal{L}, K}$ and $\EI_K$ (second row) on the refinement step 2.}
	\label{fig:quarter-annulus-example-4-local-distr-ed-md}
\end{figure}

Finally, we provide a quantitative comparison of the local distributions ${\rm e}_{{\rm d}, K}$ and 
$\overline{\rm m}^{\rm I}_{{\rm d}, K}$ as well as $|\!|\!|  e |\!|\!|_{\mathcal{L}, K}$ and $\EI_K$ 
in Figure \ref{fig:quarter-annulus-example-4-local-distr-ed-md}. The first two columns of these graphics 
expose the quantities individually, and the last column contains plots with overlapped distributions 
of the error and the error indicator. In Figure \ref{fig:quarter-annulus-example-4-local-distr-ed-md}, we see 
that $\overline{\rm m}^{\rm I}_{{\rm d}, K}$ overestimates ${\rm e}_{{\rm d}, K}$, whereas 
the local indication of $\EI_K$ is sharper w.r.t the element-wise contributions 
$|\!|\!|  e |\!|\!|_{\mathcal{L}, K}$.

{ 
\subsection{Example 5 (solution with singularity w.r.t. $x$-coordinates)}
\label{ex:l-shape-2d-t-example-5}
\rm

\begin{figure}[!t]
	\centering
	\subfloat[$u(x, y, 1)$]{
	\includegraphics[scale=0.6]{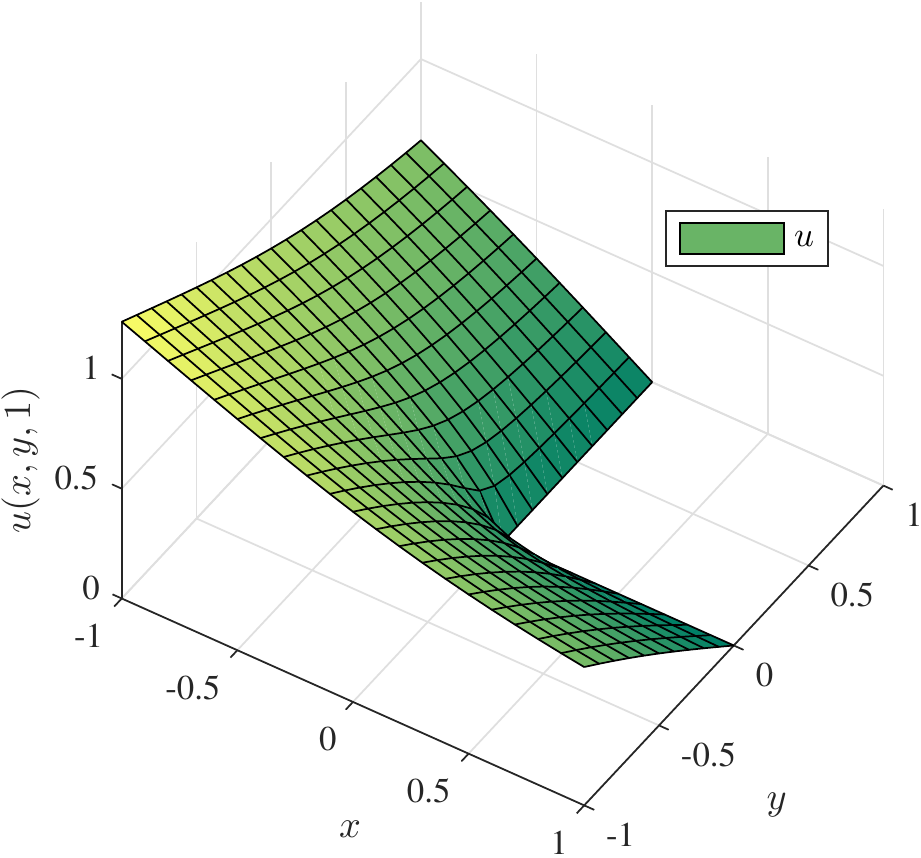}
	\label{fig:example-12-exact-solution}}
	\quad
	\subfloat[$\mathcal{K}_h$ defined on $\Omega$]{
	{\includegraphics[width=6.0cm, trim={3cm 0cm 3cm 0cm}, clip]{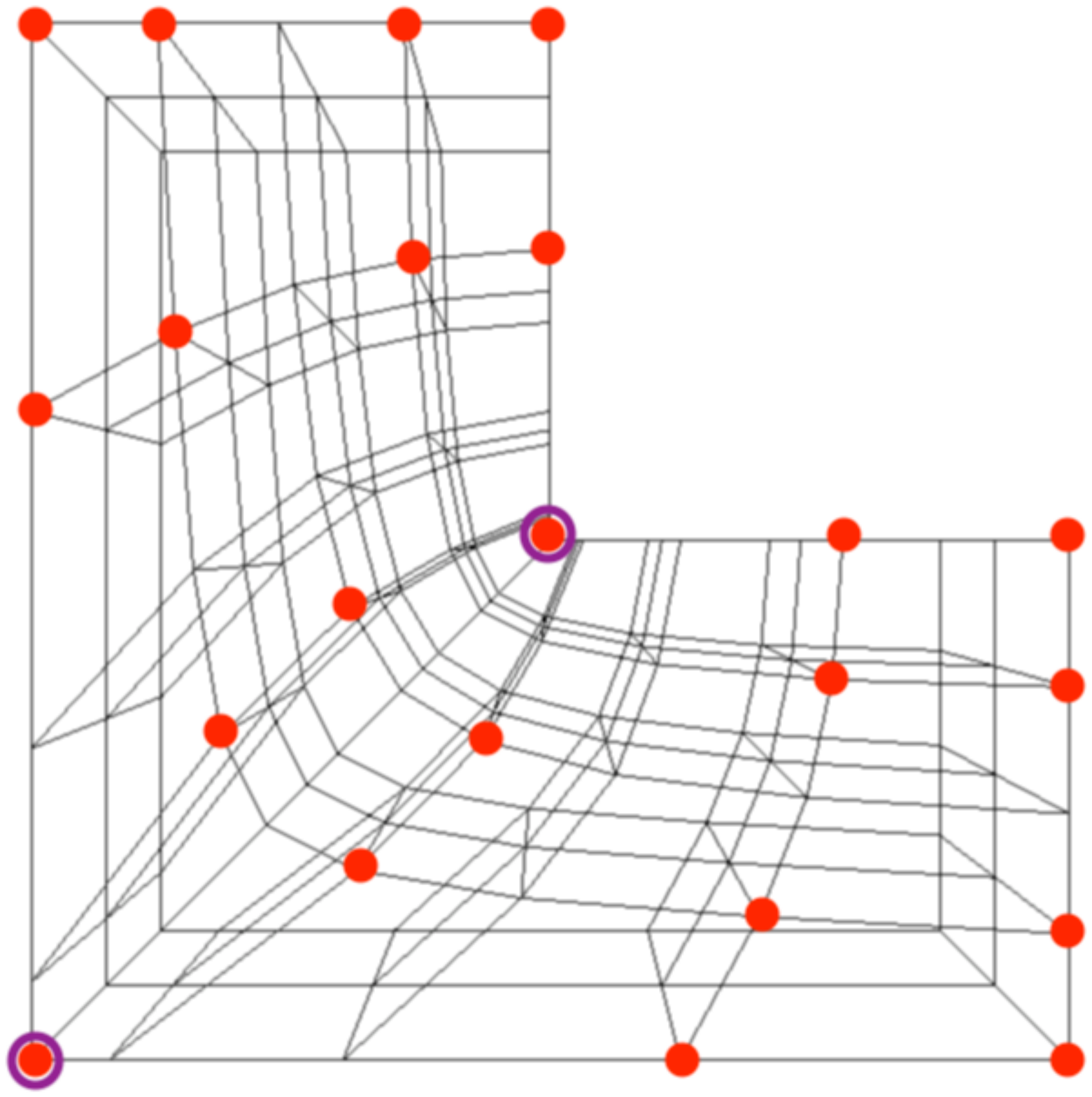}}
	\label{fig:initial-mesh-2}}	
	\caption{\small {\em Ex. \ref{ex:l-shape-2d-t-example-5}}. 
	(a) Exact solution $u = r^{1/3} \, \sin\theta \, \phi(1)$.
	(b) Initial geometry with Greville's points with double control points at the corners and 
	a corresponding mesh generated with $C^1$-continuous geometrical mapping.}
\end{figure}

In order to show how error estimates handle solutions with singularities, we consider a classical benchmark 
example on a two-dimensional L-shaped domain extended linearly in time to a cylinder $Q = \Omega \times (0, T)$, 
where 
$\Omega := \big( (-1, 1) \times (-1, 1) \big) \backslash \big( [0, 1) \times [0, 1) \big)$ 
and $T = 2$. The Dirichlet BCs are defined on $\Sigma$ by the 
{Dirichlet data}
$u_D = r^{1/3} \, \sin (\theta)$, 
where  
\begin{equation*}
r = {r(x,y) = } (x^2 + y^2) \quad \mbox{and} \quad 
\theta = {\theta(x,y) = } 
\begin{cases}
\tfrac13 \, (2\, {\rm atan2} (y, x) - \pi)\;\; \mbox{for} \;\;  y > 0, \\[-1pt]
\tfrac13 \, (2\,{\rm atan2} (y, x) + 3\,\pi)\;\; \mbox{for} \;\; y \leq 0.
\end{cases}
\end{equation*} 
The corresponding exact solution $$u(x, y, t) = r(x, y)^{1/3} \, \sin (\theta(x, y)) \, \phi(t), \quad \mbox{where} \quad 
\phi(t) = t^2 + t + 1, \quad  (x, y, t) \in \overline{Q}, $$ 
has a singularity in the point $(r, \theta) =(0, 0)$ (see 
Figure \ref{fig:example-12-exact-solution}). The RHS is given by 
$$f(x, y, t) = r(x, y)^{1/3} \, \sin (\theta(x, y)) \, \phi' (t),  \quad (x, y, t) \in Q.$$

Due to the doubled control points in the corners of the L-shape domain (denoted by red and purple circles in Figure 
\ref{fig:initial-mesh-2}), only the re-entrant corner and 
its counterpart on the other side are singular, i.e., the Jacobian of the geometry map in 
these two points is not regular. Since no integration points are placed in both corners, computational evaluation of the 
integrals remains valid. The downside of such a setting is that with the increase of refinement steps the cells near these 
corners become rather thin and lose shape-regularity. In addition, since on the functional level the requirements on the 
regularity of $\flux$ are not fulfilled, the global error estimate has rather a heuristic character, therefore, 
we only consider its performance from error indication point of view.

In Figure \ref{tab:l-shape-2d-t-example-5-meshes-on-phys-domains-different-angles}, we illustrate an evolution of 
adaptive meshes discretising physical domain (on LHS) and corresponding wireframes of meshes $\mathcal{K}_h$ (on RHS). 
L-shaped meshes (extended in time) on the right are illustrated from the point of view placed at the zero azimuth angle 
(located on the $O_z$-axis) 
in order to better see the nested THB-Splines. From both left and right columns, we can see that the refinement is 
localised in the area close to the singular point and no superfluous refinement is performed otherwise. 
Since the solution does not change in time drastically, 
the main refinement is concentrated in the area close to $(0, 0)$. The efficiency of the 
studied error bounds is also confirmed in Table \ref{tab:lshape-example-5-estimates-v-2-y-3-adaptive-ref}, 
which illustrates the decay of majorants and error identity w.r.t. refinement steps. Overall,  
$\overline{\rm M}^{\rm I}$ performs rather realistic, but even $\overline{\rm M}^{\rm I\!I}$ improves the first 
upper bound by approximately $1.5$ times. The e.o.c. (illustrated in the last column of Table 
\ref{tab:lshape-example-5-estimates-v-2-y-3-adaptive-ref}) is recovered due to the adaptive procedure dictated 
by $\overline{\rm M}^{\rm I}$.
 
\begin{table}[!t]
\scriptsize
\centering
\newcolumntype{g}{>{\columncolor{gainsboro}}c} 	
\newcolumntype{k}{>{\columncolor{lightgray}}c} 	
\newcolumntype{s}{>{\columncolor{silver}}c} 
\newcolumntype{a}{>{\columncolor{ashgrey}}c}
\newcolumntype{b}{>{\columncolor{battleshipgrey}}c}
\begin{tabular}{c|ckga|ck|c}
\parbox[c]{0.8cm}{\centering \# ref. } & 
\parbox[c]{1.4cm}{\centering  $\| \nabla_x e \|_Q$}   & 	  
\parbox[c]{1.4cm}{\centering $\overline{\rm M}^{\rm I}$ } &    
\parbox[c]{1.4cm}{\centering $\Ieff (\overline{\rm M}^{\rm I\!I})$ } & 
\parbox[c]{1.0cm}{\centering  $|\!|\!|  e |\!|\!|_{s, h}$ }   & 	  
\parbox[c]{1.4cm}{\centering $\Ieff (\overline{\rm M}^{\rm I}_{s, h})$ } & 
\parbox[c]{2cm}{\centering e.o.c. ($|\!|\!|  e |\!|\!|_{s, h}$)} \\
\bottomrule
\multicolumn{7}{l}{ \rule{0pt}{3ex}   
adaptive refinement, $\theta = 0.4$} \\[2pt]
\toprule
   2 &     1.7985e-01 &         2.54 &         1.68 &     1.7993e-01 &         2.54 &         2.24 \\
   3 &     1.2011e-01 &         2.45 &         1.64 &     1.2020e-01 &         2.45 &         1.47 \\
   4 &     7.2490e-02 &         2.57 &         1.64 &     7.2494e-02 &         2.57 &         1.49 \\
   5 &     4.8546e-02 &         2.53 &         1.62 &     4.8550e-02 &         2.53 &         1.47 \\
\end{tabular}
\caption{{\em Example 5}. 
Efficiency of $\overline{\rm M}^{\rm I}$, $\overline{\rm M}^{\rm I\!I}$, and 
$\overline{\rm M}^{\rm I}_{s, h}$ for $u_h \in S^{2}_{h}$,
$\flux_h \in \oplus^3 S_{h}^{3}$, and $w_h \in S^{3}_{h}$, 
w.r.t. adaptive refinement steps.}
\label{tab:lshape-example-5-estimates-v-2-y-3-adaptive-ref}
\end{table}

\begin{figure}[!t]
	\centering
	\captionsetup[subfigure]{oneside, labelformat=empty}
	\subfloat[ref. \# 2: \quad ${\Omega}$ and ${\mathcal{K}}_h$]{
	{\includegraphics[width=6cm, trim={2cm 0cm 1cm 4cm}, clip]{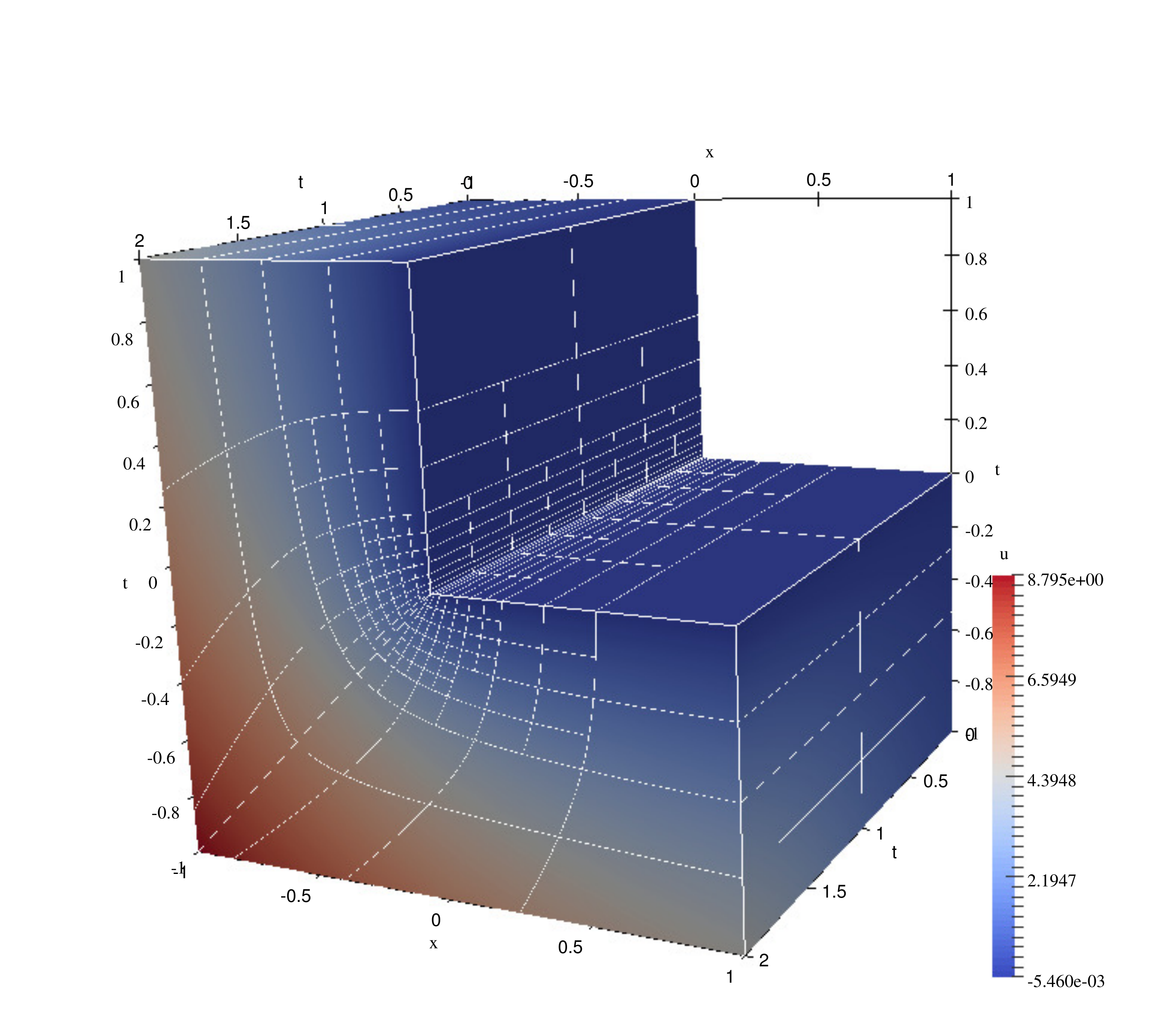}}
	}
	\qquad
	\subfloat[ref. \# 2: \quad ${\mathcal{K}}_h$, $O_z$-axis]{
	{\includegraphics[width=5cm, trim={3cm 5cm 3cm 5cm}, clip]{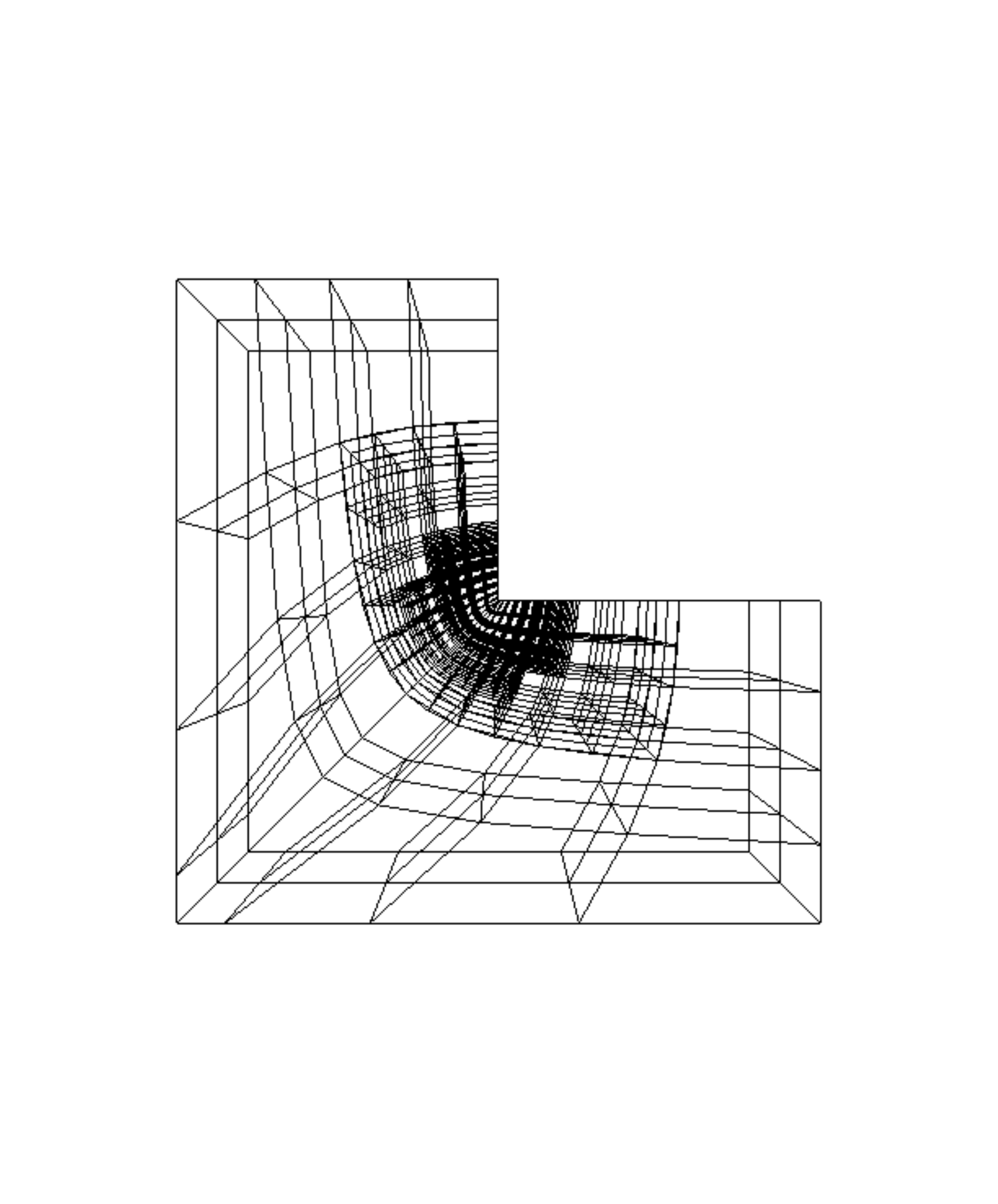}}
	} 
	\\[-1pt]
	\subfloat[ref. \# 4: \quad ${\Omega}$ and ${\mathcal{K}}_h$]{
	\includegraphics[width=6cm,  trim={2cm 0cm 1cm 4cm}, clip]{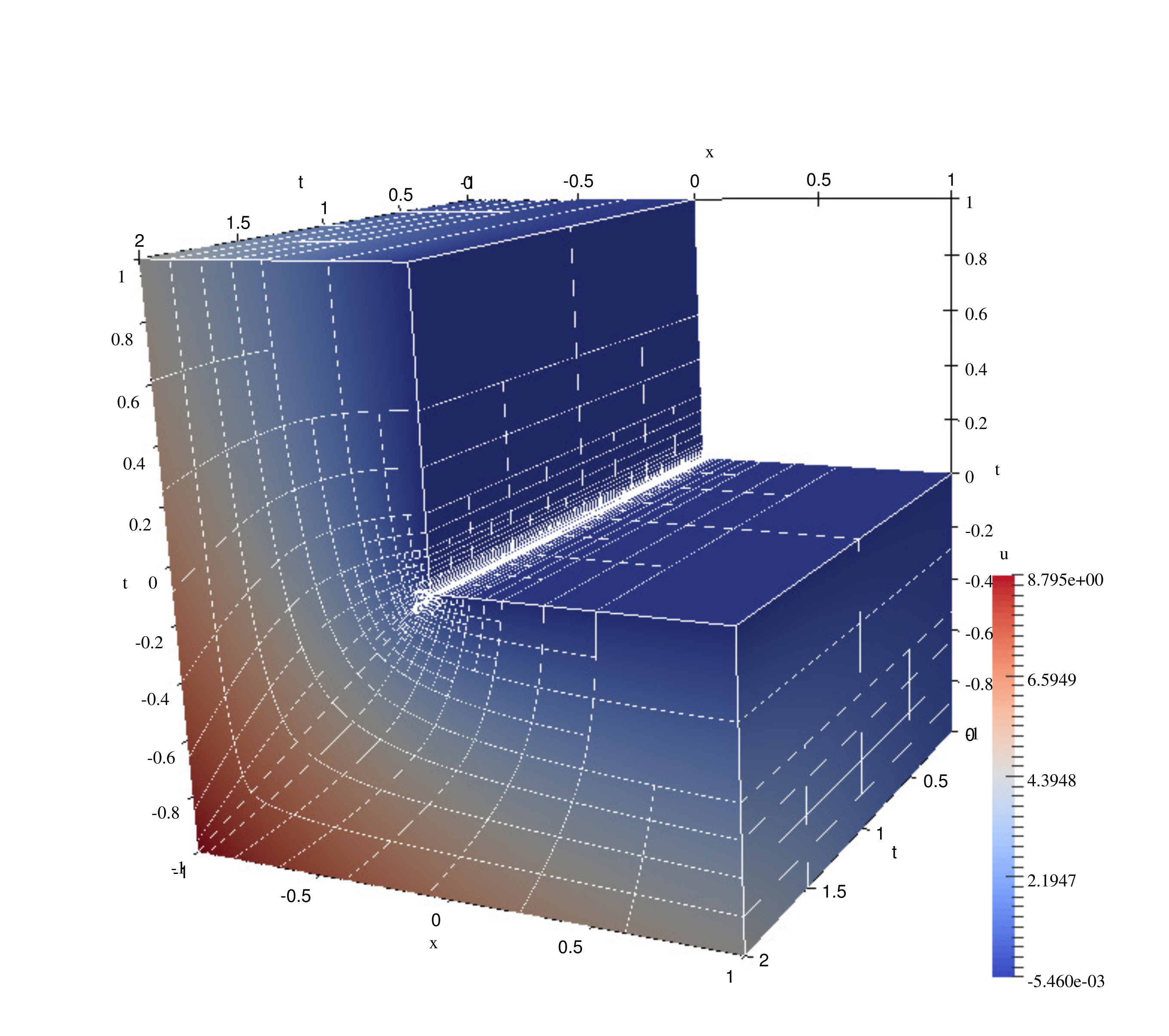}
	}
	\qquad
	\subfloat[ref. \# 4: \quad ${\mathcal{K}}_h$, $O_z$-axis]{
	\includegraphics[width=5cm, trim={3cm 5cm 3cm 5cm}, clip]{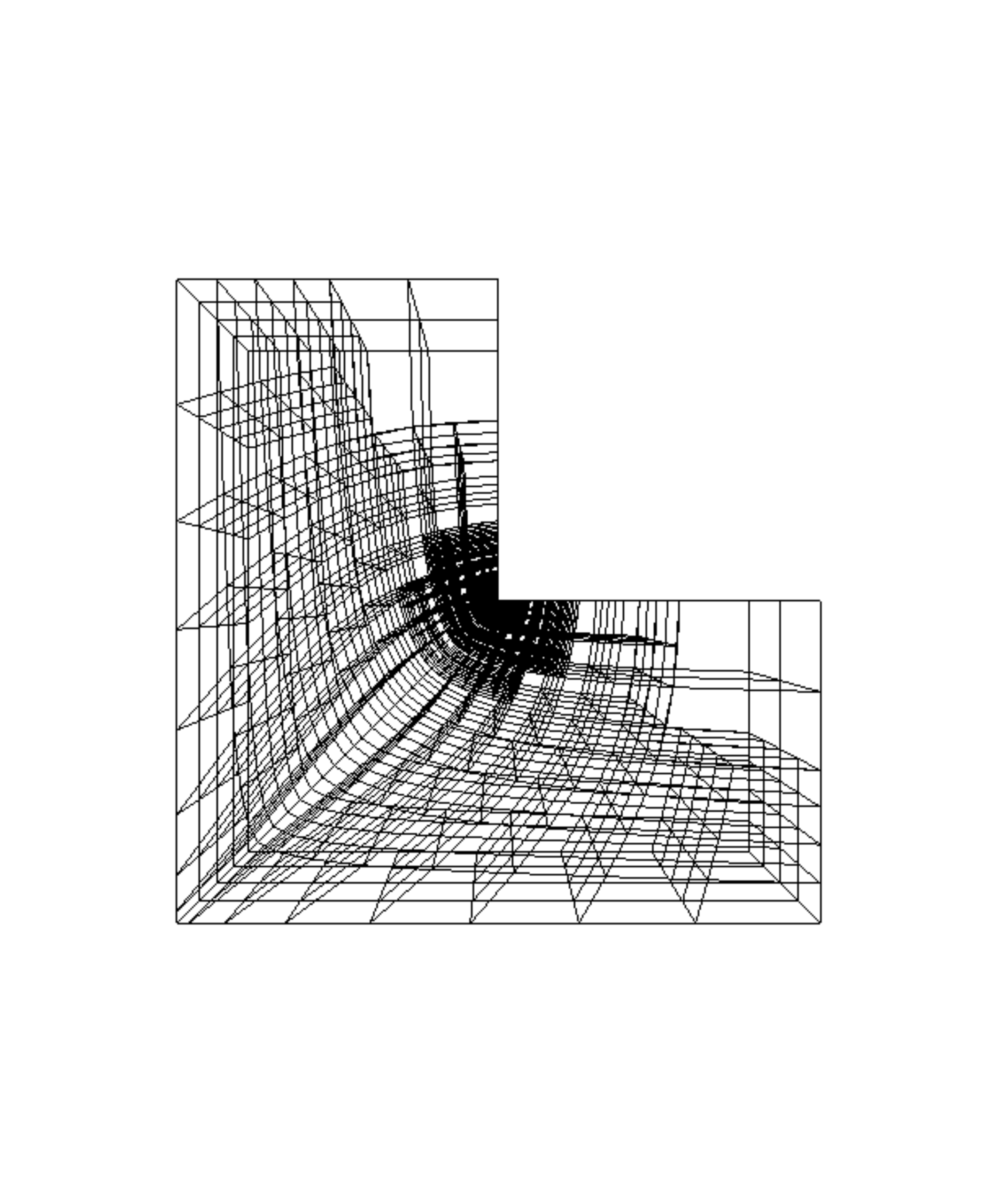}
	} 
	\caption{{\em Example 5}. Comparison of meshes on the physical domains w.r.t. adaptive refinement steps,
	 criterion ${\mathds{M}}_{{\rm \bf BULK}}(0.4)$.}
	\label{tab:l-shape-2d-t-example-5-meshes-on-phys-domains-different-angles}
\end{figure}	

Next, we assume that $u(x, y, t)$ has a more complicated dependence on $t$ and set 
$$\phi(t) = 10 \, (10\, t - \tfrac{1}{2}) \, (t - 1) \, (t - \tfrac{7}{4}).$$ Then, the 
space-time approach of solving the problem provides mesh-refinement in space and time 
automatically. Figure \ref{tab:l-shape-2d-t-example-5-2-meshes-on-phys-domains-different-angles} illustrates how 
$2$-dimensional spatial slices of the $3$-dimensional meshes involve w.r.t. time. 
Table \ref{tab:lshape-example-5-2-estimates-v-2-y-3-adaptive-ref} provides numerical evidence  
of how various error measures are estimated by $\overline{\rm M}^{\rm I}$, $\overline{\rm M}^{\rm I\!I}$, 
and $\overline{\rm M}^{\rm I}_{s, h}$.

\begin{table}[!t]
\scriptsize
\centering
\newcolumntype{g}{>{\columncolor{gainsboro}}c} 	
\newcolumntype{k}{>{\columncolor{lightgray}}c} 	
\newcolumntype{s}{>{\columncolor{silver}}c} 
\newcolumntype{a}{>{\columncolor{ashgrey}}c}
\newcolumntype{b}{>{\columncolor{battleshipgrey}}c}
\begin{tabular}{c|ckga|ck|c}
\parbox[c]{0.8cm}{\centering \# ref. } & 
\parbox[c]{1.4cm}{\centering  $\| \nabla_x e \|_Q$}   & 	  
\parbox[c]{1.4cm}{\centering $\overline{\rm M}^{\rm I}$ } &    
\parbox[c]{1.4cm}{\centering $\Ieff (\overline{\rm M}^{\rm I\!I})$ } & 
\parbox[c]{1.0cm}{\centering  $|\!|\!|  e |\!|\!|_{s, h}$ }   & 	  
\parbox[c]{1.4cm}{\centering $\Ieff (\overline{\rm M}^{\rm I}_{s, h})$ } & 
\parbox[c]{2cm}{\centering e.o.c. ($|\!|\!|  e |\!|\!|_{s, h}$)} \\
\bottomrule
\multicolumn{7}{l}{ \rule{0pt}{3ex}   
adaptive refinement, $\theta = 0.4$} \\[2pt]
\toprule
   2 &     9.3813e+00 &         1.66 &         1.44 &     1.0444e+01 &         1.49 &     4.52 \\
   3 &     5.7807e+00 &         1.74 &         1.20 &     5.9134e+00 &         1.70 &     1.61 \\
   4 &     1.8193e+00 &         2.92 &         1.36 &     1.8689e+00 &         2.85 &     2.51 \\
   5 &     1.0706e+00 &         2.44 &         1.40 &     1.0761e+00 &         2.43 &     1.43 \\
\end{tabular}
\caption{{\em Example 5}. 
Efficiency of $\overline{\rm M}^{\rm I}$, $\overline{\rm M}^{\rm I\!I}$, and 
$\overline{\rm M}^{\rm I}_{s, h}$ for $u_h \in S^{2}_{h}$,
$\flux_h \in \oplus^3 S_{h}^{3}$, and $w_h \in S^{3}_{h}$, 
w.r.t. adaptive refinement steps.}
\label{tab:lshape-example-5-2-estimates-v-2-y-3-adaptive-ref}
\end{table}

\begin{figure}[!t]
	\centering
	\captionsetup[subfigure]{oneside, labelformat=empty}
	\subfloat[${\mathcal{K}}_h$, $t = 0$]{
	{\includegraphics[width=5cm, trim={5.5cm 1.1cm 1cm 2cm}, clip]{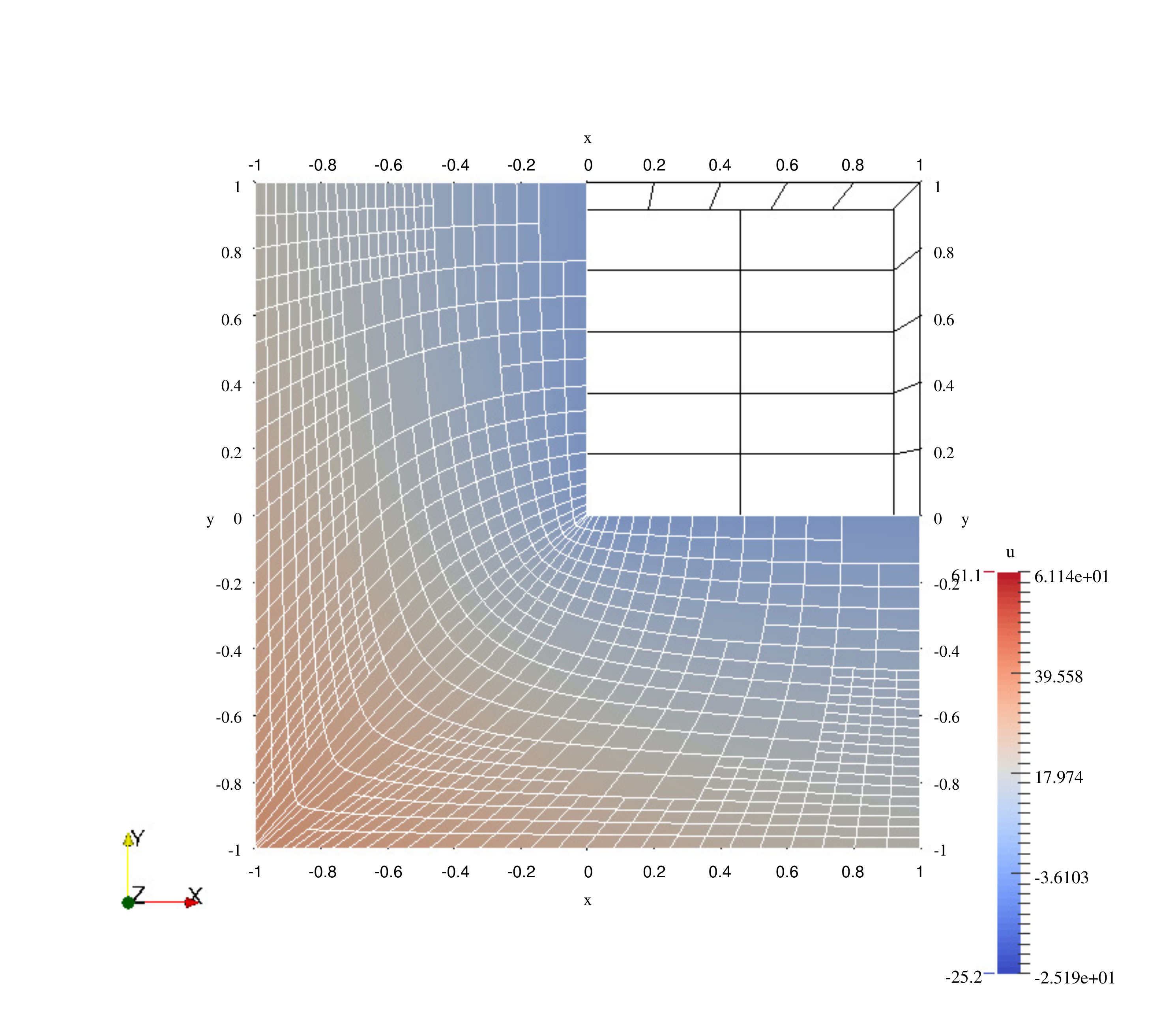}}
	} 
	\,
	\subfloat[${\mathcal{K}}_h$, $t = 1$]{
	{\includegraphics[width=5.4cm,  trim={6.5cm 0.5cm 1cm 2cm}, clip]{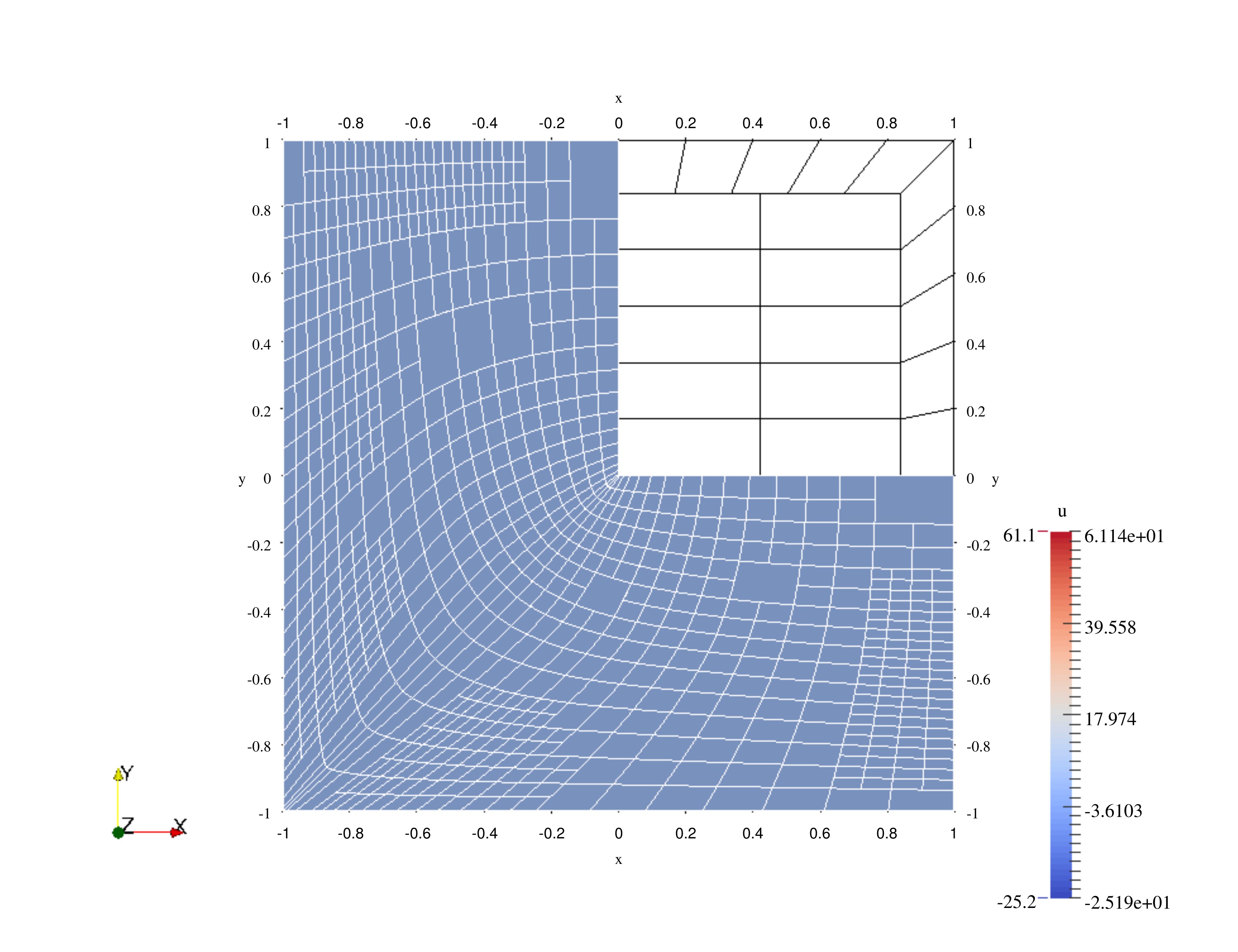}}
	}
	\,
	\subfloat[${\mathcal{K}}_h$, $t = 2$]{
	{\includegraphics[width=5.5cm, trim={6.5cm 1.1cm 1cm 2.5cm}, clip]{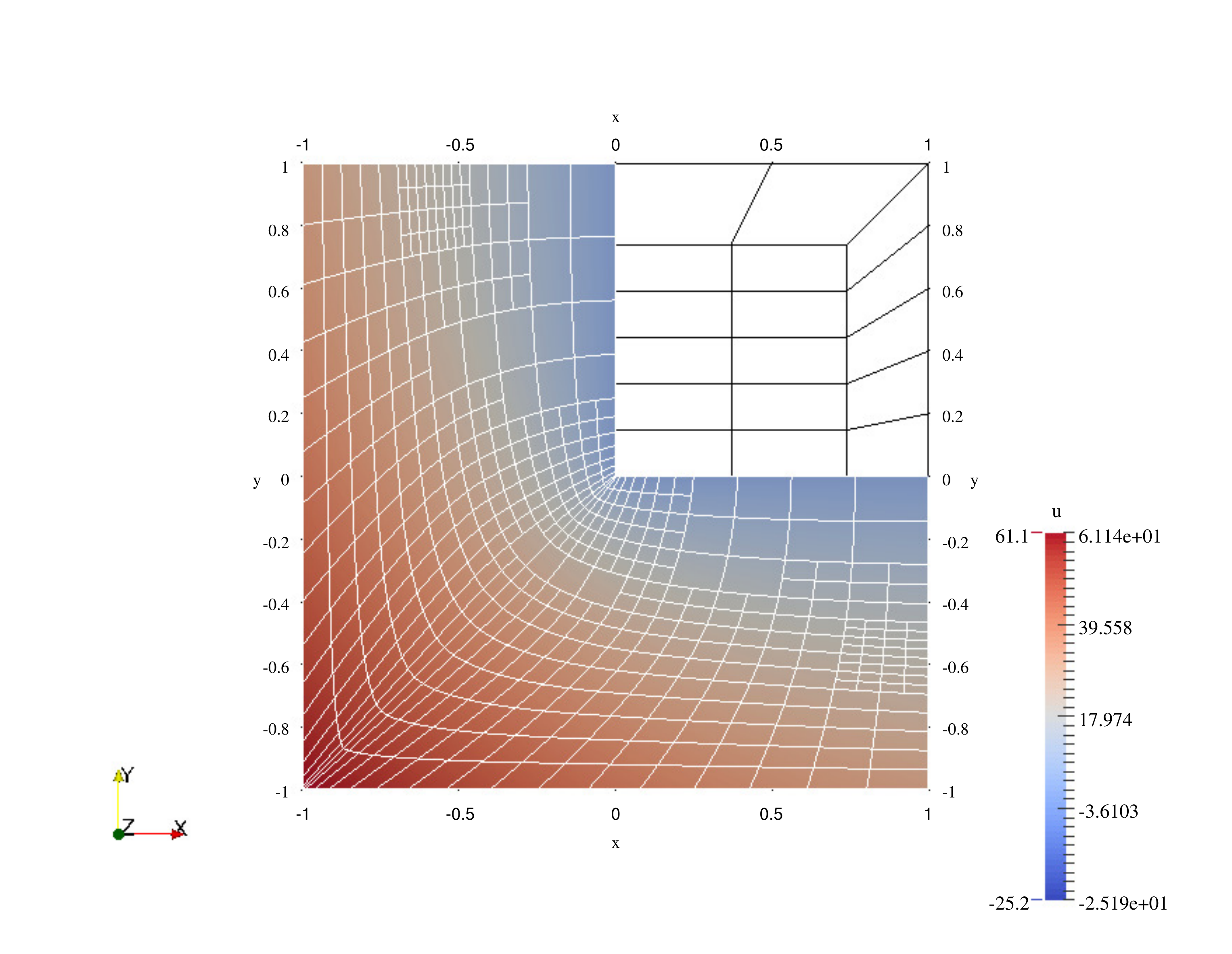}}
	} 
	\vskip -5pt
	\caption{{\em Example 5}. Comparison of meshes on the physical domains w.r.t. adaptive refinement steps,
	 criterion ${\mathds{M}}_{{\rm \bf BULK}}(0.4)$.}
	\label{tab:l-shape-2d-t-example-5-2-meshes-on-phys-domains-different-angles}
\end{figure}

\subsection{Example 6 (solution with singularity w.r.t. $t$-coordinate)}
\label{ex:time-singularity-1d-t-example-6}
\rm

For the last example, we assume that the solution has singularity w.r.t. time coordinates, i.e., we consider 
$$u(x, t) = \sin \pi x \, (1 - t)^\lambda \in H^{\ell_x, \ell_t}(Q), \ell_x, \ell_t \geq 0, \quad 
(x, t) \in \overline{Q} = (0, 1) \times (0, 2),$$
where parameter $\lambda = \Big\{ \tfrac{3}{2}, 1, \frac{1}{2} \Big\}$ 
(see Figure \ref{fig:time-singularity-1d-t-example-6} with $u$ for different $\lambda$). 
The RHS $f(x, t)$ follows from substitution of $u$ into \eqref{eq:equation}, and the Dirichlet boundary 
conditions are defined as {$u_D = u$ on $\Sigma$}.%

\begin{figure}[!t]
	\centering
	\subfloat[$\lambda = \tfrac{3}{2}$]{
	\includegraphics[scale=0.5]{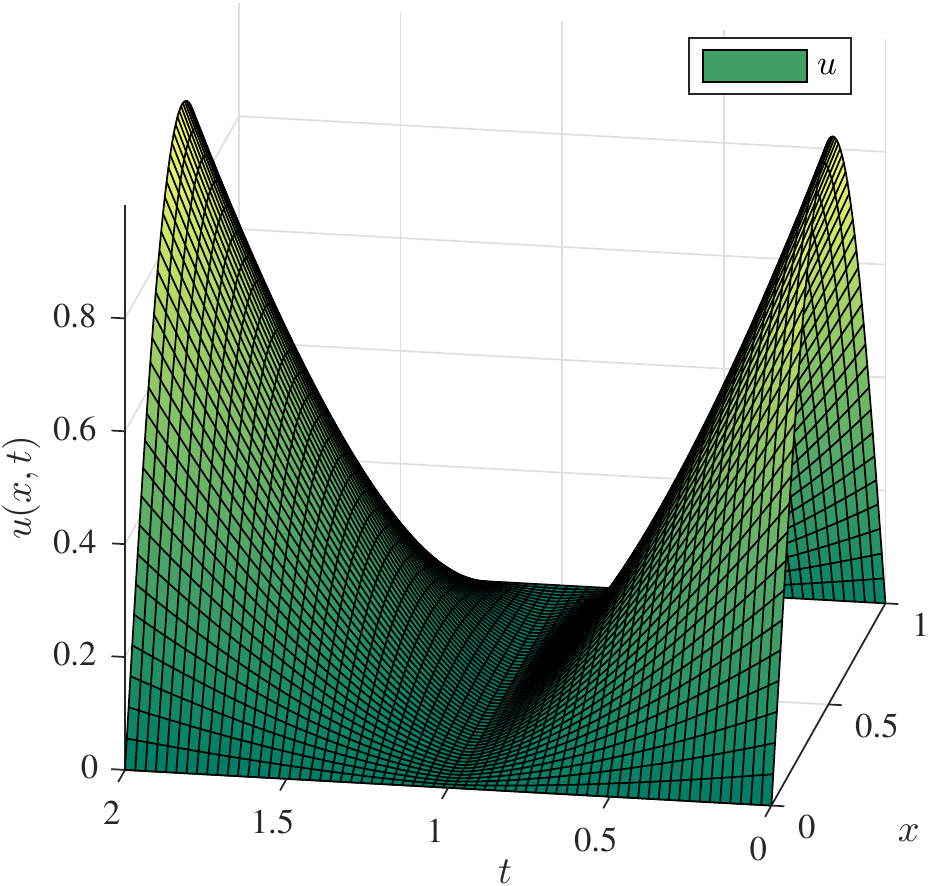}
	\label{fig:example-35-exact-solution}
	}
	\quad
	\subfloat[$\lambda = 1$]{
	\includegraphics[scale=0.5]{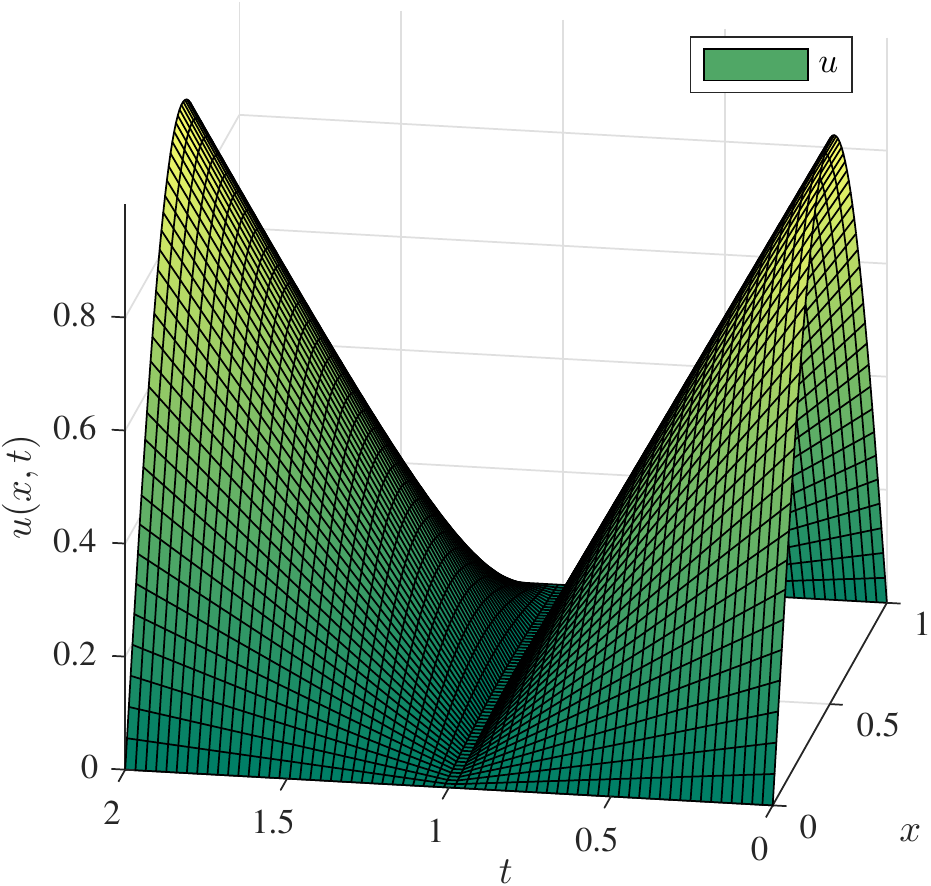}
	\label{fig:example-33-exact-solution}
	}	
	\quad
	\subfloat[$\lambda = \tfrac{1}{2}$]{
	\includegraphics[scale=0.5]{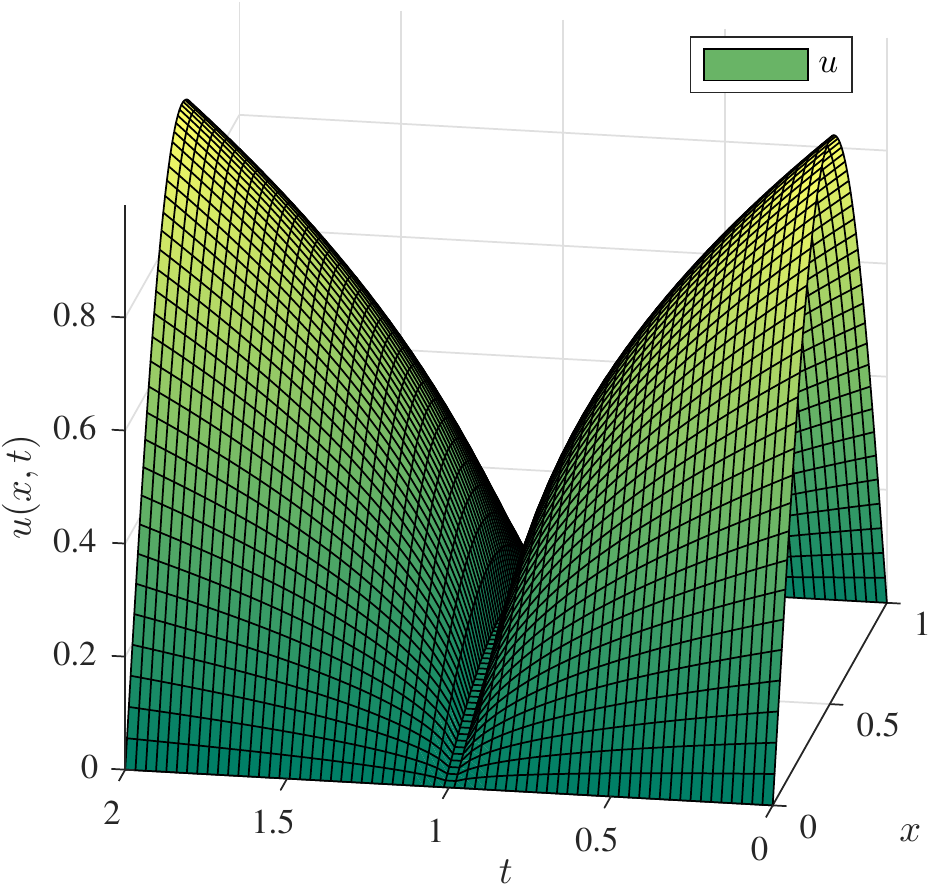}
	\label{fig:example-34-exact-solution}
	}	
	\caption{\small {\em Example 6}. 
	(a) Exact solution $u(x, t) = \sin \pi x \, (1 - t)^{\rfrac{3}{2}}$.
	(b) Exact solution $u(x, t) = \sin \pi x \, |1 - t|$. \\
	(c) Exact solution $u(x, t) = \sin \pi x \, (1 - t)^{\rfrac{1}{2}}$.}
	\label{fig:time-singularity-1d-t-example-6}
\end{figure}

The solution $u(x, t)$ is 
smooth w.r.t. to spatial coordinates, and the regularity in time depends on parameter $\lambda$. 
In particular, $\ell_t$ must satisfy the following inequality $\lambda \geq \ell_t - \tfrac{d}{p}$, 
where $d$ is the dimension of $\Omega$, and $p$ is a degree of splines used for the approximation 
of $u$ (in the current case, $d = 1$ and $p = 2$). Then, we obtain the relation $\ell_t \leq \lambda + \tfrac{1}{2}$, 
and the expected convergence in the term $h^{\rfrac{1}{2}} \, \| \partial_t (u - u_h) \|_{Q}$ is 
$(O(h^{\ell_x}) + O(h^{\ell_t - 1})) \cdot O(h^{\rfrac{1}{2}})$. These theoretical observations are 
confirmed by the numerical results presented in Table 
\ref{tab:time-singularity-1d-t-example-6-v-2-y-3-adaptive-ref}. The last column illustrates expected error 
order of convergence $O(h^{\ell_t - \rfrac{1}{2}})$. In particular, for the parameters 
$\lambda = \big\{ \frac{1}{2}, 1, \tfrac{3}{2} \big\}$, we expect $\ell_t \leq \big\{2, \tfrac{3}{2}, 1\big\}$, 
which provides approximated e.o.c. ${O (h^{\rfrac{3}{2}})}$, ${O (h^{1})}$, and 
${O (h^{\rfrac{1}{2}})}$, respectively. Table 
\ref{tab:time-singularity-1d-t-example-6-v-2-y-3-adaptive-ref} also provides the error's and estimates' decay 
in the case of adaptive refinement. The last column confirmes an improved e.o.c.

Figure \ref{fig:time-singularity-1d-t-example-6-meshes-on-phys-domains-different-lambda} presents meshes
obtained on the adaptive refinement steps 5 and 6 and reconfirms that functional error estimates detect the  
local singularities rather well. We see that for $\lambda = \tfrac{1}{2}$, the singularity at $t =1$ is captured 
and very well represented by the adaptive mesh. Moreover, for $t > 1$, where the solution is smooth, the mesh 
is not over-refined. 

\begin{table}[!t]
\scriptsize
\centering
\newcolumntype{g}{>{\columncolor{gainsboro}}c} 	
\newcolumntype{k}{>{\columncolor{lightgray}}c} 	
\newcolumntype{s}{>{\columncolor{silver}}c} 
\newcolumntype{a}{>{\columncolor{ashgrey}}c}
\newcolumntype{b}{>{\columncolor{battleshipgrey}}c}
\begin{tabular}{c|ca|ck|cc|g}
\parbox[c]{0.8cm}{\centering \# ref. } & 
\parbox[c]{1.4cm}{\centering  $\| \nabla_x e \|_Q$}   & 	  
\parbox[c]{1.4cm}{\centering $\Ieff (\overline{\rm M}^{\rm I\!I})$ } & 
\parbox[c]{1.0cm}{\centering  $|\!|\!|  e |\!|\!|_{s, h}$ }   & 	  
\parbox[c]{1.4cm}{\centering $\Ieff (\overline{\rm M}^{\rm I}_{s, h})$ } & 
\parbox[c]{1.0cm}{\centering  $|\!|\!|  e |\!|\!|_{\mathcal{L}}$ }   & 	  
\parbox[c]{1.4cm}{\centering$\Ieff ({\EI})$ } & 
\parbox[c]{1.2cm}{\centering e.o.c. ($|\!|\!|  e |\!|\!|_{s, h}$) } \\
\toprule
\multicolumn{8}{l}{$\lambda = \tfrac{3}{2}$}\\
\bottomrule
\multicolumn{7}{l}{ \rule{0pt}{3ex} uniform refinement}  & expected ${O \big(h^{\rfrac{3}{2}}\big)}$\\
\midrule
   4 &     3.7789e-03 &         
   1.13 &     8.6838e-03 &         9.05 &     2.8189e-01 &         1.00 &     2.05 \\
   5 &     9.5349e-04 &         
   1.52 &     2.8772e-03 &         9.73 &     1.4075e-01 &         1.00 &     1.74 \\
   6 &     2.4077e-04 &        
   2.28 &     9.9310e-04 &        10.14 &     7.0358e-02 &         1.00 &     1.60 \\
\midrule
\multicolumn{7}{l}{adaptive refinement, $\theta = 0.4$} &  improved e.o.c. \\
\midrule
   3 &     1.7464e-02 &        
   1.00 &     1.7464e-02 &         2.10 &     5.6998e-01 &         1.00 &     2.50  \\
   4 &     5.6537e-03 &         
   1.43 &     5.6537e-03 &         3.47 &     2.8936e-01 &         1.00 &     2.08 \\
   6 &     1.7847e-03 &         
   1.48 &     1.7847e-03 &         3.43 &     1.6984e-01 &         1.00 &     2.32 \\
   7 &     7.0591e-04 &         
   1.65 &     7.0591e-04 &         4.00 &     1.3543e-01 &         1.00 &     2.16 \\
\toprule
\multicolumn{8}{l}{$\lambda = 1$}\\
\bottomrule
\multicolumn{7}{l}{ \rule{0pt}{3ex} uniform refinement} & expected ${O(h)}$\\
\midrule
%
   4 &     1.6706e-02 &         
    1.15 &     5.0754e-02 &         2.89 &     3.7940e-01 &         1.00 &     1.32 \\
   5 &     5.9650e-03 &         
   1.71 &     2.4711e-02 &         2.88 &     2.1176e-01 &         1.00 &     1.13 \\
   6 &     2.1355e-03 &        
   2.86 &     1.2256e-02 &         3.31 &     1.2604e-01 &         1.00 &     1.06 \\
\midrule
\multicolumn{7}{l}{adaptive refinement, $\theta = 0.4$} & improved e.o.c. \\
\midrule
%
%
   3 &     5.8230e-02 &         
   1.06 &     5.8230e-02 &         2.39 &     7.4729e-01 &         0.99 &     1.73\\
   4 &     2.7175e-02 &         
   1.42 &     2.7175e-02 &         3.75 &     5.4388e-01 &         1.00 &     1.64\\
   5 &     2.5989e-02 &         
   1.35 &     2.5989e-02 &         3.27 &     5.8875e-01 &         1.00 &     0.10\\
   6 &     1.2191e-02 &         
   1.69 &     1.2191e-02 &         4.56 &     5.1917e-01 &         1.00 &     1.85\\
\toprule
\multicolumn{8}{l}{$\lambda = \tfrac{1}{2}$} \\
\bottomrule
\multicolumn{7}{l}{ \rule{0pt}{3ex} uniform refinement} & expected ${O\big(h^{\rfrac{1}{2}}\big)}$\\
\midrule
   4 &     5.3129e-02 &         
   1.40 &     2.0530e-01 &         2.00 &     9.0200e-01 &         1.00 &     0.67\\
   5 &     2.6997e-02 &        
   2.16 &     1.4296e-01 &         2.29 &     8.2282e-01 &         1.00 &     0.57\\
   6 &     1.3712e-02 &        
   3.67 &     1.0061e-01 &         2.96 &     8.0464e-01 &         1.00 &     0.53\\
\midrule
\multicolumn{7}{l}{adaptive refinement, $\theta = 0.4$} & improved e.o.c. \\
\midrule
   5 &     5.8991e-02 &         
   2.05 &     5.8991e-02 &         6.13 &     2.2131e+00 &         1.00 &     0.97\\
   6 &     3.8622e-02 &         
   3.01 &     3.8622e-02 &         9.93 &     2.3304e+00 &         1.00 &     1.11\\
   7 &     2.5883e-02 &        
   3.88 &     2.5883e-02 &        14.66 &     2.6588e+00 &         1.00 &     1.15\\
   8 &     2.0120e-02 &        
   4.24 &     2.0120e-02 &        16.98 &     3.4281e+00 &         1.00 &     0.70\\
\end{tabular}
\caption{{\em Example 6}. 
Efficiency of $\overline{\rm M}^{\rm I\!I}$, 
$\overline{\rm M}^{\rm I}_{s, h}$, and ${\EI}$ for $u_h \in S^{2}_{h}$,
$\flux_h \in \oplus^3 S_{h}^{3}$, and $w_h \in S^{3}_{h}$, 
w.r.t. uniform refinement and adaptive refinement steps.}
\label{tab:time-singularity-1d-t-example-6-v-2-y-3-adaptive-ref}
\end{table}

\begin{figure}[!t]
	\centering
	\captionsetup[subfigure]{oneside, labelformat=empty}
	%
	%
	\subfloat[ref. \# 4: \quad ${\mathcal{K}}_h$, $\lambda = \tfrac{3}{2}$]{
	{\includegraphics[width=4.6cm, trim={2cm 4cm 2cm 4cm}, clip]{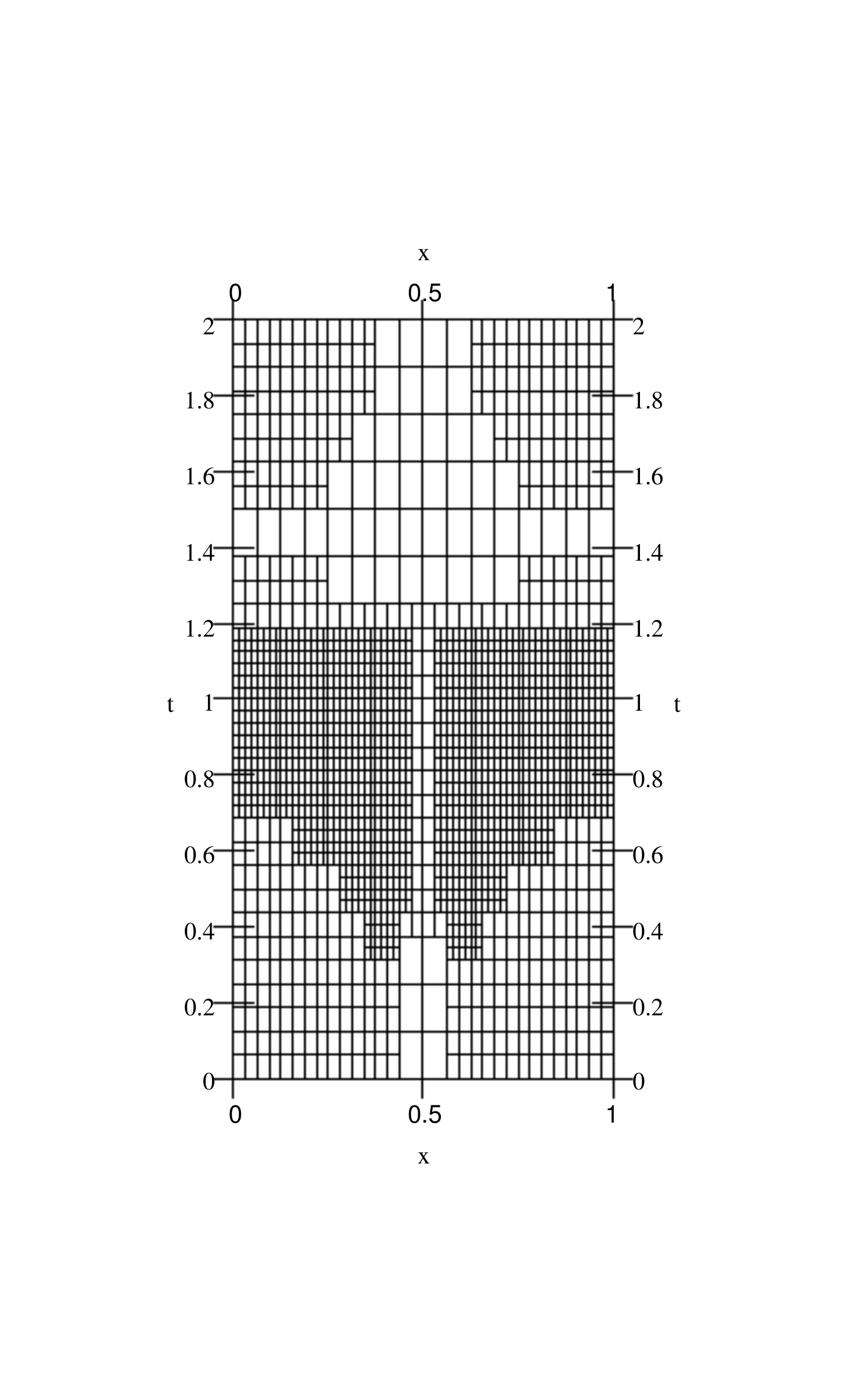}}}
	\quad
	\subfloat[ref. \# 4: \quad ${\mathcal{K}}_h$, $\lambda = 1$]{
	\includegraphics[width=4.6cm, trim={2cm 4cm 2cm 4cm}, clip]{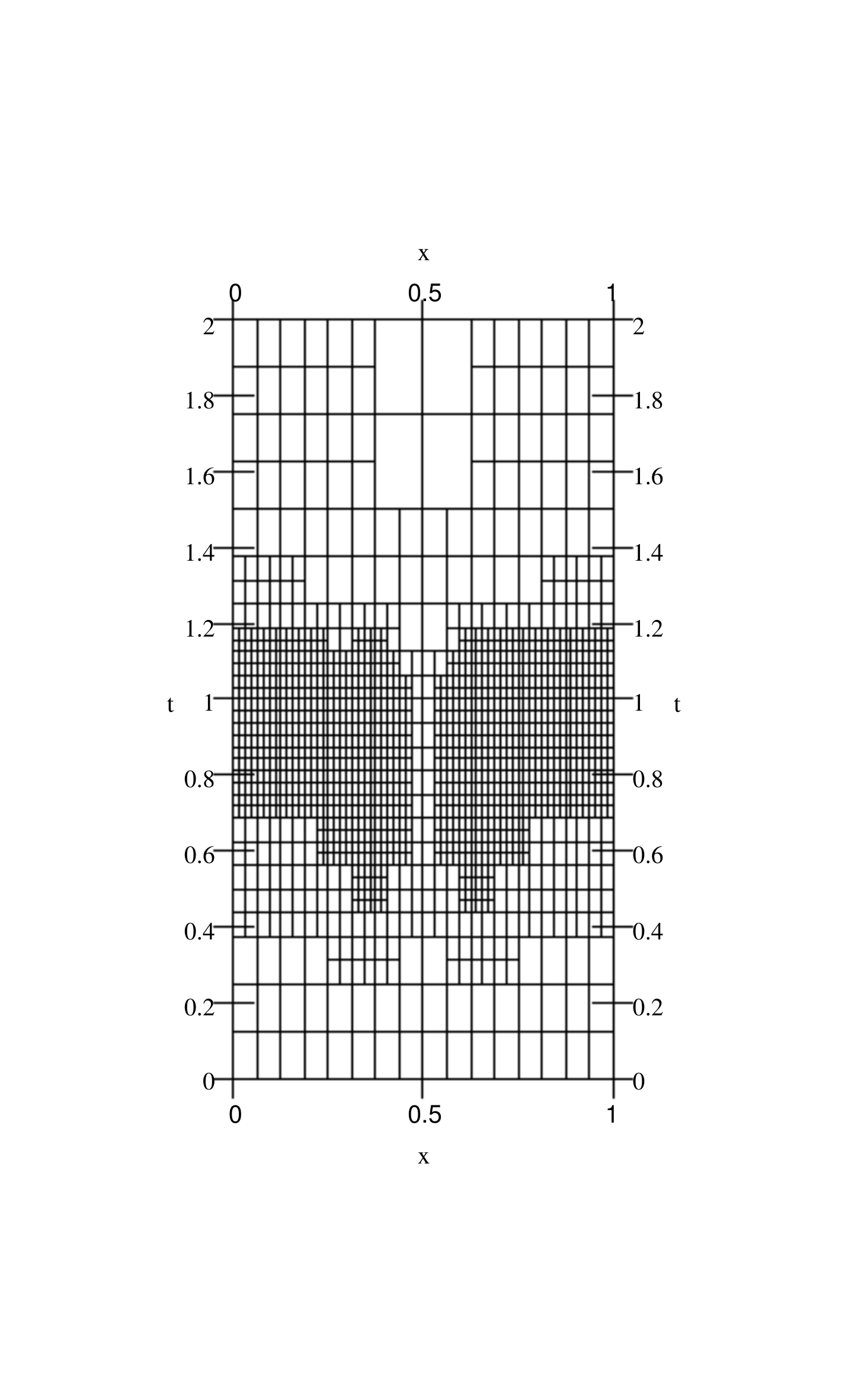}}
	\quad
	\subfloat[ref. \# 4: \quad ${\mathcal{K}}_h$, $\lambda = \tfrac{1}{2}$]{
	\includegraphics[width=4.6cm, trim={2cm 4cm 2cm 4cm}, clip]{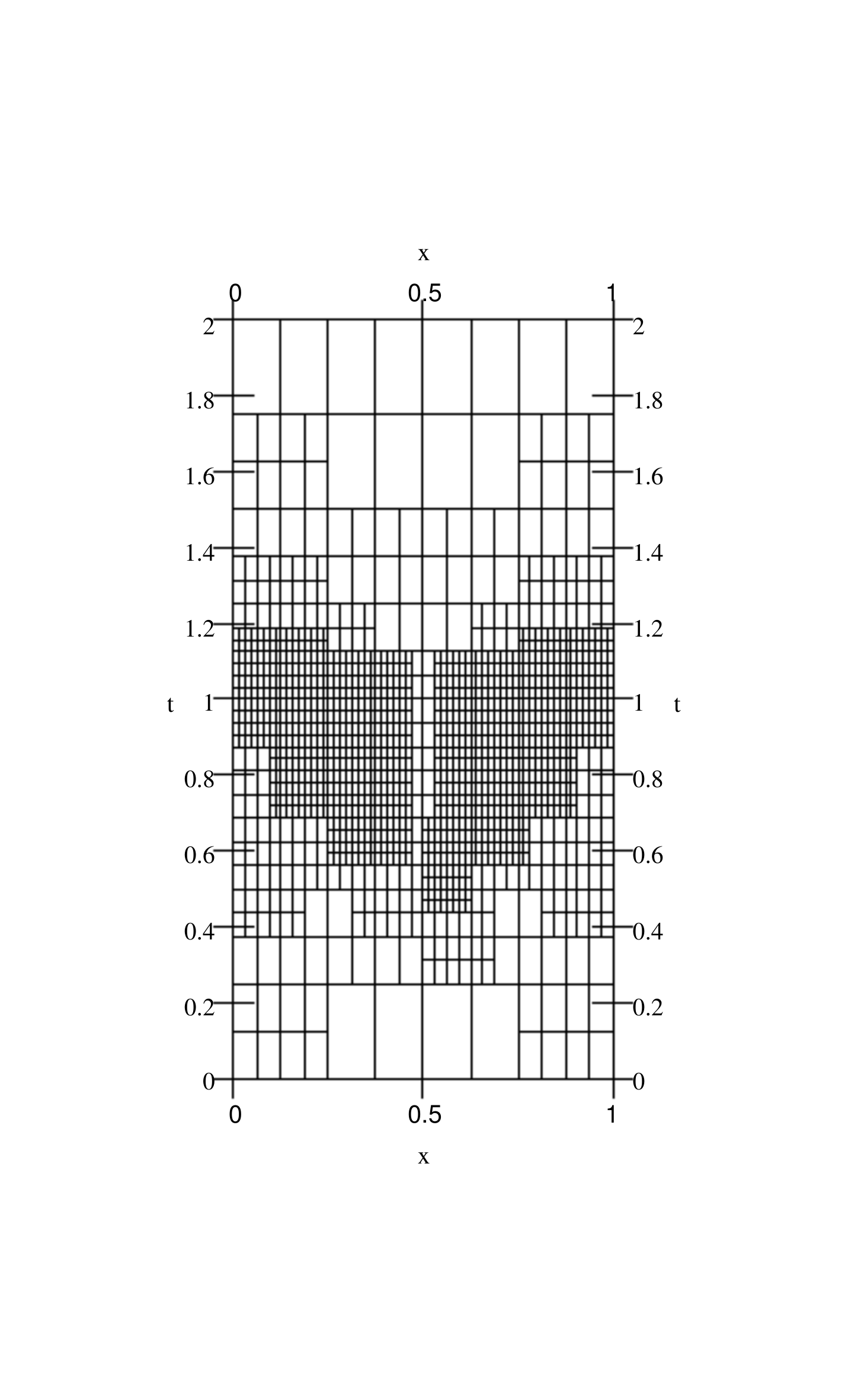}} 
	\\[-5pt]
	%
	\subfloat[ref. \# 6: \quad ${\mathcal{K}}_h$, $\lambda = \tfrac{3}{2}$]{
	{\includegraphics[width=4.6cm, trim={2cm 4cm 2cm 4cm}, clip]{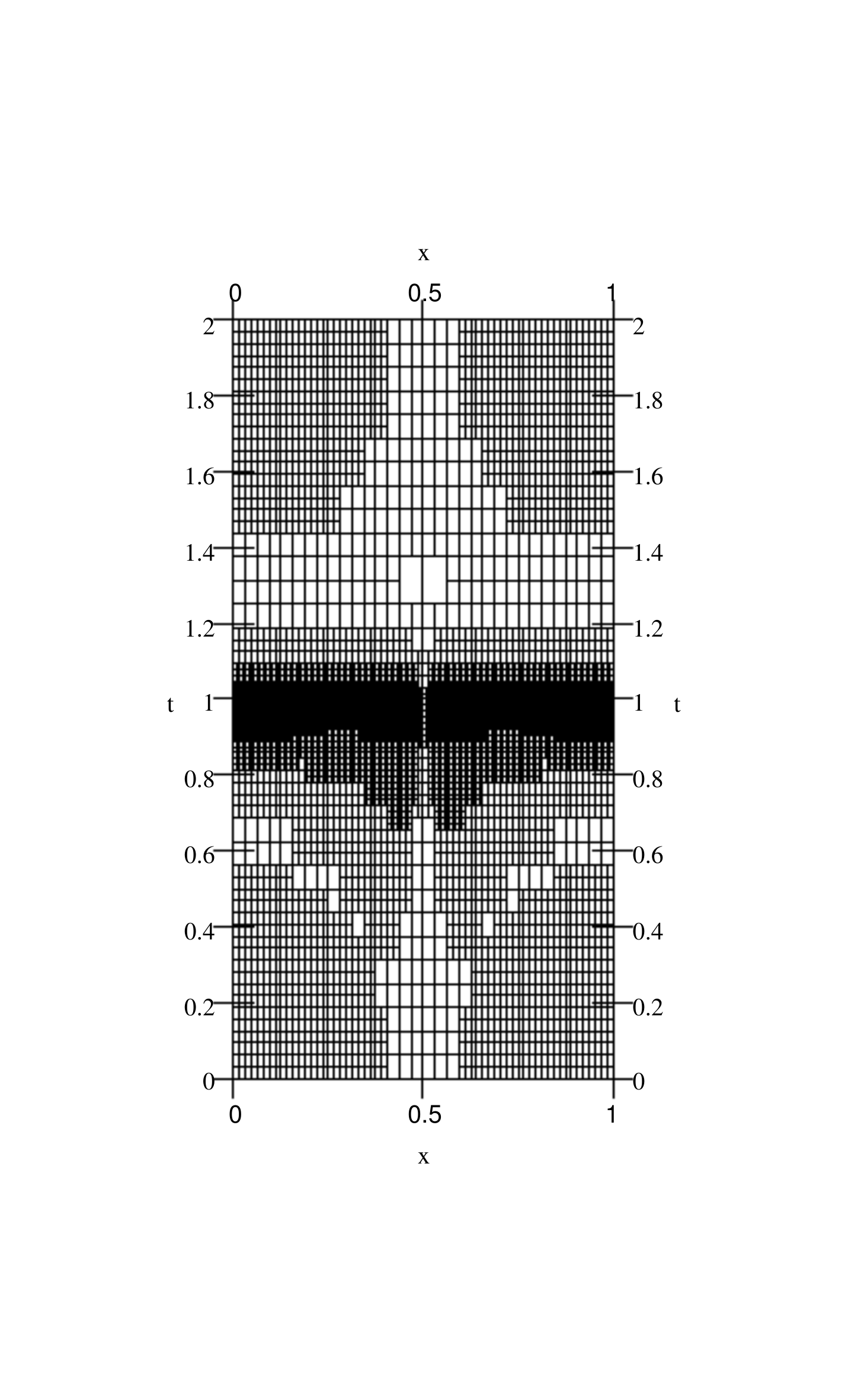}}}
	\quad
	\subfloat[ref. \# 6: \quad ${\mathcal{K}}_h$, $\lambda = 1$]{
	\includegraphics[width=4.6cm, trim={2cm 4cm 2cm 4cm}, clip]{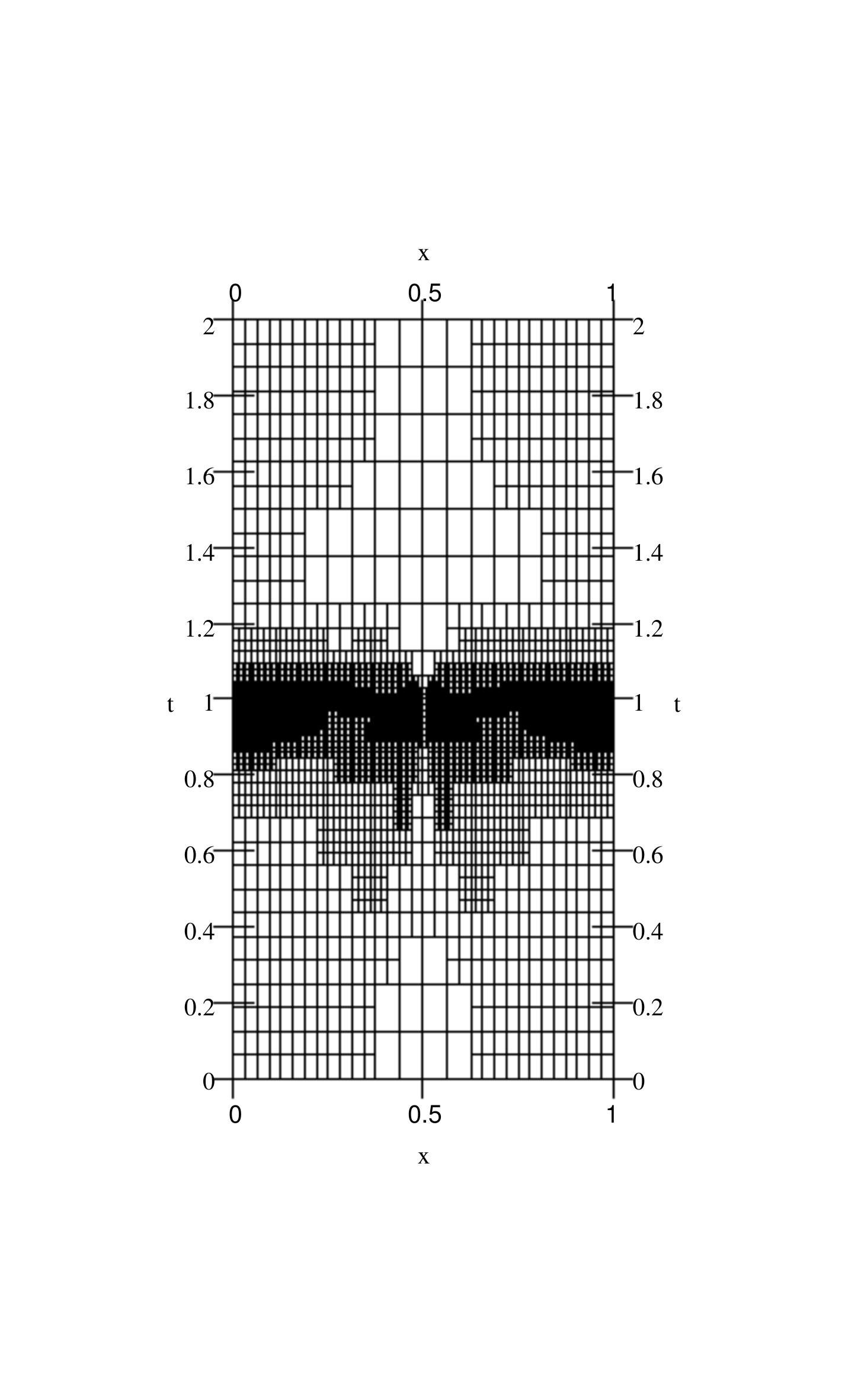}}
	\quad
	\subfloat[ref. \# 6: \quad ${\mathcal{K}}_h$, $\lambda = \tfrac{1}{2}$]{
	\includegraphics[width=4.6cm, trim={2cm 4cm 2cm 4cm}, clip]{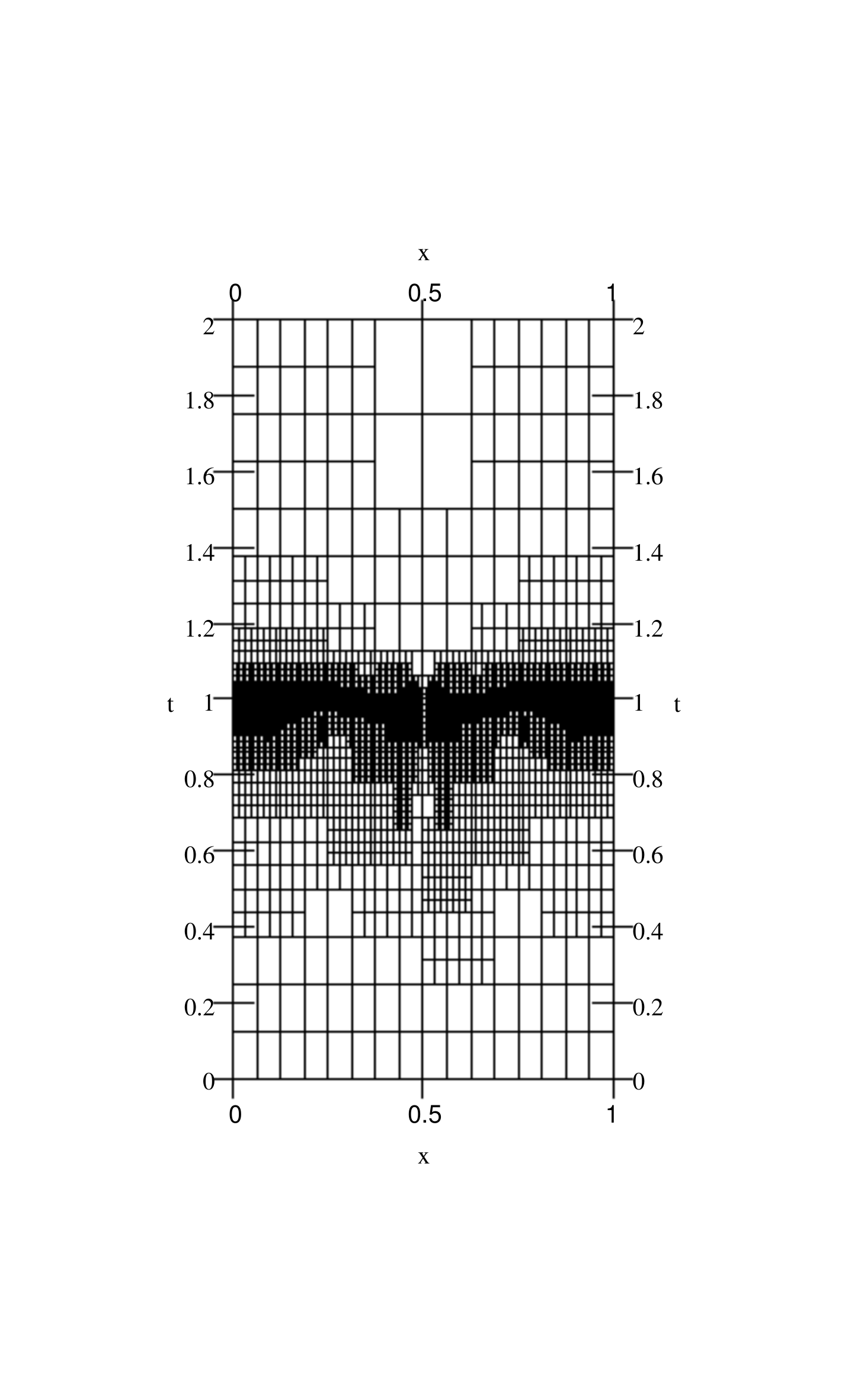}} 
	\caption{{\em Example 6}. Comparison of meshes on the physical and parametric domains w.r.t. adaptive refinement steps,
	 criterion ${\mathds{M}}_{{\rm \bf BULK}}(0.4)$.}
	\label{fig:time-singularity-1d-t-example-6-meshes-on-phys-domains-different-lambda}
\end{figure}

\section{Conclusions}
\label{sec:conclusion}

In the paper, we derived reliable space-time IgA schemes 
for parabolic initial-boundary value problems. 
In particular, we deduced {new} functional-type a posteriori error estimates 
and showed their efficient implementation in space-time IgA.
Since the derivation is based on purely functional arguments, the estimates 
are valid for any approximation from the admissible (energy) class.
They imply a posteriori error estimates
for mesh-dependent norms associated with stabilised space-time IgA schemes. 
We also proposed an efficient technique for minimising the majorant 
leading to extremely accurate guaranteed upper bounds 
of the error norm with efficiency indices close to 1.
Since this upper bound is nothing but the sum of the local contributions, 
these local contributions were be used as error indicators for mesh 
refinement. Mesh refinement in IgA is more involved than in the finite 
element method. We used THB-Splines for mesh refinement in our fully unstructured 
space-time adaptive IgA scheme.
Finally, we illustrated the reliability and efficiency of the presented a posterior 
error estimates for IgA solutions to several examples exhibiting different features
(defined on the domains with non-trivial shape, having solutions that possess singularity 
w.r.t. space and time variables).
We also reported about the cost of computing the upper bound. 
In all our examples, this was only a small portion of the cost for computing
the IgA solution. Last but not least, the numerical examples showed that 
the space-time THB-spline-based adaptive procedure works very efficiently.

}
\vskip 10pt
{\bf Acknowledgments}
The research is supported by the Austrian Science Fund (FWF) through the NFN S117-03 project. The implementation was carried out using the open-source C++ library G+Smo \cite{gismoweb}.

\bibliographystyle{plain}
\bibliography{bib/lib}

\end{document}

%% file: customized_packages_and_macros.tex
\usepackage{amsmath}
\usepackage{graphicx}								
\usepackage{subfig}
\usepackage{amsfonts} 							
\usepackage{dsfont}
\usepackage{algorithm}              
\usepackage{algorithmic}
\usepackage{booktabs} 							
\usepackage{multirow}								
\usepackage{setspace}
\usepackage{mathrsfs}
\usepackage{amssymb}
\usepackage{tikz}
\usepackage{mathtools}
\usepackage{xfrac}
\usepackage{upgreek}

\usepackage{framed}
\usepackage{afterpage}
\usepackage{capt-of}

\usepackage{geometry}
\usepackage{color}
\definecolor{gray(x11gray)}{rgb}{0.75, 0.75, 0.75}

\usepackage{tikz}	
\usetikzlibrary{arrows,chains,matrix,positioning,scopes}

\makeatletter
\tikzset{join/.code=\tikzset{after node path={%
\ifx\tikzchainprevious\pgfutil@empty\else(\tikzchainprevious)%
edge[every join]#1(\tikzchaincurrent)\fi}}}

\makeatother
\tikzset{>=stealth',every on chain/.append style={join}, every join/.style={->}}
\tikzstyle{labeled}=[execute at begin node=$\scriptstyle,   execute at end node=$]

\tikzset{
>=stealth',
help lines/.style={dashed, thick},
axis/.style={<->},
important line/.style={thick},
connection/.style={thick, dotted},
}

\def \spacetimeaxis#1 {
\begin{tikzpicture}
  	\node[anchor=south west,inner sep=0] (img) at (0,0) {#1};
	\node[below=of img, node distance=0cm, yshift=1.1cm, xshift=-0.3cm, font=\color{gray}] {$x$};
  	\node[left=of img, node distance=0cm, rotate=90, anchor=center, yshift=-0.8cm,font=\color{gray}] {$t$};
\end{tikzpicture}
}
	
\newcommand {\R}		 {\mathbf{r}}

\newcommand {\Ieff}  {I_{\rm eff}}

\newcommand {\Rd}    {{\mathds{R}}^d}

\def \EnergyNorm#1  {{\mid\!\mid\!\mid #1 \mid\!\mid\!\mid}^2 }   

\def \dvrg       {\mathrm{div}}	

\def \traspose#1 {{#1}^{rm T}}
\def \laplace    {\Delta}
\def \EI           {{\rm E \!\!\! Id}}

\newcommand {\flux}     {\boldsymbol {y}}

\newcommand {\hhat} {\hat{h}} 
\newcommand {\Khat} {\widehat{K}} 
\newcommand {\Bhat} {\widehat{B}} 
\newcommand {\Rhat} {\widehat{R}} 
\newcommand {\Qhat} {\widehat{Q}} 

\def \dxt      {\mathrm{\:d}x\mathrm{d}t}

\def\L#1{L^{#1}}
\def\H#1{H^{#1}}

\def\HD#1#2{H^{#1}_{#2}}

\def\mdI{\overline{\mathrm m}^{{\rm I}, 2}_{\mathrm{d}}}

\def\mfI{\overline{\mathrm m}^{{\rm I}, 2}_{\mathrm{eq}}}

\def\mdIK{\overline{\mathrm m}^{{\rm I}, 2}_{\mathrm{d}, K}}

\def\mdInosq{\overline{\mathrm m}^{{\rm I}}_{\mathrm{d}}}

\def\mfInosq{\overline{\mathrm m}^{{\rm I}}_{\mathrm{eq}}}

\newcommand {\CFriedrichs} {C_{{\rm F}}}

\newcommand*\rfrac[2]{{}^{#1}\!/_{#2}}

\newtheorem{theorem}{Theorem}{\bf}{\it}
\newtheorem{corollary}{Corollary}{\bf}{\it}
{\bf}{\it}
{\bf}{\it}
{\bf}{}
{\bf}{\it}

\def\ProofBegin{\noindent{\bf Proof:} \:}
\def\ProofEnd{{\hfill $\square$}}

\definecolor{formalshade}{rgb}{0.95,0.95,1}
